\documentclass[11pt]{amsart}
\usepackage{amssymb, adjustbox, enumerate, amsbsy, stmaryrd}
\usepackage{amsmath}
\usepackage[mathscr]{eucal}
\usepackage{amsthm}
\numberwithin{equation}{section}

\usepackage{geometry}

\usepackage[symbol]{footmisc}

\usepackage{amsfonts, amssymb, amscd}

\usepackage{todonotes}
  \newcommand\luo[1]{\todo[color=red!40]{#1}} %Luo
\usepackage{bm}
\usepackage{verbatim}
\usepackage{mathrsfs}
\usepackage{graphicx}
\usepackage[all]{xy}
\usepackage{tikz-cd}
\usepackage{subcaption}
\usepackage{listings}
\usepackage{subfiles}
\usepackage[toc,page]{appendix}
\usepackage{mathtools}
\usepackage{comment}
\usepackage{enumerate}
\usepackage{enumitem}
\usepackage[linesnumbered,ruled]{algorithm2e}

\usepackage{graphicx}
\graphicspath{{images/}}

\usepackage{appendix}
\usepackage{hyperref}
\hypersetup{
    colorlinks=true,
    citecolor=red,
    linkcolor=blue,
    filecolor=magenta,      
    urlcolor=red,
}
\newcommand{\bQ}{\mathbb{Q}}

\newcommand{\Cc}{\mathbb{C}}
\newcommand{\KK}{\mathbb{K}}
\newcommand{\Pp}{\mathbb{P}}
\newcommand{\Qq}{\mathbb{Q}}

\newcommand{\Rr}{\mathbb{R}}

\newcommand{\Zz}{\mathbb{Z}}

\newcommand{\aaa}{\mathfrak{a}}
\newcommand{\bbb}{\mathfrak{b}}
\newcommand{\ddd}{\mathfrak{d}}
\newcommand{\mmm}{\mathfrak{m}}
\newcommand{\rrr}{\mathfrak{r}}

\newcommand{\Span}{\operatorname{Span}}

\newcommand{\Center}{\operatorname{center}}

\newcommand{\Exc}{\operatorname{Exc}}

\newcommand{\mld}{{\rm{mld}}}

\newcommand{\fm}{\mathfrak{m}}

\newcommand{\lcm}{\operatorname{lcm}}

\newcommand{\lct}{\operatorname{lct}}
\newcommand{\ct}{\operatorname{ct}}

\newcommand{\CT}{\operatorname{CT}}

\newcommand{\Supp}{\operatorname{Supp}}

\newcommand{\mult}{\operatorname{mult}}
\newcommand{\codim}{\operatorname{codim}}
\newcommand{\Spec}{\operatorname{Spec}}

\newcommand{\Ii}{{\Gamma}}

\newcommand{\Ss}{\mathcal{S}}

\newcommand{\Zzzero}{{\mathbb Z^n_{\ge0}}}

\newtheorem{thm}{Theorem}[section]
\newtheorem{conj}[thm]{Conjecture}
\newtheorem{cor}[thm]{Corollary}
\newtheorem{lem}[thm]{Lemma}

\newtheorem{prop}[thm]{Proposition}
\newtheorem{ques}[thm]{Question}

\newtheorem{claim}[thm]{Claim}
\theoremstyle{definition}
\newtheorem{rem}[thm]{Remark}

\newtheorem{deflem}[thm]{Definition-Lemma}

\theoremstyle{definition}
\newtheorem{defn}[thm]{Definition}

\begin{document}

\title{ACC for minimal log discrepancies of terminal threefolds}
%threefold local canonical complements

%\author{Jingjun Han, Jihao Liu, and Yujie Luo}

%for terminal threefolds with DCC coefficients


\begin{abstract}
We prove that the ACC conjecture for minimal log discrepancies holds for threefolds in $[1-\delta,+\infty)$, where $\delta>0$ only depends on the coefficient set. We also study Reid's general elephant for pairs, and show Shokurov's conjecture on the existence of $(\epsilon,n)$-complements for threefolds for any $\epsilon\geq 1$. As a key important step, we prove the uniform boundedness of divisors computing minimal log discrepancies for terminal threefolds. We show the ACC for threefold canonical thresholds, and that the set of accumulation points of threefold canonical thresholds is equal to $\{0\}\cup\{\frac{1}{n}\}_{n\in\mathbb Z_{\ge 2}}$ as well. 
\end{abstract}

\author[J. ~Han]{Jingjun Han}
\address{Shanghai Center for Mathematical Sciences, Fudan University, Jiangwan Campus, Shanghai, 200438, China}
\email{hanjingjun@fudan.edu.cn}
\address{Department of Mathematics, The University of Utah, Salt Lake City, UT 84112, USA}
\email{jhan@math.utah.edu}
\address{Mathematical Sciences Research Institute, Berkeley, CA 94720, USA}
\email{jhan@msri.org}

	\author[J. ~Liu]{Jihao Liu}
\address{Department of Mathematics, Northwestern University, 2033 Sheridan Rd, Evanston, IL 60208}
\email{jliu@northwestern.edu}

\author[Y. ~Luo]{Yujie Luo}
\address{Department of Mathematics, Johns Hopkins University, Baltimore, MD 21218, USA}
\email{yluo32@jhu.edu}

\date{\today}

\maketitle

\tableofcontents

\section{Introduction}
We work over the field of complex numbers $\mathbb C$. A \emph{threefold} is a normal quasi-projective variety of dimension 3. ACC stands for the ascending chain condition whilst DCC stands for the
descending chain condition.

The minimal log discrepancy (mld) of a pair $(X,B)$, denoted by $\mld(X,B)$, is defined to be the infimum of log discrepancies of all prime divisors that are exceptional over $X$. It is one of the most basic but important invariants in birational geometry. A better understanding of this invariant will strengthen our knowledge of algebraic varieties. In particular, Shokurov \cite{Sho04a} proved that his ACC conjecture for mlds \cite[Problem 5]{Sho88} together with Ambro's lower-semicontinuity conjecture for mlds \cite[Conjecture 0.2]{Amb99} imply the termination of flips in the minimal model program (MMP).

\medskip

In this paper, we focus on the ACC conjecture for mlds. This conjecture is only known in full generality for surfaces \cite{Ale93} (see \cite{Sho94b,HL20} for other proofs), toric pairs \cite{Amb06}, and exceptional singularities \cite{HLS19}, and was widely open in general in dimension $\geq 3$ before this paper even for terminal threefold pairs (see \cite{HLS19} and references therein for other related results). The following is our first main theorem.

\begin{thm}\label{thm: intro global canonical mld acc}
Let $\Ii\subset [0,1]$ be a DCC set. Then there exists a positive real number $\delta$ depending only on $\Ii$, such that 
$$\{\mld(X,B)\mid \dim X=3, B\in \Ii\}\cap [1-\delta,+\infty)$$
satisfies the ACC, where $B\in\Ii$ means that the coefficients of $B$ belong to the set $\Ii$.
\end{thm}

Theorem \ref{thm: intro global canonical mld acc} solves the conjecture for terminal threefold pairs. All recent progress towards the ACC conjecture for mlds for threefolds \cite[Theorem 1.3]{Kaw15b}, \cite[Corollary 1.5]{Nak16}, \cite[Theorem 1.3]{Jia21} are special cases of Theorem \ref{thm: intro global canonical mld acc}. Indeed, we prove a slightly stronger version of Theorem \ref{thm: intro global canonical mld acc} for germs $(X\ni x,B)$ instead of pairs $(X,B)$, see Theorem \ref{thm:  terminal mld acc}. 

%Theorem \ref{thm: intro global canonical mld acc} not only recovers all recent progress towards the ACC conjecture for mlds for threefolds (\cite[Theorem 1.3]{Kaw15b}, \cite[Corollary 1.5]{Nak16}, \cite[Theorem 1.3]{Jia21}), but also solves the conjecture for terminal threefold pairs. Indeed, we prove a slightly stronger version of Theorem \ref{thm: intro global canonical mld acc} for germs $(X\ni x,B)$ instead of pairs $(X,B)$, see Theorem \ref{thm:  terminal mld acc}. 

 Although the MMP and the abundance conjecture are settled in dimension 3, and we even have a complete classification of terminal threefold singularities including flips as well as divisorial contractions (cf. \cite{Mor85,Rei87,KM92,Kaw01,Kaw02,Kaw03,Kaw05,Kaw12,Yam18}), the ACC conjecture for mlds for terminal threefold pairs remains open. Thus Theorem \ref{thm: intro global canonical mld acc} strengthens our grasp on terminal threefolds. Note that many important results in birational geometry were first observed and proved for terminal threefolds before generalizing to other larger classes of singularities and to higher dimensions, such as the existence of flips \cite{Mor88,Sho92}. It is our hope that Theorem \ref{thm: intro global canonical mld acc} will shed light on the study of algebraic varieties in higher dimensions.

\medskip

We emphasis that our proof of Theorem \ref{thm: intro global canonical mld acc} not only depends on the classification of threefold terminal singularities \cite{Kaw01,Kaw02,Kaw03,Kaw05,Yam18}, but also heavily relies on the state-of-the-art theory on the boundedness of complements \cite{Bir19,HLS19,Sho20}. Moreover, the proof of Theorem \ref{thm: intro global canonical mld acc} is intertwined with the proof of Theorem \ref{thm: mn intro DCC coeff}, which shows a uniform bound of $a(E,X,0)$ for some mld-computing divisor $E$.

\begin{thm}\label{thm: mn intro DCC coeff}
Let $\Ii\subset [0,1]$ be a DCC set. Then there exists a positive integer $l$ depending only on $\Ii$ satisfying the following. Assume that $(X\ni x,B)$ is a threefold pair such that $X$ is terminal, $B\in\Ii$, and $\mld(X\ni x,B)\geq 1$. Then there exists a prime divisor $E$ over $X\ni x$, such that $a(E,X,B)=\mld(X\ni x,B)$ and $a(E,X,0)\leq l$.
\end{thm}

Theorem \ref{thm: mn intro DCC coeff} generalizes a result of Kawakita \cite[Theorem 1.3(ii)]{Kaw21}, which requires $X$ to be smooth and $\Ii$ to be a finite set. When $X\ni x$ is a fixed germ and $\Ii$ is a finite set, the existence of such a uniform bound $l$ was predicted by Nakamura \cite[Conjecture 1.1]{MN18}, and it is equivalent to the ACC conjecture for mlds for fixed germs (cf. \cite[Theorem 4.6]{Kaw21}). When $\dim X=2$, Theorem \ref{thm: mn intro DCC coeff} is proved by G.~Chen and the first named author for smooth germs \cite[Theorem B.1]{CH21}, and is proved by the first and the third named authors for all lc germs \cite[Theorem 1.1]{HL20}. All the above results indicate that the uniform boundedness of $a(E,X,0)$ should hold under more general settings (see Conjecture \ref{conj: bdd mld computing divisor}). We refer the reader to Theorem~\ref{thm: toric mld bdd} for a proof of Conjecture \ref{conj: bdd mld computing divisor} for toric varieties.

\medskip

Theorem \ref{thm: intro global canonical mld acc} has many applications towards other topics on threefolds, both for local singularities and global algebraic structures. We list a few of them in the rest part of the introduction, and refer the reader to Section \ref{sec: Ideal-adic versions} for more applications.

\medskip

\noindent\textbf{Reid's general elephant for pairs and Shokurov's boundedness of complements conjecture}. For a terminal threefold singularity $x\in X$, we say that a Weil divisor $H$ is an \emph{elephant} of $x\in X$ if $H\in |-K_X|$ and $(X,H)$ is canonical near $x$. By \cite[6.4(B)]{Rei87}, elephant exists for any terminal threefold singularity. As an application of Theorem \ref{thm: intro global canonical mld acc}, we generalize Reid's general elephant theorem to the category of pairs.

\begin{thm}\label{thm: intro cc}
Let $\Ii\subset [0,1]\cap\Qq$ be a finite set. Then there exists a positive integer $N$ depending only on $\Ii$ satisfying the following.

Let $(X\ni x,B)$ be a threefold pair such that $X$ is terminal, $B\in\Ii$, and $(X,B)$ is canonical near $x$. Then on a neighborhood of $x$, there exists an element $G\in |-N(K_X+B)|$ such that $(X,B+\frac{1}{N}G)$ is canonical near $x$.
\end{thm}

We remark that in Theorem~\ref{thm: intro cc}, if $x\in X$ is a threefold terminal singularity that is not smooth, then we can choose $G\in |-N(K_X+B)|$ such that $(X,B+\frac{1}{N}G)$ is canonical near $x$ and $\mld(X\ni x, B+\frac{1}{N}G)=1$ (see Theorem~\ref{thm: N-complements for canonical pair finite rational coeff}).

It is worth mentioning that Reid's general elephant theorem is a special case of Theorem \ref{thm: intro cc} when $\Ii=\{0\}$ and $x$ is a closed point, where we can take $N=1$. We refer the reader to Koll\'ar, Mori, Prohokorv, Kawakita's previous works  \cite{KM92,Kaw02,MP08,MP09,MP21} and reference therein for other results on general elephant for terminal threefolds. 

Theorem \ref{thm: intro cc} is closely related to Theorem \ref{thm: intro ecc}, which gives an affirmative answer to Shokurov's conjecture on the boundedness of $(\epsilon,N)$-complements (\cite[Conjecture 1.1]{CH21}; see \cite[Conjecture]{Sho04b}, \cite[Conjectures 1.3, 1.4]{Bir04} for some embryonic forms) for terminal threefold germs. We refer the reader to Subsection~\ref{section: complements} for basic notation on complements.

\begin{thm}\label{thm: intro ecc}
Let $\epsilon\geq 1$ be a real number and $\Ii\subset [0,1]$ a DCC set. Then there exists a positive integer $N$ depending only on $\epsilon$ and $\Ii$ satisfying the following. 

Assume that $(X\ni x,B)$ is a threefold pair, such that $X$ is terminal, $B\in\Ii$, and $\mld(X\ni x,B)\geq\epsilon$. Then there exists an $N$-complement  $(X\ni x,B^+)$ of $(X\ni x,B)$ such that $\mld(X\ni x,B^+)\geq\epsilon$.
\end{thm}

We refer the reader to Theorem \ref{ref: main epsilon comp threefolds p version} for a more detailed version of Theorem \ref{thm: intro ecc}.

We remark that the boundedness of $(0,N)$-complements proved in \cite{Bir19,HLS19} plays an important role in several breakthroughs in birational geometry including the proof of Birkar-Borisov-Alexeev-Borisov Theorem and the
openness of K-semistability in families of log Fano pairs (cf. \cite{Bir21,Xu20}), while the conjecture on the boundedness of $(\epsilon,N)$-complements was only known for surfaces (\cite[Main Theorem 1.6]{Bir04}, \cite[Theorem 1.6]{CH21}) before. 

Shokurov suggested that Theorem \ref{thm: intro ecc} may hold without assuming $\Ii$ is a DCC set as Theorem \ref{thm: intro ecc} implies this holds if the number of components of $B$ is bounded from above (cf. \cite[Proof of Corollary 1.8]{CH21}). 

\medskip

A special case of Theorem \ref{thm: intro ecc} gives an affirmative answer on Shokurov's index conjecture (cf. \cite[Conjecture 7.3]{CH21}, \cite[Question 5.2]{Kaw15a}) for terminal threefolds.

\begin{thm}\label{thm: index conjecture DCC threefold}
Let $\epsilon\geq 1$ be a real number and $\Ii\subset [0,1]\cap\Qq$ a DCC set. Then there exists a positive integer $I$ depending only on $\epsilon$ and $\Ii$ satisfying the following. 

Let $(X\ni x, B)$ be a threefold pair such that $X$ is terminal, $B\in\Ii$, and $\mld(X\ni x,B)=\epsilon$. Then $I(K_X+B)$ is Cartier near $x$.
\end{thm}
Kawakita showed that for any canonical threefold singularity $x\in X$ with $\mld(X\ni x)=1$, $IK_X$ is Cartier for some $I\leq 6$ \cite[Theorem 1.1]{Kaw15a}, hence \cite[Theorem 1.1]{Kaw15a} can be viewed as a complementary result to Theorem \ref{thm: index conjecture DCC threefold} when $\epsilon=1$. We refer the reader to Theorem \ref{thm: index conjecture threefold} for an explicit bound of $I$ when $\Ii$ is a finite set.

\medskip

%\textit{ACC for $a$-lc thresholds and canonical thresholds}.
\noindent\textbf{Other Applications}. As an application of our main theorems, we show the ACC for $a$-lc thresholds (a generalization of lc thresholds, see Definition \ref{defn: alct local}) for terminal threefolds when $a\geq 1$:

\begin{thm}[ACC for $a$-lc thresholds for terminal threefolds]\label{thm: alct acc terminal threefold}
Let $a\geq 1$ be a real number, and $\Ii\subset [0,1],~\Ii'\subset [0,+\infty)$ two DCC sets. Then the set of $a$-lc thresholds,
   $$\{a\text{-}\lct(X\ni x,B;D)\mid \dim X=3, X\text{ is terminal}, B\in\Ii, D\in\Ii'\},$$
     satisfies the ACC.
\end{thm}

Theorem \ref{thm: alct acc terminal threefold} implies the ACC for canonical thresholds in dimension $3$:
\begin{thm}\label{thm: 3fold acc ct}
Let $\Ii\subset[0,1]$ and $\Ii'\subset[0,+\infty)$ be two DCC sets. Then the set
$$\CT(3,\Ii,\Ii'):=\{\ct(X,B;D)\mid\dim X=3, B\in \Ii,D\in\Ii'\}$$
satisfies the ACC.
\end{thm}

It is worth to mention that the canonical thresholds in dimension $3$ is deeply related to Sarkisov links in dimension $3$ (cf. \cite{Cor95,Pro18}). Moreover, we have a precise description of the accumulation points of $\CT(3,\{0\},\Zz_{\ge 1})$:

\begin{thm}\label{thm: Accumulation points}
The set of accumulation points of $\CT(3,\{0\},\Zz_{\ge 1})$ is $\{0\}\cup\{\frac{1}{m}\mid m\in\Zz_{\ge 2}\}$.
\end{thm}

Theorem \ref{thm: Accumulation points} plays a crucial role in the proof of  Theorem \ref{thm: intro global canonical mld acc}. We refer the reader to Theorem \ref{thm: explicite accumulation points of canonical thresholds} for a more detailed version of Theorem \ref{thm: Accumulation points}. 

\cite[Theorem 1.7]{Ste11} proved Theorem \ref{thm: 3fold acc ct} when $\Ii=\{0\}$, $\Ii'=\Zz_{\ge1}$, and $X$ is smooth, and \cite[Theorem 1.2]{Che19} proved Theorem \ref{thm: 3fold acc ct} when $\Ii=\{0\}$ and $\Ii'=\Zz_{\ge1}$. \cite[Theorem 1.3]{Che19} proved that $\frac{1}{2}$ is the largest accumulation point of $\CT(3,\{0\},\Zz_{\ge1})$. We refer the reader to \cite{Shr06,Pro08} for other related results. 

After we completed this paper, we were informed by J.J. Chen and J.A. Chen that J.J. Chen has independently obtained Theorems \ref{thm: 3fold acc ct} and \ref{thm: Accumulation points} \cite[Theorems 1.1, 1.2]{Che22}.

\medskip

Theorem \ref{thm: non-canonical klt CY threefold bdd flop} is another application of our main theorems:

\begin{thm}\label{thm: non-canonical klt CY threefold bdd flop}
Let $\Ii\subset[0,+\infty)$ be a DCC set. Then the set of non-canonical klt threefold log Calabi-Yau pairs $(X,B)$ with $B\in \Ii$ forms a bounded family modulo flops.
\end{thm}
Theorem \ref{thm: non-canonical klt CY threefold bdd flop} is a generalization of \cite[Theorem 1.4]{BDS20} for threefolds. Jiang proved Theorem \ref{thm: non-canonical klt CY threefold bdd flop} for the case when $\Ii=\{0\}$ \cite[Theorem 1.6]{Jia21}.  

%\han{maybe delete for submit version}

We also remark that the assumption ``non-canonical klt'' is natural and necessary as rationally connected Calabi-Yau varieties are not canonical, and the set of $(X:=Y\times \Pp^1,F_1+F_2)$ is not birationally bounded, where $Y$ takes all K3 surfaces and $F_1,F_2$ are two fibers of $X\to \Pp^1$.

\medskip

The proofs of our main theorems depend on many different theories, such as the classification of threefold divisorial contractions, the theory of lc complements, the theory of uniform rational polytopes, and singular Riemann-Roch formula. Table \ref{tbl: flowchart} below provides a flowchart of the structure of the paper.

\begin{table}[ht]
\caption{Flowchart of the structure of the paper}\label{tbl: flowchart}
\begin{adjustbox}{width=0.8\textwidth,center}
%\begin{center}
$\xymatrix{
& *+[F]\txt{Special cases$^{(\star)}$ of Theorem \ref{thm: mn intro DCC coeff}\\ (Lemmas \ref{lem: bounded index Mustata nakamura conjecture}, \ref{lem: can extract mld place smooth case}, \ref{lem: can extract mld place cA/n case}, Theorem \ref{thm: Nakamura Mustata Conjecture})}\ar@{->}[d]\ar@/^7pc/[ddd]\\
& *+[F]\txt{ACC for cts\\ (Theorems \ref{thm: 3fold acc ct}, \ref{thm: ACC ct local with base terminal})}\ar@{->}[dl]_{(\star)}\ar@{->}[d]_{(\dagger)}\ar@/^4pc/[dd]\\
*+[F]\txt{Accumulation points\\ (Theorem \ref{thm: Accumulation points})}\ar@{->}[dr] & *+[F]\txt{$1$-gap for mlds\\ (Theorem \ref{thm: 1-gap pair})}\ar@{->}[d] \\
& *+[F]\txt{ACC for mlds$^{(\dagger)}$\\ (Theorems \ref{thm: intro global canonical mld acc}, \ref{thm:  terminal mld acc}),\\ ACC for $a$-lcts\\ Theorem \ref{thm: alct acc terminal threefold}) }\ar@/^1pc/[ddl]_{(\ddagger)}\ar@{->}[dl]\ar@{->}[d]_{(\amalg)}\ar@/^5pc/[dd] \\
*+[F]\txt{Log Calabi-Yau\\ birational boundedness\\ (Theorem \ref{thm: non-canonical klt CY threefold bdd flop})} & *+[F]\txt{Index conjecture\\ (Theorems \ref{thm: index conjecture DCC threefold}, \ref{thm: index conjecture threefold})}\ar@{->}[d] \\
*+[F]\txt{Uniform boundedness\\ (Theorem \ref{thm: mn intro DCC coeff})} & *+[F]\txt{General elephants$^{(\dagger)}$ and local complements$^{(\ddagger)}$\\ (Theorems \ref{thm: intro cc}, \ref{thm: intro ecc})}\\
}$ 
%\end{center}
\end{adjustbox}
\begin{flushleft}
$(\star)$: $3$-fold divisorial contraction classification \cite{Kaw01,Kaw02,Kaw03,Kaw05,Yam18}. $(\dagger)$: The boundedness of lc $n$-complements \cite{Bir19,HLS19}. $(\ddagger)$: The theory of uniform rational polytopes \cite{HLS19,CH21}. $(\amalg)$: Singular Riemann-Roch formula (Subsection \ref{subsection: rr}).
\end{flushleft}
\end{table}

\medskip

%start with the proof of Theorem \ref{thm: intro global canonical mld acc} and
\noindent\textbf{Sketch of the proof of Theorem \ref{thm: intro global canonical mld acc}}. As we mentioned earlier, the proof is intertwined with the proof of Theorem \ref{thm: mn intro DCC coeff}. For simplicity, we only deal with the interval $[1,+\infty)$. Suppose on the contrary, there exists a sequence of threefold germs $\{(X_i\ni x_i, B_i)\}_{i=1}^{\infty}$, where $X_i$ is terminal and $B_i\in\Ii$ for each $i$, such that $\{\mld(X_i\ni x_i,B_i)\}_{i=1}^{\infty}\subset(1,+\infty)$ is strictly increasing. Then the index of $X_i\ni x_i$ is bounded from above. We may assume that there exists a finite set $\Ii_0$, such that $\lim_{i\to +\infty}B_i=\bar{B}_i$ and $\bar{B}_i\in\Ii_0$. By the ACC for threefold canonical thresholds (Theorem \ref{thm: 3fold acc ct}) and \cite[Corollary 1.3]{Nak16}, we may assume that $\mld(X_i\ni x_i,\bar{B}_i)=\alpha\ge 1$. for some constant $\alpha$. In order to derive a contradiction, it suffices to show a special case of Theorem \ref{thm: mn intro DCC coeff}, that is, there exists a prime divisor $\bar{E}_i$ over $X_i\ni x_i$, such that 
$a(\bar{E}_i,X_i,\bar{B}_i)=\mld(X_i\ni x_i,\bar{B}_i)=\alpha$, and $a(\bar{E}_i,X_i,0)\le l$ for some constant number $l$, see Step 3 of the proof of Theorem \ref{thm:  terminal mld acc}. If $\alpha>1$, and $x_i\in X_i$ is neither smooth nor of $cA/n$ type for each $i$, then we show the special case by the uniform canonical rational polytopes (Theorem \ref{thm: uniform canonical polytope}) and the accumulation points of the set of canonical thresholds in dimension 3 (Theorem \ref{thm: Accumulation points}). Otherwise, we show the following key fact: there exists a divisorial contraction $Y_i\to X_i$ from a terminal variety $Y_i$ which extracts a prime divisor $E_i'$ over $X_i\ni x_i$ such that $a(E_i',X_i,\bar{B}_i)=\mld(X_i\ni x_i,\bar{B}_i)$. Note that when $\alpha=1$, we may show the fact by standard tie breaking trick even in higher dimensions, see Lemma \ref{lem: can extract divisor computing ct that is terminal strong version}. The case when $\alpha>1$ and either $x_i\in X_i$ is smooth or of type $cA/n$, which is one of our key observations, depends on the proofs of the classification of divisorial contractions for terminal threefolds \cite{Kaw01,Kaw02,Kaw03,Kaw05,Yam18}, see Lemmas \ref{lem: can extract mld place smooth case}, \ref{lem: can extract mld place cA/n case}. Finally, Theorem \ref{thm: intro global canonical mld acc} follows from the key fact and Lemma \ref{lem: bounded index Mustata nakamura conjecture}. We remark that Lemma \ref{lem: bounded index Mustata nakamura conjecture} implies the ACC for threefold canonical thresholds (Theorem \ref{thm: 3fold acc ct}).

\medskip
\noindent \textbf{Acknowledgement.} JH would like to thank Vyacheslav V. Shokurov who teaches him the theory of complements, and Caucher Birkar for suggesting him the question on the boundedness of complements for terminal threefolds in August 2019. Part of this work was done when the YL visited Fudan University in August 2021. YL would like to thank their hospitality.

We would like to thank Caucher Birkar, Guodu Chen, Jungkai Alfred Chen, Jheng-jie Chen, Christopher D. Hacon, Yong Hu, Chen Jiang, Junpeng Jiao, Masayuki Kawakita, Yuchen Liu, Lu Qi, Vyacheslav V. Shokurov, Charles Stibitz, Lingyao Xie, Chenyang Xu, Qingyuan Xue, and Ziquan Zhuang for helpful discussions and answering our questions. This work was supported by a grant from the Simons Foundation (Grant Number 814268, MSRI). YL would like to thank his advisor Chenyang Xu for constant support, encouragement, and numerous inspiring conversations.

\section{Preliminaries}

We adopt the standard notation and definitions in \cite{KM98,BCHM10} and will freely use them. All varieties are assumed to be normal quasi-projective and all birational morphisms are assumed to be projective. We denote by $\xi_n$ the $n$-th root of unity $e^{\frac{2\pi i}{n}}$, and denote by $\Cc[x_1,\dots,x_d]$ (resp. $\Cc\{x_1,\dots,x_d\}$, $\Cc[[x_1,\dots,x_d]]$) the ring of power series (resp. analytic power series, formal power series) with the coordinates $x_1,\dots,x_d$.

\begin{comment}
Let $\KK=\Qq$ or $\Rr$ be either the rational number field $\Qq$ or the real number field $\Rr$.
 Let $X$ be a normal variety. A {\it $\KK$-divisor} is a finite $\KK$-linear combination $D=\sum d_{i} D_{i}$ of prime Weil divisors $D_{i}$, and $d_{i}$ denotes the {\it coefficient} of $D_i$ in $D$. A {\it $\KK$-Cartier divisor} is a $\KK$-linear combination of Cartier divisors. 
 
  We use $\sim_{\KK}$ to denote the $\KK$-linear equivalence between $\KK$-divisors. For a projective morphism $X\to Z$, we use $\sim_{\KK,Z}$ to denote the relative $\KK$-linear equivalence. % and use $\equiv_{Z}$ to denote the relative numerical equivalence. 
\end{comment}

\subsection{Pairs and singularities}\label{the definitions of singularities and pairs}

\begin{defn}\label{defn contraction}
A \emph{contraction} is a projective morphism $f: Y\rightarrow X$ such that $f_*\mathcal{O}_Y=\mathcal{O}_X$. In particular, $f$ is surjective and has connected fibers.
\end{defn}

\begin{defn}\label{defn: divisorial contraction}
Let $f: Y\rightarrow X$ be a birational morphism, and $\Exc(f)$ the exceptional locus of $f$. We say that $f$ is a \emph{divisorial contraction} of a prime divisor $E$ if $\Exc(f)=E$ and $-E$ is $f$-ample.
\end{defn}
%$f$ is a contraction,

\begin{defn}[Pairs, {cf. \cite[Definition 3.2]{CH21}}] \label{defn sing}
A \emph{pair} $(X/Z\ni z, B)$ consists of a contraction $\pi: X\rightarrow Z$, a (not necessarily closed) point $z\in Z$, and an $\mathbb{R}$-divisor $B\geq 0$ on $X$, such that $K_X+B$ is $\Rr$-Cartier over a neighborhood of $z$ and $\dim z<\dim X$. If $\pi$ is the identity map and $z=x$, then we may use $(X\ni x, B)$ instead of $(X/Z\ni z,B)$. In addition, if $B=0$, then we use $X\ni x$ instead of $(X\ni x,0)$. When we consider a pair $(X\ni x, \sum_{i} b_iB_i)$, where $B_i$ are distinct prime divisors and $b_i>0$, we always assume that $x\in \Supp B_i$ for each $i$.

If $(X\ni x,B)$ is a pair for any codimension $\geq 1$ point $x\in X$, then we call $(X,B)$ a pair. A pair $(X\ni x, B)$ is called a \emph{germ} if $x$ is a closed point. We say $x\in X$ is a \emph{singularity} if $X\ni x$ is a germ.
\end{defn}

\begin{defn}[Singularities of pairs]\label{defn: relative mld}
 Let $(X/Z\ni z,B)$ be a pair associated with the contraction $\pi: X\to Z$, and let $E$ be a prime divisor over $X$ such that $z\in\pi(\Center_X E)$. Let $f: Y\rightarrow X$ be a log resolution of $(X,B)$ such that $\Center_Y E$ is a divisor, and suppose that $K_Y+B_Y=f^*(K_X+B)$ over a neighborhood of $z$. We define $a(E,X,B):=1-\mult_EB_Y$ to be the \emph{log discrepancy} of $E$ with respect to $(X,B)$. 
 
 For any prime divisor $E$ over $X$, we say that $E$ is \emph{over} $X/Z\ni z$ if $\pi(\Center_X E)=\bar z$. If $\pi$ is the identity map and $z=x$, then we say that $E$
 is \emph{over} $X\ni x$. We define
 $$\mld(X/Z\ni z,B):=\inf\{a(E,X,B)\mid E\text{ is over }Z\ni z\}$$
 to be the \emph{minimal log discrepancy} (\emph{mld}) of $(X/Z\ni z,B)$.
 
 Let $\epsilon$ be a non-negative real number. We say that $(X/Z\ni z,B)$ is lc (resp. klt, $\epsilon$-lc,$\epsilon$-klt) if $\mld(X/Z\ni z,B)\geq 0$ (resp. $>0$, $\geq\epsilon$, $>\epsilon$). We say that $(X,B)$ is lc (resp. klt, $\epsilon$-lc, $\epsilon$-klt) if $(X\ni x,B)$ is lc (resp. klt, $\epsilon$-lc, $\epsilon$-klt) for any codimension $\geq 1$ point $x\in X$. 
 
 We say that $(X,B)$ is \emph{canonical} (resp. \emph{terminal}, \emph{plt}) if $(X\ni x,B)$ is $1$-lc (resp. $1$-klt, klt) for any codimension $\geq 2$ point $x\in X$.% We say that $(X/Z\ni z,B)$ is \emph{canonical} (resp. \emph{terminal}) if $(X/Z\ni z,B)$ is $1$-lc (resp. $1$-klt).
 
 For any (not necessarily closed) point $x\in X$, we say that $(X,B)$ is lc (resp. klt, $\epsilon$-lc, $\epsilon$-klt, canonical, terminal) near $x$ if $(X,B)$ is lc (resp. klt, $\epsilon$-lc, $\epsilon$-klt, canonical, terminal) in a neighborhood of $x$. If $X$ is (resp. klt, $\epsilon$-lc, $\epsilon$-klt, canonical, terminal) near a closed point $x$, then we say that $x\in X$ is an lc (resp. klt, $\epsilon$-lc, $\epsilon$-klt, canonical, terminal) singularity. We remark that if $(X\ni x, B)$ is lc, then $(X,B)$ is lc near $x$.
\end{defn}

\begin{defn}\label{defn: alct local}  Let $a$ be a non-negative real number, $(X\ni x,B)$ (resp. $(X,B)$) a pair, and $D\geq 0$ an $\Rr$-Cartier $\Rr$-divisor on $X$. We define
$$a\text{-}\lct(X\ni x,B;D):=\sup\{-\infty,t\mid t\geq 0, (X\ni x,B+tD)\text{ is }a\text{-}lc\}$$
$$\text{(resp. }a\text{-}\lct(X,B;D):=\sup\{-\infty,t\mid t\geq 0, (X,B+tD)\text{ is }a\text{-}lc\} \text{)}$$
to be the \emph{$a$-lc threshold} of $D$ with respect to $(X\ni x,B)$ (resp. $(X,B)$). We define
$$\ct(X\ni x,B;D):=1\text{-}\lct(X\ni x,B;D)$$
$$\text{(resp. }\ct(X,B;D):=\sup\{-\infty,t\mid t\geq 0, (X,B+tD)\text{ is canonical}\} \text{)}$$
to be the \emph{canonical threshold} of $D$ with respect to $(X\ni x,B)$ (resp. $(X,B)$). We define $\lct(X\ni x,B;D):=0\text{-}\lct(X\ni x,B;D)$ (resp. $\lct(X,B;D):=0\text{-}\lct(X,B;D)$) to be the \emph{lc threshold} of $D$ with respect to $(X\ni x,B)$ (resp. $(X,B)$).
\end{defn}

\begin{lem}\label{lem: Terminal blow up}
Let $(X\ni x,B)$ be a pair such that $X$ is terminal and $\dim x=\dim X-2$. Let $E_1$ be the exceptional divisor obtained by blowing up $x\in X$. If $\mult_x B\leq 1$, then $$\mld(X\ni x, B)=a(E_1,X,B)=2-\mult_x B\geq 1.$$
Moreover, $\mld(X\ni x,B)\geq 1$ if and only if $\mult_x B\leq 1$.
\end{lem}
%\han{for submitted version, we may delete the proof.}
\begin{proof}
Since $X$ is terminal, by \cite[Corollary 5.18]{KM98}, $X$ is smooth in codimension $2$. Since $\dim x=\dim X-2$, possibly shrinking $X$ to a neighborhood of $x$, we may assume that $X$ is smooth. By \cite[Lemma 2.45]{KM98}, there exists a sequence of blow-ups
$$X_n\xrightarrow{f_n} X_{n-1}\xrightarrow{f_{n-1}}\cdots\xrightarrow{f_2} X_1\xrightarrow{f_1} X_0:=X,$$
such that for any $1\leq i\leq n$,
 \begin{itemize}
     \item  $f_i$ is a blow-up of $X_{i-1}$ at a point $x_{i-1}$ of codimension at least $2$ with the exceptional divisor $E_i$,
     \item $X_i$ is smooth, $\bar{x}_{i-1}:=\Center_{X_{i-1}}E_n$, and
     \item $a(E_n,X,B)=\mld(X\ni x,B)$.
 \end{itemize}
In particular, $x_{i-1}$ dominates $x$, $x_0=x$, and $\dim x_{i-1}=\dim X-2$ for $1\leq i\leq n$. For any $0\leq i\leq n$, we let $B_{i}$ be the strict transform of $B$ on $X_i$. For any $1\leq i\leq n$, we have $$f_i^*{B_{i-1}}=B_{i}+(\mult_{x_{i-1}} B_{i-1}) E_i.$$ Let $U_{i-1}$ be an open neighborhood of $x_{i-1}$ such that $U_{i-1}$ and $\bar{x}_{i-1}|_{U_{i-1}}$ are both smooth. Then $f_{i}^{-1}(U_{i-1})\cap E_i$ is covered by smooth rational curves that are contracted by $f_i$ and whose intersection numbers with $E_i$ are all equal to $-1$ (cf. \cite[\S2,~Theorem~8.24(c)]{Har77}), from which we may choose a general curve and denote it by $C_i$, such that $C_i\not\subset \Supp B_i$ and $C_i\cap\bar x_i\not=\emptyset$ when $i\not=n$. Thus
 $$0=f_i^*{B_{i-1}}\cdot C_i=(B_i+(\mult_{x_{i-1}} B_{i-1}) E_i)\cdot C_i=B_{i}\cdot C_i-\mult_{x_{i-1}} B_{i-1},$$ 
 which implies that $\mult_{x_{i-1}}B_{i-1}= B_i\cdot C_i$. Since $x_i$ dominates $x_{i-1}$ and $\dim x_i=\dim X-2$ for $1\leq i\leq n-1$, we may choose $C_i$ so that $C_i\not\subset\bar x_i$. For any $1\leq i\leq n-1$, let $\bar C_i$ be the birational transform of $C_i$ on $X_{i+1}$. We have $\bar C_i\not\subset E_{i+1}$. By the projection formula,
$$\mult_{x_{i-1}} B_{i-1}=B_{i}\cdot C_i= f_{i+1}^*B_{i}\cdot\bar C_{i}\geq (\mult_{x_i}B_{i})E_{i+1}\cdot\bar C_{i} \geq  \mult_{x_i}B_{i}.$$
By induction on $i$, we have $1\geq \mult_{x} B\geq \mult_{x_{i}} B_{i}$ for any $0\leq i\leq n-1$, thus
\begin{equation}\label{equ: 2-mult equ 1}
    a(E_1,X,B)=2-\mult_xB\leq 2-\mult_{x_{n-1}}B_{n-1}=a(E_n,X_{n-1},B_{n-1}).
\end{equation}
Moreover, since
$$K_{X_i}+B_i=f_i^*(K_{X_{i-1}}+B_{i-1})+(1-\mult_{x_{i-1}B_{i-1}})E_i\geq f_i^*(K_{X_{i-1}}+B_{i-1})$$ for $1\leq i\leq n$,
by induction, we have $K_{X_{n-1}}+B_{n-1}\geq (f_{1}\circ\dots\circ f_{n-1})^*(K_X+B)$, hence
\begin{equation}\label{equ: 2-mult equ 2}
a(E_n,X_{n-1},B_{n-1})\leq a(E_n,X,B)=\mld(X\ni x,B)\leq a(E_1,X,B).
\end{equation}
Lemma \ref{lem: Terminal blow up} now follows from Inequalities (\ref{equ: 2-mult equ 1}) and (\ref{equ: 2-mult equ 2}).
\end{proof}

\begin{comment}
\begin{thm}[{\cite[Theorem 18.22]{Kol+92}}]\label{thm: number of coefficients local}
Let $(X\ni x,\sum_{i=1}^m b_iB_i)$ be an lc pair such that $K_X$ and $B_i\geq 0$ are $\Qq$-Cartier Weil divisors near $x$ and $x\in\Supp B_i$ for each $i$. Then $\sum_{i=1}^m b_i\le \dim X$.
\end{thm}
\end{comment}

\subsection{Complements}\label{section: complements}
\begin{defn}\label{defn: complement}
Let $n$ be a positive integer, $\epsilon$ a non-negative real number, $\Ii_0\subset (0,1]$ a finite set, and $(X/Z\ni z,B)$ and $(X/Z\ni z,B^+)$ two pairs. We say that $(X/Z\ni z,B^+)$ is an \emph{$(\epsilon,\Rr)$-complement} of $(X/Z\ni z,B)$ if 
\begin{itemize}
    \item $(X/Z\ni z,B^+)$ is $\epsilon$-lc,
    \item $B^+\geq B$, and
    \item $K_X+B^+\sim_{\Rr}0$ over a neighborhood of $z$.
\end{itemize}
We say that $(X/Z\ni z,B^+)$ is an \emph{$(\epsilon,n)$-complement} of $(X/Z\ni z,B)$ if
\begin{itemize}
\item $(X/Z\ni z,B^+)$ is $\epsilon$-lc,
\item $nB^+\geq \lfloor (n+1)\{B\}\rfloor+n\lfloor B\rfloor$, and
\item $n(K_X+B^+)\sim 0$ over a neighborhood of $z$.
\end{itemize}
A $(0,\Rr)$-complement is also called an \emph{$\Rr$-complement}, and a $(0,n)$-complement is also called an \emph{$n$-complement}. We say that $(X/Z\ni z,B)$ is $(\epsilon,\Rr)$-complementary (resp. $(\epsilon,n)$-complementary, $\Rr$-complementary, $n$-complementary) if $(X/Z\ni z,B)$ has an $(\epsilon,\Rr)$-complement (resp. $(\epsilon,n)$-complement, $\Rr$-complement, $n$-complement). 

 We say that $(X/Z\ni z,B^+)$ is a \emph{monotonic $(\epsilon,n)$-complement} of $(X/Z\ni z,B)$ if $(X/Z\ni z,B^+)$ is an $(\epsilon,n)$-complement of $(X/Z\ni z,B)$ and $B^+\geq B$.%A  monotonic $(0,n)$-complement is also called a \emph{monotonic $n$-complement}. 
 
 We say that $(X/Z\ni z,B^+)$ is an \emph{$(n,\Ii_0)$-decomposable $\Rr$-complement} of $(X/Z\ni z,B)$ if there exist a positive integer $k$, $a_1,\dots,a_k\in\Ii_0$, and  $\Qq$-divisors $B_1^+,\dots,B_k^+$ on $X$, such that
\begin{itemize}
\item $\sum_{i=1}^ka_i=1$ and  $\sum_{i=1}^ka_iB_i^+=B^+$,
\item $(X/Z\ni z,B^+)$ is an $\Rr$-complement of $(X/Z\ni z,B)$, and
\item  $(X/Z\ni z,B_i^+)$ is an $n$-complement of itself for each $i$.
\end{itemize}
\begin{comment}
We say that $(X/Z\ni z,B^+)$ is a \emph{canonical $n$-complement} of $(X/Z\ni z,B)$ if $(X/Z\ni z,B^+)$ is an $n$-complement of $(X/Z\ni z,B)$ and $(X/Z\ni z,B^+)$ is canonical.
\end{comment} 
\end{defn}

\begin{thm}[{\cite[Theorem 1.10]{HLS19}}]\label{thm: ni decomposable complement}
Let $d$ be a positive integer and $\Ii\subset [0,1]$ a DCC set. Then there exists a positive integer $n$ and a finite set $\Ii_0\subset (0,1]$ depending only on $d$ and $\Ii$ and satisfy the following. 

Assume that $(X/Z\ni z,B)$ is a pair of dimension $d$ and $B\in\Ii$, such that $X$ is of Fano type over $Z$ and $(X/Z\ni z,B)$ is $\Rr$-complementary. Then $(X/Z\ni z,B)$ has an $(n,\Ii_0)$-decomposable $\Rr$-complement. Moreover, if $\bar\Ii\subset\Qq$, then $(X/Z\ni z,B)$ has a monotonic $n$-complement.
\end{thm}

\subsection{Index of canonical threefolds}

\begin{defn}\label{def: definition of index}
Let $(X\ni x,B)$ be a pair such that $B\in\Qq$, and $(X^{\rm an}\ni x,B^{\rm an})$ the corresponding analytic pair. The \emph{index} (resp. \emph{analytic index}) of $(X\ni x,B)$ is the minimal positive integer $I$ such that $I(K_X+B)$ is  (resp. $I(K_{X^{\rm an}}+B^{\rm an})$) is Cartier near $x$. 
\end{defn}

The following lemma indicates that the index of $X\ni x$ coincides with the analytic index of $X\ni x$. Hence we will not distinguish the index and the analytic index in our paper.

\begin{lem}[{\cite[Lemma~1.10]{Kaw88}}]\label{lem: analytic index and algebraic index}
Let $X$ be a variety and $D$ a Weil divisor on $X$. Let $X^{\rm an}$ be the underlying analytic space of $X$ and $D^{\rm an}$ the underlying analytic Weil divisor of $D$ on $X^{\rm an}$. Then $D^{\rm an}$ is Cartier on $X^{\rm an}$ if and only if $D$ is Cartier on $X$.
\end{lem}

\begin{comment}
\begin{lem}[{\cite[Lemma~5.1]{Kaw88}}]\label{lem: local index is local cartier index}
Let $x\in X$ be a terminal threefold singularity and $I$ a positive integer such that $IK_X$ is Cartier near $x$. Then $ID$ is Cartier near $x$ for any $\Qq$-Cartier Weil divisor $D$ on $X$.
\end{lem}
\end{comment}

\begin{thm}[cf. {\cite[Theorem 1.1]{Kaw15a}}]\label{thm: canonical index}
Let $X$ be a canonical threefold and $x\in X$ a (not necessarily closed) point such that $\mld(X\ni x)=1$. Then $IK_X$ is Cartier near $x$ for some positive integer $I\leq 6$.
\end{thm}
\begin{proof}
If $\dim x=2$ then $K_X$ is Cartier near $x$. If $\dim x=0$, then the theorem follows from \cite[Theorem 1.1]{Kaw15a}. If $\dim x=1$, then we let $f: Y\rightarrow X$ be the terminalization of $X\ni x$. By \cite[Theorem 4.5]{KM98}, $Y$ is smooth over a neighborhood of $x$. Since $K_Y=f^*K_X$, $K_X$ is Cartier near $x$ by the cone theorem.
\end{proof}

\begin{lem}\label{lem: canonical threshold attain mld=1}
Let $(X\ni x, B)$ be a threefold germ such that $\mld(X\ni x, B)\geq 1$. Let $D\ge0$ be an $\Rr$-Cartier $\Rr$-divisor on $X$ and $t:=\ct(X\ni x, B;D)$. Then $\mld(X\ni x, B+tD)=1$ if one of the following holds:
\begin{enumerate}
    \item $\mult_S (B+tD)<1$ for any prime divisor $S\subset \Supp D$,
    \item $X$ is $\Qq$-factorial terminal near $x$ and $x\in X$ is not smooth, and
    \item $X$ is $\Qq$-Gorenstein, $x\in X$ is not smooth, and $D$ is a $\Qq$-Cartier prime divisor.
\end{enumerate}
\end{lem}

\begin{proof}
If $\mld(X\ni x, B+tD)>1$, then $t=\lct(X\ni x, B;D)$. For (1), since $\mult_S (B+tD)<1$ for any prime divisor $S\subset \Supp D$, there exists a curve $C$ passing through $x$, such that $\mld(X\ni \eta_C, B+tD)=0$, where $\eta_C$ is the generic point of $C$. By \cite[Theorem~0.1]{Amb99}, $\mld(X\ni x, B+tD)\leq 1+\mld(X\ni \eta_C, B+tD)=1$, a contradiction.

For (2) and (3), by (1), we may assume that there exists a $\Qq$-Cartier prime divisor $S\subset \Supp D$ such that $\mult_S (B+tD)=1$. By \cite[Appendix, Theorem]{Sho92} and \cite[Theorem 0.1]{Mar96}, there exists a divisor $E$ over $X\ni x$ such that $a(E,X,0)=1+\frac{1}{I}$, where $I$ is the index of $x\in X$. Since $\mult_E (B+tD)\geq \mult_E S\geq \frac{1}{I}$, $a(E,X,B+tD)=a(E,X,0)-\mult_E(B+tD)\leq 1$, hence $\mld(X\ni x, B+tD)=1$.
\end{proof}
% (see Definition~\ref{def: definition of index})

\begin{thm}\label{thm: 3-dim terminal number of coefficients local}
Let $(X\ni x,B:=\sum_{i=1}^m b_iB_i)$ be a threefold germ such that $\mld(X\ni x,B)\geq 1$, $X$ is terminal, and each $B_i\geq 0$ is a $\Qq$-Cartier Weil divisor. Then we have the following:
\begin{enumerate}
    \item If $X\ni x$ is smooth, then $\sum_{i=1}^m b_i\leq 2$.
    \item If $X\ni x$ is not smooth, then $\sum_{i=1}^m b_i\le 1$. Moreover, if $\sum_{i=1}^m b_i=1$, then $\mld(X\ni x, B)=1$.
\end{enumerate}
\end{thm}
\begin{proof}
If $X\ni x$ is smooth, then let $E$ be the exceptional divisor of the blowing-up of $X$ at $x$. Since $X\ni x$ is smooth, $\mult_EB_i\geq 1$ for each $i$. Thus
$$1\leq \mld(X\ni x,B)\leq a(E,X,B)=3-\mult_EB=3-\sum_{i=1}^mb_i\mult_EB_i\leq 3-\sum_{i=1}^mb_i,$$
and we get (1).

If $X\ni x$ is not smooth, then let $I$ be the index of $X\ni x$. By \cite[Appendix, Theorem]{Sho92} and \cite[Theorem~0.1]{Mar96}, there exists a prime divisor $E$ over $X\ni x$ such that $a(E,X,0)=1+\frac{1}{I}$. Moreover, by \cite[Lemma~5.1]{Kaw88}, $IB_i$ is Cartier near $x$, and $I\mult_EB_i\geq 1$ for each $i$. Thus
\begin{align*}
    1&\leq\mld(X\ni x,B)\leq a(E,X,B)=a(E,X,0)-\mult_EB\\
    &=1+\frac{1}{I}-\sum_{i=1}^mb_i\mult_EB_i\leq 1+\frac{1}{I}-\frac{1}{I}\sum_{i=1}^mb_i,
\end{align*}
which implies (2).
\end{proof}

\subsection{Singular Riemann-Roch formula and Reid basket}\label{section: RR-formula}

\begin{defn}
Let $X\ni x$ be a smooth germ such that $\dim X=d$, and $\mathfrak{m}_x$ (resp. $\mathfrak{m}_x^{\rm an}$) the maximal ideal of the local ring $\mathcal{O}_{X,x}$ (resp. analytic local ring $\mathcal{O}_{X,x}^{\rm an}$). We say that $x_1,\dots,x_d\in \mathfrak{m}_x$ (resp. $x_1,\dots,x_d\in\mathfrak{m}_x^{\rm an}$) is a local coordinate system (resp. analytic local coordinate system) of $x\in X$ if $x_1,\dots,x_d$ span the linear space $\mathfrak{m}_x/\mathfrak{m}_x^2$ (resp. $\mathfrak{m}_x^{\rm an}/(\mathfrak{m}_x^{\rm an})^2$). We also call $x_1,\dots,x_d$ local coordinates (resp. analytic local coordinates) of $x\in X$.
\end{defn}

\begin{defn}[Cyclic quotient singularities]
Let $d$ and $n$ be two positive integers, and $a_1,\dots,a_d$ integers. A \emph{cyclic quotient singularity of type} $\frac{1}{n}(a_1,\dots,a_d)$ is the cyclic quotient singularity $(o\in\mathbb C^d)/\bm{\mu}$ given by the action
$$\bm{\mu}:(x_1,\dots,x_d)\rightarrow (\xi_n^{a_1}x_1,\xi_n^{a_2}x_2,\dots,\xi_n^{a_d}x_d)$$
on $\mathbb C^d$, where $x_1,\dots,x_d$ are the local coordinates of $\mathbb C^d\ni o$. We may also use $(o\in\mathbb C^d)/\frac{1}{n}(a_1,\dots,a_d)$ to represent the singularity $(o\in\mathbb C^d)/\bm{\mu}$, and use $(\mathbb C^d\ni o)/\frac{1}{n}(a_1,\dots,a_d)$ to represent the germ $(\mathbb C^d\ni o)/\bm{\mu}$.
\end{defn}

By the terminal lemma (cf. {\cite[Corollary 1.4]{MS84}}), if a cyclic quotient threefold singularity $x\in X$ is terminal, then $$(x\in X)\cong (o\in \mathbb C^3)/\frac{1}{n}(1,-1,b)$$ for some positive integers $b,n$ such that $\gcd(b,n)=1$. We say that the terminal (cyclic quotient threefold) singularity $x\in X$ is of type $\frac{1}{n}(1,-1,b)$ in this case.

\begin{defn}[{\cite[Theorem~10.2(2)]{Rei87}}]\label{defn: generalized sum}
For any integers $1\leq u<v$ and real numbers $s_u,\dots,s_v$, we define
$\sum_{i=v}^u s_i:=-\sum_{i=u+1}^{v-1}s_i$
when $v\geq u+2$, and define $\sum_{i=v}^us_i:=0$ when $v=u+1$.

Let $n$ be a positive integer and $m$ a real number. We define
$$\overline{(m)}_n:=m-\lfloor\frac{m}{n}\rfloor n.$$

Let $b,n$ be two positive integers such that $\gcd(b,n)=1$. Let $x\in X$ be a terminal cyclic quotient singularity of type $\frac{1}{n}(1,-1,b)$ and $D$ a Weil divisor on $X$, such that $\mathcal{O}_{X}(D)\cong\mathcal{O}_X(iK_X)$ for some integer $i$ near $x$. We define
$$c_x(D):=-i\frac{n^2-1}{12n}+\sum_{j=1}^{i-1}\frac{\overline{(jb)}_n(n-\overline{(jb)}_n)}{2n}.$$
We remark that $c_x(D)$ is independent of the choices of $i$ and $b$ by construction.
\end{defn}

\begin{comment}
\begin{rem}\label{rem: generalized summation form}
We introduce a generalized summation here, for any integers $a,b$, we define $$\sum_{i=a}^b :=\sum_{i=a}^{\max\{a,b\}}-\sum_{i=b+1}^{\max\{a,b\}}.$$ Under this notation, $$c_Q(D):=-i_Q'\frac{r_Q^2-1}{12r_Q}+\sum_{j=1}^{i_Q'-1}\frac{\overline{jb_Q}(r_Q-\overline{jb_Q})}{2r_Q}$$ for any integer $i_Q'$ such that $\overline{i_Q'}=i_Q$.
\end{rem}
\end{comment}

\begin{deflem}[{\cite[(6.4)]{Rei87}}]\label{def: definition of general cP}
Let $x\in X$ be a terminal threefold singularity. By the classification of threefold terminal singularities (cf. \cite[(6.1) Theorem]{Rei87}, \cite[Theorems 12,23,25]{Mor85}), we have an analytic isomorphism $(x\in X)\cong (y\in Y)/\bm{\mu}$ for some isolated cDV singularity $(y\in Y)\subset(o\in\mathbb C^4)$ and cyclic group action $\bm{\mu}$ on $o\in\mathbb C^4$. Moreover, $y\in Y$ is defined by an equation $(f=0)\subset (o\in\mathbb C^4)$ with analytic local coordinates $x_1,x_2,x_3,x_4$, and there exists $1\leq i\leq 4$ such that $f/x_i$ is a rational function that is invariant under the $\bm\mu$-action. Now we consider the $1$-parameter deformation ${Y_\lambda}$ of $Y$, such that $Y_{\lambda}$ is given by $(f+\lambda x_i=0)$. Then the deformation ${Y_\lambda}$ is compatible with the action $\bm{\mu}$, and we let $X_{\lambda}:=Y_{\lambda}/\bm{\mu}$ for each $\lambda$. For a general $\lambda\in\mathbb C$, the singularities of $X_{\lambda}$ are terminal cyclic quotient singularities $Q_1,\dots,Q_m$ for some positive integer $m$. We have the following.
\begin{enumerate}
    \item $Q_1,\dots,Q_m$ only rely on $x\in X$ and are independent of the choice of $\lambda$, and we define the set of \emph{fictitious singularities} of $x\in X$ to be $I_x:=\{Q_1,\dots,Q_m\}$.
    \item For any $\Qq$-Cartier Weil divisor $D$ on $X$, $D$ is deformed to a Weil divisor $D_{\lambda}$ on $X_{\lambda}$. We define $c_x(D):=\sum_{j=1}^mc_{Q_j}(D_{\lambda})$.
\end{enumerate}
For such a general $\lambda$, $X_{\lambda}$ is called a \emph{$Q$-smoothing} of $x\in X$.

\end{deflem}

The following theorem indicates that $c_x(D)$ is well-defined.
 
 \begin{thm}[{\cite[Theorem~10.2(1)]{Rei87}}]\label{thm: cx(D) only analytic type}
 Let $x\in X$ be a terminal threefold singularity and $D$ a $\Qq$-Cartier Weil divisor on $X$. Then $c_x(D)$ depends only on the analytic type of $x\in X$ and $\mathcal{O}_X(D)$ near $x$. Moreover, if $x\in X$ is smooth, then $c_x(D)=0$.
\end{thm}

\begin{thm}[{\cite[Theorem~10.2]{Rei87}}]\label{thm: singular rr}
Let $X$ be a projective terminal threefold, and $D$ a $\Qq$-Cartier Weil divisor on $X$. Then $$\chi(\mathcal{O}_X(D))=\chi(\mathcal{O}_X)+\frac{1}{12}D(D-K_X)(2D-K_X)+\frac{1}{12}D\cdot c_2(X)+\sum_{x\text{ is a closed point}} c_x(D).$$ 
\end{thm}

\begin{defn}[Reid basket for divisorial contractions]\label{def: Reid basket for divisorial contractions}
Let $f: Y\to X$ be a divisorial contraction of a prime divisor $F$ such that $Y$ is a terminal threefold.

For any closed point $y\in F$, consider a $Q$-smoothing of $y\in F\subset Y$ as in Definition-Lemma~\ref{def: definition of general cP}, and let $I_y$ be the corresponding set of fictitious singularities. For each $Q_y\in I_y$, let $Y_{Q_y}$ be the deformed variety on which $Q_y\in Y_{Q_y}$ is a cyclic quotient terminal threefold singularity of type $\frac{1}{r_{Q_y}}(1,-1,b_{Q_y})$, and $F_{Q_y}\subset Y_{Q_y}$ the deformed divisor of $F\subset Y$. Let $f_{Q_y}$ be the smallest non-negative integer such that $F_{Q_y}\sim f_{Q_y}K_{Y_{Q_y}}$ near $Q_y$. Possibly replacing $b_{Q_y}$ with $r_{Q_y}-b_{Q_y}$, we may assume that $v_{Q_y}:=\overline{(f_{Q_y}b_{Q_y})}_{r_{Q_y}}\leq \frac{r_{Q_y}}{2}$. The \emph{Reid basket} for the divisorial contraction $f:Y\to X$ with the exceptional divisor $F$ is defined as $$J:=\{(r_{Q_y},v_{Q_y})\mid y\in F, Q_y\in I_y, v_{Q_y}\neq 0\}.$$
\end{defn}

\subsection{Weighted blow-ups over quotient of complete intersection singularities}

%We adopt the notation and definitions as in \cite[Section 2]{Che19}.

\begin{defn}
A \emph{weight} is a vector $w\in\mathbb Q_{>0}^d$ for some positive integer $d$.%\luo{I think it should be a function instead of a vector, anyway}
\end{defn}

\begin{defn}[Weights of monomials and polynomials]
Let $d$ be a positive integer and $w=(w_1,\dots,w_d)\in \mathbb{Q}^d_{>0}$ a weight. For any vector $\bm{\alpha}=(\alpha_1,\dots, \alpha_d)\in \Zz_{\geq 0}^d$, we define $\bm{x}^{\bm{\alpha}} :=x_1^{\alpha_1}\dots x_d^{\alpha_d}$, and $$w(\bm{x}^{\bm{\alpha}}):=\sum_{i=1}^dw_i\alpha_i$$
to be \emph{the weight of $\bm{x}^{\bm{\alpha}}$ with respect to $w$}.
For any analytic power series $0\neq h:=\sum_{\bm{\alpha}\in \Zz^d_{\geq 0}} a_{\bm{\alpha}}\bm{x}^{\bm{\alpha}}$, we define $$w(h):=\min\{w(\bm{x}^{\bm{\alpha}})\mid a_{\bm{\alpha}}\neq 0\}$$
to be \emph{the weight of $h$ with respect to $w$}. If $h=0$, then we define $w(h):=+\infty$.
\end{defn}

\begin{defn}
Let $h\in \mathbb{C}\{x_1,\dots,x_d\}$ be an analytic power series and $G$ a group which acts on $\Cc\{x_1,\dots,x_d\}$. We say that $h$ is \emph{semi-invariant} with respect to the group action $G$ if for any $g\in G$, $\frac{g(h)}{h}\in \mathbb{C}$. If the group action is clear from the context, then we simply say that $h$ is semi-invariant.
\end{defn}

\begin{defn}\label{def: Weights of threefolds terminal singularities}
Let $(X\ni x,B:=\sum_{i=1}^k b_iB_i)$ be a germ such that $X$ is terminal and $B_i\geq 0$ are $\Qq$-Cartier Weil divisors on $X$. Let $d,n$ and $m<d$ be positive integers such that
$$(X\ni x)\cong (\phi_1=\cdots =\phi_m=0)\subset (\mathbb C^d\ni o)/\frac{1}{n}(a_1,\dots,a_d)$$ for some semi-invariant analytic power series $\phi_1\dots,\phi_m\in \mathbb{C}\{x_1,\dots,x_d\}$. By \cite[Lemma~5.1]{Kaw88}, $B_i$ can be identified with $\big((h_i=0)\subset \mathbb (\mathbb{C}^d\ni o)/\frac{1}{n}(a_1,\dots,a_d)\big)|_X$ for some semi-invariant analytic power series $h_i\in \mathbb{C}\{x_1,\dots,x_d\}$ near $x\in X$, and we say that $B_i$ is defined by $(h_i=0)$ near $x$ or $B_i$ is locally defined by $(h_i=0)$ for simplicity. We define the set of \emph{admissible weights} of $X\ni x$ to be $$\{\frac{1}{n}(w_1,\dots,w_d)\in\frac{1}{n}\mathbb Z^d_{>0}\mid \text{ there exists $b\in\mathbb{Z}$ such that $w_i\equiv b a_i~\mathrm{mod}~n$, $1\leq i\leq d$}\}.$$ For any admissible weight $w=\frac{1}{n}(w_1,\dots,w_d)$, 
we define 
$$w(X\ni x):=\frac{1}{n}\sum_{i=1}^d w_i-\sum_{i=1}^m w(\phi_i)-1, \text{ and } w(B):=\sum_{i=1}^k b_iw(h_i).$$
By construction, $w(B)$ is independent of the choices of $b_i$ and $B_i$.
\end{defn}

\begin{comment}
%\han{for submitted version we may delete the proof}
\begin{lem}\label{lem: well-definess of w(B)} Let $d_1,\dots,d_m,d_1',\dots,d_{m'}'$ be real numbers and $D_1,\dots,D_m,D_1',\dots,D_{m'}'$ $\Qq$-Cartier Weil divisors such that $$\sum_{i=1}^m d_iD_i=\sum_{i=1}^{m'} d_i'D_i'.$$
Then $\sum_{i=1}^m d_iw(D_i)=\sum_{i=1}^{m'} d_i'w(D_i')$.
\end{lem}

\begin{proof}
Let $r_1,\dots,r_n\in \Rr$ be a basis for the $\Qq$-linear space spanned by the real numbers $\{d_1,\dots,d_m,d_1',\dots,d_{m'}'\}$, we may write
\begin{align*}
    \sum_{i=1}^m d_iD_i=\sum_{j=1}^nr_j\sum_{i=1}^{m}d_{i,j}D_i,~and~
    \sum_{i=1}^{m'} d_i'D_i'&=\sum_{j=1}^nr_j\sum_{i=1}^{m'}d_{i,j}'D_i'
\end{align*}
for some rational numbers $\{d_{i,j}\}_{1\leq i\leq m,1\leq j\leq n}$ and $\{d_{i,j}'\}_{1\leq i\leq m',1\leq j\leq n}$. For each prime divisor $D$, we have $$\mult_D (\sum_{i=1}^m d_iD_i)=\sum_{j=1}^nr_j\mult_D(\sum_{i=1}^{m}d_{i,j}D_i)=\sum_{j=1}^nr_j\mult_D(\sum_{i=1}^{m'}d_{i,j}'D_i')=\mult_D (\sum_{i=1}^{m'} d_i'D_i'),$$ hence $\mult_D(\sum_{i=1}^{m}d_{i,j}D_i)=\mult_D(\sum_{i=1}^{m'}d_{i,j}'D_i')$ for each $j$, this implies that for each $j$, we have $$\sum_{i=1}^{m}d_{i,j}D_i=\sum_{i=1}^{m'}d_{i,j}'D_i'.$$ Hence 
\begin{align*}
    \sum_{i=1}^m d_iw(D_i)&=\sum_{j=1}^nr_j\sum_{i=1}^{m}d_{i,j}w(D_i)=\sum_{j=1}^nr_j\sum_{i=1}^{m'}d_{i,j}'w(D_i')=\sum_{i=1}^{m'} d_i'w(D_i').
\end{align*}
\end{proof}
\end{comment}

\begin{defn}\label{defn: compare two weights}
Let $d$ be a positive integer, $\mu$ a real number, and $w:=(w_1,\dots,w_d)\in\mathbb Q^d_{>0}$, $w':=(w'_1,\dots,w'_d)\in\mathbb Q^d_{>0}$ two weights. If $w_i\geq w_i'$ for each $i$, then we write $w\succeq w'$, and if $w_i=\mu w_i'$ for each $i$, then we write $w=\mu w'$.
\end{defn}

%In the following, we construct the weighted blow-ups for a subvariety of $\Cc^d/\frac{1}{n}(a_1,\dots,a_d)$ (cf. \cite[\S 3]{Hay99}).

\begin{deflem}\label{defn: exceptional divisor and weighted blow-up}
Let $f': W\to (\Cc^d\ni o)/\frac{1}{n}(a_1,\dots,a_d)$ be the weighted blow-up at $o$ with the (admissible) weight $w:=\frac{1}{n}(w_1,\dots,w_d)$ with respect to the coordinates $x_1,\dots,x_d$ (cf. \cite[\S 10]{KM92} and \cite[\S 3.2]{Hay99}). The exceptional locus for $f'$, denoted by $E'$, is isomorphic to the cyclic quotient of the weighted projective space $\mathbf{P}(w_1,\dots,w_d)/\bm \eta$, where the cyclic group action is given by $$\bm \eta: [x_1:\dots:x_d]_w\to [\xi_n^{a_1}x_1:\dots,\xi_n^{a_d}x_d]_w,$$ and $[x_1:\dots:x_d]_w$ denotes the image of $(x_1,\dots,x_d)\in \Cc^d\setminus \{o\}$ under the natural quotient morphism $ \Cc^d\setminus\{o\}\to \mathbf{P}(w_1,\dots,w_d)$. We remark that if the admissible weight $w$ satisfies $w_i\equiv b a_i\mod n$ for $1\leq i\leq d$ and some integer $b$ such that $\gcd(b,n)=1$, then $E'\cong \mathbf{P}(w_1,\dots,w_d)$ (cf. \cite[\S 3.2]{Hay99}).

Now we have an induced morphism $f: Y\to X$ by restricting $f'$ to $Y$, which is the strict transform of $X$ under $f'$. We call $f: Y\to X$ the \emph{weighted blow-up with weight $w$ at $x\in X$}, and $E:=E'|_{Y}$ the \emph{exceptional divisor of the weighted blow-up $f: Y\to X$ with the weight $w$ at $x\in X$} (cf. \cite[\S 3.7]{Hay99}). If $E$ is an integral scheme, then we also say $f$ \emph{extracts a prime divisor}.
\end{deflem}

\begin{rem}
The weighted blow-up constructed as above depends on the choice of local coordinates. However, we will not mention them later in this paper when the local coordinates for a weighted blow-up are clear from the context.
\end{rem}

We will use the following well-known lemma frequently: %\han{maybe prove it in full version?}
\begin{lem}[cf. {\cite[the proof of Theorem 2]{Mor85} and \cite[\S 3.9]{Hay99}}]\label{lem: weighted blowup log discrepancies}
Let $d,n$ and $m<d$ be positive integers and $(X\ni x,B)$ a $\Qq$-Gorenstein germ. Under suitable analytic local coordinates $x_1,\dots,x_d$,
$$(X\ni x)\cong (\phi_1=\cdots=\phi_m=0)\subset (\mathbb C^d\ni o)/\frac{1}{n}(a_1,\dots,a_d)$$ for some semi-invariant analytic power series $\phi_1,\dots, \phi_m$. For any admissible weight $w$ of $X\ni x$, let $E$ be the exceptional divisor of the corresponding weighted blow-up $f: Y\rightarrow X$ at $x$ (cf. Definition-Lemma~\ref{defn: exceptional divisor and weighted blow-up}). Then 
$$K_Y=f^*K_X+w(X\ni x)E, \text{ and }~f^*B=B_Y+w(B)E,$$
where $B_Y$ is the strict transform of $B$ on $Y$. In particular, if $E$ is a prime divisor, then $a(E,X,B)=1+w(X\ni x)-w(B)$.
\end{lem}

\begin{defn}[Newton polytopes]
Let $d$ be a positive integer. For any analytic power series $h=\sum_{\bm{\alpha}\in\Zz_{\geq 0}^d} a_{\bm{\alpha}}x^{\bm{\alpha}}\in\mathbb C\{x_1,\dots,x_d\}$, the \emph{Newton polytope} $\mathcal{N}(h)$ of $h$ is a Newton polytope in $\mathbb Z^d_{\ge 0}$ (cf. Definition \ref{defn: newton polytope zd}) defined by 
$$\mathcal{N}(h):=\{\bm{\alpha}+\mathbb Z^d_{\ge 0}\mid \bm{\alpha}\in\mathbb Z^d_{\ge 0},~a_{\bm{\alpha}}\not=0\}.$$
\end{defn}

\subsection{Divisorial contractions between terminal threefold singularities} Let $f: Y\to X$ be a divisorial contraction of a prime divisor $E$ between two terminal threefolds such that $f(E)$ is a closed point on $X$. Then $f$ is classified into two types: the \emph{ordinary type} and the \emph{exceptional type}. Moreover, in the ordinary type case, any non-Gorenstein singularity
on $Y$ which contributes to the Reid basket of $f$ is a cyclic quotient terminal singularity (see the paragraph after \cite[Theorem~1.1]{Kaw05}).

\begin{defn}\label{defn: powerseries contains a monomial term notation}
Let $h\in\Cc\{x_1,\dots,x_d\}$ be an analytic power series. Let $a_{\bm{\alpha}}\bm{x}^{\bm{\alpha}}$ be a monomial for some $a_{\bm{\alpha}}\in \Cc$ and $\bm{\alpha}\in \Zz^d_{\geq 0}$. By $a_{\bm{\alpha}}\bm{x}^{\bm{\alpha}}\in h$, we mean the monomial term $\bm{x}^{\bm{\alpha}}$ appears in the analytic power series $h$ with the coefficient $a_{\bm{\alpha}}$.
\end{defn}

\begin{thm}[{\cite[Theorem~1.2]{Kaw05}}]\label{thm: kaw05 1.2 strenghthened}
In the statement of this theorem, $d,r,r_1,r_2,\alpha$ are assumed to be positive integers, $\lambda,\mu$ are assumed to be complex numbers, and $g,p,q$ are assumed to be polynomials with no non-zero constant terms.

Let $x\in X$ be a terminal singularity that is not smooth, $f: Y\rightarrow X$ a divisorial contraction of a prime divisor $E$ over $X\ni x$ (see Definition~\ref{defn: divisorial contraction}), such that $Y$ is terminal over a neighborhood of $x$. Let $n$ be the index of $X\ni x$. We may write
$$K_Y=f^*K_X+\frac{a}{n}E$$
for some positive integer $a$. If $f: Y\to X$ is of ordinary type, then one of the following holds:
\begin{enumerate}

    \item $x\in X$ is a $cA/n$ type singularity. Moreover, under suitable analytic local coordinates $x_1,x_2,x_3,x_4$,
    \begin{enumerate}
    \item we have an analytic identification
    $$(X\ni x)\cong (\phi:= x_1x_2+g(x_3^n,x_4)=0)\subset (\mathbb C^4\ni o)/\frac{1}{n}(1,-1,b,0),$$
    \item $f$ is a weighted blow-up with the weight $w=\frac{1}{n}(r_1,r_2,a,n),$
        \item $x_3^{dn}\in g(x_3^n,x_4)$,
        \item $nw(\phi)=r_1+r_2=adn$, and
        \item $a\equiv br_1\mod n$.
    \end{enumerate}
    \item $x\in X$ is a $cD$ type singularity. In this case, one of the following holds:
    \begin{itemize}
        \item[(2.1)] Under suitable analytic local coordinates $x_1,x_2,x_3,x_4$,
        \begin{enumerate}
        \item  we have an analytic identification $$(X\ni x)\cong (\phi:= x_1^2+x_1q(x_3,x_4)+x_2^2x_4+\lambda x_2x_3^2+\mu x_3^3+p(x_2,x_3,x_4)=0)\subset (\mathbb C^4\ni o),$$ where $p(x_2,x_3,x_4)\in(x_2,x_3,x_4)^4$,
            \item $f$ is a weighted blow-up with the weight $w=(r+1,r,a,1),$ where $a$ is an odd integer,
            \item $\mu' x_3^d\in \phi$ for some $\mu'\neq 0$ and an odd integer $d\geq 3$, and if $d=3$, then $\mu'=\mu$,
            \item $w(\phi)=w(x_2^2x_4)=w(x_3^d)=2r+1=ad$,
            \item if $q(x_3,x_4)\not=0$, then $w(x_1q(x_3,x_4))=2r+1$, and
            \item if $d>3$, then $\mu=\lambda=0$.
        \end{enumerate}
        \item[(2.2)] Under suitable analytic local coordinates $x_1,x_2,x_3,x_4,x_5$, 
        \begin{enumerate}
        \item  we have an analytic identification  $$(X\ni x)\cong \binom{\phi_1:= x_1^2+x_2x_5+p(x_2,x_3,x_4)=0}{\phi_2:= x_2x_4+x_3^d+q(x_3,x_4)x_4+x_5=0}\subset(\mathbb C^5\ni o),$$where $p(x_2,x_3,x_4)\in(x_2,x_3,x_4)^4$,
        \item $f$ is a weighted blow-up with the weight
        $w=(r+1,r,a,1,r+2),$
            \item $r+1=ad$ and $d\geq 2$,
            \item $w(\phi_1)=2(r+1)$, 
            \item $w(\phi_2)=r+1$, and
            \item if $q(x_3,x_4)\not=0$, then $w(q(x_3,x_4)x_4)=r+1$.
        \end{enumerate}
    \end{itemize}
    \item $x\in X$ is a $cD/2$ type singularity. In this case, one of the following holds:
    \begin{itemize}
        \item[(3.1)]  Under suitable analytic local coordinates $x_1,x_2,x_3,x_4$, 
        \begin{enumerate}
        \item  we have an analytic identification  $$(X\ni x)\cong (\phi:= x_1^2+x_1x_3q(x_3^2,x_4)+x_2^2x_4+\lambda x_2x_3^{2\alpha-1}+p(x_3^2,x_4)=0)
            \subset (\mathbb C^4\ni o)/\frac{1}{2}(1,1,1,0),$$
            \item $f$ is a weighted blow-up with the weight $w=\frac{1}{2}(r+2,r,a,2),$
            \item $w(\phi)=w(x_2^2x_4)=r+1=ad$, where $a,r$ are odd integers,
            \item if $q(x^2_3,x_4)\not=0$, then $w(x_1x_3q(x_3^2,x_4))=r+1$, and
            \item $x_3^{2d}\in p(x_3^2,x_4)$.
        \end{enumerate}
        \item[(3.2)] Under suitable analytic local coordinates $x_1,x_2,x_3,x_4,x_5$,
        \begin{enumerate}
        \item we have an analytic identification $$(X\ni x)\cong \binom{\phi_1:= x_1^2+x_2x_5+p(x_3^2,x_4)=0}{\phi_2:= x_2x_4+x_3^d+q(x_3^2,x_4)x_3x_4+x_5=0}\subset(\mathbb C^5\ni o)/\frac{1}{2}(1,1,1,0,1),$$
        \item $f$ is a weighted blow-up with the weight
        $w=\frac{1}{2}(r+2,r,a,2,r+4),$
            \item $r+2=ad$ and $d$ is an odd integer,
            \item $w(\phi_1)=r+2$,
            \item $w(\phi_2)=\frac{r+2}{2}$, and
            \item if $q(x_3^2,x_4)\not=0$, then $w(q(x_3^2,x_4)x_3x_4)=\frac{r+2}{2}$.
        \end{enumerate}
    \end{itemize}
\end{enumerate}
Moreover, if $a\geq 5$, then $f$ is of ordinary type. The cases are summarized in Table \ref{tbl:MT}:

\begin{table}[ht]
\caption{A summary of Theorem \ref{thm: kaw05 1.2 strenghthened}}\label{tbl:MT}
\begin{tabular}{|c|c|c|c|c|}
\hline
\multicolumn{1}{|c|}{Case}  &\multicolumn{1}{c|}{Type} & Local coordinates &
\multicolumn{1}{c|}{$w$} & \multicolumn{1}{c|}{$w(\phi)$ or $w(\phi_i)$}\\
\hline
(1) & $cA/n$  & $(\phi=0)\subset\mathbb C^4/\frac{1}{n}(1,-1,b,0)$ & $\frac{1}{n}(r_1,r_2,a,n)$ & $\frac{r_1+r_2}{n}$\\
\hline
(2.1) & $cD$  & $(\phi=0)\subset\mathbb C^4$&
$(r+1,r,a,1)$ & $2r+1$\\
\hline
(2.2) & $cD$  & $(\phi_1=\phi_2=0)\subset\mathbb C^5$&
$(r+1,r,a,1,r+2)$ & $2(r+1)$ \text{and} $r+1$\\
\hline
(3.1) & $cD/2$  & $(\phi=0)\subset\mathbb C^4/\frac{1}{2}(1,1,1,0)$ &
$\frac{1}{2}(r+2,r,a,2)$ & $r+1$\\
\hline
(3.2) & $cD/2$  & $(\phi_1=\phi_2=0)\subset\mathbb C^4/\frac{1}{2}(1,1,1,0,1)$&
$\frac{1}{2}(r+1,r,a,2,r+4)$ & $r+2$ \text{and} $\frac{r+2}{2}$\\
\hline
\end{tabular}
\end{table}
\end{thm}

\begin{proof}
Most part of this theorem are identical to \cite[Theorem 1.2]{Kaw05} but with small differences for further applications. For the reader's convenience, we give a proof here.

Since $f: Y\to X$ is a divisorial contraction of ordinary type and $x\in X$ is not smooth, by \cite[Theorem 1.2]{Kaw05}, we have the following possible cases.

\medskip

$x\in X$ is of type $cA/n$, then we are in case (1). (1.a) and (1.b) follow from \cite[Theorem 1.2(i)]{Kaw05}. By \cite[Theorem 1.2(i.a)]{Kaw05}, (1.e) holds, and we may pick a positive integer $d$ such that $r_1+r_2=adn$. By \cite[Theorem 1.2(i.c)]{Kaw05} and \cite[Theorem 1.2(i.d)]{Kaw05}, (1.c) and (1.d) hold.

\medskip

$x\in X$ is of type $cD$ or $cD/2$, then we are in either case (2) or case (3). Now (2.1.a), (2.1.b), (3.1.a), (3.1.b) follow directly from \cite[Theorems 1.2(ii.a) and 1.2(ii.a.1)]{Kaw05}, and (2.2.a), (2.2.b), (3.2.a), (3.2.b) follow directly from \cite[Theorem~1.2(ii.b)]{Kaw05}. (2.1.d), (3.1.c) follow from \cite[Theorems 1.2(ii.a.1) and 1.2(ii.a.2)]{Kaw05}, and (2.2.c), (3.2.c) follow from \cite[Theorem~1.3(ii.b.1)]{Kaw05}. (2.1.e), (3.1.d) follow from \cite[Theorem~1.2(ii.a.2)]{Kaw05}, and (2.2.d), (2.2.e), (2.2.f), (3.2.d), (3.2.e), (3.2.f) follow from \cite[Theorem~1.2(ii.b.2)]{Kaw05}. For (2.1.f), if $d>3$, by (2.1.d) we have $\mu=0$. Assume that $\lambda\neq 0$. By (2.1.d), $w(x_2x_3^2)=r+2a\geq 2r+1$, hence $2a\geq r+1$ and $2r+1=ad>\frac{1}{4}(2r+1)d$. It follows that $d<4$, a contradiction.

\medskip

(2.1.c), (3.1.e) are not contained in the statement of \cite[Theorem~1.2]{Kaw05}, however, they are implied by the proofs of the corresponding results. To be more specific, (2.1.c) is stated in \cite[\S 4, Case Ic]{Che15} and (3.1.e) is stated in \cite[Page 13, Line 10]{CH11}.

\medskip

For the moreover part, the Theorem follows directly from \cite[Theorem 1.3]{Kaw05}.
\end{proof}

\subsection{Terminal blow-ups}

\begin{defn}\label{defn: canonical place}
Let $(X\ni x,B)$ be an lc pair. We say that $x$ is a \emph{canonical center} of $(X,B)$ if $\mld(X\ni x,B)=1$ and $\dim x\leq\dim X-2$. A prime divisor (resp. An analytic prime divisor) $E$ that is exceptional over $X$ is called a \emph{canonical place} of $(X,B)$ if $a(E,X,B)=1$. Moreover, if $\Center_X E=\bar{x}$, then $E$ is called a \emph{canonical place of $(X\ni x,B)$}.
\end{defn}

\begin{lem}\label{lem: tie breaking canonical pair}
Let $(X\ni x,B)$ be a germ such that $X$ is terminal and $\mld(X\ni x, B)=1$. Then there exists a pair $(X,B')$ that is klt near $x$, such that
\begin{itemize}
    \item $x$ is the only canonical center of $(X,B')$,
    \item there exists exactly one canonical place $E$ of $(X\ni x,B')$, and
    \item $a(E,X,B)=1$.
\end{itemize}
\end{lem}

%\han{submiited version, just say by tie-breaking, cf. cot04}

\begin{proof}
Possibly shrinking $(X,B)$ to a neighborhood of $x$, we may assume that $(X,B)$ is an lc pair. We play the so-called ``tie-breaking trick'' and follow the proof of \cite[Proposition 8.7.1]{Kol07}.

\medskip

\noindent\textbf{Step 1}. Let $f:W\to X$ be a log resolution of $(X,B)$. We may write 
$$K_W=f^{*}K_X+\sum_{i\in \mathfrak{I}} a_iE_i, \text{ and } f^{*}B=B_W+\sum_{i\in \mathfrak{I}} b_iE_i,$$
where $B_W:=f_*^{-1}B$ is smooth, and $\{E_i\}_{i\in \mathfrak{I}}$ is the set of $f$-exceptional divisors. Let 
$$\mathfrak{I}_x:=\{i\in \mathfrak{I}\mid \Center_X E_i=x\}, \text{ and } \mathfrak{I}_{x,0}:=\{i\in \mathfrak{I}_x\mid a_i=b_i\}\neq\emptyset.$$ Then $a_i>0$ for any $i\in \mathfrak{I}$, and $b_i>0$, $a_i\ge b_i$ for each $i\in \mathfrak{I}_x$.

Let $C$ be a very ample Cartier divisor such that $x\in\Supp C$. Possibly replacing $C$ with $C'\in H^0(\mathcal{O}_X(kC)\otimes \mathfrak{m}_x)$ for some $k\gg 1$ and some irreducible $C'$ which is sufficiently general, we may assume that $C$ is a prime divisor whose support does not contain any $\Center_XE_i$ $i\in \mathfrak{I}\setminus \mathfrak{I}_x$, or any canonical center of $(X,B)$ except $x$. We may write $f^{*}C=C_W+\sum_{i\in \mathfrak{I}_x}c_iE_i$, where $C_W$ is the strict transform of $C$ on $W$, and $c_i>0$ for each $i\in \mathfrak{I}_x$. Now we may choose a real number $0<\epsilon\ll 1$, such that
$$t:=\min_{i\in \mathfrak{I}_{x,0}}\frac{a_i-(1-\epsilon)b_i}{c_i}=\epsilon\min_{i\in \mathfrak{I}_{x,0}}\frac{ a_i}{c_i}<\min_{i\in 
\mathfrak{I}_x\setminus \mathfrak{I}_{x,0}}\frac{a_i-b_i}{c_i}<\min_{i\in \mathfrak{I}_x\setminus \mathfrak{I}_{x,0}}\frac{a_i-(1-\epsilon)b_i}{c_i},$$
and $(1-\epsilon)B+tC\in[0,1)$. Let $K_W+D_{\epsilon}:=f^*(K_X+(1-\epsilon)B+tC)$ and $K_W+D:=f^*(K_X+B)$. We have $$D_{\epsilon}-D=\epsilon(\min_{i\in \mathfrak{I}_{x,0}}\frac{ a_i}{c_i}f^*C-f^*B).$$ 
Consider the finite sets
\begin{align*}
    \mathfrak{J}_{1,\epsilon}&:=\{\codim y-\mult_y D_{\epsilon}\mid y\in W, f(y)\neq x,\codim f(y)\geq 2\},~\text{and}\\
    \mathfrak{J}_{2,\epsilon}&:=\{\codim y-\mult_y D_{\epsilon}\mid y\in W, f(y)=x\}.
\end{align*}
If $F$ is an $f$-exceptional divisor over $X\ni x'$ for some $x'\neq x$ such that $a(F,X,B)=1$, then $a(F,X,(1-\epsilon)B+tC)>1$. By taking $\epsilon$ sufficiently small, we may assume that $1\notin \mathfrak{J}_{1,\epsilon}$. Since $B_W$ is smooth and $\mult_{E_i}D_{\epsilon}\le 0$ for any $i\in \mathfrak{I}_x$, $1\in\mathfrak{J}_{2,\epsilon}\subset [1,+\infty)$. By \cite[Lemma~3.3]{CH21}, $x$ is the only canonical center of $(X,(1-\epsilon)B+tC)$. Note that $1\in \mathfrak{J}_{2,\epsilon}$.

Possibly replacing $B$ with $(1-\epsilon)B+tC$, we may assume that $x$ is the only canonical center of $(X,B)$ and $(X,B)$ is klt.

\medskip

\noindent\textbf{Step 2}. We may assume that $f$ is a composition of blow-ups of centers of codimension at least 2, hence there exists an $f$-ample $\Qq$-divisor $-\sum_{i\in \mathfrak{I}} e_iE_i$. Note that $e_i>0$ for each $i$. Then there exists $i_0\in \mathfrak{I}_{x,0}$, such that possibly replacing $e_{i_0}$ with a bigger positive rational number, we may assume that $\lambda:=\frac{a_{i_0}}{e_{i_0}}<\frac{a_i}{e_i}$ for any $i\in \mathfrak{I}_{x,0}\backslash \{i_0\}$. Moreover, there exists a positive real number $\epsilon'_0<1$, such that $\lambda':=\epsilon'_0\lambda<\frac{a_i-b_i}{e_i}<\frac{a_i-(1-\epsilon'_0)b_i}{e_i}$ for any $i\in \mathfrak{I}_x\backslash \mathfrak{I}_{x,0}$. Let $E:=E_{i_0}$. Let $H$ be an ample $\Qq$-Cartier $\Qq$-divisor on $X$ such that $f^*H-\sum_{i\in \mathfrak{I}} e_iE_i$ is ample. Now $a_i-(1-\epsilon'_0)b_i-\lambda' e_i>0$ for any $i\in \mathfrak{I}_x\setminus \{i_0\}$, and $a_{i_0}-(1-\epsilon'_0)b_{i_0}-\lambda' e_{i_0}=0$. We have
$$K_W+(1-\epsilon'_0)B_W+\lambda' (f^{*}H-\sum_{i\in \mathfrak{I}} e_iE_i)=f^{*}(K_X+(1-\epsilon'_0)B+\lambda' H)+\sum_{i\in \mathfrak{I}}(a_i-(1-\epsilon'_0)b_i-\lambda' e_i)E_i.$$
Let $A_W\ge0$ be a $\Qq$-divisor such that $A_W\sim_{\Qq} f^{*}H-\sum_{i\in \mathfrak{I}} e_iE_i$, the coefficients of the prime components of $A_W$ are sufficiently small, and $\Supp A_W\cup \Supp B_W\cup_{i\in \mathfrak{I}} \Supp E_i$ has simple normal crossings. Consider the pair $(X\ni x,B':=(1-\epsilon'_0)B+\lambda' A)$, where $A$ is the strict transform of $A_W$ on $X$. We may write $$K_W+D'_{\epsilon'_0}:=f^{*}(K_X+(1-\epsilon')B+\lambda' A).$$ 
For any $\epsilon'\ge 0$, let
\begin{align*}
    \mathfrak{J}'_{1,\epsilon'}&:=\{\codim y-\mult_y D'_{\epsilon'}\mid y\in W, f(y)\neq x,\codim f(y)\geq 2\},~\text{and}\\
    \mathfrak{J}'_{2,\epsilon',\ge2}&:=\{\codim y-\mult_y D'_{\epsilon'}\mid y\in W, f(y)=x,\codim y\geq 2\}.
\end{align*}
Since $1\notin \mathfrak{J}'_{1,0}$, by taking $\epsilon'$ small enough, we may assume that $1\notin \mathfrak{J}'_{1,\epsilon'}$. Thus $x$ is the only canonical center of $(X,B')$. By our choices of $\epsilon'_0$, $\lambda$, and $A_W$, $\mathcal{J}'_{2,\epsilon'_0,\ge2}\subset (1,+\infty)$. Hence by \cite[Lemma~3.3]{CH21}, $E$ is the only canonical place of $(X\ni x,B')$. The pair $(X\ni x, B')$ satisfies all the requirements.
\end{proof}

\begin{defn}\label{defn: terminal blow-up}
Let $(X\ni x,B)$ be an lc germ. A \emph{terminal blow-up} of $(X\ni x,B)$ is a birational morphism $f: Y\rightarrow X$ which extracts a prime divisor $E$ over $X\ni x$, such that $Y$ is terminal, $a(E,X,B)=\mld(X\ni x,B)$, and $-E$ is $f$-ample. 
\end{defn}

\begin{lem}\label{lem: can extract divisor computing ct that is terminal strong version}
Let $(X\ni x, B)$ be a germ such that $X$ is terminal and $\mld(X\ni x, B)=1$. Then there exists a terminal blow-up $f: Y\to X$ of $(X\ni x,B)$. Moreover, if $X$ is $\Qq$-factorial, then $Y$ is $\Qq$-factorial.
\end{lem}

\begin{proof}
By Lemma \ref{lem: tie breaking canonical pair}, possibly shrinking $X$ to a neighborhood of $x$,  we may assume that $(X,B)$ is klt, and there exists exactly one canonical place $E_W$ of $(X\ni x,B)$. By \cite[Corollary~1.4.3]{BCHM10}, there exists a birational morphism $g: W\rightarrow X$ of $(X,B)$ such that $E_W$ is the only $g$-exceptional divisor. We may write $K_W+B_W:=g^*(K_X+B)$, where $B_W$ is the strict transform of $B$ on $W$. Since $(W,B_W)$ is klt, $(W,(1+\epsilon)B_W)$ is klt for some positive real number $\epsilon$. Let $\phi:W\dashrightarrow Y$ be the lc model of $(W,(1+\epsilon)B_W)$ over $X$. Since $X$ is terminal and $E_W$ is a canonical place of $(X\ni x,B)$,
$$K_W+(1+\epsilon)B_W=g^*(K_X+(1+\epsilon)B)-eE_W$$
 for some positive real number $e$. It follows that $-E$ is ample over $X$, where $E$ is the strict transform of $E_W$ on $Y$. Thus $f:Y\to X$ is an isomorphism over $X\backslash\{x\}$, and $\Exc(f)=\Supp E$.
 
It suffices to show that $Y$ is terminal. Let $F$ be any prime divisor that is exceptional over $Y$. If $\Center_{X}F=x$, then $a(F,Y,0)\ge a(F,X,B)>1$ as $E$ is the only canonical place of $(X\ni x,B)$. If $\Center_{X}F\neq x$, then $a(F,Y,0)=a(F,X,0)>1$ as $f$ is an isomorphism over $X\backslash\{x\}$ and $X$ is terminal. 

Suppose that $X$ is $\Qq$-factorial. Then for any prime divisor $D_Y\neq E$ on $Y$, $f^{*}f_{*}D_Y-D_Y=(\mult_{E_Y}f_{*}D_Y)E$. It follows that $D_Y$ is $\Qq$-Cartier, and $Y$ is $\Qq$-factorial.
\end{proof}

\section{ACC for threefold canonical thresholds}\label{section: ACC for threefold canonical thresholds}

The main goal for this section is to show Theorem~\ref{thm: 3fold acc ct}.

\subsection{Weak uniform boundedness of divisors computing mlds}

\begin{lem}\label{lem: bounded index Mustata nakamura conjecture}
Let $I$ be a positive integer, $\alpha \geq 1$ a real number,  and $\Ii\subset [0,1]$ a DCC set. Then there exists a positive integer $l$ depending only on $I,\alpha$ and $\Ii$ satisfying the following. Assume that
\begin{enumerate}
    \item $(X\ni x,B:=\sum_{i} b_iB_i')$ is a threefold germ,
    \item $X$ is terminal,
    \item each $b_i\in\Ii$ and each $B_i'\geq 0$ is a $\Qq$-Cartier Weil divisor,
    \item $\mld(X\ni x,B)= \alpha$,
    \item $IK_X$ is Cartier near $x$, and
    \item there exists a terminal blow-up (see Definition~\ref{defn: terminal blow-up}) $f: Y\rightarrow X$ of $(X\ni x,B)$.
\end{enumerate}
Then there exists a prime divisor $E$ over $X\ni x$, such that $a(E,X,B)=\mld(X\ni x,B)$ and $a(E,X,0)=1+\frac{a}{I}$ for some positive integer $a\leq l$.
\end{lem}
\begin{proof}

Suppose that the lemma does not hold. Then there exists a sequence of threefold germs $\{(X_i\ni x_i, B_i:=\sum_{j=1}^{p_i} b_{i,j}B_{i,j}'))\}_{i=1}^{\infty}$ and terminal blow-ups $f_i: Y_i\rightarrow X_i$ corresponding to $(X\ni x,B:=\sum_{i} b_iB_i')$ and $f:Y\to X$ as in (1)-(6), such that
\begin{itemize}
\item $f_i$ extracts a prime divisor $E_i$,
\item $K_{Y_i}=f_i^*K_{X_i}+\frac{a_i}{I}E_i$ for some positive real number $a_i$, and
\item the following sequence of non-negative integers$$A_i:=\inf\{a_i'\mid F_i\text{ is over }X_i\ni x_i, a(F_i,X_i,B_i)=\alpha, a(F_i,X_i,0)=1+\frac{a_i'}{I}\}$$
    is strictly increasing, and in particular, $\lim_{i\rightarrow+\infty}A_i=+\infty$.
\end{itemize}
Since $f_i$ is a terminal blow-up, $a(E_i,Y_i,B_i)=\alpha$, hence $a_i\geq A_i$. Thus $\lim_{i\rightarrow+\infty}a_i=+\infty$. Possibly passing to a subsequence, we may assume that $a_i$ is strictly increasing and $a_i>5I$ for each $i$. By \cite[Theorem~1.1]{Kaw01} and Theorem \ref{thm: kaw05 1.2 strenghthened}, analytically locally, we have an embedding $$(X_i\ni x_i)\hookrightarrow (\mathbb{C}^{m_i}\ni o)/\frac{1}{n_i}(\alpha_{1,i},\dots,\alpha_{m,i})$$ for each $i$, where $n_i$ is the index of $X_i\ni x_i$, $n_i\mid I$, $m_i\in \{3,4,5\}$, $\alpha_{1,i},\dots,\alpha_{m,i}\in \Zz\cap[1,I]$, and $f_i$ is an admissible weighted blow-up with the weight $w_i\in \frac{1}{n_i}\Zz_{>0}^{m_i}$. Moreover, for each $i,j$, we may assume that $B_{i,j}'$ is defined by $(h_{i,j}=0)$ for some semi-invariant analytic power series $h_{i,j}$ near $x_i$. 

By Theorems \ref{thm: 3-dim terminal number of coefficients local} and \ref{thm: acc newton polytope}, possibly passing to a subsequence, we may assume that
\begin{itemize}
    \item there exist positive integers $n,m,\alpha_1,\dots,\alpha_m$ and a non-negative integer $p$, such that $n_i=n$, $m_i=m$, $(\alpha_{1,i},\dots,\alpha_{m,i})=(\alpha_1,\dots,\alpha_m)$, and $p_i=p$ for each $i$, 
    \item $b_{i,j}$ is increasing for any fixed $j$, and
    \item $\mathcal{N}(h_{i,j})\subset \mathcal{N}(h_{i',j})$ for any $i>i'$ and any $j$.
\end{itemize}
By Lemma \ref{lem: weighted blowup log discrepancies},
$$\frac{a_i}{I}=w_i(X_i\ni x_i)=w_i(B_i)+\alpha-1=\sum_{j=1}^{p} b_{i,j}w_i(h_{i,j})+\alpha-1.$$

We will show the following claim:

\begin{claim}\label{claim: auxiliary weights}

Possibly passing to a subsequence, we may assume that $(X_1\ni x_1, B_1)$ and $(X_2\ni x_2, B_2)$ satisfy the following:
\begin{enumerate}
    \item $A_2>a_1$, and 
    \item there exists an admissible weighted blow-up $f'$ of $(X_2\ni x_2)\subset(\mathbb{C}^m\ni o)/\frac{1}{n}(\alpha_1,\dots,\alpha_m)$ with the weight $w'\in\frac{1}{n}\Zz_{\ge 1}^m$, such that $w_1(B_2)\leq w'(B_2)$, $w'(X_2\ni x_2)=w_1(X_1\ni x_1)$, and the exceptional divisor $E'$ of $f'$ is an analytic prime divisor. 
\end{enumerate}
\end{claim}

%We remark that, even though $X_1\ni x_1$ and $X_2\ni x_2$ are embedded into the space $(\Cc^m\ni o)/\frac{1}{n}(\alpha_1,\dots,\alpha_m)$ possibly with different analytic local coordinates, we can still view $w_1,w_2$ and $w'$ as admissible weights respect to those analytic local coordinates simultaneously.

We proceed the proof assuming Claim~\ref{claim: auxiliary weights}. By Lemma~\ref{lem: algebraic approximation of weighted blow-up}, we may assume that $E'$ is a prime divisor over $X_2$. Since $$w_1(X_1\ni x_1)=w_1(B_1)+\alpha-1\leq w_1(B_2)+\alpha-1\leq w'(B_2)+\alpha-1\leq w'(X_2\ni x_2)=w_1(X_1\ni x_1),$$ $w'(X_2\ni x_2)=w'(B_2)+\alpha-1$ and $a(E',X_2,B_2)=\alpha=\mld(X_2\ni x_2,B_2)$. Since $a(E',X_2,0)=1+w'(X_2\ni x_2)=1+w_1(X_1\ni x_1)=1+\frac{a_1}{I}$, it follows that $a_1<A_2\leq a_1$, a contradiction. This finishes the proof.
\end{proof}
\begin{proof}[Proof of Claim~\ref{claim: auxiliary weights}] Claim~\ref{claim: auxiliary weights}(1) follows from the fact that the sequence $\{A_i\}_{i=1}^{\infty}$ is strictly increasing. We now prove Claim~\ref{claim: auxiliary weights}(2) case by case. By \cite[Theorem~1.1]{Kaw01} and Theorem~\ref{thm: kaw05 1.2 strenghthened}, we only need to consider the following cases. 
\medskip

\noindent {\bf Case 1}. $X_i\ni x_i$ ($i=1,2$) are all smooth. Then $m=3$, and analytically locally, $(X_i\ni x_i)\cong (\mathbb{C}^3\ni o)$. We may take $w'=w_1$ in this case.

\medskip

\noindent {\bf Case 2}. $X_i\ni x_i$ ($i=1,2$) are all of type $cD/n$ for $n=1$ or $2$. By Theorem~\ref{thm: kaw05 1.2 strenghthened}(2-3), we only need to consider the following two subcases:

\medskip

\noindent {\bf Case 2.1}. $m=4$. $f_i: Y_i\to X_i$ are divisorial contractions as in Theorem~\ref{thm: kaw05 1.2 strenghthened}(2.1) when $n=1$ and as in Theorem~\ref{thm: kaw05 1.2 strenghthened}(3.1) when $n=2$. In particular, there exist positive integers $d_i$ and $r_i$, such that  $2r_i+n=na_id_i$, and analytically locally, $$(X_i\ni x_i)\cong (\phi_i=0)\subset (\mathbb{C}^4\ni o)/\frac{1}{n}(1,1,1,0)$$ for some semi-invariant analytic power series $\phi_i$, and each $f_i$ is a weighted blow-up with the weight $w_i:=\frac{1}{n}(r_i+n,r_i,a_i,n)$. Possibly passing to a subsequence, we may assume that $d_1\leq d_2$. Let $s_2:=\frac{n}{2}(a_1d_2-1)$. Since $\frac{1}{n}(2s_2+n)=a_1d_2$ and $5\leq a_1< a_2$, by Lemmas~\ref{lem: cD irreducible case 1}(1) and \ref{lem: cD/2 irreducible case 1}, the weighted blow-up at $x_2\in X_2$ with the weight $ w':=\frac{1}{n}(s_2+n,s_2,a_1,n)$ extracts an analytic prime divisor over $X_2\ni x_2$, and $$w'(X_2\ni x_2)=\frac{a_1}{n}=w_1(X_1\ni x_1).$$ Since $d_2\geq d_1$ and $2r_1+n=na_1d_1$, $s_2=\frac{n}{2}(a_1d_2-1)\geq \frac{n}{2}(a_1d_1-1)= r_1$, hence $w_1(B_2)\leq w'(B_2)$.

\medskip

\noindent {\bf Case 2.2}. $m=5$. $f_i: Y_i\to X_i$ are divisorial contractions as in Theorem~\ref{thm: kaw05 1.2 strenghthened}(2.2) when $n=1$ and as in Theorem~\ref{thm: kaw05 1.2 strenghthened}(3.2) when $n=2$. In particular, there exist positive integers $d_i$ and $r_i$, such that $r_i+n=a_id_i$, and analytically locally, $$(X_i\ni x_i)\cong (\phi_{i,1}=\phi_{i,2}=0)\subset (\mathbb C^5\ni o)/\frac{1}{n}(1,1,1,0,1)$$ for some semi-invariant analytic power series $\phi_{i,1},\phi_{i,2}$, and each $f_i$ is a weighted blow-up with the weight $w_i:=\frac{1}{n}(r_i+n,r_i,a_i,n,r_i+2n)$. Possibly passing to a subsequence, we may assume that $d_1\leq d_2$. Let $s_2:=a_1d_2-n$. Since $s_2+n=a_1d_2$ and $5\leq a_1< a_2$, by Lemmas~\ref{lem: cD irreducible case 2}(1) and \ref{lem: cD/2 irreducible case 2}, the weighted blow-up at $x_2\in X_2$ with the weight $ w':=\frac{1}{n}(s_2+n,s_2,a_1,n,s_2+2n)$ extracts an analytic prime divisor over $X_2\ni x_2$, and
$$w'(X_2\ni x_2)=\frac{a_1}{n}=w_1(X_1\ni x_1).$$ Since $d_2\geq d_1$ and $r_1+n=a_1d_1$, $s_2=a_1d_2-n\geq a_1d_1-n= r_1$, hence $w_1(B_2)\leq w'(B_2)$.

\medskip

\noindent {\bf Case 3}. $m=4$, $X_i\ni x_i$ ($i=1,2$) are of type $cA/n$ and $f_i: Y_i\to X_i$ are divisorial contractions as in Theorem~\ref{thm: kaw05 1.2 strenghthened}(1). In particular, possibly passing to a subsequence, there exist positive integers $d_i,r_{1,i},r_{2,i},b$, such that $r_{1,i}+r_{2,i}=a_id_in$, $b\in [1,n-1]$, $\gcd(b,n)=1$, and analytically locally, $$(X_i\ni x_i)\cong (\phi_i=0)\subset (\mathbb C^4\ni o)/\frac{1}{n}(1,-1,b,0)$$for some semi-invariant analytic power series $\phi_i$, and each $f_i$ is a weighted blow-up with the weight $w_i:=\frac{1}{n}(r_{i,1},r_{i,2},a_i,n)$. Possibly passing to a subsequence, we may assume that $d_1\leq d_2$. Let $s_{2,1}:=r_{1,1}$ and $s_{2,2}:=a_1d_2n-r_{1,1}$. Since $s_{2,1}+s_{2,2}=a_1d_2n$ and $5\leq a_1<a_2$, by Lemma~\ref{lem: irreducibility of exceptional divisors extracted by weighted blow up}, the weighted blow-up at $x_2\in X_2$ with the weight $ w':=\frac{1}{n}(s_{2,1},s_{2,2},a_1,n)$ extracts an analytic prime divisor over $X_2\ni x_2$, and $$ w'(X_2\ni x_2)=\frac{a_1}{n}=w_1(X_1\ni x_1).$$ Since $d_1\leq d_2$ and $r_{1,1}+r_{1,2}=a_1d_1n$, $s_{2,2}=a_1d_2n-r_{1,1}\geq a_1d_1n-r_{1,1}=r_{1,2}$, hence $w_1(B_2)\leq w'(B_2)$.
\end{proof}

Now we prove some boundedness results on divisors computing mlds when the germ is either smooth or a terminal singularity of type $cA/n$.

\begin{lem}\label{lem: general sections map surjectively to O(1)}
Let $X$ be a smooth variety of dimension $n$ for some $n\in \Zz_{\geq 2}$ and $x\in X$ a closed point. Let $\pi: \tilde{X}\to X$ be the blow-up of $X$ at $x$ with the exceptional divisor $E$. For any hyperplane section $\tilde{H}\in |\mathcal{O}_{E}(1)|$ on $E$, there exists a Cartier divisor $H$ on $X$, such that $x\in \Supp H$, $\mult_x H=1$ and $\pi_*^{-1}H|_E=\tilde{H}$.
\end{lem}

\begin{proof}
Let $\fm_x$ be the maximal ideal of the local ring $\mathcal{O}_{X,x}$. Then we have a canonical isomorphism (cf. \cite[\S2,~Theorem~8.24(b)]{Har77}) $$E\cong \mathbf{P}^{n-1}=\mathrm{Proj}_{\Cc} \oplus_{i=0}^{+\infty} \fm_x^i/\fm_x^{i+1},$$ where $\fm_x^0:=\mathcal{O}_{X,x}$. Thus there exists a nonzero element $\tilde{h}\in\fm_x/\fm_x^2$ such that $\tilde{H}$ is defined by $(\tilde{h}=0)$ on $E$. Let $h\in \fm_x$ be a preimage of $\tilde{h}$ under the morphism $\fm_x\to \fm_x/\fm_x^2$, and $H$ the Cartier divisor locally defined by $(h=0)$ near $x$. We have $\mult_x H=1$ and $\pi_*^{-1}H|_E=\tilde{H}$ (cf. \cite[Exercise~III-29]{EH00} and \cite[Exercise~IV-24]{EH00}).
\end{proof}

\begin{lem}\label{lem: can extract mld place smooth case}
Let $(X\ni x, B)$ be a threefold germ such that $X$ is smooth and $\mld(X\ni x, B)\geq 1$. Then there exists a terminal blow-up (see Definition~\ref{defn: terminal blow-up}) $f: Y\to X$ of $(X\ni x,B)$, such that $Y$ is $\Qq$-factorial. 
\end{lem}

\begin{proof}
Let $g_1: X_1\to X_0:=X$ be the blow-up at $x\in X$ and $F_1$ the $g_1$-exceptional divisor. When $\mult_x B\leq 1$, by \cite[Proposition~6(i)]{Kaw17}, $a(F_1,X,B)=\mld(X\ni x, B)$, we may take $Y=X_1$, and $f=g_1$ in this case. From now on, we may assume that $\mult_x B>1$.

Let $g: W\to X_1$ be a birational morphism which consists of a sequence of blow-ups at points with codimension at least two, such that the induced morphism $h: W\to X$ is a log resolution of $(X,B)$. By Lemma~\ref{lem: general sections map surjectively to O(1)}, there exists a Cartier divisor $H$ on $X$ passing through $x$, such that $\mult_x H=1$, $h^*(H+B)$ is an snc divisor on $W$, and $H_{X_1}:=(g_1^{-1})_*H$ does not contain the center of any $g$-exceptional divisor on $X_1$.

Let $t:=\ct(X\ni x, B;H)$. Since $\mult_x B> 1$ and $1\leq a(E_1,X,B+tH)=3-\mult_x B-t$, $t<1$. By Lemma~\ref{lem: canonical threshold attain mld=1}(1), $\mld(X\ni x,B+tH)=1$. By Lemma~\ref{lem: can extract divisor computing ct that is terminal strong version}, there exists a terminal blow-up $f: Y\to X$ of $(X\ni x, B+tH)$ which extracts a prime divisor $E$ over $X\ni x$ such that $Y$ is $\Qq$-factorial. In particular, $a(E,X,B+tH)=\mld(X\ni x, B+tH)$.

It suffices to show that $a(E,X,B)=\mld(X\ni x,B)$. By \cite[Theorem 1.1]{Kaw01}, under suitable analytic local coordinates $(x_1,x_2,x_3)$, $f$ is the weighted blow-up of $X$ with the weight $(1,a,b)$ for some coprime positive integers $a$ and $b$. By \cite[Proof of Proposition~3.6, line 6]{Kaw01}, $\mult_E F_1=1$\footnote{We recall that $F_n$ and $E$ in \cite[Proposition~3.6]{Kaw01} are the same divisorial valuation, see \cite[Remark~3.3]{Kaw01}, and we use the same notion of $F_1$ and $E$ as in \cite[\S 3]{Kaw01}, see \cite[Construnction 3.1]{Kaw01}.}. By construction, $\mult_E H_X=\mult_E(H_{X_1}+F_1)=\mult_E F_1=1$. Let $F$ be any prime divisor over $X\ni x$. We have $$a(F,X,B)-t\mult_F H=a(F,X,B+tH)\geq a(E,X,B+tH)=a(E,X,B)-t\mult_E H.$$ Since $\mult_F H\geq 1=\mult_E H,$ $a(F,X,B)\geq a(E,X,B)$. It follows that $a(E,X,B)=\mld(X\ni x, B)$, and $f: Y\to X$ is the desired terminal blow-up of $(X\ni x, B)$.
\end{proof}

\begin{lem}\label{lem: can extract mld place cA/n case}
Let $x\in X$ be a threefold terminal singularity of type $cA/n$. Then there exists a Cartier divisor $C$ near $x$ satisfying the following.

Let $(X\ni x, B)$ be a pair such that $\mld(X\ni x, B)\geq 1$. Then there exists a terminal blow-up of $(X\ni x, B)$ (see Definition~\ref{defn: terminal blow-up}) which extracts a prime divisor $E$ over $X\ni x$, such that $a(E,X,B+tC)=1$ and $\mult_E C=1$, where $t:=\ct(X\ni x, B;C)$.
\end{lem}

\begin{proof}
By \cite[(6.1) Theorem]{Rei87} (cf. \cite[Theorems 12,23,25]{Mor85}), analytically locally, $$(X\ni x)\cong (\phi:= x_1x_2+g(x_3^n,x_4)=0)\subset (\mathbb C^4\ni o)/\frac{1}{n}(1,-1,b,0),$$ such that $b\in [1,n-1]\cap \Zz$, $\gcd(b,n)=1$, $\phi$ is a semi-invariant analytic power series, and $x_3^{dn}\in g(x_3^n,x_4)$ for some positive integer $d$. When $n=1$, $x\in X$ is a terminal singularity of type $cA_m$ (see \cite[Page 333, Line 11]{Kaw03}) for some positive integer $m$. 

\begin{claim}\label{claim: choosing a Cartier divisor C with mult 1}
There exist a Cartier divisor $C$ on $X$ and an integer $k\in \{1,4\}$ depending only on $x\in X$ that satisfy the following.
\begin{enumerate}
    \item $x_k$ is invariant under the $\xi_n$-action on $\Cc^4$, and $C$ is a Cartier divisor locally defined by an invariant analytic power series $(x_k+h=0)$, where $h\in (\mathfrak{m}_o^{\rm an})^2$.
    \item Let $t:=\ct(X\ni x,B;C)$. Then $\mld(X\ni x, B+tC)=1$, and there exists a terminal blow-up $f$ of $(X\ni x, B+tC)$ which extracts a prime divisor $E$ over $X\ni x$.
    \item Possibly choosing a new local analytic coordinates $x_1',x_2',x_3',x_4'$, where $x_k=x_k'+p_k'$ for some $p_k'\in \mathfrak{m}_o^{\rm an}$ such that $\lambda x_k'\notin p_k'$ for any $\lambda\in \Cc^*$, analytically locally, $f$ is a weighted blow-up with the weight $w=(w(x_1'),w(x_2'),w(x_3'),w(x_4'))$ such that $w(x_k')=1$.
\end{enumerate}
\end{claim}

We proceed the proof assuming Claim~\ref{claim: choosing a Cartier divisor C with mult 1}. By Claim~\ref{claim: choosing a Cartier divisor C with mult 1}(1,3), $C$ is locally defined by $(x_k'+h_k'=0)$ for some $h_k'\in \mathfrak{m}_o^{\rm an}$ such that $\lambda x_k'\notin h_k'$ for any $\lambda\in \Cc^*$ under the new coordinates $x_1', x_2', x_3', x_4'$, hence $1\leq \mult_E C=w(x_k'+h')\leq w(x_k')=1$, and $\mult_E C=1$. Let $F\neq E$ be any prime divisor over $X\ni x$. We have $$a(F,X,B)-t\mult_F C=a(F,X,B+tC)\geq a(E,X,B+tC)=a(E,X,B)-t\mult_E C.$$ Since $\mult_F C\geq 1=\mult_E C$, $a(F,X,B)\geq a(E,X,B)$. It follows that $a(E,X,B)=\mld(X\ni x, B)$, hence $f$ is a terminal blow-up of $(X\ni x, B)$.
\end{proof}

\begin{proof}[Proof of Claim~\ref{claim: choosing a Cartier divisor C with mult 1}]
We have an analytic isomorphism $$\psi: \tilde{X}\ni \tilde{x}\to \tilde{Y}:=(\phi: x_1x_2+g(x_3^n,x_4)=0)\subset (\mathbb C^4\ni o),$$ where $\pi:\tilde{X}\ni \tilde{x} \to X\ni x$ is the index one cover (cf. \cite[Definition~5.19]{KM98}). Under the analytic isomorphism $\psi$, the $\xi_n$-action on $\tilde{Y}$ induces the cyclic group action on $\tilde{X}\ni \tilde{x}$ which corresponds to $\pi$. By Lemma~\ref{lem: formal approximation divisors}, we can find a Cartier divisor $\tilde{C}$ on $\tilde{X}$ whose image under $\psi$ is locally defined by $(x_k+h=0)$ for some $h\in (\fm_o^{\rm an})^2$, and $x_k+h$ is invariant under $\xi_n$-action. Set $C:=\pi(\tilde{C})$, we finish the proof of Claim~\ref{claim: choosing a Cartier divisor C with mult 1}(1). For Claim~\ref{claim: choosing a Cartier divisor C with mult 1}(2), since $C$ is a prime divisor that is Cartier, by Lemma~\ref{lem: canonical threshold attain mld=1}(3), $\mld(X\ni x, B+tC)=1$. By Lemma~\ref{lem: can extract divisor computing ct that is terminal strong version}, there exists a terminal blow-up of $(X\ni x, B+tC)$ which extracts a prime divisor $E$ over $X\ni x$. Now we prove Claim~\ref{claim: choosing a Cartier divisor C with mult 1}(3) case by case. 

\medskip

\noindent {\bf Case 1}. $n\geq 2$. By \cite[Theorem~1.3]{Kaw05}, $f$ is a divisorial contraction of ordinary type. By Theorem~\ref{thm: kaw05 1.2 strenghthened}(1) and \cite[Lemmas~6.1,~6.2~and~6.5]{Kaw05}, there exist analytic local coordinates $x_1',x_2',x_3',x_4'$, such that analytically locally, $f$ is a weighted blow-up with the weight $w:=\frac{1}{n}(r_1',r_2',a,n)$, where $r_1',r_2',a$ are positive integers such that $an\mid r_1'+r_2'$. Moreover, by \cite[Proof of Lemma 6.3, Line 7]{Kaw05}, $x_4'=x_4+x_1p$ for some $p\in \mathfrak{m}_o^{\rm an}$. In this case, we take $k=4$.

\medskip

\noindent {\bf Case 2}. $n=1$ and $m\geq 2$. By \cite[Theorem~1.13]{Kaw03} and \cite[Theorem~2.6]{Yam18}, there exist analytic local coordinates $x_1',x_2',x_3',x_4'$, such that analytically locally, $f$ is a weighted blow-up with the weight $w=(r_1,r_2,a,1)$ for some positive integers $r_1,r_2,a$. Moreover, by \cite[Proof of Lemma 6.1, Page 309, Line 5]{Kaw03}, the coordinates change relation for $x_4$ is given by $x_4=x_4'+c'x_i'$ for some $1\leq i\leq 3$ and $c'\in \Cc$. Thus we may take $k=4$.

\medskip

\noindent {\bf Case 3}. $n=1$ and $m=1$. By \cite[Theorem~1.1]{Kaw02}, there exist analytic local coordinates $x_1',x_2',x_3',x_4'$, such that analytically locally, $f$ is a weighted blow-up with the weight $w$, where either $w=(s,2a-s,a,1)$ for some positive integers $s$ and $a$, or $w=(1,5,3,2)$. Moreover, by \cite[Claim 6.13]{Kaw02}, the change of coordinates relations for $x_i$ is given by $x_i=x_i'+p_i'(x_4')$ and $x_4=x_4'$ for $p_i'\in {\fm}_o^{\rm an}$ and $1\leq i\leq 3$. Thus we may take $k=1$ or $4$.
\end{proof}

\begin{lem}\label{lem: finiteness of description of terminal threefold singularities}
Let $\Ii\subset [0,1]$ be a set such that $\gamma_0:=\inf \{b\mid b\in\Ii\setminus\{0\}\}>0$. Let $(X\ni x,B)$ be a threefold germ, such that 
\begin{itemize}
    \item $x\in X$ is a terminal singularity of type $cA/n$ for some $n>N:=\lceil\frac{3}{\gamma_0}\rceil$,
    \item $B:=\sum_{i} b_iB_i$ for some $b_i\in \Ii$, where $B_i\geq 0$ are $\Qq$-Cartier Weil divisors, %\han{if you use some $m$ here then other place you also need, maybe delete it}
    \item $\mld(X\ni x,B)\geq 1$, and
    \item there exists a terminal blow-up (see Definition \ref{defn: terminal blow-up}) $f: Y\to X$ of $(X\ni x,B)$ which extracts a prime divisor $E$ over $X\ni x$, such that $a(E,X,0)=1+\frac{a}{n}$ for some positive integer $a\geq 3$.
\end{itemize}
Then there exists a prime divisor $\bar E$ over $X\ni x$ such that $a(\bar E,X,B)=\mld(X\ni x, B)$ and $a(\bar E,X,0)=1+\frac{3}{n}$. Moreover,
\begin{enumerate}
\item if $\mld(X\ni x, B)>1$, then $a=3$, and 
\item if $\mld(X\ni x,B)=1$ and $\Gamma$ is either a DCC set or an ACC set, then $b_i\in\Gamma_0$ for some  finite set $\Gamma_0$ depending only on $\Gamma$.
\end{enumerate}
\end{lem}

\begin{proof}
Since $n>1$, by \cite[Theorem~1.3]{Kaw05}, $f: Y\to X$ is a divisorial contraction of ordinary type as in Theorem~\ref{thm: kaw05 1.2 strenghthened}(1). In particular, under suitable analytic local coordinates $x_1,x_2,x_3,x_4$, there exist positive integers $r_1,r_2,b,d$ such that $r_1+r_2=adn$, $b\in \Zz\cap [1,n-1]$, $\gcd(b,n)=1$, $a\equiv br_1\mod n$, and analytically locally, $$(X\ni x)\cong (\phi(x_1,x_2,x_3,x_4)=0)\subset (\mathbb C^4\ni o)/\frac{1}{n}(1,-1,b,0)$$
for some invariant analytic power series $\phi$, and $f: Y\to X$ is a weighted blow-up with the weight $w:=\frac{1}{n}(r_1,r_2,a,n)$. Assume that each $B_{i}$ is locally defined by $(h_i=0)$ for some semi-invariant analytic power series $h_{i}$.

\medskip
    
Since $n>N\geq \frac{3}{\gamma_0}$, we can pick positive integers $s_1,s_2$, such that
\begin{itemize}
    \item $s_1+s_2=3dn$,
    \item $3\equiv bs_1\mod n$, and
    \item $s_1,s_2>n$.
\end{itemize}
Let $\bar w:=\frac{1}{n}(s_1,s_2,3,n)$. Since $a\geq 3$, by Lemma~\ref{lem: irreducibility of exceptional divisors extracted by weighted blow up}, the weighted blow-up with the weight $\bar w$ extracts an analytic prime divisor $\bar E$ over $X\ni x$, such that $a(\bar E,X,0)=1+\bar w(X\ni x)=1+\frac{3}{n}$. By Lemma~\ref{lem: algebraic approximation of weighted blow-up}, we may assume that $\bar E$ is a prime divisor over $X\ni x$. Since $a(\bar{E},X,B)=1+\bar w(X\ni x)-\bar{w}(B)\geq \mld(X\ni x, B)\geq 1$, $$\gamma_0>\frac{3}{n}=\bar w(X\ni x)\geq \bar w(B)= \sum_{i} b_i\bar w(B_i)\geq\gamma_0\sum_{i}\bar w(B_i),$$
which implies that $\bar w(h_i)=\bar w(B_i)<1$ for each $i$.
Since $\bar w(x_1)=\frac{s_1}{n}>1$, $\bar w(x_2)=\frac{s_2}{n}>1$, and $\bar w(x_4)=1$, for each $i$, there exists a positive integer $l_i$, such that up to a scaling of $h_i$, $x_3^{l_i}\in h_i$ for each $i$, and $\bar{w}(B_i)=\bar w(h_i)=\bar w(x_3^{l_i})$. In particular, 
$$\mult_{\bar E} B=\bar w(B)=\sum_{i} b_il_i\bar w(x_3)=\frac{3}{n}\sum_{i} b_il_i,$$ 
and $1+\frac{3}{n}\geq\bar w(B)+\mld(X\ni x,B)=\frac{3}{n}\sum_{i} b_il_i+\mld(X\ni x,B)$. This implies that
\begin{align}\label{eqn: the sum relation for coeff 1}
\mld(X\ni x,B)-1\leq \frac{3}{n}(1-\sum_{i} b_il_i).
\end{align}
On the other hand, 
$$\mult_E B=w(B)=\sum_{i} b_iw(B_i)\leq\sum_{i} b_iw(x_3^{l_i})=\frac{a}{n}\sum_{i} b_il_i,$$ 
and
$$\frac{a}{n}-\mld(X\ni x,B)+1=w(X\ni x)-\mld(X\ni x,B)+1=w(B)\leq\frac{a}{n}\sum_{i} b_il_i.$$
Combining with \eqref{eqn: the sum relation for coeff 1}, we have
\begin{align}\label{eqn: the sum relation for coeff 2}
\frac{a}{n}(1-\sum_{i} b_il_i)\leq \mld(X\ni x,B)-1\leq \frac{3}{n}(1-\sum_{i} b_il_i).
\end{align}
If $\mld(X\ni x, B)>1$, then by \eqref{eqn: the sum relation for coeff 2}, $a\leq 3$, hence $a=3$. It follows that $a(\bar E,X,B)=1+\frac{3}{n}-\bar w(B)=\mld(X\ni x, B)$ in this case. If $\mld(X\ni x, B)=1$, then by \eqref{eqn: the sum relation for coeff 2}, $\sum_{i} b_il_i=1$. In particular, $\bar{w}(B)=\frac{3}{n}=\bar{w}(X\ni x)$, hence $\mld(X\ni x,B)=a(\bar{E},X,B)=1$.

\medskip

When $\mld(X\ni x, B)=1$ and $\Gamma$ is a DCC set or an ACC set, the equality $\sum_{i} b_il_i=1$ implies that $B$ belongs to a finite subset $\Gamma_0\subset \Gamma$.
\end{proof}

\begin{lem}\label{lem: bdd mld computing dcc cA/n type}
Let $\gamma_0$ be a positive real number. Let $(X\ni x,B:=\sum_{i}b_i B_i)$ be a threefold germ, where $x\in X$ is a terminal singularity of type $cA/n$ for some $n>\lceil\frac{3}{\gamma_0}\rceil$, $b_i\geq\gamma_0$ and $\Qq$-Cartier Weil divisors $B_i\geq 0$, such that $\mld(X\ni x,B)\ge 1$. Then there exists a prime divisor $E$ over $X\ni x$, such that $a(E,X,B)=\mld(X\ni x,B)$ and $a(E,X,0)\leq 1+\frac{3}{n}$.
\end{lem}

\begin{proof}
This follows from Lemmas~\ref{lem: can extract mld place cA/n case} and \ref{lem: finiteness of description of terminal threefold singularities}.
\end{proof}

\begin{thm}\label{thm: Nakamura Mustata Conjecture}
Let $\Ii \subset [0,1]$ be a DCC set. Then there exists a positive integer $l$ depending only on $\Ii$ satisfying the following.

Let $(X\ni x,B)$ be a threefold pair such that $X$ is terminal, $B\in \Ii$, and $\mld(X\ni x,B)=1$. Then there exists a prime divisor $E$ over $X\ni x$, such that $a(E,X,B)=1$ and $a(E,X,0)\le 1+\frac{l}{I}$, where $I$ is the index of $X\ni x$. In particular,  $a(E,X,0)\le 1+l$.
\end{thm}

\begin{proof}
Let $Y$ be a small $\Qq$-factorialization of $X$, and let $K_Y+B_Y:=f^*(K_X+B)$. There exists a point $y\in Y$ such that $f(y)=x$ and $\mld(Y\ni y,B_Y)=\mld(X\ni x, B)=1$. Moreover, the index of $Y\ni y$ divides the index of $X\ni x$. Possibly replacing $(X\ni x, B)$ with $(Y\ni y, B_Y)$, we may assume that $X$ is $\Qq$-factorial.

If $\dim x=2$, then the theorem is trivial as we can take $l=0$. If $\dim x=1$, then $X$ is smooth near $x$ and $I=1$. By Lemma~\ref{lem: Terminal blow up}, if $E$ is the exceptional divisor of the blow-up at $x\in X$, then $a(E,X,B)=\mld(X\ni x, B)$. Since $a(E,X,0)=2$, we may take $l=1$ in this case. 

Now we may assume that $\dim x=0$. By Lemma \ref{lem: can extract divisor computing ct that is terminal strong version}, there exists a terminal blow-up (see Definition \ref{defn: terminal blow-up}) $f: Y\to X$ of $(X\ni x,B)$ which exactly extracts a prime divisor $E$ over $X\ni x$. We may write $$K_Y-\frac{a}{I}E=f^*K_X$$
for some positive integer $a$. Moreover, we may assume that $a\geq 5$.
 
By Theorem~\ref{thm: kaw05 1.2 strenghthened}, $f$ is a divisorial contraction of ordinary type. If $x\in X$ is a terminal singularity of type other than $cA/n$, then $I\leq 2$, and the theorem follows from Lemma~\ref{lem: bounded index Mustata nakamura conjecture}. If $x\in X$ is a terminal singularity of type $cA/n$, then by Lemma~\ref{lem: finiteness of description of terminal threefold singularities}, there exists an integer $N'$ depending only on $\Ii'$, such that if $n=I\geq N'$, then there exists a prime divisor $\bar E$ over $X\ni x$ with $a(\bar E,X,B)=1$ and $a(\bar E,X,0)=1+\frac{3}{I}$. Hence when $I\geq N'$, we may take $l=3$, and when $I<N'$, the theorem follows from Lemma~\ref{lem: bounded index Mustata nakamura conjecture}.
\end{proof}

\subsection{Proof of Theorem~\ref{thm: 3fold acc ct}}
%\subsection{The terminal germ case}

\begin{thm}\label{thm: 3fold acc ct general terminal germ case}
Let $\Ii\subset [0,1]$ be a DCC set. Then there exists a finite set $\Ii_0\subset\Ii$ depending only on $\Ii$ satisfying the following. Assume that
\begin{itemize}
    \item $(X\ni x,B:=\sum_ib_iB_i)$ is a threefold pair, 
    \item $X$ is terminal,
    \item  $b_i\in \Ii$ and $B_i\geq 0$ are $\Qq$-Cartier Weil divisors, and
    \item $\mld(X\ni x,B)=1$.
\end{itemize}
Then $b_i\in\Ii_0$ for all $i$.
\end{thm}

\begin{proof}
We may assume that $\dim x\leq 1$. If $\dim x=1$, then $X$ is smooth near $x$. By Lemma~\ref{lem: Terminal blow up}, $\mld(X\ni x, B)=2-\sum_i b_i\mult_x B_i=1$ and $\mult_xB_i\in\mathbb Z_{>0}$ for each $i$, hence $b_i$ belongs to a finite set $\Ii_0\subset \Ii$ depending only on $\Ii$. If $\dim x=0$, then we let $n$ be the index of $X\ni x$. By Theorem~\ref{thm: Nakamura Mustata Conjecture}, there exists a prime divisor $E$ over $X\ni x$ such that $a(E,X,B)=a(E,X,0)-\mult_E B=1$ and $a(E,X,0)=1+\frac{a}{n}$ for some $a\leq 
l$, where $l$ is a positive integer depending only on $\Ii$. By \cite[Lemma~5.1]{Kaw88}, $\mult_E B_i=\frac{1}{n}c_i$ for some positive integers $c_i$. It follows that $\frac{a}{n}=\frac{1}{n}\sum_{i=1}^m c_ib_i$. Thus $b_i\in\Ii_0\subset\Ii$ for all $i$ for some finite set $\Ii_0\subset \Ii$ depending only on $\Ii$.
\end{proof}

As a consequence of Theorem~\ref{thm: 3fold acc ct general terminal germ case}, we show the ACC for $\ct(X\ni x, B;D)$ for terminal threefold singularities $x\in X$. 

\begin{thm}\label{thm: ACC ct local with base terminal}
Let $\Ii\subset [0,1],~\Ii'\subset [0,+\infty)$ be two DCC sets. Then the set 
$$\{\ct(X\ni x, B;D)\mid \dim X=3,~X \text{ is terminal},~B\in \Ii,~D\in \Ii'\}$$
satisfies the ACC.
\end{thm}

\begin{proof}
Pick $t\in\{\ct(X\ni x, B;D)\mid \dim X=3,~X \text{ is terminal},~B\in \Ii,~D\in \Ii'\}$. Then there exists a threefold pair $(X\ni x,B)$ and an $\Rr$-Cartier $\Rr$-divisor $D\in\Ii'\backslash\{0\}$ on $X$, such that $X$ is terminal, $B\in\Ii$, and $t=\ct(X\ni x, B;D)$. 

We only need to show that $t$ belongs to an ACC set depending only on $\Ii$ and $\Ii'$. By \cite[Theorem 1.1]{HMX14}, we may assume that $\mld(X\ni x,B+tD)=1$, and there exists a prime divisor $E$ over $X\ni x$ such that $a(E,X,B+tD)=1$. Possibly replacing $X$ with a small $\Qq$-factorialization $X'$ and replacing $x$ with the generic point of $\Center_{X'}E$, we may assume that $X$ is $\Qq$-factorial. By Theorem~\ref{thm: 3fold acc ct general terminal germ case}, $B+tD$ belongs to a finite set depending only on $\Ii$ and $\Ii'$, hence $t$ belongs to an ACC set depending only on $\Ii$ and $\Ii'$.
\end{proof}

%\han{for submitted version, we may just say replace it by terminalization.}
\begin{proof}[Proof of Theorem~\ref{thm: 3fold acc ct}] 
Let $(X,B)$ be a canonical threefold pair and $D\geq 0$ a non-zero $\Rr$-Cartier $\Rr$-divisor on $X$, such that $B\in \Ii$ and $D\in\Ii'$. Let $t:=\ct(X,B;D)$. We only need to show that $t$ belongs to an ACC set.

We may assume that $t>0$. In particular, $(X,B)$ is canonical. By \cite[Corollary~1.4.3]{BCHM10}, there exists a birational morphism $f: Y\to X$ that exactly extracts all the exceptional divisors $E$ over $X$ such that $a(E,X,0)=1$. Since $(X,B)$ is canonical, $X$ is canonical, hence $Y$ is terminal and $a(E,X,B+tD)=1$ for any $f$-exceptional divisor $E$ such that $a(E,X,0)=1$. We have $K_Y+B_Y+tD_Y=f^*(K_X+B+tD)$, where $B_Y,D_Y$ are the strict transforms of $B,D$ on $Y$ respectively. Possibly replacing $(X,B)$ and $D$ with $(Y,B_Y)$ and $D_Y$ respectively, we may assume that $X$ is terminal.

Now there exists a point $x$ on $X$ such that $t=\ct(X\ni x,B;D)$. Theorem~\ref{thm: 3fold acc ct} follows from Theorem \ref{thm: ACC ct local with base terminal}.
\end{proof}

\section{ACC for minimal log discrepancies on \texorpdfstring{$[1,+\infty)$}{Lg}}\label{section: ACC for minimal log discrepancies on 1,oo}

In this section, we prove the following theorem:

\begin{thm}\label{thm:  terminal mld acc}
Let $\Ii\subset [0,1]$ be a DCC set. Then the set
$$\{\mld(X\ni x,B)\mid \dim X=3, X\text{ is terminal near $x$},B\in\Ii\}\cap [1,+\infty)$$
satisfies the ACC.
\end{thm}

\subsection{Uniform canonical rational polytopes}\label{subsection: uniform can rational polytope}

\cite{HLS19} established a general theory to show the boundedness of complements for DCC coefficients from the boundedness of complements for finite rational coefficients. Some embryonic forms of this theory can be found in \cite{Kaw14,Nak16,Liu18}. For other related results we refer the reader to \cite{Che20,HLQ20,CH21,HLQ21}.
We will follow this theory in our paper. As the key step, we need to show the existence of uniform canonical rational polytopes in this section. Recall that the proof of the uniform lc rational polytopes \cite{HLS19} is based on some ideas in the proof of accumulation points of lc
thresholds \cite{HMX14}, which relies on applying the adjunction formula to the lc places. Our proof is quite different from \cite{HLS19} as we could not apply the adjunction formula to canonical places.

%\han{maybe only prove for canonical pairs, the rest you only use for canonical but not $\epsilon$-lc? or write together}

\begin{lem}\label{lem: fixed germ uniform canonical polytope}
Let $I,c,m$ be three non-negative integers, $r_1,\dots,r_c$ real numbers such that $1,r_1,\dots,r_c$ are linearly independent over $\Qq$, and $s_1,\dots,s_m:\mathbb R^{c+1}\rightarrow\mathbb R$ $\Qq$-linear functions. Let $\bm{r}:=(r_1,\dots,r_c)$. Then there exists an open subset $U\subset\mathbb R^c$ depending only on $I,\bm{r}$ and $s_1,\dots,s_m$, such that $U\ni\bm{r}$ satisfies the following.

Let $x\in X$ be a terminal threefold singularity such that $IK_X$ is Cartier near $x$, $B_1,\dots,B_m\geq 0$ Weil divisors on $X$ such that $(X\ni x,B:=B(\bm{r}))$ is lc and $\mld(X\ni x,B)\geq 1$, where $B(\bm{v}):=\sum_{j=1}^ms_j(1,\bm{v})B_j$ for any $\bm{v}\in\Rr^c$. Then $(X\ni x,B(\bm{v}))$ is lc and $\mld(X\ni x,B(\bm v))\geq 1$ for any $\bm{v}\in U$.
\end{lem}

\begin{proof}
By \cite[Theorem~5.6]{HLS19}, we may pick an open subset $U_0\subset\mathbb R^c$ such that $\bm{r}\in U_0$ and $(X\ni x,B(\bm{v}))$ is lc for any $\bm{v}\in U_0$. By \cite[Lemma~5.1]{Kaw88}, $IB_j$ is Cartier near $x$ for $1\leq j\leq m$, we may write $B=\sum_{j=1}^m \frac{s_j(1,\bm r)}{I} IB_j$. By \cite[Theorem~1.2]{Nak16}, $\{a(E,X,B)\mid \Center_{X}E=x\}$ belongs to a discrete set depending only on $I,\bm{r}$ and $s_1,\ldots,s_m$. In particular, we may let
$$\alpha:=\min\{a(E,X,B)\mid \Center_{X}E=x,\, a(E,X,B)>1\}.$$
Now we let
$$U:=\{\frac{1}{\alpha}\bm{r}+\frac{\alpha-1}{\alpha}\bm{v}_0\mid \bm{v}_0\in U_0\}.$$
We show that $U$ satisfies our requirements. For any prime divisor $E$ over $X\ni x$, if $a(E,X,B)=1$, then $a(E,X,B(\bm{v}))=1$ for any $\bm{v}\in U$ as $r_1,\dots,r_c$ are linearly independent over $\Qq$. If $a(E,X,B)>1$, then $a(E,X,B)\geq \alpha$. By the construction of $U$, for any $\bm{v}\in U$, there exists $\bm{v}_0\in U_0$ such that $\bm{v}=\frac{1}{\alpha}\bm{r}+\frac{\alpha-1}{\alpha}\bm{v}_0$. Hence
$$a(E,X,B(\bm{v}))=\frac{1}{\alpha}a(E,X,B)+\frac{\alpha-1}{\alpha}a(E,X,B(\bm{v}_0))\geq 1.$$
It follows that $\mld(X\ni x,B(\bm{v}))\geq 1$ for any $\bm{v}\in U$.
\end{proof}

\begin{lem}\label{lem: index conjecture threefold non-klt case}
Let $(X\ni x, B)$ be a threefold pair such that $X$ is $\Qq$-factorial and $x\in X$ is a terminal singularity of type $cA/n$. Assume that $\mld(X\ni x, B)\geq 1$ and $\lfloor B\rfloor\neq \emptyset$. Then $B=\lfloor B\rfloor$ is a prime divisor, $K_X+B$ is Cartier near $x$, and $\mld(X\ni x, B)=1$.
\end{lem}

\begin{proof}
Let $S\subset\lfloor B\rfloor$ be a prime divisor. By Theorem~\ref{thm: 3-dim terminal number of coefficients local}(2), $B=\lfloor B\rfloor=S$ and $\mld(X\ni x, B)=1$. If $n=1$, by \cite[Lemma 5.1]{Kaw88}, $K_X+B$ is Cartier near $x$. We may assume that $n\geq 2$. By \cite[Theorem 1.3]{Kaw05}, Theorem \ref{thm: kaw05 1.2 strenghthened}(2), and Lemmas \ref{lem: algebraic approximation of weighted blow-up} and \ref{lem: irreducibility of exceptional divisors extracted by weighted blow up}, there exist analytic local coordinates $x_1,x_2,x_3,x_4$ and a positive integer $d$, such that analytically locally, $x\in X$ is a hyperquotient singularity of the form $$(X\ni x)\cong (\phi:= x_1x_2+g(x_3^n,x_4)=0)\subset (\mathbb C^4\ni o)/\frac{1}{n}(1,-1,b,0),$$
where $b\in [1,n-1]\cap \Zz$ and $\gcd(b,n)=1$. Moreover, the weighted blow-up with the weight $w:=\frac{1}{n}(s_1,dn-s_1,1,n)$ extracts a prime divisor $E$ such that $a(E,X,0)=1+\frac{1}{n}$, where $s_1\in [1,n-1]\cap \Zz$ and $bs_1\equiv 1\mod n$. Let $h$ be a semi-invariant analytic power series which defines $S$. Since $a(E,X,S)=a(E,X,0)-\mult_E S=1+\frac{1}{n}-\mult_E S\geq 1$, $w(h)=\mult_E S=\frac{1}{n}$, and if $s_1=1$ (resp. $dn-s_1=1$), then either $x_3\in h$ or $x_1\in h$ (resp. $x_2\in h$) up to a scaling of $h$. If $s_1=1$ (resp. $dn-s_1=1$), then $b=1$ (resp. $b=-1$), and the analytic Cartier divisor $(x_3=0)$ is linearly equivalent to $(x_1=0)$ (resp. $(x_2=0)$), hence the $\Qq$-Cartier divisor $(h=0)$ is linear equivalent to the $\Qq$-Cartier divisor $(x_3=0)$. By \cite[(6.4)(B.1)]{Rei87} and Lemma~\ref{lem: analytic index and algebraic index}, $K_X+S$ is Cartier near $x$.
\end{proof}

\begin{thm}\label{thm: uniform canonical polytope}
Let $c,m$ be two non-negative integers, $r_1,\dots,r_c$ real numbers such that $1,r_1,\dots,r_c$ are linearly independent over $\Qq$, and $s_1,\dots,s_m:\mathbb R^{c+1}\rightarrow\mathbb R$ $\Qq$-linear functions. Let $\bm{r}:=(r_1,\dots,r_c)$. Then there exists an open subset $U\subset\mathbb R^c$ depending only on $\bm{r}$ and $s_1,\dots,s_m$, such that $U\ni\bm{r}$ satisfies the following.

Let $X$ be a terminal threefold, $x\in X$ a point, $B_1,\dots,B_m\geq 0$ distinct Weil divisors on $X$, and $B(\bm{v}):=\sum_{j=1}^ms_j(1,\bm{v})B_j$ for any $\bm{v}\in\mathbb R^c$. Assume that $(X\ni x,B:=B(\bm{r}))$ is lc and $\mld(X\ni x, B)\geq 1$. Then $(X\ni x,B(\bm{v}))$ is lc and $\mld(X\ni x, B(\bm{v}))\geq 1$ for any $\bm{v}\in U$. Moreover, if $\mld(X\ni x, B)> 1$, then we may choose $U$ so that $\mld(X\ni x, B(\bm{v}))> 1$ for any $\bm{v}\in U$.
\end{thm}

\begin{proof}
Possibly replacing $X$ with a small $\Qq$-factorialization, we may assume that $X$ is $\Qq$-factorial.

By construction, we may assume that $s_j(1,\bm r)>0$ for each $j$. If $\dim x=2$, then the theorem is trivial. If $\dim x=1$, then $X$ is smooth near $x$. By Lemma~\ref{lem: Terminal blow up}, $\mld(X\ni x, B)=2-\mult_x B> 1$, thus $\mult_x B=\sum_{j=1}^m s_j(1,\bm r)\mult_x B_j<1$, where $\mult_x B_j$ are non-negative integers. Hence $\sum_{j=1}^m s_j(1,\bm r)\mult_x B_j\leq 1-\epsilon_0$ for some $\epsilon_0\in(0,1)$ depending only on $\bm r$ and $s_1,\dots,s_m$. By Lemma~\ref{lem: Terminal blow up}, $\mld(X\ni x, \frac{1}{1-\epsilon_0} B)\geq 1$, hence we can take $U:=\{\bm v\mid 0<s_j(1,\bm v)< \frac{1}{1-\epsilon_0}s_j(1,\bm r)\text{ for each $j$}\}$ in this case. Hence we may assume that $\dim x=0$.

If $\lfloor B\rfloor\not=0$, by \cite[(6.1) Theorem]{Rei87} and Lemma \ref{lem: index conjecture threefold non-klt case}, we may assume that $12K_X$ is Cartier near $x$, and the theorem follows from Lemma \ref{lem: fixed germ uniform canonical polytope}. Thus we may assume that $\lfloor B \rfloor=0$. By \cite[Corollary~1.4.3]{BCHM10}, there exists a birational morphism $f: Y\to X$ from a $\Qq$-factorial variety $Y$ that exactly extracts all the exceptional divisors $F$ over $X\ni x$ such that $a(F,X,B)=1$. In particular, $a(F,X,B(\bm v))=1$ for all $\bm v\in \Rr^c$. It follows that $f^*(K_X+B(\bm v))=K_Y+f_*^{-1}B(\bm v)$ for all $\bm v\in \Rr^c$. Hence it suffices to prove the theorem for all pairs $(Y\ni y, f_*^{-1}B)$, where $y\in f^{-1}(x)$ is a closed point. From now on, we may assume that $\mld(X\ni x, B)>1$.

By \cite[(6.1) Theorem]{Rei87}, if $x\in X$ is a terminal singularity of types other than $cA/n$, then the index of $X\ni x$ is $\leq 4$, and the theorem holds by Lemma~\ref{lem: fixed germ uniform canonical polytope}. From now on, we may assume that $x\in X$ is of type $cA/n$.

%Consider the following set $$\Ii:=\{\ct(X\ni x,0;\sum_{j=1}^m s_j(1,\mathbf{r})B_j)\mid~x\in X \text{ is of type $cA/n$},~ B_1,\dots,B_m> 0\text{ Weil divisors}\}.$$ By Theorem~\ref{thm: ACC ct local with base terminal}, $\Gamma$ is an ACC set.

\begin{claim}\label{claim: ct lower bound to 1}
There exist a positive integer $N$ and a positive real number $\epsilon$ depending only on $\bm r$ and $s_1,\dots,s_m$ satisfying the following.

For any terminal threefold singularity $x\in X$ of type $cA/n$ and $B:=\sum_{j=1}^m s_j(1,\mathbf{r})B_j$, where $B_j\geq 0$ are $\Qq$-Cartier Weil divisors on $X$ and $\lfloor B\rfloor=0$, if $\mld(X\ni x,B)>1$ and $n>N$, then $t:=\ct(X\ni x,0;B)>1+\epsilon$.
\end{claim}

We proceed the proof assuming Claim~\ref{claim: ct lower bound to 1}. By Lemma~\ref{lem: fixed germ uniform canonical polytope}, we may assume that $n>N$. By Claim~\ref{claim: ct lower bound to 1}, $t> 1+\epsilon$, hence we can take $U:=\{\bm v\mid 0< s_j(1,\bm v)< (1+\epsilon)s_j(1,\bm r)\text{ for each $j$}\}$ in this case. Moreover, if $\mld(X\ni x, B)>1$, then possibly replacing $U$ with $\{\frac{1}{2}\bm v+\frac{1}{2}\bm r\mid \bm v\in U\}$, we have $\mld(X\ni x, B(\bm v))>1$ for all $\bm v\in U$.
\end{proof}

\begin{proof}[Proof of Claim~\ref{claim: ct lower bound to 1}] Since $\mld(X\ni x,B)>1$, $t>1$. Since $\lfloor B\rfloor=0$, if $\lfloor tB\rfloor\neq  0$, then $t>1+\epsilon$ for some $\epsilon>0$ depending only on $\bm r$ and $s_1,\dots,s_m$. Thus we may assume that $\lfloor tB\rfloor=0$. By Lemma~\ref{lem: canonical threshold attain mld=1}(1), $\mld(X\ni x, tB)=1$. Since $t>1$, by Lemma~\ref{lem: bdd mld computing dcc cA/n type}, there exists a positive integer $N$ depending only on $\bm r$ and $s_1,\dots,s_m$, such that if $n>N$, then there exists a prime divisor $\bar{E}$ over $X\ni x$, such that $a(\bar E,X,tB)=1$ and $a(\bar E,X,0)=1+\frac{a}{n}$ for some positive integer $a\leq 3$. Since $$a(\bar E,X,tB)=a(\bar E,X,0)-\mult_{\bar E} tB=1+\frac{a}{n}-\frac{t}{n}\sum_{j=1}^m l_js_j(1,\bm r)=1,$$ 
where $l_j:=n\mult_{\bar E}B_j\in \Zz_{>0}$ for each $j$, we have $t\sum_{j=1}^m l_js_j(1,\bm r)=a$. Since $a,s_j(1,\bm{r})$ belong to a finite set of positive real numbers for any $j$, and $l_j$ belongs to a discrete set of positive real numbers, $t$ belongs to a set whose only accumulation point is $0$. Since $t>1$, there exists a positive real number $\epsilon$ depending only on $\bm r$ and $s_1,\dots,s_m$, such that $t>1+\epsilon$.
\end{proof}

\subsection{Accumulation points of canonical thresholds}

In this section, we prove Theorem~\ref{thm: Accumulation points}. %We introduce some notation first.

%Given two admissible weights (cf. Definition~\ref{def: Weights of threefolds terminal singularities}) $w=\frac{1}{n}(w_1,\dots,w_d)$ and $w'=\frac{1}{n}(w_1',\dots,w_d')$ of the quotient space $X\ni x\cong (\phi_1=\cdots=\phi_m=0)\ni o\subset\mathbb C^d/\frac{1}{n}(b_1,\dots,b_d)$, where $\{\phi_i\}_{1\leq i\leq m}$ are semi-invariant analytic power serieses (cf. Definition~\ref{def: Weights of threefolds terminal singularities}). Recall that $w'\succeq\mu w$ for some $\mu\in \Rr_{\geq 0}$ means that $w_i'\geq \mu w_i$ for each $1\leq i\leq d$. Let $f: Y\to X$ (resp. $f': Y'\to X$) be the weighted blow-up with the weight $w$ (resp. $w'$) at $x$. We assume that $f$ (resp. $f'$) extracts an analytic prime divisor $E$ (resp. $E'$) over $X\ni x$.

%For any $\Qq$-Cartier Weil divisor $D$ on $X$, we may write $K_{Y}=f^*K_X+\frac{a}{n}E$ and $f^*D=D_Y+\frac{m}{n}E$ (resp. $K_{Y'}=f'^*K_X+\frac{a'}{n}E'$ and $f'^*D=D_{Y'}+\frac{m'}{n}E'$) for some positive integers $a,a',m,m'$. We say that the weight $w$ \emph{calculates the canonical threshold $\ct(X\ni x,0; D)$} if $E$ is a canonical place (see Definition~\ref{defn: canonical place}) of $(X\ni x, \ct(X\ni x,0; D)D)$. Recall that by Definition~\ref{def: Weights of threefolds terminal singularities}, $a=w(X\ni x)$ and $a'=w'(X\ni x)$. We state the following simple lemma for the reader's convenience.

\begin{lem}[{Lemma~2.1]{Che19}}]\label{lem: weighted blow up comparison inequality}
Let $(X\ni x)\cong(\phi_1=\dots=\phi_m=0)\subset(\mathbb C^d\ni o)/\frac{1}{n}(b_1,\dots,b_d)$ be a germ, where $\phi_1,\dots,\phi_m$ are semi-invariant analytic power series. Let $w,w'\in\frac{1}{n}\mathbb Z^d_{>0}$ be two weights and $f: Y\rightarrow X, f': Y'\rightarrow X$ weighted blow-ups with the weights $w,w'$ at $x\in X$ respectively, such that $f$ extracts an analytic prime divisor $E$ and $f'$ extracts an analytic prime divisor $E'$ respectively.

Let $B\geq 0$ be a $\Qq$-Cartier Weil divisor on $X$ such that $1=a(E,X,\ct(X\ni x,0;D)D)$, $m:=n\mult_ED$, and $m':=n\mult_{E'}D$. Then for any real number $\mu\geq 0$ such that $w'\succeq\mu w$ (see Definition \ref{defn: compare two weights}),
$$\lceil\mu m\rceil \leq m'\leq\lfloor \frac{w'(X\ni x)}{w(X\ni x)}m\rfloor.$$
\end{lem}

\begin{lem}\label{lem: The accumulation point of ct reduce to 1-index cover}
Let $\mathcal{TS}$ be a set of terminal threefold singularities, and $$\mathcal{TS}_{1}:=\{(\tilde{x}\in \tilde{X})\mid \tilde{x}\in \tilde{X} \text{ is an index one cover of }(x\in X)\in \mathcal{TS}\}.$$ Then the set of accumulation points of $$\{\ct(X\ni x,0;D)\mid (x\in X)\in \mathcal{TS}, D\in \Zz_{> 0}\}$$ is a subset of the set of accumulation points of $$\{\ct(X\ni x,0;D)\mid (x\in X)\in \mathcal{TS}_1, D\in \Zz_{> 0}\}.$$
\end{lem}

\begin{proof}
Let $\{(X_i\ni x_i, D_i)\}_{i=1}^{\infty}$ be a sequence of germs, such that $(x_i\in X_i)\in \mathcal{TS}$ and $D_i\geq 0$ are non-zero $\Qq$-Cartier Weil divisors. Let $n_i$ be the index of the terminal singularity $x_i\in X_i$ for each $i$. By Theorem~\ref{thm: ACC ct local with base terminal}, we may assume that the sequence $\{c_i:=\ct(X_i\ni x_i,0;D_i)\}$ is strictly decreasing with the limit point $c\geq 0$. It suffices to show that $c$ is an accumulation point of $\{\ct(X\ni x,0;D)\mid (x\in X)\in \mathcal{TS}_1, D\in \Zz_{> 0}\}$. We may assume that $c>0$.

We may assume that $1>c_i$ for each $i$, and by Lemma~\ref{lem: canonical threshold attain mld=1}(1), $\mld(X_i\ni x_i,c_iD_i)=1$. For each $i$, consider the pair $(X_i\ni x_i, c_iD_i)$, by Lemma~\ref{lem: can extract divisor computing ct that is terminal strong version}, there exists a terminal blow-up (see Definition~\ref{defn: terminal blow-up}) of $(X_i\ni x_i,c_iD_i)$ which extracts a prime divisor $E_i$ over $X_i\ni x_i$. We may write $a(E_i,X_i,0)=1+\frac{a_i}{n_i}$ and $\mult_{E_i}D_i=\frac{m_i}{n_i}$ for some positive integers $a_i,m_i$. Set $t_i:=\lct(X_i,0; D_i)$ for each $i$, then we have $c_i=\frac{a_i}{m_i}$ and $t_i\le \frac{a_i+n_i}{m_i}$. Since $\{c_i\}_{i=1}^{\infty}$ is strictly decreasing and $c>0$, $\lim_{i\to +\infty} m_i=+\infty$ and $\lim_{i\to +\infty}a_i=+\infty$. In particular, possibly passing to a subsequence, we may assume that $a_i\geq 3$ for all $i$. Since $c_i>c$, by Lemma~\ref{lem: finiteness of description of terminal threefold singularities}(2) and \cite[(6.1) Theorem]{Rei87}, $n_i\leq \max\{4,\frac{3}{c}\}$ for all $i$.

Possibly shrinking $X_i$ to a neighborhood of $x_i$, we may assume that $(X_i,c_iD_i)$ is lc for each $i$. It follows that  $$c=\lim_{i\to +\infty} \frac{a_i}{m_i}\leq \lim_{i\to +\infty} t_i\leq \lim_{i\to +\infty} \frac{a_i+n_i}{m_i}=c,$$ hence $c=\lim_{i\to +\infty} t_i$. 

For each $i$, let $\pi_i: (\widetilde{X}_i\ni \tilde{x}_i)\to (X_i\ni x_i)$ be the index one cover of $x_i\in X_i$. Set $\widetilde{D}_i:=\pi_i^{-1}D_i$, $\widetilde{c}_i:=\ct(\widetilde{X}_i\ni \widetilde{x}_i,0;\widetilde{D}_i)$ and $\widetilde{t}_i:=\lct(\widetilde{X}_i,0;\widetilde{D}_i)$. Possibly shrinking $X_i$ to a neighborhood of $x_i$ again, we may assume that $(\tilde{X}_i,\tilde{c}_i\tilde{D}_i)$ is lc. By \cite[Proposition~5.20]{KM98}, $\widetilde{t}_i=t_i$ and $\widetilde{c}_i\geq c_i$. Now $$c=\lim_{i\to +\infty} t_i=\lim_{i\to +\infty} \widetilde{t}_i\geq \lim_{i\to +\infty} \widetilde{c}_i\geq \lim_{i\to +\infty} c_i=c,$$ which implies that $c=\lim_{i\to +\infty} \widetilde{c}_i$.
\end{proof}

\begin{thm}\label{thm: explicite accumulation points of canonical thresholds}
Let $\mathfrak{T}$ be the set of all terminal threefold singularities. Then the set of accumulation points of $$\mathcal{CT}_t:=\{\ct(X\ni x,0;D)\mid (x\in X)\in \mathfrak{T},\, D\in \Zz_{> 0}\}$$ is $\{0\}\cup\{\frac{1}{k}\mid k\in \Zz_{\geq 2}\}$. Moreover, $0$ is the only accumulation point of
      \[
		\mathcal{CT}_{t,\neq{ sm,cA/n}}:=\left\{\ct(X\ni x,0;D)\left| \begin{array}{l} \text{ $(x\in X)\in \mathfrak{T}$, $(x\in X)$ is neither smooth nor}\\ 
		\text{of type $cA/n$ for any $n\in\Zz_{>0}$}, D\in \Zz_{>0}, 
		\end{array}\right.
		\right\}.
     \]
\end{thm}

\begin{proof}
\noindent {\bf Step 0}. By \cite[Theorem~3.6]{Ste11}, $\{0\}\cup\{\frac{1}{k}\mid k\in \Zz_{\geq 2}\}$ is a subset of the set of accumulation points of $\mathcal{CT}_t$. For any $(x\in X)\in\mathfrak{T}$, let $D\geq 0$ be a non-zero $\Qq$-Cartier Weil divisor on of $X$. Then $\ct(X\ni x,0,kD)=\frac{1}{k}\ct(X\ni x, 0;D)$ for any positive integer $k$, hence $0$ is an accumulation point of $\mathcal{CT}_{t,\neq{ sm,cA/n}}$. 

It suffices to show the corresponding reverse inclusions. Let $c>0$ be an accumulation point of $\mathcal{CT}_t$. We will finish the proof by showing that $c$ is not an accumulation point of $\mathcal{CT}_{t,\neq{ sm,cA/n}}$ in {\bf Step 2} and $c=\frac{1}{k+1}$ for some positive integer $k$ in {\bf Step 3}.

\medskip

\noindent {\bf Step 1}.  By Theorem~\ref{thm: ACC ct local with base terminal}, $c<1$. Let $k$ be a positive integer such that $\frac{1}{k+1}\leq c<\frac{1}{k}$. Consider the set $$I_k:=\{\frac{p}{q}\mid p,q\in \Zz_{>0},p\leq 16(k+1)^2\},$$ which is discrete away from 0. By Theorem~\ref{thm: ACC ct local with base terminal}, there exists a positive real number $\epsilon$, such that for any $c'\in \mathcal{CT}_t$, if $0<|c-c'|< \epsilon$, then $c<c'<\frac{1}{k}$, and $c'\notin I_k$.

By Lemma \ref{lem: The accumulation point of ct reduce to 1-index cover}, there exists a Gorenstein terminal threefold singularity $x\in X$ and a non-zero $\Qq$-Cartier Weil divisor $D\geq 0$ on $X$, such that $0<|c-\ct(X\ni x,0; D)|<\epsilon$.  We have $c<\ct(X\ni x, 0; D)<\frac{1}{k}$ and $\ct(X\ni x, 0; D)\notin I_k$. Let $D_0:=\ct(X\ni x, 0; D)D$. By Lemma~\ref{lem: canonical threshold attain mld=1}(1), $\mld(X\ni x, D_0)=1$. By Lemma~\ref{lem: can extract divisor computing ct that is terminal strong version}, there exists a terminal blow-up $f: Y\to X$ of $(X\ni x,D_0)$ which extracts a prime divisor $E$ over $X\ni x$, and $K_Y=f^*K_X+aE$, $f^*D=f_*^{-1}D+mE$ for some positive integers $a,m$. We have $\ct(X\ni x, 0;D)=\frac{a}{m}$. If $a\leq 4$ or $a\mid m$, then $\frac{a}{m}\in I_k$, a contradiction. Hence $a\geq 5$ and $a\nmid m$. By Theorem~\ref{thm: kaw05 1.2 strenghthened}, $f$ is a divisorial contraction of ordinary type as in Theorem~\ref{thm: kaw05 1.2 strenghthened}(1-3) when $x\in X$ is not smooth or as in \cite[Theorem~1.1]{Kaw01} when $x\in X$ is smooth.

\medskip

\noindent {\bf Step 2}. We show that $x\in X$ is either smooth or of type $cA$ in this step. In particular, by Theorem~\ref{lem: The accumulation point of ct reduce to 1-index cover}, $c$ is not an accumulation point of $\mathcal{CT}_{t,\neq sm,cA/n}$.

Otherwise, by Theorem~\ref{thm: kaw05 1.2 strenghthened}, $X$ is of type $cD$, and there are two cases:

\medskip

\noindent {\bf Case 2.1}. $f$ is a divisorial contraction as in Theorem~\ref{thm: kaw05 1.2 strenghthened}(2.1). In particular, under suitable analytic local coordinates $x_1,x_2,x_3,x_4$, we have $$(X\ni x)\cong (\phi:= x_1^2+x_1q(x_3,x_4)+x_2^2x_4+\lambda x_2x_3^2+\mu x_3^3+p(x_2,x_3,x_4)=0)\subset (\mathbb C^4\ni o)$$ for some analytic power series $\phi$ as in Theorem~\ref{thm: kaw05 1.2 strenghthened}(2.1), and $f$ is a weighted blow-up with the weight $w:=(r+1,r,a,1)$, where $r$ is a positive integer, $2r+1=ad$ for some integer $d\geq 3$, and $a$ is an odd number. We have $w(X\ni x)=a$ and $w(D)=m$.

Since $a\nmid m$, $\ct(X\ni x,0;D)\in (\frac{1}{k+1},\frac{1}{k})$, and we have $$\frac{2r+1}{md}=\frac{a}{m}=\ct(X\ni x,0; D)\in (\frac{1}{k+1},\frac{1}{k}),$$ hence $k(2r+1)<dm<(k+1)(2r+1)$. Consider the weighted blow-up $f': Y'\to X$ (resp. $f'': Y''\to X$) at $x\in X$ with the weight $w':=(d,d,2,1)$ (resp. $w'':=(1+r-d,r-d,a-2,1)$). Since $a-2\geq 3$, by Lemma~\ref{lem: cD irreducible case 1}(1), $f''$ extracts an analytic prime divisor, and $w''(X\ni x)=a-2$. Since $a\geq 5$, by Lemma~\ref{lem: cD irreducible case 1}(2), $f'$ extracts an analytic prime divisor, and $w'(X\ni x)=2$. Since $w'\succeq \frac{d}{r+1}w$ and $w''\succeq \frac{r-d}{r}w$, by Lemma~\ref{lem: weighted blow up comparison inequality}, $\lfloor\frac{2}{a}m\rfloor\geq\lceil\frac{d}{r+1}m\rceil$ and $\lfloor\frac{a-2}{a}m\rfloor\geq\lceil\frac{r-d}{r}m\rceil$. Thus $$m-1=\lfloor\frac{2}{a}m\rfloor+\lfloor\frac{a-2}{a}m\rfloor\geq \lceil\frac{d}{r+1}m\rceil+\lceil\frac{r-d}{r}m\rceil\geq \lceil m-\frac{dm}{r(r+1)}\rceil,$$ where the first equality follows from $a\nmid 2m$ as $a\nmid m$ and $a$ is an odd number. It follows that $\frac{dm}{r(r+1)}\geq 1$. Hence $(k+1)(2r+1)>dm\geq r(r+1)>r(r+\frac{1}{2})$, and $r<2(k+1)$. Since $2r+1=ad$, $a\leq 4k+4$. Therefore, $\ct(X\ni x,0; D)=\frac{a}{m}\in I_k$, a contradiction.

\medskip
        
\noindent {\bf Case 2.2}. $f$ is a divisorial contraction as in Theorem \ref{thm: kaw05 1.2 strenghthened}(3.1). In particular, under suitable analytic local coordinates $x_1,x_2,x_3,x_4,x_5$, we have $$(X\ni x)\cong \binom{\phi_1:= x_1^2+x_2x_5+p(x_2,x_3,x_4)=0}{\phi_2:= x_2x_4+x_3^d+q(x_3,x_4)x_4+x_5=0}\subset(\mathbb C^5\ni o)$$ for some analytic power series $\phi_1,\phi_2$ as in Theorem~\ref{thm: kaw05 1.2 strenghthened}(3.1), and $f$ is a weighted blow-up with the weight $w:=(r+1,r,a,1,r+2)$, where $r$ is a positive integer such that $r+1=ad$ and $d\geq 2$ is an integer. We have $w(X\ni x)=a$ and $w(D)=m$.
        
Since $a\nmid m$, $\ct(X\ni x,0;D)\in (\frac{1}{k+1},\frac{1}{k})$, and we have $$\frac{r+1}{dm}=\frac{a}{m}=\ct(X\ni x,0; D)\in (\frac{1}{k+1},\frac{1}{k}),$$ hence $k(r+1)<dm<(k+1)(r+1)$. Consider the weighted blow-up $f': Y'\to X$ (resp. $f'': Y''\to X$) at $x\in X$ with the weight $w':=(d,d,1,1,d)$ (resp. $w'':=(r+1-d,r-d,a-1,1,r+2-d)$). Since $a-1\geq 4$, by Lemma~\ref{lem: cD irreducible case 2}(1), $f''$ extracts an analytic prime divisor, and $w''(X\ni x)=a-1$. Since $a\geq 5$, by Lemma~\ref{lem: cD irreducible case 2}(2), $f'$ extracts an analytic prime divisor, and $w'(X\ni x)=1$. Since $w'\succeq \frac{d}{r+2}w$ and $w''\succeq \frac{r-d}{r}w$, by Lemma~\ref{lem: weighted blow up comparison inequality}, $\lfloor\frac{1}{a}m\rfloor\geq\lceil\frac{d}{r+2}m\rceil$ and $\lfloor\frac{a-1}{a}m\rfloor\geq\lceil\frac{r-d}{r}m\rceil$. Thus $$m-1=\lfloor\frac{1}{a}m\rfloor+\lfloor\frac{a-1}{a}m\rfloor\geq \lceil\frac{d}{r+2}m\rceil+\lceil\frac{r-d}{r}m\rceil\geq \lceil m-\frac{2dm}{r(r+2)}\rceil,$$ where the first equality follows from $a\nmid m$. This implies that $\frac{2dm}{r(r+2)}\geq 1$. Hence $(k+1)(r+1)>dm\geq \frac{1}{2}r(r+2)>\frac{1}{2}r(r+1)$, and $r<2(k+1)$. Since $r+1=ad$, $a\leq 2k+2$. Therefore, $\ct(X\ni x,0; D)=\frac{a}{m}\in I_k$, a contradiction.

\medskip

\noindent {\bf Step 3}. We show that $c=\frac{1}{k+1}$ in this step. 

By \noindent {\bf Step 1}, \noindent {\bf Step 2}, \cite[Theorem~1.1]{Kaw01}, and Theorem~\ref{thm: kaw05 1.2 strenghthened}, there are two cases. 

\medskip

\noindent {\bf Case 3.1}. $x\in X$ is smooth, and under suitable analytic local coordinates $x_1,x_2,x_3$, $f$ is a weighted blow-up with the weight $w:=(1,r_1,r_2)$ for some positive integers $r_1,r_2$, such that $\gcd(r_1,r_2)=1$. Now $n=1$, $a=r_1+r_2$, and $\ct(X\ni x,0;D)=\frac{r_1+r_2}{m}$, such that $r_1+r_2\nmid m$. Possibly switching $x_2,x_3$, we may assume that $r_1\leq r_2$. By \cite[Proposition~3.3(1)]{Che19}, $$\ct(X\ni x,0;D)=\frac{r_1+r_2}{m}\leq \frac{1}{r_1}+\frac{1}{r_2}$$
when $r_1\geq 2$. When $r_1=1$, we have
$$\ct(X\ni x,0;D)<1<\frac{1}{r_1}+\frac{1}{r_2}.$$
Since $a\nmid m$, $\ct(X\ni x,0;D)\in (\frac{1}{k+1},\frac{1}{k})$, and we have $\frac{1}{r_1}+\frac{1}{r_2}>\frac{1}{k+1}$, and $r_1<2(k+1)$. If $k+2\leq r_1$, then $r_2< (k+1)(k+2)$ and $a=r_1+r_2\leq 16(k+1)^2$. It follows that $\ct(X\ni x,0;D)=\frac{r_1+r_2}{m}\in I_k$, a contradiction. Hence $1\leq r_1\leq k+1$.

Consider the weighted blow-up with the weight $w':=(1,r_1,r_2-1)$. This weighted blow-up extracts an analytic prime divisor $E'$ that is isomorphic to $\mathbf{P}(1,r_1,r_2-1)$, and $w'(X\ni x)=r_1+r_2-1$. Since $w'\succeq \frac{r_2-1}{r_2}w$, by Lemma~\ref{lem: weighted blow up comparison inequality}, $$m-(k+1)=\lfloor\frac{r_1+r_2-1}{r_1+r_2}m\rfloor\geq \lceil\frac{r_2-1}{r_2}m\rceil,$$ where the equality follows from $\frac{m}{r_1+r_2}\in (k,k+1)$. It follows that $\frac{m}{r_2}\geq k+1$ and $(k+1)r_2\leq m< (k+1)(r_1+r_2)$. Thus $\ct(X\ni x,0;D)$ belongs to the set $$\{\frac{r_1+r_2}{m}\mid r_1,r_2,m\in \mathbb{Z}, 1\leq r_1\leq k+1, r_1\leq r_2, (k+1)r_2\leq m< (k+1)(r_1+r_2)\},$$ which has only one accumulation point $\frac{1}{k+1}$.

\medskip

\noindent {\bf Case 3.2}. $x\in X$ is of type $cA$ and $f: Y\to X$ is a divisorial contraction of ordinary type as in Theorem~\ref{thm: kaw05 1.2 strenghthened}(1). In particular, under suitable analytic local coordinates $x_1,x_2,x_3,x_4$, we have $$(X\ni x)\cong (\phi:= x_1x_2+g(x_3,x_4)=0)\subset (\mathbb C^4\ni 0)$$ for some analytic power series $\phi$ as in Theorem~\ref{thm: kaw05 1.2 strenghthened}(1), and $f$ is a weighted blow-up with the weight $w:=(r_1,r_2,a,1)$, where $r_1,r_2,d$ are positive integers such that $r_1+r_2=ad$. We have $w(X\ni x)=a$ and $w(D)=m$. By \cite[Proposition 4.2]{Che19}, $$\ct(X\ni x,0;D)=\frac{a}{m}=\frac{r_1+r_2}{dm}\leq \frac{1}{r_1}+\frac{1}{r_2}.$$

Possibly switching $x_1,x_2$, we may assume that $r_1\leq r_2$. Since $a\mid m$, $\ct(X\ni x,0;D)\in (\frac{1}{k+1},\frac{1}{k})$, and we have $\frac{1}{r_1}+\frac{1}{r_2}>\frac{1}{k+1}$, hence $r_1<2(k+1)$. If $k+2\leq r_1$, then $r_2< (k+1)(k+2)$, hence $a\leq r_1+r_2\leq 16(k+1)^2$, which implies that $\ct(X\ni x,0;D)=\frac{a}{m}\in I_k$, a contradiction. Hence $1\leq r_1\leq k+1$.

Consider the weight $w':=(r_1,r_2-d,a-1,1)$. By Lemma~\ref{lem: irreducibility of exceptional divisors extracted by weighted blow up}, the weighted blow-up with the weight $w'$ extracts an analytic prime divisor, and $w'(X\ni x)=a-1$. Since $w'\succeq \frac{r_2-d}{r_2}w$, by Lemma~\ref{lem: weighted blow up comparison inequality}, $$m-(k+1)=\lfloor\frac{a-1}{a}m\rfloor\geq \lceil\frac{r_2-d}{r_2}m\rceil,$$ where the equality follows from $\frac{m}{a}\in (k,k+1)$. It follows that $\frac{dm}{r_2}\geq k+1$ and $(k+1)r_2\leq dm< (k+1)(r_1+r_2)$. Thus $\ct(X\ni x,0;D)$ belongs to the set $$\{\frac{r_1+r_2}{dm}\mid r_1,r_2,d,m\in \mathbb{Z}_{>0}, 1\leq r_1\leq k+1, r_1\leq r_2, (k+1)r_2\leq dm< (k+1)(r_1+r_2)\},$$ which has only one accumulation point $\frac{1}{k+1}$.
%To sum up, $\ct(X\ni x, 0;D)$ belongs to a set whose only accumulation point is $\frac{1}{k+1}$, hence $c=\frac{1}{k+1}$. 
\end{proof}

%Now we prove Theorem~\ref{thm: Accumulation points}.

\begin{proof}[Proof of Theorem~\ref{thm: Accumulation points}]
Let $X$ be a canonical threefold and $D\geq 0$ a non-zero $\Qq$-Cartier Weil divisor on $X$. Consider the pair $(X, D_0:=\ct(X,0; D)D)$. For any exceptional prime divisor $F$ over $X$ such that $a(F,X,0)=1$, we have $a(F,X,D_0)=1$. By \cite[Corollary~1.4.3]{BCHM10}, there exists a $\Qq$-factorial variety $X'$ and a birational morphism $g: X'\to X$ that exactly extracts all exceptional divisors $F$ such that $a(F,X,0)=1$. By construction, $X'$ is terminal and $K_{X'}+g_*^{-1}D_0=g^*(K_X+D_0)$, hence $\ct(X,0;D)=\ct(X',0;g_*^{-1}D)$. Possibly replacing $(X,D)$ with $(X',g_*^{-1}D)$, we may assume that $X$ is terminal.

Now either $\ct(X,0; D)=1$ or $t:=\ct(X,0; D)=\ct(X\ni x,0;D)<1$ for some point $x\in X$ of codimension $\geq 2$. If $\dim x=1$, then by Lemma~\ref{lem: Terminal blow up}, $\mld(X\ni x, tD)=2-\mult_x tD=1$, hence $t=\frac{1}{\mult_x D}\in \{\frac{1}{m}\mid m\in \Zz_{>0}\}$. If $\dim x=0$, by Theorem~\ref{thm: explicite accumulation points of canonical thresholds}, we are done.
\end{proof}

\subsection{Proof of Theorem~\ref{thm:  terminal mld acc}}
%In this section, we apply Lemmas \ref{lem: bounded index Mustata nakamura conjecture}, \ref{lem: can extract mld place smooth case}, \ref{lem: can extract mld place cA/n case}, and Theorems~\ref{thm: Nakamura Mustata Conjecture}, \ref{thm: ACC ct local with base terminal}, \ref{thm: uniform canonical polytope}, \ref{thm: explicite accumulation points of canonical thresholds} to prove that the ACC conjecture for mlds holds in the interval $[1,+\infty)$ for threefolds.

%, which can be viewed as a weak version of Musta\c{t}\v{a}-Nakamura's \han{gen MN conj proposed by HL20?}Conjecture on the divisors calculating the minimal log discrepancies (see \cite[Conjecture~1.1]{MN18}),
\begin{proof}[Proof of Theorem~\ref{thm:  terminal mld acc}]
%\noindent\textbf{Step 1}. Let $(X\ni x,B)$ be a threefold pair, such that $X$ is terminal and $B\in\Ii$. When $\dim x=1$, by Lemma~\ref{lem: Terminal blow up}, we have $\mld(X\ni x, B)=2-\mult_x B$, which belongs to an ACC set. Hence we may assume that $x\in X$ is a closed point.

%\medskip

\noindent\textbf{Step 1}. Suppose that Theorem~\ref{thm:  terminal mld acc} does not hold, then there exists a sequence of threefold pairs $\{(X_i\ni x_i, B_i)\}_{i=1}^{\infty}$, where $X_i$ is terminal and $B_i\in\Ii$ for each $i$, such that $\{\mld(X_i\ni x_i,B_i)\}_{i=1}^{\infty}\subset(1,+\infty)$ is strictly increasing. Possibly replacing $X_i$ with a small $\Qq$-factorialization, we may assume that $X_i$ is $\Qq$-factorial. If $\dim x_i=1$, then by Lemma~\ref{lem: Terminal blow up}, $\mld(X_i\ni x_i, B_i)=2-\mult_{x_i} B_i$, which belongs to an ACC set. Possibly passing to a subsequence, we may assume that $\dim x_i=0$ for each $i$. \cite[Theorem 0.1]{Amb99}, we may let $\beta:=\lim_{i\to+\infty}\mld(X_i\ni x_i,B_i)$. By Theorem \ref{thm: 3-dim terminal number of coefficients local}, possibly passing to a subsequence, there exists a non-negative integer $p$, such that $B_i:=\sum_{j=1}^{p}b_{i,j}B_{i,j}$ for each $i$, where $B_{i,j}$ are distinct prime divisors. Set $b_j:=\lim_{i\to +\infty}b_{i,j}$ for $1\leq j\leq p$ and $\bar B_{i}:=\sum_{j=1}^p b_j B_{i,j}$ for each $i$.

Let $n_i$ be the index of $X_i\ni x_i$. By \cite[Appendix, Theorem]{Sho92}, if $n_i\geq 2$, then there exists a prime divisor $F_i$ over $X_i\ni x_i$, such that $a(F_i,X_i,0)=1+\frac{1}{n_i}$. Thus $$1+\frac{1}{n_i}\geq a(F_i,X_i,B_i)\geq \mld(X_i\ni x_i,B_i)\geq \mld(X_1\ni x_1,B_1)>1,$$ and $n_i\le \frac{1}{\mld(X_1\ni x_1,B_1)-1}$. Hence, possibly passing to a subsequence, we may assume that there exists a positive integer $n$ such that $n_i=n$ for all $i$. By \cite[Lemma~5.1]{Kaw88}, $nD_i$ is Cartier near $x_i$ for any $\Qq$-Cartier Weil divisor $D_i$ on $X_i$ and for each $i$. 

By \cite[Theorem~0.1]{Amb99} and Theorem~\ref{thm: ACC ct local with base terminal}, $1\leq\mld(X_i\ni x_i)\le 3$. By \cite[Corollary~1.3]{Nak16}, $\{\mld(X_i\ni x_i, \bar B_i)\mid i\in \Zz_{>0}\}\subset [1,3]$ is a finite set. Possibly passing to a subsequence, we may assume that there exists a positive real number $\alpha\geq 1$, such that $\mld(X_i\ni x_i,\bar B_i)=\alpha<\beta$ for all $i$. %By Theorem~\ref{thm: ACC ct local with base terminal}, $\alpha\ge 1$.

\medskip

\noindent\textbf{Step 2}. In this step, we show that for each $i$, there exists a prime divisor $\bar{E}_i$ over $X_i\ni x_i$, such that 
$a(\bar{E}_i,X_i,\bar{B}_i)=\mld(X_i\ni x_i,\bar{B}_i)=\alpha$, and $a(\bar{E}_i,X_i,0)\le l$ for some positive real number $l$ depending only on $\{{b}_j\}_{j=1}^p$.

By Lemmas \ref{lem: bounded index Mustata nakamura conjecture}, \ref{lem: can extract mld place smooth case}, \ref{lem: can extract mld place cA/n case}, and Theorem \ref{thm: Nakamura Mustata Conjecture}, it suffices to show that if $\alpha>1$, then either $x_i\in X_i$ is smooth or is of $cA/n$ type for all but finitely many $i$.

Otherwise, $\alpha>1$, and possibly passsing to a subsequence, we may assume that $x_i\in X_i$ is neither smooth nor of $cA/n$ type for each $i$. Let $E_i$ be a prime divisor over $X_i\ni x_i$, such that $a(E_i,X_i,B_i)=\mld(X_i\ni x_i,B_i)$. By Theorem~\ref{thm: uniform canonical polytope}, there exist a positive integer $m$ depending only on $\{{b}_j\}_{j=1}^p$, and $\Qq$-Cartier $\Qq$-divisors $\bar B_i'>0$ on $X_i$, such that for each $i$,
\begin{enumerate}
    \item $m\bar B_{i}'$ is a Weil divisor,
    \item $\mld(X_i\ni x_i, \bar B_{i}')> 1$, and
    \item $a(E_i,X_i,\bar B_{i}')\le a(E_i,X_i,\bar{B}_i)<\mld(X_i\ni x_i,B_i)$.
\end{enumerate}
%In particular, $a(E_i,X_i,\bar B_{i}')\le a(E_i,X_i,\bar{B_i})\le a(E_i,X_i,B_i)=\mld(X_i\ni x_i,B_i)\le \beta$.\han{do we need more details?}
Since $\lim_{i\to +\infty }\mld(X_i\ni x_i,B_i)=\beta$, by \cite[Theorem~1.2]{Nak16}, possibly passing to a subsequence, we may assume that there exists a positive real number $\gamma$, such that $\alpha\le a(E_i,X_i,\bar{B}_i)=\gamma<\mld(X_1\ni x_1,B_1)$ for each $i$.

By Lemma \ref{lem: canonical threshold attain mld=1}(2), $t_i':=\ct(X_i\ni x_i,0;\bar B_{i}')>1$ for each $i$. We have
$$\gamma-(t_i'-1)\mult_{E_i}\bar B_i'\ge a(E_i,X_i,\bar B_i')-(t_i'-1)\mult_{E_i}\bar B_i'=a(E_i,X_i,t_i'\bar B_i')\ge 1,$$
which implies that 
$t_i'-1\le \frac{\gamma-1}{\mult_{E_i}\bar B_i'}$. Since 
$$\mult_{E_i}(\bar B_i-B_i)=a(E_i,X_i,B_i)-a(E_i,X_i,\bar B_i)
\ge\mld(X_1\ni x_1,B_1)-\gamma>0,$$
$\lim_{i\to +\infty}\mult_{E_i}\bar B_i=+\infty$. Thus $\lim_{i\to +\infty}\mult_{E_i}\bar B_i'=+\infty$ as $a(E_i,X_i,\bar B_i')\le a(E_i,X_i,\bar B_i)$. It follows that $\lim_{i\to +\infty} t_i'=1$. Hence $\ct(X_i\ni x_i,0;m\bar B_i')>\frac{1}{m}$, and $\lim_{i\to +\infty}\ct(X_i\ni x_i,0;m\bar B_i')=\frac{1}{m}$, which contradicts Theorem~\ref{thm: explicite accumulation points of canonical thresholds}. 

\medskip

\noindent\textbf{Step 3}.
By \textbf{Step 2}, $\mult_{\bar{E_i}}\bar {B_i}=a(\bar{E_i},X_i,0)-a(\bar{E_i},X_i,B_i)\le l-\beta$. Thus $\lim_{i\to +\infty}\mult_{\bar{E_i}}(\bar {B_i}-B_i)=0$. Hence
\begin{align*}
\alpha=&\mld(X_i\ni x_i,\bar{B_i})+\lim_{i\to+\infty}\mult_{\bar{E_i}}(\bar {B_i}-B_i)=\lim_{i\to +\infty}a(\bar{E_i},X_i,B_i)\\
\ge &\lim_{i\to +\infty}\mld(X_i\ni x_i,B_i)=\beta,
\end{align*}
a contradiction.
\end{proof}

\subsection{Proof of Theorem~\ref{thm: alct acc terminal threefold}}

\begin{proof}[Proof of Theorem~\ref{thm: alct acc terminal threefold}]
We follow the argument in \cite{BS10,HLS19}. Suppose that the theorem does not hold. Then there exist a sequence of threefold $a$-lc germs $(X_i\ni x_i,B_i)$ such that $X_i$ is terminal and $B_i\in\Ii$, and a strictly increasing sequence of positive real numbers $t_i$, such that for every $i$, there exists an $\Rr$-Cartier $\Rr$-divisor $D_i$ on $X_i$, such that $D_i\in\Ii'$ and $t_i=a$-$\lct(X_i\ni x_i,B_i;D_i)$.	It is clear that $t:=\lim_{i\rightarrow+\infty}t_i<+\infty$. Let $a_i:=\mld(X_i\ni x_i,B_i+tD_i)$. By Theorem \ref{thm: ACC ct local with base terminal}, possibly passing to a subsequence, we may assume that $a_i\ge 1$. Let $\{\epsilon_i\}_{i=1}^{\infty}$ be a strictly decreasing sequence which converges to $0$, such that $0<\epsilon_i<1$ and $t_i':=t_i+\epsilon_i(t-t_i)\in(t_i,t_{i+1})$ for any $i$. Then all the coefficients of $B_i+t_i'D_i$ belong to a DCC set. By Theorem \ref{thm:  terminal mld acc}, the sequence $\{\mld(X_i\ni x_i,B_i+t_i'D_i)\}_{i=1}^{\infty}$ satisfies the ACC. By the convexity of minimal log discrepancies, we have
	\begin{align*}
	a&> \mld(X_i\ni x_i,B_i+t_i'D_i)\\
	&=\mld(X_i\ni x_i,\frac{t_i'-t_i}{t-t_i}(B_i+tD_i)+\frac{t-t_i'}{t-t_i}(B_i+t_iD_i))\\
	&\geq\frac{t_i'-t_i}{t-t_i}\mld(X_i\ni x_i,B_i+tD_i)+\frac{t-t_i'}{t-t_i}\mld(X_i\ni x_i,B_i+t_iD_i)\\
	&\geq\frac{t_i'-t_i}{t-t_i}a_i+\frac{t-t_i'}{t-t_i}a=a-\frac{(t_i'-t_i)(a-a_i)}{t-t_i}\\
	&=a-\epsilon_i(a-a_i)\geq (1-\epsilon_i)a.
	\end{align*}
	Therefore, possibly passing to a subsequence, we may assume that $\mld(X_i\ni x_i,B_i+t_i'D_i)$ is strictly increasing and converges to $a$, which contradicts Theorem \ref{thm:  terminal mld acc}.
\end{proof}

\section{Boundedness of complements}\label{section: Boundedness of complements}

\subsection{Boundedness of indices for strictly canonical germs}\label{subsection: rr}

\begin{defn}
Let $(X\ni x,B)$ be pair. We say that $(X\ni x,B)$ is \emph{strictly canonical} if $\mld(X\ni x,B)=1$. We say that $(X,B)$ is \emph{strictly canonical} if $\mld(X,B)=1$.
\end{defn}

\begin{lem}\label{lem: vanishing of the high direct image}
Let $(X\ni x,B)$ be a strictly canonical germ such that $B$ is a $\Qq$-divisor. Let $f: Y\rightarrow X$ be a birational morphism which extracts a prime divisor $E$ over $X\ni x$ such that $a(E,X,B)=1$ and $-E$ is nef over $X$. Then the following holds.

 Let $\mathfrak{m}_x$ be the maximal ideal sheaf for $x\in X$, $m$ the smallest positive integer such that $mB$ is a Weil divisor, and $r$ the smallest positive integer such that $rm(K_X+B)$ is Cartier near $x$. Then for any $i\in \mathbb{Z}$, $$f_*\mathcal{O}_Y(im(K_Y+B_Y)-E)=\begin{cases}\mathfrak{m}_x\mathcal{O}_X(im(K_X+B)) & \text{if } r\mid i,\\
    \mathcal{O}_X(im(K_X+B)) & \text{if } r\nmid i,
    \end{cases}$$ where $K_Y+B_Y:=f^*(K_X+B)$.
\end{lem}

\begin{proof}
If $r\mid i$, then by the projection formula, $$f_*\mathcal{O}_Y(im(K_Y+B_Y)-E)=f_*\mathcal{O}_X(-E)\otimes \mathcal{O}_X(im(K_X+B))=\mathfrak{m}_x\mathcal{O}_X(im(K_X+B)),$$ where the last equality follows from $$f_*\mathcal{O}_Y(-E)(U)=\{u\in K(X)\mid ((u)-E)|_{f^{-1}(U)}\geq 0\}=\{u\in \mathcal{O}_X(U)\mid \mult_x (u)>0\},$$ where $U$ is an arbitrary open neighborhood of $x\in X$, $K(X)$ is the field of rational functions of $X$, and $(u)$ is the Cartier divisor defined by the rational function $u$.

If $r\nmid i$, then $f_*\mathcal{O}_Y(im(K_Y+B_Y)-E)(U)\subset\mathcal{O}_X(im(K_X+B))(U)$ for any open set $U\subset X$ as
\begin{align*}
    f_*\mathcal{O}_Y(im(K_Y+B_Y)-E)(U)&=\{u\in K(X)\mid ((u)+imf^*(K_X+B)-E)|_{f^{-1}(U)}\geq 0\}, \text{ and}\\
    \mathcal{O}_X(im(K_X+B))(U)&=\{u\in K(X)\mid ((u)+im(K_X+B))|_U\geq 0\}.
\end{align*}
  Suppose that $u\in K(X)$ satisfies $((u)+im(K_X+B))|_U\geq 0$. Since $r\nmid i$, $(u)+im(K_X+B)$ is not Cartier at $x$, so there exists a Weil divisor $D$ passing through $x$ such that $((u)+im(K_X+B)-D)|_U\geq 0$, which implies that $((u)+im(K_Y+B_Y)-f^*D)|_{f^{-1}(U)}\geq 0$. Since $E\subset \Supp (f^*D)$, we obtain $((u)+im(K_Y+B_Y)-E)|_{f^{-1}(U)}\geq 0$. Thus $f_*\mathcal{O}_Y(im(K_Y+B_Y)-E)=\mathcal{O}_X(im(K_X+B))$ in this case.
\end{proof}

\medskip

\noindent{\bf Notation} $(\star).$ Let $f: Y\to X$ be a divisorial contraction of a prime divisor $E$ over $X\ni x$ as in Theorem~\ref{thm: kaw05 1.2 strenghthened}(1) (see also \cite[Theorem~1.2(1)]{Kaw05}). Recall that in this case, $x\in X$ is a terminal singularity of type $cA/n$. In particular, under suitable analytic local coordinates $x_1,x_2,x_3,x_4$, $$(X\ni x)\cong (\phi:= x_1x_2+g(x_3^n,x_4)=0)\subset (\mathbb C^4\ni o)/\frac{1}{n}(1,-1,b,0),$$ where $b\in [1,n-1]\cap\Zz$ such that $\gcd(b,n)=1$ and $f$ is a weighted blow-up with the weight   $w=\frac{1}{n}(r_1,r_2,a,n)$ for some positive integers $a,r_1,r_2$, such that $an\mid r_1+r_2$ and $a\equiv br_1\mod n$.

Let $J'$ be the Reid basket for $f: Y\to X$ (see Definition~\ref{def: Reid basket for divisorial contractions}). By \cite[Theorem~1.2]{Kaw05}, we have three cases: $J'=\emptyset$, $J'=\{(r'_{Q'},1)_{Q'}\}$, or $J'=\{(r'_{Q'_1},1)_{Q'_1},(r'_{Q'_2},1)_{Q'_2}\}$, where $r'_{Q'},r'_{Q'_1},r'_{Q'_2}\in \mathbb{Z}_{\geq 2}$, and $Q'$, $Q'_1$, $Q'_2$ are fictitious singularities (see Definition-Lemma~\ref{def: definition of general cP}). In the case when $J'=\{(r'_{Q'_1},1)_{Q'_1},(r'_{Q'_2},1)_{Q'_2}\}$, $Q'_1,Q'_2$ come from two different non-Gorenstein points on $Y$. In the following, we introduce the set $J:=\{(r_{Q_1},1)_{Q_1},(r_{Q_2},1)_{Q_2}\}$ for $f$, here $Q_1,Q_2$ may not be fictitious cyclic quotient singularities any more as they could be smooth points. We let $Q_1,Q_2$ be any smooth points on $Y$ and $(r_{Q_1},r_{Q_2}):=(1,1)$ when $J'=\emptyset$, $Q_1$ any smooth closed point on $Y$, $Q_2=Q'$ and $(r_{Q_1},r_{Q_2}):=(1,r'_{Q'})$ when $J'=\{(r_{Q'},1)_{Q'}\}$, and $Q_1:=Q_1',Q_2:=Q_2'$, $(r_{Q_1},r_{Q_2}):=(r'_{Q'_1},r'_{Q'_2})$ when $J'=\{(r'_{Q'_1},1)_{Q'_1},(r'_{Q'_2},1)_{Q'_2}\}$.

\begin{lem}\label{lem: the blow up index coincide with the reid basket index}
With \noindent{\bf Notation} $(\star).$ Up to a permutation, we have $r_{Q_1}=r_1$ and $r_{Q_2}=r_2$. Moreover, $Q_1$, $Q_2$ are indeed singularities (possibly smooth) on $Y$.
\end{lem}

\begin{proof}
By \cite[Theorem~6.5, Page 112, Line 12-14]{Kaw05} and \cite[Proposition~2.15, Page 9, Line 15]{CH11}, there are two cyclic quotient terminal singularities $P_1,P_2\in Y$ of type $\frac{1}{r_1}(1,-1,b_1)$, $\frac{1}{r_2}(1,-1,b_2)$ respectively and possibly a $cA/n$ type singularity $P_3\in Y$. By \cite[Theorem~4.3]{Kaw05}, possibly changing the order of the indices, $Q_1, Q_2$ are $P_1,P_2$ on Y respectively. It follows that $r_{Q_1}=r_1$, and $r_{Q_2}=r_2$.
\end{proof}

\begin{lem}\label{lem: unboundness of r1 and r2}
Let $(X\ni x, B)$ be a threefold germ and $B$ a $\Qq$-divisor, such that $X$ is terminal and $\mld(X\ni x, B)=1$. Let $f: Y\to X$ be a divisorial contraction of a prime divisor $E$ over $X\ni x$ as in \noindent{\bf Notation} $(\star)$, such that $a(E,X,B)=1$. Let $m$ be the smallest positive integer such that $mB$ is a Weil divisor near $x$, and $r$ the smallest positive integer such that $rm(K_X+B)$ is Cartier near $x$. Then $r\mid \gcd(r_1,r_2)$.
\end{lem}

\begin{proof}
Possibly shrinking and compactifying $X$, we may assume that $X$ is projective and terminal. In particular, $Y$ is also projective and terminal. For each $r\in \Zz_{>0}$ and $i\in \Zz$, we define $$\delta_r(i):=\begin{cases}1 & \text{if } r\mid i,\\
0 & \text{if } r\nmid i.
\end{cases}$$ Let $D_{i,m}:=im(K_Y+B_Y)=imf^*(K_X+B)$ for each $i\in \Zz$. Since $E$ is $\Qq$-Cartier, by \cite[Proposition 5.26]{KM98}, we have the following short exact sequence
$$0\to\mathcal{O}_Y(D_{i,m}-E) \to\mathcal{O}_Y(D_{i,m})\to\mathcal{O}_E(D_{i,m}|_E) \to 0.$$

Since $D_{i,m}-E$ and $D_{i,m}$ are both $f$-big and $f$-nef, by the Kawamata-Viehweg vanishing theorem \cite[Theorem~1.2.5]{KMM87}, $R^jf_*\mathcal{O}_Y(D_{i,m}-E)=R^jf_*\mathcal{O}_Y(D_{i,m})=0$ for all $j\in\Zz_{>0}$. It follows that $h^j(\mathcal{O}_E(D_{i,m}|_E))=0$ for all $j\in\Zz_{>0}$. By Lemma~\ref{lem: vanishing of the high direct image}, 
\begin{align*}
    \delta_r(i)&=h^0(\mathcal{O}_X(im(K_X+B))/f_*\mathcal{O}_Y(im(K_Y+B_Y))-E))\\
    &=h^0(\mathcal{O}_E(D_{i,m}|_E))=\chi(\mathcal{O}_E(D_{i,m}|_E))\\
    &=\chi(\mathcal{O}_Y(D_{i,m}))-\chi(\mathcal{O}_Y(D_{i,m}-E)).
\end{align*}

For each fictitious singularity $Q\in Y_Q$ of some closed point on $Y$, $$(D_{1,m})_Q\sim d_QK_{Y_Q} \text{  and  } E_Q\sim f_QK_{Y_Q}$$ near $Q\in Y_Q$ for some integers $f_Q, d_Q\in [1,r_Q]$, where $Y_Q$ is the deformed variety on which $Q$ appears as a cyclic quotient terminal singularity, and $E_Q, (D_{i,m})_Q$ are the corresponding deformed divisors on $Y_Q$ (see Definition-Lemma~\ref{def: definition of general cP}). By Theorem~\ref{thm: singular rr},
$$
\delta_r(i)=\chi(\mathcal{O}_Y(D_{i,m}))-\chi(\mathcal{O}_Y(D_{i,m}-E))=\Delta_1+\Delta_2+\frac{1}{12}E\cdot c_2(X)$$
with 
$$\Delta_1=T(D_{i,m})-T(D_{i,m}-E),\Delta_2=\sum_{y\in Y,\dim y=0}(c_y(D_{i,m})-c_y(D_{i,m}-E)),$$
and
$$T(D)=\frac{1}{12}D(D-K_Y)(2D-K_Y).$$
Since $D_{i,m}\cdot E^2=D_{i,m}^2\cdot E=D_{i,m}\cdot E\cdot K_Y=0$, $\Delta_1=\frac{1}{6}E^3+\frac{1}{4}E^2\cdot K_Y$. For any fictitious point $Q$, if $E_Q$ is Cartier, then $c_{Q}((D_{i,m})_{Q})=c_{Q}((D_{i,m}-E)_{Q})$. By Definition-Lemma \ref{def: definition of general cP} and Definition~\ref{def: Reid basket for divisorial contractions},
\begin{align*}
    \Delta_2&=\sum_{Q\in J} c_Q(D_{i,m})-c_Q(D_{i,m}-E)=\sum_{Q\in J}(A_Q(id_Q)-A_Q(id_Q-f_Q)),
\end{align*}
where $J$ is defined as in \noindent{\bf Notation} $(\star)$ for $f:Y\to X$, and $$A_Q(i):=-i\frac{r^2_Q-1}{12r_Q}+\sum_{j=1}^{i-1}\frac{\overline{(jb_Q)}_{r_Q}(r_Q-\overline{(jb_Q)}_{r_Q})}{2r_Q}.$$ By $Q\in J$, we mean $Q$ is a fictitious singularity that contributes to $J$.
Here we allow $i<0$ if we adopt the notation of generalized summation (see Definition~\ref{defn: generalized sum}). It is worthwhile to mention that $A_Q(i)=A_Q(\overline{(i)}_{r_Q})$ as $\gcd(b_Q,r_Q)=1$.
Now
\begin{equation}\label{eqn: RR for delta}
   \delta_r(i+1)-\delta_r(i)=\sum_{Q\in J}\sum_{j=id_Q}^{(i+1)d_Q-1}(B_Q(jb_Q)-B_Q(jb_Q-v_Q)).
\end{equation}
Here, for each $Q\in J$, $B_{Q}(i)$ is an even periodic function with period $r_{Q}$ defined by $$B_{Q}(i):=\frac{\overline{(i)}_{r_Q}(r_Q-\overline{(i)}_{r_Q})}{2r_Q}.$$ 

By Lemma~\ref{lem: the blow up index coincide with the reid basket index}, $J=\{(r_1,1)_{Q_1},(r_2,1)_{Q_2}\}$, where $Q_1,Q_2$ are cyclic quotient terminal singularities (might be smooth) on $Y$. It follows that
\begin{equation}\label{eqn: subscription formula for delta function}
   \delta_r(i+1)-\delta_r(i)=\sum_{k=1,2}\sum_{j=id_{Q_k}}^{(i+1)d_{Q_k}-1}(B_{Q_k}(jb_{Q_k})-B_{Q_k}(jb_{Q_k}-1)).
\end{equation}

\begin{claim}\label{claim: description of the index $r$}
We have the following equality:
$$r=\mathrm{lcm}\{ \frac{r_1}{\gcd(r_1,d_{Q_1})}, \frac{r_2}{\gcd(r_2,d_{Q_2})}\}.$$
\end{claim}

We proceed the proof assuming Claim~\ref{claim: description of the index $r$}. Let $l_1:=\frac{r_1}{\gcd(r_1,r_2)}$ and $l_2:=\frac{r_2}{\gcd(r_1,r_2)}$. We have $\gcd(l_1l_2,l_1+l_2)=1$ as $\gcd(l_1,l_2)=1$. By Claim~\ref{claim: description of the index $r$}, $r\mid  \lcm(r_1,r_2)=\gcd(r_1,r_2)l_1l_2$. Since $r\mid n$ and $n\mid r_1+r_2$, $r\mid \gcd(r_1,r_2)(l_1+l_2)$. Hence $r\mid \gcd(r_1,r_2)$.
\end{proof}

\begin{proof}[Proof of Claim~\ref{claim: description of the index $r$}]
Let $\lambda\in \Zz$ such that $r_{Q_k}\mid \lambda d_{Q_k}$ for $k=1,2$. By (\ref{eqn: subscription formula for delta function}),
\begin{align*}
    \delta_r(\lambda+1)-\delta_r(\lambda)&=\sum_{k=1,2}\sum_{j=\lambda d_{Q_k}}^{(\lambda+1)d_{Q_k}-1}(B_{Q_k}(jb_{Q_k})-B_{Q_k}(jb_{Q_k}-1))\\
    &=\sum_{k=1,2}\sum_{j=0}^{d_{Q_k}-1}(B_{Q_k}(jb_{Q_k})-B_{Q_k}(jb_{Q_k}-1))=\delta_r(1)-\delta_r(0).
\end{align*}
Thus $r\mid \lambda$. By Lemma~\ref{lem: the blow up index coincide with the reid basket index}, $r\mid \mathrm{lcm}\{ \frac{r_1}{\gcd(r_1,d_{Q_1})}, \frac{r_2}{\gcd(r_2,d_{Q_2})}\}$. 

Since $rD_{1,m}=rm(K_Y+B_Y)$ is Cartier, $r_{k}\mid rd_{Q_k}$ for $k=1,2$. Thus $\frac{r_{k}}{\gcd(r_{k},d_{Q_k})}\mid r$ for $k=1,2$, and the claim is proved.
\end{proof}

%\han{submit version may be removed}
\begin{rem}
Let $n=r(4r^2-2r-1)$,
$a=r$, $b=4r^2+2r-1$,
$r_1=r_{Q_1}=(2r-1)^2r^2,d_{Q_1}=(2r-1)^2r^2,b_{Q_1}=4r^3-r+1$, $r_2=r_{Q_2}=2r^2(r-1)$, $d_{Q_2}=2r(r-1)$, and $b_{Q_2}=2r^2-1$. Then $(n,a,b,r_{Q_1},d_{Q_1},b_{Q_1},r_{Q_2},d_{Q_2},b_{Q_2})$ satisfies both \eqref{eqn: RR for delta} and Claim \ref{claim: description of the index $r$}. Moreover, as $r\mid n$, $\gcd(b,n)=1$, $n\mid a-br_{Q_1}$, $an\mid r_{Q_1}+r_{Q_2}$, and $\gcd(\frac{a-br_{Q_1}}{n},r_{Q_1})=1$, $(n,a,b,r_{Q_1},d_{Q_1},b_{Q_1},r_{Q_2},d_{Q_2},b_{Q_2})$ also satisfies the restrictions proved in \cite[Theorem~1.2(1)]{Kaw05}. Hence we could not show Theorem \ref{thm: index conjecture threefold} by simply applying singular Riemann-Roch formula for terminal threefold as \cite{Kaw15a} did for the case when $B=0$, $X$ is canonical and $x$ is an isolated canonical center of $X$.
\end{rem}

\begin{thm}\label{thm: index conjecture threefold}
Let $\Gamma\subset[0,1]$ be a set and $m$ a positive integer such that $m\Gamma\subset \mathbb{Z}$. Then the positive integer $N:=12m^2$ satisfies the following.

Let $(X\ni x, B)$ be a threefold germ such that $X$ is terminal, $B\in \Ii$, and $\mld(X\ni x, B)=1$. Then $I(K_X+B)$ is Cartier near $x$ for some positive integer $I\le N$. 
\end{thm}

\begin{proof}
Let $r$ be the smallest positive integer such that $rm(K_X+B)$ is Cartier near $x$.

\medskip

By \cite[(6.1) Theorem]{Rei87}, if $x\in X$ is a terminal singularity of types other than $cA/n$, then the index of $x\in X$ divides $12$. By \cite[Lemma~5.1]{Kaw88}, $12m(K_X+B)$ is Cartier near $x$. From now on, we may assume that $x\in X$ is a terminal singularity of type $cA/n$.

\medskip

By Lemma~\ref{lem: can extract divisor computing ct that is terminal strong version}, there exists a terminal blow-up (see Definition \ref{defn: terminal blow-up}) $f: Y\to X$ of $(X\ni x,B)$ which extracts a prime divisor $E$ over $X\ni x$. By \cite[Theorem~1.1]{Kaw05}, $f$ is either of ordinary type or of exceptional type.

If $f$ is of exceptional type, then by \cite[Theorem~1.3]{Kaw05}, $x\in X$ is a terminal singularity of type $cA$. Hence by \cite[Lemma~5.1]{Kaw88}, $m(K_X+B)$ is Cartier.

We may now assume that $f$ is of ordinary type, and we write $f^*K_X+\frac{a}{n}E=K_Y$ for some positive integer $a\geq 1$. Now $f: Y\to X$ is a divisorial contraction of ordinary type as in Theorem~\ref{thm: kaw05 1.2 strenghthened}(1). In particular, under suitable analytic local coordinates $x_1,x_2,x_3,x_4$, $$(X\ni x)\cong (\phi:= x_1x_2+g(x_3^n,x_4)=0)\subset (\mathbb C^4\ni o)/\frac{1}{n}(1,-1,b,0),$$ where $b\in [1,n-1]\cap \Zz$, $\gcd(b,n)=1$, and $f$ is a weighted blow-up at $x\in X$ with the weight $w:=\frac{1}{n}(r_1,r_2,a,n)$. Now $mB$ is a Weil divisor locally defined by a semi-invariant analytic power series $(h(x_1,x_2,x_3,x_4)=0)$.

\medskip

\begin{claim}\label{claim: control the index or finding an elephant in B}
Either $r\leq 3m$, or $x_3^m\in h$ (up to a scaling of $h$).
\end{claim}

We proceed the proof assuming Claim~\ref{claim: control the index or finding an elephant in B}. If $r\leq 3m$, then $Im(K_X+B)$ is Cartier for some $I\leq 3m$. Otherwise, $r>3m$. By Claim~\ref{claim: control the index or finding an elephant in B}, up to a scaling of $h$, we have $h=x_3^m+p$ for some analytic power series $p$ such that $\lambda x_3^m\notin p$ for any $\lambda\in \Cc^*$. Recall that $\xi_n$ is the primitive $n$-th root of unity. Since $h=x_3^m+p$ is semi-invariant with respect to the $\xi_n$-action: $(x_1,x_2,x_3,x_4)\rightarrow (\xi_nx_1,\xi_n^{-1}x_2,\xi_n^{b}x_3,x_4)$, $\xi_n(h)/h=\xi_n^{mb}$. Since $\xi_n(x_3^m)/x_3^m=\xi_n^{mb}$, $\xi_n(\frac{h}{x_3^m})=\frac{h}{x_3^m}$, and $(\frac{h}{x_3^m})$ is $\xi_n$-invariant, hence $(\frac{h}{x_3^m})$ is a rational function on $X$ which defines a principle divisor. Now $mB=(x_3^m+p=0)\sim (x_3^m=0)$ near $x$. Let $S$ be the analytic Cartier divisor locally defined by $(x_3=0)$ on $X$. By \cite[(6.4)(B.1)]{Rei87} and Lemma~\ref{lem: analytic index and algebraic index}, $K_X+S$ is Cartier near $x$. It follows that $m(K_X+B)\sim m(K_X+S)$ is Cartier near $x$.
\end{proof}

\begin{proof}[Proof of Claim~\ref{claim: control the index or finding an elephant in B}]
Assume that $r>3m$. By Lemma~\ref{lem: unboundness of r1 and r2}, $r_1> 3m$ and $r_2> 3m$. Note also that $n\geq r>3m$. When $a\leq 2$,  $$w(h)=mw(B)=mw(X\ni x)=\frac{am}{n}< \frac{3m}{n}.$$ Since $w(x_k)=\frac{r_k}{n}> \frac{3m}{n}$ for $k=1,2$ and $w(x_4)=1> \frac{3m}{n}$, up to a scaling of $h$, $x_3^l\in h$ for some $l\in \Zz_{>0}$ and $w(h)=\frac{la}{n}$. Thus $l=m$, and the claim follows in this case.

\medskip

When $a\geq 3$, we can pick positive integers $s_1,s_2$ such that
\begin{itemize}
    \item $s_1+s_2=3dn$,
    \item $3\equiv bs_1\mod n$, and
    \item $s_1,s_2>n$.
\end{itemize}
Let $\bar w:=\frac{1}{n}(s_1,s_2,3,n)$. Since $a\geq 3$, by Lemma~\ref{lem: irreducibility of exceptional divisors extracted by weighted blow up}, the weighted blow-up with the weight $\bar w$ extracts a prime analytic divisor $\bar E$ such that $\bar w(X\ni x)=\frac{3}{n}$. By Lemma~\ref{lem: algebraic approximation of weighted blow-up}, we may assume that $E$ is a prime divisor over $X\ni x$. Since $\mld(X\ni x,B)=1$, $a(\bar{E},X,B)=1+\bar{w}(X\ni x)-\bar{w}(B)\geq 1$, thus
$$1>\frac{3m}{n}=m\bar w(X\ni x)\geq m\bar{w}(B)=\bar w(h).$$
Since $\bar w(x_1)=\frac{s_1}{n}>1$, $\bar w(x_2)=\frac{s_2}{n}>1$, and $\bar w(x_4)=1$, there exists a positive integer $l$, such that up to a scaling of $h$, $x_3^l\in h$ and $\bar{w}(h)=\bar w(x_3^l)=\frac{3l}{n}$. Since
$\frac{3m}{n}\geq\bar w(h)=\frac{3l}{n}$, $l\leq m$. On the other hand,
$$\frac{am}{n}=mw(X\ni x)=w(mB)=w(h)\leq w(x_3^l)=\frac{al}{n},$$
which implies that $l\geq m$. Thus $l=m$, and $x_3^m\in h$ up to a scaling of $h$.
\end{proof}

\begin{proof}[Proof of Theorem~\ref{thm: index conjecture DCC threefold}]
If $\dim x=2$, then the theorem is trivial. If $\dim x=1$, then $X$ is smooth near $x$. By Lemma~\ref{lem: Terminal blow up}, $\mld(X\ni x, B)=2-\mult_{x} B=\epsilon$, hence the coefficients of $B$ belong to a finite set of rational numbers depending only on $\epsilon$ and $\Ii$, and the theorem holds in this case. Thus we may assume that $\dim x=0$ and $X$ is not smooth at $x$.

Let $B_1$ be any component of $B$ with coefficient $b_1$. Let $f:Y\to X$ be a small $\Qq$-factorialization, and $B_{1,Y}$ the strict transform of $B_1$ on $Y$. We may write $K_Y+B_Y:=f^{*}(K_X+B)$. Let $Y\dashrightarrow Z$ be the canonical model of $(Y,B_Y-b_1B_{1,Y})$ over $X$. Let $B_{Z}$ and $B_{1,Z}$ be the strict transforms of $B_Y$ and $B_{1,Y}$ on $Z$ respectively. Since $-B_{1,Z}$ is ample over $X$, $\Supp B_{1,Z}$ contains $g^{-1}(x)$, where $g:Z\to X$ is the natural induced morphism. Moreover, since $g$ is small, $K_Z=g^{*}K_X$ and $Z$ is terminal.

Since $\mld(X\ni x,B)=\epsilon$, there exists a point $z\in g^{-1}(x)$ such that $\mld(Z\ni z,B_Z)=\epsilon$. Then $$b_1=\epsilon\text{-}\lct(Z\ni z,B_Z-b_1B_{1,Z};B_{1,Z})\in \Ii.$$ By Theorem~\ref{thm: alct acc terminal threefold}, $b_1$ belongs to a finite set depending only on $\Ii$. Hence we may assume that the coefficients of $B$ belong to a finite set $\Ii'$ depending only on $\Ii$.

If $\epsilon=1$, then Theorem~\ref{thm: index conjecture DCC threefold} follows from Theorem~\ref{thm: index conjecture threefold}. If $\epsilon>1$, then we let $n$ be the index of $X\ni x$. By \cite[Appendix, Theorem]{Sho92} and \cite[Theorem~0.1]{Mar96}, there exists a prime divisor $E$ over $X\ni x$ such that $a(E,X,0)=1+\frac{1}{n}$. Since $1+\frac{1}{n}\geq a(E,X,B)\geq \epsilon$, $n\leq \frac{1}{\epsilon-1}$, and the theorem follows from \cite[Lemma 5.1]{Kaw88}.
\end{proof}

\subsection{Boundedness of complements for finite rational coefficients}
We prove Theorem~\ref{thm: N-complements for canonical pair finite rational coeff} in this section, and Theorem~\ref{thm: intro cc} follows as a direct corollary. Also, for any $\epsilon\geq 1$, we prove the existence of $(\epsilon,N)$-complements for terminal pairs (see Theorem~\ref{thm: epsilon-N-complements for canonical pair finite rational coeff}).

%In this subsection, when we consider a singularity $x\in X$, we always assume that $X$ is an affine variety. Let $D>0$ be a Weil divisor on $X$. By $|D|$ we mean the set of all effective Weil divisors on $X$ that is linearly equivalent to $D$.\luo{affine, add}

% such that $x\in \Supp D$
\begin{lem}\label{lem: linear system base locus x}
Let $x\in X$ be an isolated singularity such that $X$ is affine. Let $D\geq 0$ be a Weil divisor on $X$ and $E$ a prime divisor over $X\ni x$. Then there exists a finite dimensional linear system ${{\mathfrak {d}}}\subset |D|$, such that
\begin{enumerate}
    \item ${{\mathfrak {d}}}$ contains $D$,
    \item the base locus of ${{\mathfrak {d}}}$ is $x$, and
    \item if $D$ is $\Qq$-Cartier, then $\mult_E D'\geq \mult_E D$ for any $D'\in \mathfrak{d}$.
\end{enumerate}
In particular, if $x\in X$ is a terminal threefold singularity such that $X$ is affine, then there exists a finite dimensional linear system $\mathfrak{d}\subset |-K_X|$ such that $\mathfrak{d}$ contains an elephant (cf. \cite[(6.4)(B)]{Rei87}) of $x\in X$, and the base locus of ${{\mathfrak {d}}}$ is $x$.
\end{lem}
\begin{proof}
Let $\mathcal{I}_{E}$ be the ideal sheaf on $X$ such that for any open set $x\in U\subset X$,
$$\mathcal{I}_{E}(U)=\begin{cases}\{u\in \mathcal{O}_X(U)\mid \mult_E (u)>0\}& \text{if $D$ is not $\Qq$-Cartier},\\
\{u\in \mathcal{O}_X(U)\mid \mult_E (u)\geq \mult_ED\} & \text{if $D$ is $\Qq$-Cartier},
\end{cases}$$
where $(u)$ is the Cartier divisor defined by the rational function $u$, and $\mathcal{I}_E(U)=\mathcal{O}_X(U)$ for any open set $x\notin U\subset X$. Since $X$ is affine, the coherent sheaf $\mathcal{O}_X(D)\otimes \mathcal{I}_{E}$ is globally generated. Since $\mathcal{O}_X(D)\otimes \mathcal{I}_{E}|_{X\backslash\{x\}}=\mathcal{O}_X(D)|_{X\backslash\{x\}}$ is an invertible sheaf, there exist finitely many sections $s_1,\ldots,s_m\in H^0(X,\mathcal{O}_X(D)\otimes \mathcal{I}_{E})$, such that the linear system defined by $s_1,\ldots,s_m$ is base point free on $X\backslash\{x\}$. Pick $s_0\in H^0(X,\mathcal{O}_X(D)\otimes \mathcal{I}_{E})$ such that $(s_0)=0$. Then the linear system ${{\mathfrak {d}}}$ defined by $s_0,\ldots,s_m$ satisfies our requirements.
\end{proof}

\begin{lem}\label{lem: ct is attained by elephant at x}
Let $(X\ni x,B)$ be a canonical threefold germ such that $X$ is affine terminal but not smooth, $m$ a positive integer such that $mB\in \Zz$, and $\mathfrak{d}$ a finite dimensional linear system whose base locus is $x$. Then for any integer $N>m$ and any general element $(D_1,\dots,D_N)\in \mathfrak{d}^N$, the divisor $D:=\sum_{i=1}^N D_i$ satisfies that $\ct(X,B;D)=\ct(X\ni x, B;D)$.
\end{lem}

\begin{proof}
Let $f:Y\to X$ be a log resolution of $(X,B)$ such that
\begin{itemize}
    \item $f^{*}|{{\mathfrak {d}}}|=F+|M|$, and
    \item $\Supp B_Y\cup \Supp F\cup \Supp(\Exc(f))$ is snc,
\end{itemize}
where $B_Y$ is the strict transform of $B$ on $Y$, $F$ is the fixed part of $f^{*}|{{\mathfrak {d}}}|$, and $M$ is a base point free Cartier divisor. Note that $\Supp F\subset f^{-1}(x)$ as ${{\mathfrak {d}}}$ is base point free on $X\setminus\{x\}$. Let $M':=\sum_{i=1}^N M_i$, where $M_1,\cdots,M_N$ are $N$ general elements in $|M|$. Set $D:=f_{*}(NF+M')$. Then $D=\sum_{i=1}^N D_i$, where $D_i:=f_*(F+M_i)$ for each $i$, and $(D_1,\dots,D_N)\in \mathfrak{d}^N$ is a general element. Since $D$ has $N$ distinct components, by Theorem~\ref{thm: 3-dim terminal number of coefficients local}(2), $t:=\ct(X,B;D)\leq \frac{1}{N}$.

If $\lfloor B\rfloor\neq 0$, then $B=\lfloor B\rfloor$ is a $\Qq$-Cartier prime divisor and $t=0$. By Theorem~\ref{thm: 3-dim terminal number of coefficients local}(2), $\mld(X\ni x, B)=1$. Hence $\ct(X\ni x, B;D)=\ct(X,B;D)=0$ in this case.

If $\lfloor B\rfloor= 0$, then $B\in (0,1-\frac{1}{m}]$. Since $K_Y+B_Y+tM'+tNF=f^*(K_X+B+tD)+G$ for some $\Qq$-divisor $G\geq 0$, $\Supp M' \cup \Supp B_Y$ is snc, and $t\leq \frac{1}{N}<\frac{1}{m}$, we have $(Y,B_Y+tM'+tNF)$ is terminal on $Y\setminus\Supp F$. Hence $\ct(X,B;D)=\ct(X\ni x, B;D)$ as $\Supp F$ contains at least one canonical place of $(X,B+tD)$.
\end{proof}

\begin{lem}\label{lem: bounded index strict canonical complement}
Let $I$ be a positive integer and $\Gamma\subset [0,1]\cap \Qq$ a finite set. Then there exists a positive integer $N$ depending only on $I$ and $\Ii$ satisfying the following.

Let $(X\ni x,B)$ be a threefold germ such that $X$ is terminal, $B\in \Ii$, $\mld(X\ni x, B)\geq 1$, and $IK_X$ is Cartier. Then there exists a monotonic $N$-complement $(X\ni x, B^+)$ of $(X\ni x, B)$ such that $\mld(X\ni x, B^+)=1$, and if $(X\ni x,B)$ is canonical near $x$ and $x\in X$ is not smooth, then $(X\ni x, B^+)$ is canonical near $x$.
\end{lem}

\begin{proof}
Possibly shrinking $X$ near $x$, we may assume that $(X,B)$ is lc and $X$ is affine. For any positive real number $\epsilon<1$, since $(X\ni x,(1-\epsilon)B)$ is a klt germ, by \cite[Lemmas 3.12 and 3.13]{HLS19}, there exists a $\Qq$-factorial weak plt blow-up $f_{\epsilon}:Y_{\epsilon}\to X$ of $(X\ni x,(1-\epsilon)B)$, such that $(Y_{\epsilon},(1-\epsilon)B_{Y_{\epsilon}}+E_{Y_{\epsilon}})$ is $\Qq$-factorial plt, where $B_{Y_{\epsilon}}$ is the strict transform of $B$ on $Y_{\epsilon}$, and $E_{Y_{\epsilon}}$ is the reduced exceptional divisor of $f_{\epsilon}$. By \cite[Theorem 1.1]{HMX14}, we may choose $\epsilon<1$ such that $(Y_{\epsilon},B_{Y_{\epsilon}}+E_{Y_{\epsilon}})$ is lc. Let $Y:=Y_{\epsilon},B_Y:=B_{Y_{\epsilon}}$, $E:=E_{Y_{\epsilon}}$, and $f:=f_{\epsilon}$. By \cite[Corollary 3.20]{HLS19}, there exists a $\Qq$-divisor $G_Y\ge 0$ on $Y$, such that $N'G_Y$ is a Weil divisor, and $(Y/X\ni x,B_Y+E+G_Y)$ is an $N'$-complement of $(Y/X\ni x,B_Y+E)$ for some positive integer $N'$ depending only on $\Ii$. Then $K_Y+B_Y+G_Y+E=f^*(K_X+B+G)$, $(X\ni x,B+G)$ is an $N'$-complement of $(X\ni x,B)$, and $\mld(X\ni x,B+G)=0$, where $G:=f_{*}G_Y$.

Let $m$ be a positive integer such that $m\Ii\subset \Zz$ and $N'':=m(m+1)N'$. Since $mN'G\in |-mN'(K_X+B)|$ near $x$, by Lemma~\ref{lem: linear system base locus x}, there exists a finite dimensional linear system $\mathfrak{d}\subset |-mN'(K_X+B)|$ such that $\mathfrak{d}$ contains $mN'G$, the base locus of ${{\mathfrak {d}}}$ is $x$, and $\mult_E G'\geq \mult_E mN'G$ for any $G'\in \mathfrak{d}$. By Lemma~\ref{lem: ct is attained by elephant at x}, when $(X,B)$ is canonical and $x\in X$ is not smooth, the divisor $G'':=G_1'+\dots+G_{m+1}'\in |-N''(K_X+B)|$ satisfies that $t:=\ct(X,B;\frac{1}{N''}G'')=\ct(X\ni x,B;\frac{1}{N''}G'')$, where $(G_1',\dots,G_{m+1}')\in \mathfrak{d}^{m+1}$ is a general element. By construction, $(X\ni x,B+\frac{1}{N''}G'')$ is an $N''$-complement of $(X\ni x,B)$. Since $\mult_E \frac{1}{N''}G''\geq \mult_E G$, $\mld(X\ni x,B+\frac{1}{N''}G'')=0$. In particular, $t<1$.

By Theorem~\ref{thm: index conjecture threefold}, it suffices to show that $t$ belongs to a finite set of rational numbers depending only on  $I$ and $\Ii$. By Lemma \ref{lem: canonical threshold attain mld=1}(1), there exists a prime divisor $F$ over $X\ni x$ such that $a(F,X,B+\frac{t}{N''}G'')=1$. Then $a(F,X,B+\frac{1}{N''}G'')=\frac{i}{N''}$ for some non-negative integer $i<N''$. We have $$a(F,X,B+\frac{t}{N''}G'')=a(F,X,B+\frac{1}{N''}G'')+(1-t)\mult_F \frac{1}{N''}G''=\frac{i}{N''}+(1-t)\mult_F \frac{1}{N''}G'',$$ and $\mult_F \frac{1}{N''}G''=\frac{1}{1-t}(1-\frac{i}{N''})$. By Theorem~\ref{thm: 3fold acc ct general terminal germ case}, $\delta\le 1-t$ for some positive real number $\delta$ depending only on $\Ii$. Thus $\mult_F \frac{1}{N''}G''\leq \frac{1}{\delta}(1-\frac{i}{N''})$. By \cite[Lemma~5.1]{Kaw88}, $IG''$ is Cartier near $x$. It follows that $\mult_F \frac{1}{N''}G''$ belongs to a finite set depending only on $I$ and $\Ii$. Hence $t=1-(1-\frac{i}{N''})\frac{1}{\mult_F \frac{1}{N''}G''}$ belongs to a finite set of rational numbers depending only on $I$ and $\Ii$.
\end{proof}

\begin{thm}\label{thm: N-complements for canonical pair finite rational coeff}
Let $\Gamma\subset [0,1]\cap \Qq$ be a finite set. Then there exists a positive integer $N$ depending only on $\Ii$ satisfying the following.

Let $(X\ni x,B)$ be a threefold germ such that $X$ is a terminal, $B\in \Ii$, and $\mld(X\ni x, B)\geq 1$. Then there exists a monotonic $N$-complement $(X\ni x, B^+)$ of $(X\ni x, B)$ such that $\mld(X\ni x, B^+)=1$. Moreover, if $(X ,B)$ is canonical near $x$ and $x\in X$ is not smooth, then $(X, B^+)$ is canonical near $x$.
\end{thm}

\begin{proof}
Let $m$ be a positive integer such that $m\Ii\subset \Zz$. 

By Theorem~\ref{thm: index conjecture threefold}, we may assume that $\mld(X\ni x,B)>1$. Let $n$ be the index of $X\ni x$. By Lemma~\ref{lem: bounded index strict canonical complement}, we may assume that $n>4m$. By \cite[(6.1) Theorem]{Rei87} (cf. \cite[Theorems 12,23,25]{Mor85}),  $x\in X$ is of type $cA/n$ for some $n>4m$.

Possibly shrinking $X$ near $x$, we may assume that $X$ is affine. By Lemma~\ref{lem: linear system base locus x}, there exists a finite dimensional linear system $\mathfrak{d}\subset |-K_X|$ such that $\mathfrak{d}$ contains an elephant of $x\in X$ and the base locus of ${{\mathfrak {d}}}$ is $x$. By Lemma~\ref{lem: ct is attained by elephant at x}, if $(X,B)$ is canonical and $x\in X$ is not smooth, then the divisor $D=D_1+\dots+D_{m+1}$ satisfies $t:=\ct(X,B;\frac{1}{m+1}D)=\ct(X\ni x, B;\frac{1}{m+1}D)$, where $(D_1,\dots,D_{m+1})\in \mathfrak{d}^{m+1}$ is a general element. Since $(X, \frac{1}{m+1}D)$ is canonical near $x$, by Theorem~\ref{thm: 3-dim terminal number of coefficients local}(2), $\mld(X\ni x,\frac{1}{m+1}D)=1$.

By Lemma~\ref{lem: can extract divisor computing ct that is terminal strong version}, there exists a terminal blow-up (see Definition \ref{defn: terminal blow-up}) $f: Y\to X$ of $(X\ni x, B+\frac{t}{m+1}D)$ which extracts a prime divisor $E$ over $X\ni x$. Since $x\in X$ is a terminal singularity of type $cA/n$ for some $n> 4m>1$, by \cite[Theorem~1.3]{Kaw05} and Theorem~\ref{thm: kaw05 1.2 strenghthened}, $f: Y\to X$ is a divisorial contraction of ordinary type as in Theorem~\ref{thm: kaw05 1.2 strenghthened}(1). We may write $K_Y=f^*K_X+\frac{a}{n}E$ for some positive integer $a$. %In particular, under suitable analytic local coordinates $x_1,x_2,x_3,x_4$, there exist positive integers $r_1,r_2,b,d$, where $r_1+r_2=adn$ and $a\equiv br_1\mod n$, such that analytically locally, $$(X\ni x)\cong (\phi(x_1,x_2,x_3,x_4)=0)\subset (\mathbb C^4\ni o)/\frac{1}{n}(1,-1,b,0)$$ 
%for some invariant analytic power series $\phi$, and $f: Y\to X$ is a weighted blow-up at $x\in X$ with the weight $\frac{1}{n}(r_1,r_2,a,n)$. Assume that $mB$ is locally defined by $(h(x_1,x_2,x_3,x_4)=0)$ for some semi-invariant analytic power series $h$.

Since $(X,\frac{1}{m+1}D)$ is canonical near $x$, $a(E,X,\frac{1}{m+1}D)=a(E,X,0)-\mult_E \frac{1}{m+1}D=1+\frac{a}{n}-\mult_E \frac{1}{m+1}D\geq 1$. It follows that $\mult_E \frac{1}{m+1}D\leq \frac{a}{n}$. Since $B=\frac{1}{m}mB$, $n>3m$, and $\mld(X\ni x, B)>1$, by Lemma~\ref{lem: finiteness of description of terminal threefold singularities}(1), $a\leq 3$, and $\mult_E \frac{1}{m+1}D\in\{\frac{i}{(m+1)n}\mid i\in\Zz\cap [1,3(m+1)]\}$. Since $nm(K_X+B)$ is Cartier, $a(E,X,B)=1+\frac{k}{nm}$ for some positive integer $k$. Since $a(E,X,B+\frac{t}{m+1}D)=1$, $t\mult_E \frac{1}{m+1}D=a(E,X,B)-a(E,X,B+\frac{t}{m+1}D)=\frac{k}{nm}$, which implies that $t=\frac{k}{nm\mult_E \frac{1}{m+1}D}\in \frac{1}{(3m+3)!m}\Zz_{>0}.$ Now the coefficients of $B+\frac{t}{m+1}D$ belong to $\frac{1}{(3m+3)!(m+1)m}\Zz\cap[0,1]$. By Theorem~\ref{thm: index conjecture DCC threefold}, $(X\ni x,B+\frac{t}{m+1}D)$ is a monotonic $N$-complement of $(X\ni x,B)$ for some positive integer $N$ depending only on $\Ii$ satisfying all the required properties.
\end{proof}

\begin{proof}[Proof of Theorem~\ref{thm: intro cc}]
When $X$ is smooth near $x$, in particular, when $\dim x\geq 1$, we may take $G=0$, and $(X,B)$ is an $m$-complement of itself, where $m$ is a positive integer such that $m\Ii\subset \Zz$. When $x\in X$ is a closed point that is not smooth, by Theorem~\ref{thm: N-complements for canonical pair finite rational coeff}, we are done.
\end{proof}

\begin{thm}\label{thm: epsilon-N-complements for canonical pair finite rational coeff}
Let $\epsilon\geq 1$ be a rational number and $\Gamma\subset [0,1]\cap \Qq$ a finite set. Then there exists a positive integer $N$ depending only on $\epsilon$ and $\Ii$ satisfying the following.

Let $(X\ni x,B)$ be a threefold $\epsilon$-lc pair such that  $X$ is terminal and $B\in \Ii$. Then there exists a monotonic $(\epsilon,N)$-complement $(X\ni x, B^+)$ of $(X\ni x, B)$.
\end{thm}

\begin{proof}
If $\dim x=2$, then the theorem is trivial. When $X$ is smooth near $x$, in particular, when $\dim x=1$, we may take $G=0$, and $(X,B)$ is a monotonic $(\epsilon,m)$-complement of itself, where $m$ is a positive integer such that $m\Ii\subset \Zz$. From now on, we may assume that $x\in X$ is a closed point that is not smooth. 

When $\epsilon=1$, the theorem follows from Theorems~\ref{thm: N-complements for canonical pair finite rational coeff}. When $\epsilon>1$, let $n$ be the index of the terminal singularity $X\ni x$. By \cite[Appendix, Theorem]{Sho92} and \cite[Theorem~0.1]{Mar96}, there exists a prime divisor $E$ over $X\ni x$ such that $a(E,X,0)=1+\frac{1}{n}$, hence $1+\frac{1}{n}\geq a(E,X,B)\geq \epsilon,$ and $n\leq \lfloor\frac{1}{\epsilon-1}\rfloor$. It follows that $(X\ni x, B)$ is a monotonic $(\epsilon,\lfloor\frac{1}{\epsilon-1}\rfloor!m)$-complement of itself.
\end{proof}

\subsection{Boundedness of complements for DCC coefficients}

\begin{defn}\label{defn: definition of norm of divisor}
For any $\bm v=(v_1,\dots,v_m)\in \Rr^m$, we defined $||\bm v||:=\max_i\{v_i\}$. For an $\Rr$-divisor $B=\sum b_iB_i$, where $B_i$ are the distinct prime divisors of $\Supp B$, we define $||B||:=\max_i\{b_i\}$.
\end{defn}

\begin{thm}\label{thm: linearmld terminal threefold}

Let $m$ be a positive integer, $\epsilon\geq 1$ a real number,  and $\bm{v}=(v_1^0,\ldots,v_m^0)\in\Rr^m$ a point. Then there exist a rational polytope $\bm{v}\in P\subset\Rr^m$ with vertices $\bm{v}_j=(v_1^j,\ldots,v_m^j)$, positive real numbers $a_j$, and positive real numbers $\epsilon_j$ depending only on $m,\epsilon$ and $\bm{v}$ satisfying the following. 
	\begin{enumerate}
	 \item $\sum_{j} a_j=1,\sum_{j} a_j\bm{v}_j=\bm{v}$, and $\sum_{j} a_j\epsilon_j\ge\epsilon$.
	 \item Assume that $(X\ni x,B:=\sum_{i=1}^m v_i^0B_{i})$ is a threefold germ such that $X$ is terminal, $(X\ni x,B)$ is $\epsilon$-lc, and $B_1,\ldots,B_m\ge0$ are Weil divisors. Then for any $j$,
  	 $$\mld(X\ni x,\sum_{i=1}^m v_i^jB_{i})\ge \epsilon_j.$$
	\end{enumerate}
	Moreover, if $\epsilon>1$, then the function $P\to \Rr$ defined by
   	 $$(v_1,\ldots,v_m)\mapsto \mld(X\ni x,\sum_{i=1}^m v_iB_{i})$$
  	 is a linear function; if $\epsilon\in \Qq$, then we may pick $\epsilon_j=\epsilon$ for any $j$.
\end{thm}
\begin{proof}
%Possibly replacing $X$ by its small $\Qq$-factorization, we may assume that $X$ is $\Qq$-factorial.\luo{So, just take $\Qq$-factorization without any explanation?}
\noindent\textbf{Step 1}. There exist $\Qq$-linearly independent real numbers $r_0=1,r_1,\ldots,r_c$ for some $0\le c\le m$, and $\Qq$-affine functions $s_i:\Rr^{c}\to \Rr$ such that $s_i(\bm{r}_0)=v_i^0$ for any $1\le i\le m$, where $\bm{r}_0:=(r_1,\dots,r_c)$. Note that the map $\Rr^c\to V$ defined by $$\bm{r}\mapsto (s_1(\bm{r}),\dots,s_m(\bm{r}))$$ is one-to-one, where $V\subset\Rr^m$ is the rational envelope of $\bm{v}$.
      
If $c=0$, then $P=V=\{\bm{v}\}$, and there is nothing to prove. Suppose that $c\ge1$. Let $B(\bm{r}):=\sum_{i=1}^m s_i(\bm{r})B_i$. Then $B(\bm{r}_0)=\sum_{i=1}^m v_i^0B_{i}=B$. By \cite[Lemma 5.4(1)]{HLS19}, $K_X+B(\bm{r})$ is Cartier near $x$ for any $\bm{r}\in\Rr^c$. 

\medskip

\noindent\textbf{Step 2}. We will show that there exist a positive real number $\delta$ and a $\Qq$-affine function $f(\bm{r})$ depending only on $m,\epsilon,c,\bm{r}_0,$ $s_1,\dots,s_m$ such that $f(\bm{r}_0)\ge\epsilon$, and for any $\bm{r}\in\Rr^c$ satisfying $||\bm{r}-\bm{r}_0||\le\delta$, $\mld(X\ni x,B(\bm{r}))\ge f(\bm{r})$, moreover, when $\epsilon>1$, $$\mld(X\ni x,B(\bm{r}))=a(E,X,B(\bm{r}))\ge f(\bm{r})$$ for some prime divisor $E$ over $X\ni x$. 

When $\epsilon=1$, we may take $f(\bm r)=1$, and the assertion follows from Theorem~\ref{thm: uniform canonical polytope}. When $\epsilon>1$, by \cite[Appendix, Theorem]{Sho92}, for all germs $(X\ni x,B)$ which is $\epsilon$-lc, the index of $X\ni x$ is bounded from above by $I_0:=\lfloor \frac{1}{\epsilon-1}\rfloor$. Note that by \cite[Lemma~5.1]{Kaw88}, $I_0!D$ is Cartier for any Weil divisor $D$ on $X$. Also note that $(X\ni x,B(\bm{r}_0))$ is $\epsilon$-lc, and $\mld(X\ni x)\le 3$ (cf. \cite[Theorem~0.1]{Amb99}). The existence of $\delta$ and $f(\bm{r})$ in this case follows from \cite[Lemma 4.7]{CH21}. 
      
\medskip

\noindent \textbf{Step 3}. We finish the proof in this step. It follows from the same line of the proof of \cite[Theorem 7.15]{CH21}. 

Note that if $\epsilon\in\Qq$, then $f(\bm{r})=\epsilon$ for any $\bm{r}\in\Rr^c$. We may find $2^c$ positive rational numbers $r_{i,1},r_{i,2}$ such that $r_{i,1}<r_i<r_{i,2}$ and $\max\{r_i-r_{i,1},r_{i,2}-r_i\}\le\delta$ for any $1\le i\le c$. By our choice of $\delta$, the function $\Rr^c\to \Rr$ defined by $$\bm{r}\mapsto \mld(X\ni x,B(\bm{r}))$$ is a linear function on $\bm{r}\in U_c:=[r_{1,1},r_{1,2}]\times\cdots\times[r_{c,1},r_{c,2}]$. 
      
Let $\bm{r}_j$ be the vertices of $U_c$. Set $\epsilon_j:=f(\bm{r}_j),v_i^j:=s_i(\bm{r}_j)$ for any $i,j$. Note that if $\epsilon\in\Qq$, then $\epsilon_j=f(\bm{r}_j)=\epsilon$ for any $j$. Let $P:=\{(s_1(\bm{r}),\dots,s_m(\bm{r}))\mid \bm{r}\in U_c\}\subset V$. Then the function $P\to \Rr:\,(v_1,\ldots,v_m)\mapsto \mld(X\ni x,\sum_{i=1}^m v_iB_{i})$ is linear, $(v_1^j,\dots,v_m^j)$ are vertices of $P,$ and $$\mld(X\ni x,\sum_{i=1}^m v_i^jB_i)=\mld(X\ni x,B(\bm{r}_j))\ge f(\bm{r}_j)\ge\epsilon_j$$ for any $j$. 
      
Finally, we may find positive real numbers $a_j$ such that $\sum_{j} a_j=1$ and $\sum_{j} a_j\bm{r}_j=\bm{r}_0$. Then $\sum_{j} a_j\bm{v}_j=\bm{v}$ and $\sum_{j} a_j\epsilon_j\ge\epsilon$ as $\sum_j a_j v_i^j=\sum_j a_js_i(\bm{r}_j)=s_i(\sum_j a_j\bm{r}_j)=s_i(\bm{r}_0)=v_i^0$ for any $1\le i\le m,$ and $\sum_{j} a_j\epsilon_j=\sum_{j} a_jf(\bm{r}_j)=f(\sum_{j} a_j\bm{r}_j)=f(\bm{r}_0)\ge\epsilon.$
\end{proof}

\begin{thm}\label{ref: main epsilon comp threefolds finite set p version}
Let $p$ be a positive integer, $\epsilon\ge 1$ a real number, and $\Ii\subset [0,1]$ a finite set. Then there exists a positive integer $N$ depending only on $\epsilon,p$ and $\Ii$, such that $p\mid N$ and $N$ satisfies the following. 

Let $(X\ni x,B)$ be a pair such that $X$ is a terminal threefold, $B\in\Ii$, and $\mld(X\ni x,B)\geq \epsilon$. Then there exists an $N$-complement $(X\ni x,B^+)$ of $(X\ni x,B)$ such that $\mld(X\ni x,B^+)\geq \epsilon$. Moreover, if $\Span_{\Qq_{\ge0}}(\bar{\Ii}\cup\{\epsilon\}\backslash\Qq)\cap (\Qq\backslash\{0\})=\emptyset$, then we may pick $B^+\ge B.$ 
\end{thm}

%\han{proof maybe delete in the submit version}
\begin{proof}
 By Theorem \ref{thm: linearmld terminal threefold}, there exist three finite sets $\Ii_1\subset(0,1]$, $\Ii_2\subset[0,1]\cap \Qq$ and $\mathcal{M}$ of non-negative rational numbers depending only on $\epsilon,\Ii$, such that 
   \begin{itemize}
       %\item $\sum a_i=1,$
       \item $\sum a_i\epsilon_i     \ge\epsilon$,
       \item  $K_X+B=\sum a_i(K_X+B^i),$ and
       \item $(X\ni x,B^i)$ is $\epsilon_i$-lc at $x$ for any $i$, 
       %\item $-I(K_X+B^i)$ is Cartier near $x$ for any $i$,
       
       %\item If $\epsilon\in\Qq$, then  $\epsilon=\epsilon_i$ for any $i$, and if $\epsilon\neq 1$, then $\epsilon_i\neq1$ for any $i$.
     \end{itemize}
     for some $a_i\in\Ii_1,B^i\in\Ii_2$ and $\epsilon_i\in\mathcal{M}$. By Theorem~\ref{thm: epsilon-N-complements for canonical pair finite rational coeff}, there exists a positive integer $n_0$ which only depends on $\Gamma_2$ and $\mathcal{M}$, such that $({X}\ni x,{B}^i)$ has an $(\epsilon_i,n_0)$-complement $({X}\ni x,{B}^i+{G}^i)$ for some $\Qq$-Cartier divisor ${G}^i\ge0$ for any $i$.
     Let ${G}:=\sum a_i{G}^i$.
     
     By \cite[Lemma 6.2]{CH21}, there exists a positive integer $n$ depending only on $\epsilon,p,n_0,\Ii$, $\Ii_1,\Ii_2,\mathcal{M}$, such that there exist positive rational numbers $a_i'$ with the following properties:
     \begin{itemize}
         \item $pn_0|n$,
         \item $\sum a_i'=1$,
         \item $\sum a_i'\epsilon_i\ge\epsilon$,
         \item $na_i'\in n_0\Zz$ for any $i$, and
         \item $nB'\ge n\lfloor B\rfloor+\lfloor (n+1)\{B\}\rfloor$, where $B':=\sum a_i'{B}^i$.
     \end{itemize}
     Let ${G}':=\sum a_i'{G}^i$. Then
     $$
     n(K_{{X}}+B'+G')=n\sum a_i'(K_{X}+{B}^i+{G}^i)
     =\sum \frac{a_i'n}{n_0}\cdot n_0(K_{X}+{B}^i+{G}^i)\sim_Z0
     $$
     and
     $$a(E,X,B'+G')=\sum a_i'(E,X,B^i+G^i)\ge a_i'\epsilon_i\ge\epsilon$$
     for any prime divisor $E$ over $X\ni x$. Hence $(X\ni x,B'+G')$ is an $(\epsilon,n)$-complement of $(X\ni x,B)$.

     Moreover, if $\Span_{\Qq_{\ge0}}(\bar{\Ii}\cup\{\epsilon\}\backslash\Qq)\cap (\Qq\backslash\{0\})=\emptyset$, then ${B}'\ge B$ by \cite[Lemmas 6.2, 6.4]{CH21}.
     \end{proof}

Proposition \ref{prop: I R-Cartier terminal threefold} and Theorem \ref{thm: R-Cartier terminal threefold} study the inversion of stability property for $\Rr$-Cartier divisors, and give a positive answer to \cite[Conjecture 7.8]{HL20} in some special cases.

\begin{prop}\label{prop: I R-Cartier terminal threefold}
	Let $I$ be a positive integer and $\Ii\subset [0,1]$ a finite set. Then there exists a positive real number $\tau$ depending only on $I$ and $\Ii$ satisfying the following. 
	
	Let $x\in X$ be a terminal threefold singularity, and $B\ge0,B'\ge0$ two $\Rr$-divisors on $X$, such that 
	\begin{enumerate}
	\item $IK_X$ is Cartier near $x$,
	    \item $B\geq B'$,  $||B-B'||<\tau,B\in\Ii$,
	    \item $\mld(X\ni x,B')\ge 1$,
	    and
	    \item $K_X+B'$ is $\Rr$-Cartier.
	\end{enumerate}
	Then $K_X+B$ is $\Rr$-Cartier.
\end{prop}
\begin{comment}
a sequence of positive real numbers $\tau_i$, a sequence of terminal threefold singularities $x_i\in X_i$, and $\Rr$-divisors $B_i,B_i'$ on $X_i$, such that 
	\begin{itemize}
	\item $\lim_{i\to +\infty}\tau_i=0$,
	\item $IK_{X_i}$ is Cartier near $x_i$,
	    \item $B_i\geq B_i'$,  $||B_i-B_i'||<\tau_i,B_i\in\Ii, B_i'\in \Ii'$,
	    \item $\mld(X_i\ni x_i,B_i')\ge 1$,
	    \item $K_{X_i}+B_i'$ is $\Rr$-Cartier, and
	    \item $K_{X_i}+B_i$ is not $\Rr$-Cartier.
	\end{itemize}
\end{comment}
\begin{proof}%\han{add details?}
Suppose that the proposition does not hold, then there exist $X_i\ni x_i,B_i,B_i',\tau_i$ corresponding to $X\ni x,B,B',\tau$ as in the assumptions, and a DCC set $\Ii'$, such that
\begin{itemize}
    \item $\lim_{i\rightarrow+\infty}\tau_i=0$, 
    \item $B_i'\in\Ii'$, and
    \item $K_{X_i}+B_i$ is not $\Rr$-Cartier.
\end{itemize}
Let $f_i:Y_i\to X_i$ be a small $\Qq$-factorialization of $X_i$. Let $B_{Y_i}$ be the strict transform of $B_i$ on $Y_i$. Possibly replacing $Y_i$ with a minimal model of $(Y_i,B_{Y_i})$ over $X_i$, we may assume that $K_{Y_i}+B_{Y_i}$ is big and nef over $X$. We may write
$K_{Y_i}+B_{Y_i}':=f_i^{*}(K_{X_i}+B_{i}')$. Since $\mld(Y_i/X_i\ni x_i,B_{Y_i}')=\mld(X_i\ni x_i,B_i')\ge 1$, by Theorem \ref{thm: 3fold acc ct}, possibly passing to a subsequence, we may assume that $(Y_i/X_i\ni x_i,B_{Y_i})$ is $1$-lc over $x_i$ for any $i$. Since the Cartier index of any Weil divisor on $Y_i$ is bounded from above by $I$, by \cite[Theorem~1.2]{Nak16} and \cite[Theorem~0.1]{Amb99}, $\{\mld(Y_i/X_i\ni x_i,B_{Y_i})\}_{i=1}^{\infty}$ belongs to a finite set. Thus possibly passing to a subsequence, we may assume that there exists a real number $\epsilon\ge 1$, such that $\epsilon:=\mld(Y_i/X_i\ni x_i,B_{Y_i})$ for any $i$. Since $K_{Y_i}+B_{Y_i}'\le K_{Y_i}+B_{Y_i}$, $(Y_i/X_i\ni x_i,B_{Y_i}')$ is an $(\epsilon,\Rr)$-complement of itself. 

Note that $Y_i$ is of Fano type over $X_i$. Let $Y_i'$ be a minimal model of $-(K_{Y_i}+B_{Y_i})$ over $X_i$. Then $(Y_i'/X_i\ni x_i,B_{Y_i'}')$ is an $(\epsilon,\Rr)$-complement of itself, where $B_{Y_i'}'$ is the strict transform of $B_{Y_i}'$ on $Y_i'$. In particular, $\mld(Y_i'/X_i\ni x_i,B_{Y_i'}')\ge \epsilon$. By Theorem \ref{thm: alct acc terminal threefold}, possibly passing to a subsequence, we may assume that $(Y_i'/X_i\ni x_i,B_{Y_i'})$ is $\epsilon$-lc over $x_i$, where $B_{Y_i'}$ is the strict transform of $B_{Y_i}$ on $Y_i'$. Thus $(Y_i'/X_i\ni x_i,B_{Y_i'})$ is $(\epsilon,\Rr)$-complementary as $-(K_{Y_i'}+B_{Y_i'})$ is big and nef over $X_i$. By \cite[Lemma 3.13]{CH21}, $(Y_i/X_i\ni x_i,B_{Y_i})$ has an $(\epsilon,\Rr)$-complement $(Y_i/X_i\ni x_i,B_{Y_i}+G_{Y_i})$ for some $\Rr$-divisor $G_{Y_i}\ge0$.   

Let $Y_i\to Z_i$ be the canonical model of $({Y_i},B_{Y_i})$ over $X$ and $B_{Z_i}$ the strict transform of $B_{Y_i}$ on $Z_i$. Then $-G_{Z_i}$ is ample over $X$, where $G_{Z_i}$ is the strict transform of $G_{Y_i}$ on $Z_i$. Since $K_{X_i}+B_{i}$ is not $\Rr$-Cartier, the natural induced morphism $g_i:Z_i\to X_i$ is not the identity, and $G_{Z_i}\neq 0$. It follows that $\Supp G_{Z_i}$ contains $g_i^{-1}(x_i)$. Thus $$\epsilon=\mld(Y_i/X_i\ni x_i,B_{Y_i})=\mld(Z_i/X_i\ni x_i,B_{Z_i})>\mld(Z_i/X_i\ni x_i,B_{Z_i}+G_{Z_i})\ge\epsilon,$$ a contradiction.  
\end{proof}
\begin{rem}
Note that on any fixed potential klt variety $X$, the Cartier index of any Weil $\Qq$-Cartier divisor is bounded from above (cf. \cite[Lemma 7.14]{CH21}). Thus the proof of Proposition \ref{prop: I R-Cartier terminal threefold} also works for any fixed potential klt variety $X$ by assuming the ACC conjecture for minimal log discrepancies. It would be interesting to ask if it is necessary to assume $X$ is fixed in higher dimensional cases. 
\end{rem}

\begin{thm}\label{thm: R-Cartier terminal threefold}
	Let $\Ii\subset [0,1]$ be a finite set. Then there exists a positive real number $\tau$ depending only on $\Ii$ satisfying the following. 
	
	Let $x\in X$ be a terminal threefold singularity, and $B\ge0,B'\ge0$ two $\Rr$-divisors on $X$, such that 
	\begin{enumerate}
	    \item $B'\le B$,  $||B-B'||<\tau,B\in\Ii$,
	    \item $\mld(X\ni x,B')\ge 1$,
	    and
	    \item $K_X+B'$ is $\Rr$-Cartier.
	\end{enumerate}
	Then $K_X+B$ is $\Rr$-Cartier.
\end{thm}
\begin{proof}
Let $\tau$ be the positive real number constructed in Proposition \ref{prop: I R-Cartier terminal threefold} which only depends on $\Ii$ and $I:=1$.

Let $f:Y\to X$ be an index one cover of $K_X$. We may write $K_Y+B_{Y}':=f^{*}(K_X+B')$, and $K_Y+B_{Y}:=f^{*}(K_X+B)$. Then $K_Y$ is Cartier, $B_{Y}'\le B_{Y}$, $||B_{Y}-B_{Y}'||<\tau$, and $B_Y\in\Ii$. Moreover, by \cite[Proposition 5.20]{KM98}, $Y$ is terminal, and $\mld(Y\ni y,B_{Y}')\ge \mld(X\ni x,B')\ge 1$ for any point $y\in f^{-1}(x)$. Thus by Proposition \ref{prop: I R-Cartier terminal threefold}, $K_Y+B_{Y}$ is $\Rr$-Cartier. We conclude that $K_X+B$ is $\Rr$-Cartier as $f$ is a finite morphism. %\han{maybe use lemma in sec 4}
\end{proof}

\begin{thm}\label{ref: main epsilon comp threefolds p version}
Let $p$ be a positive integer, $\epsilon\ge 1$ a real number, and $\Ii\subset [0,1]$ a finite set. Then there exists a positive integer $N$ depending only on $\epsilon,p$ and $\Ii$, such that $p\mid N$ and $N$ satisfies the following. 

Let $(X\ni x,B)$ be a pair such that $X$ is a terminal threefold, $B\in\Ii$ and $\mld(X\ni x,B)\geq \epsilon$. Then there exists an $N$-complement $(X\ni x,B^+)$ of $(X\ni x,B)$ such that $\mld(X\ni x,B^+)\geq \epsilon$. Moreover, if $\Span_{\Qq_{\ge0}}(\bar{\Ii}\cup\{\epsilon\}\backslash\Qq)\cap (\Qq\backslash\{0\})=\emptyset$, then we may pick $B^+\ge B.$ 
\end{thm}

\begin{proof}
By Theorems \ref{thm: R-Cartier terminal threefold}, \ref{thm: alct acc terminal threefold}, \cite[Lemma 5.17]{HLS19} (see also \cite[Lemma 5.5]{CH21}) and follow the same lines of the proof of \cite[Theorem 5.6]{CH21} (see also \cite[Theorem 5.18]{HLS19}), possibly replacing $\Ii$ by a finite subset of $\bar{\Ii}$, we may assume that $\Ii$ is a finite set. Now the theorem follows from Theorem \ref{ref: main epsilon comp threefolds finite set p version}.
\end{proof}

\begin{proof}[Proof of Theorem \ref{thm: intro ecc}] This is a special case of Theorem \ref{ref: main epsilon comp threefolds p version}.
\end{proof}

\section{Proofs of the main results}\label{section: Proofs of the main results}

\subsection{Proof of Theorem \ref{thm: intro global canonical mld acc}}

In this subsection, we prove the following theorem:

\begin{thm}\label{thm: 1-gap pair}
Let $\Ii\subset [0,1]$ be a DCC set. Then $1$ is not an accumulation point of $$\{\mld(X,B)\mid \dim X=3,B\in\Ii\}$$ from below. 
\end{thm}

We give a proof of Theorem~\ref{thm: 1-gap pair} in this subsection.

%\han{relative setting?}
\begin{defn}
Let $(X,B)$ be a pair. We say that $(X,B)$ is \emph{extremely non-canonical} if $\mld(X,B)<1$, and the set
$$\{E\mid E\text{ is exceptional over }X, a(E,X,B)\leq 1\}$$
contains a unique element. In particular, any extremely non-canonical pair is klt.

A pair $(X\ni x,B)$ is called \emph{extremely non-canonical} if $(X,B)$ is extremely non-canonical near $x$ and $\mld(X\ni x,B)=\mld(X,B)<1$.
\end{defn}

\begin{lem}\label{lem: refduction to extremely canonical}
Let $d$ be a positive integer and $\Ii\subset [0,1]$ a set. Let $(X,B)$ be a klt pair of dimension $d$ such that $B\in\Ii$ and $\mld(X,B)<1$. Then there exists a $\Qq$-factorial extremely non-canonical klt pair $(Y,B_Y)$ of dimension $d$, such that $B_Y\in \Gamma$ and $\mld(X,B)\leq\mld(Y,B_Y)$.
\end{lem}

\begin{proof}
Since $(X,B)$ is klt, by \cite[Corollary 1.4.3]{BCHM10}, there exists a birational morphism $f: W\rightarrow X$ from a $\Qq$-factorial variety $W$ which extracts exactly all the exceptional divisors $E$ over $X$ such that $a(E,X,B)=1$. Let $K_W+B_W:=f^*(K_X+B)$. Possibly replacing $(X,B)$ with $(W,B_W)$, we may assume that $a(E,X,B)\neq1$ for any prime divisor $E$ that is exceptional over $X$.

Since $(X,B)$ is klt and $\mld(X,B)<1$, there exist prime divisors $E_1,\dots,E_k$ that are exceptional over $X$, such that
$$\{E_1,\dots,E_k\}=\{E\mid E\text{ is exceptional over }X,\, a(E,X,B)<1\}.$$
Let $\alpha_j:=1-a(E_j,X,B)$ for each $j$. By \cite[Lemma 5.3]{Liu18}, there exists $i\in\{1,2,\dots,k\}$ and a birational morphism $h: Y\rightarrow X$ from a $\mathbb{Q}$-factorial variety $Y$, such that
\begin{itemize}
    \item $f$ exactly extracts $E_1,\dots,E_{i-1},E_{i+1},\dots,E_k$, and
    \item $\mult_{E_i}\sum_{j\not=i}\alpha_jE_{j,Y}<\alpha_i$, where $E_{j,Y}=\Center_YE_j$ for each $j\neq i$.% denotes the center of $E_j$ on $Y$ when viewed as divisorial valuations over $X$.
\end{itemize}
Let $B_Y:=h_*^{-1}B$. Then
\begin{align*}
    &\mld(Y,B_Y)\leq a(E_i,Y,B_Y)=a(E_i,Y,B_Y+\sum_{j\neq i}\alpha_j E_{j,Y})+\mult_{E_i} \sum_{j\neq i}\alpha_j E_{j,Y}\\
   <&a(E_i,X,B)+\alpha_i=1,
\end{align*}
and for any prime divisor $F\not=E_i$ that is exceptional over $Y$,
$$a(F,Y,B_Y)\geq a(F,Y,B_Y+\sum_{j\not=i}\alpha_jE_{j,Y})=a(F,X,B)>1.$$
Thus $(Y,B_Y)$ satisfies our requirements.
\end{proof}

%$\begin{thm}[{\cite[Theorem 1.4]{LX21}, see also \cite[Theorem 1.3]{Jia21}}]\label{thm: 12/13}Let $X$ be a $\Qq$-Gorenstein threefold. If $\mld(X)<1$, then $\mld(X)\leq\frac{12}{13}$.\end{thm} \han{maybe let Jiang happy, Jiang first}
%\han{try to use decomposable complements in this paper.}

\begin{lem}\label{lem: low bound lct by complement} Let $d$ be a positive integer and $\Ii\subset [0,1]$ a DCC set. Then there exists a positive real number $t$ depending only on $d$ and $\Ii$ satisfying the following.

Let $(X,B)$ be a klt pair of dimension $d$, $E$ a prime exceptional divisor over $X$ such that $a(E,X,B)<1$, and $x$ the generic point of $\Center_XE$. Let $f: Y\rightarrow X$ be a birational morphism which only extracts $E$. Then $(Y,B_Y+tE)$ is lc over a neighborhood of $x$.
\end{lem}
\begin{proof}
By Theorem \ref{thm: ni decomposable complement}, there exist a positive integer $n$ and a finite set $\Ii_0\subset (0,1]$, such that $(X\ni x,B)$ has an $(n,\Ii_0)$-decomposable $\Rr$-complement $(X\ni x,B^+)$ of $(X\ni x,B)$. In particular, there exist real numbers $a_1,\dots,a_k\in\Ii_0$ and lc pairs $(X\ni x,B_i^+)$, such that $\sum_{i=1}^ka_i=1$, $\sum_{i=1}^ka_iB_i^+=B^+$, and each $(X\ni x,B_i^+)$ is an $n$-complement of itself. Let
$$\Ii_0':=\{\sum_{i=1}^k s_ia_i\mid ns_i\in\Zz_{\geq 0}\}.$$
Then $\Ii_0'\subset [0,+\infty)$ is a discrete set, and we may let
$$\gamma_0:=\max\{\gamma\in\Ii_0'\mid \gamma<1\}.$$
Since $na(E,X,B_i^+)\in\Zz_{\ge0}$ for every $i$, $$1>a(E,X,B)\geq a(E,X,B^+)=\sum_{i=1}^ka_ia(E,X,B_i^+)\in\Ii_0',$$
so $a(E,X,B^+)\leq \gamma_0$. Thus $(Y,B_Y+(1-\gamma_0)E)$ is lc over a neighborhood of $x$. We may take $t:=1-\gamma_0$.
\end{proof}

\begin{lem}\label{lem: 1-gam canonical extreme case}
Let $\Ii\subset [0,1]$ be a DCC set. Then there exists a positive real number $\epsilon$ depending only on $\Ii$ satisfying the following.

Let $(X,B)$ be a $\Qq$-factorial extremely non-canonical threefold pair such that $X$ is strictly canonical. Then $\mld(X,B)\leq 1-\epsilon$. 
\end{lem}

\begin{proof}
Since $(X,B)$ is extremely non-canonical, $(X,B)$ is klt. Let $E$ be the unique prime divisor that is exceptional over $X$ such that $a(E,X,B)\leq 1$. Then $a(E,X,B)<1$. Moreover, for any prime divisor $F\not=E$ that is exceptional over $X$, $a(F,X,0)\geq a(F,X,B)>1$. Since $X$ is strictly canonical, $a(E,X,0)=1$. In particular, $\mult_E B>0$.

Let $\gamma_0:=\min\{\gamma\in\Ii\mid \gamma>0\}$ and let $m:=\lceil\frac{1}{\gamma_0}\rceil$. Then $B\geq\frac{1}{m}\Supp B$, 
$$a(E,X,B)\leq a(E,X,\frac{1}{m}\Supp B)<a(E,X,0)=1,$$
and
$$1<a(F,X,B)\leq a(F,X,\frac{1}{m}\Supp B)$$
for any prime divisor $F\not=E$ that is exceptional over $X$. Thus possibly replacing $\Ii$ with $\{0,\frac{1}{m}\}$ and $B$ with $\frac{1}{m}\Supp B$, we may assume that $\Ii$ is a finite set of rational numbers.

Let $f: Y\rightarrow X$ be a birational morphism which extracts $E$, and let $B_Y$ be the strict transform of $B$ on $Y$. Then $K_Y=f^{*}K_X$, and $$K_Y+B_Y+(1-a(E,X,B))E=f^*(K_X+B).$$ Let $x$ be the generic point of $\Center_XE$. If $\dim x=1$, by taking general hyperplane sections, the lemma follows from  \cite[Theorem~3.8]{Ale93} (see also \cite{Sho94b} and \cite[Theorem~1.5]{HL20}). Therefore, we may assume that $x$ is a closed point. 

By Theorem~\ref{thm: canonical index}, $60K_Y$ is Cartier over a neighborhood of $x$. By our construction, $Y$ is terminal. By \cite[Lemma~5.1]{Kaw88}, $60D$ is Cartier over an neighborhood of $x$ for any Weil divisor $D$ on $Y$.

By Lemma~\ref{lem: low bound lct by complement}, there exists a positive integer $n$ depending only on $\Ii$, such that $n\Ii\subset\Zz_{\geq 0}$ and $(Y,B_Y+\frac{1}{n}E)$ is lc over a neighborhood of $x$. 

If $a(E,X,B)\leq 1-\frac{1}{2n}$, then we are done. Thus we may assume that $a(E,X,B)>1-\frac{1}{2n}$. In this case, by the boundedness of length of extremal rays (cf. \cite[Theorem~4.5.2(5)]{Fuj17}), there exists a $(K_Y+B_Y+\frac{1}{n}E)$-negative extremal ray $R$ in $\overline{NE}(Y/X)$ which is generated by a rational curve $C$, such that 
$$0>(K_Y+B_Y+\frac{1}{n}E)\cdot C\geq -6.$$
Since $(K_Y+B_Y+(1-a(E,X,B))E)\cdot C=0$, we have
$$0< (a(E,X,B)-1+\frac{1}{n})(-E\cdot C)\leq 6.$$
Hence $0<(-E\cdot C)<12n$ as $a(E,X,B)>1-\frac{1}{2n}$. Since  $60n(K_Y+B_Y)$ is Cartier over a neighborhood of $x$, we have
$$a(E,X,B)=1-\frac{60n(K_Y+B_Y)\cdot C}{60n(-E\cdot C)},$$
so $a(E,X,B)<1-\frac{1}{720n^2}$ and we are done.
\end{proof}

\begin{lem}\label{lem: 1-gap X terminal extreme case}
Let $\Ii\subset [0,1]$ be a DCC set. Then there exists a positive real number $\epsilon$ depending only on $\Ii$ satisfying the following.

Let $(X,B)$ be a $\Qq$-factorial extremely non-canonical threefold pair such that $X$ is terminal. Then $\mld(X,B)\leq 1-\epsilon$.
\end{lem}
\begin{proof}
Let $E$ be the unique divisor that is exceptional over $X$ such that $a(E,X,B)\leq 1$. Then $0<a(E,X,B)<1$. Since $X$ is terminal, $a(E,X,0)>1$. Thus $\mult_EB>0$.%Moreover, for any prime divisor $F\not=E$ that is exceptional over $X$, we have $a(F,X,0)\geq a(F,X,B)>1$.

Let $x$ be the generic point of $\Center_XE$. If $\dim x=1$, by taking general hyperplane sections, the lemma follows from  \cite[Theorem~3.8]{Ale93} (see also \cite{Sho94b} and \cite[Theorem~1.5]{HL20}). Therefore, we may assume that $x$ is a closed point. 

\medskip

Let $t:=\ct(X,0;B)$. Since $X$ is terminal and $(X,B)$ is extremely non-canonical, we have $a(E,X,cB)=1$ and $t=\ct(X\ni x,0;B)<1$. By Theorem \ref{thm: 3fold acc ct}, there exists a real number $\delta\in (0,1)$ depending only on $\Ii$ such that $t\leq 1-\delta$. Let $\gamma_0:=\min\{\gamma\in\Ii\mid \gamma>0\}$, $m:=\lceil\frac{1}{\delta\gamma_0}\rceil$, and $B':=\frac{1}{m}\lfloor mB\rfloor$. Then $||B-B'||<\delta\gamma_0$. Hence $B\geq B'\geq tB$ and $\Supp(B'-tB)=\Supp B$. Since $\mult_EB>0$,
$$a(E,X,B)\leq a(E,X,B')<a(E,X,tB)=1,$$
and
$$1<a(F,X,B)\leq a(F,X,B')$$
for any prime divisor $F\not=E$ that is exceptional over $X$. Thus possibly replacing $\Ii$ with $\frac{1}{m}\Zz_{\geq 0}\cap [0,1]$, $B$ with $B'$, and $t$ with $\ct(X,0;B')$ respectively, we may assume that $\Ii\subset \frac{1}{m}\mathbb{Z}\cap[0,1]$.

\medskip

If $x\in X$ is a terminal singularity of types other than $cA/n$ or of type $cA/n$ with $n\le 2$, then by \cite[(6.1) Theorem]{Rei87}, the index of $X\ni x$ divides $12$. By \cite[Lemma~5.1]{Kaw88}, $12m(K_X+B)$ is Cartier, and we may take $\epsilon=\frac{1}{12m}$ in this case. Thus we may assume that $x\in X$ is a terminal singularity of type $cA/n$ for some $n\geq 3$. 

\medskip

By construction, $(X,tB)$ is extremely non-canonical. By Lemma~\ref{lem: can extract divisor computing ct that is terminal strong version}, there exists a terminal blow-up $f: Y\to X$ of $(X\ni x,tB)$ which extracts $E$. Since $n\geq 3$, by \cite[Theorem 1.3]{Kaw05}, $f$ is of ordinary type. Let $a:=a(E,X,B)+1$. By Theorem \ref{thm: kaw05 1.2 strenghthened}(1), under suitable analytic local coordinates $x_1,x_2,x_3,x_4$, there exist positive integers $r_1,r_2,b,d$, where $\gcd(b,n)=1$, $r_1+r_2=adn$ and $a\equiv br_1\mod n$, such that analytically locally, $$(X\ni x)\cong (\phi(x_1,x_2,x_3,x_4)=0)\subset (\mathbb C^4\ni o)/\frac{1}{n}(1,-1,b,0)$$
for some invariant analytic power series $\phi$, and $f: Y\to X$ is a weighted blow-up at $x\in X$ with the weight $w:=\frac{1}{n}(r_1,r_2,a,n)$. Assume that $mB$ is locally defined by $(h(x_1,x_2,x_3,x_4)=0)$ for some semi-invariant analytic power series $h$.

\begin{claim}\label{claim: bound n for extremaly canonical variety}
    If either $d\geq 4$ or $a\geq 4$, then $n\leq 3m$.
\end{claim}

We proceed the proof assuming Claim~\ref{claim: bound n for extremaly canonical variety}. If either $a\geq 4$ or $d\geq 4$, then by Claim~\ref{claim: bound n for extremaly canonical variety}, $n\leq 3m\leq \frac{3m}{t}$. Since $a(E,X,tB)=a(E,X,0)-t\mult_E B=1+\frac{a}{n}-t\mult_E B$,
\begin{equation}\label{eqn: relation of c and n}
    \frac{a}{n}=t\mult_E B.
\end{equation}
It follows that
$$\mult_E B=\frac{a}{tn}\geq \frac{a}{3m}\geq  \frac{1}{3m}.$$ Thus $a(E,X,B)=a(E,X,tB)-(1-t)\mult_E B\leq 1-\frac{\delta}{3m}$. We can take $\epsilon=\frac{\delta}{3m}$ in this case.

\medskip

We may now assume that $a\leq 3$ and $d\leq 3$. Since $a\equiv br_1\mod n$, $\gcd(r_1,n)\mid 6$. Since $r_1+r_2=adn$, $\gcd(r_2,n)\mid 6$, and $\gcd(r_1,r_2)\mid adn$. This implies that $\gcd(r_1,r_2)\mid 216$. Let $m'$ be the smallest positive integer such that $m'tB$ is an integral divisor and $r$ the smallest positive integer such that $rm'(K_X+tB)$ is Cartier. By Lemma~\ref{lem: unboundness of r1 and r2}, $r\mid \gcd(r_1,r_2)$. Thus $r\mid 216$. By (\ref{eqn: relation of c and n}), $t=\frac{a}{n\mult_E B}=\frac{a}{N}$, where $N=n\mult_E B$ is a positive integer. We may write $tB=\frac{a}{mN} mB$, then $216mN(K_X+tB)$ is Cartier.

\medskip

By \cite[4.8 Corollary]{Sho94a}, there exists a prime divisor $E_1\neq E$ over $X\ni x$ such that $a(E_1,X,0)= 1+\frac{a_1}{n}$ for some positive integer $a_1\leq 2$.

\medskip

Since $a(E_1,X,tB)>1$, $$1+\frac{a_1}{n}=a(E_1,X,0)\geq a(E_1,X,tB)\geq 1+\frac{1}{216mN},$$ hence $n\leq 432mN=\frac{432am}{t}$. It follows that $tn\leq 432am\leq 1296m$. By (\ref{eqn: relation of c and n}),
$$\mult_E B=\frac{a}{tn}\geq \frac{1}{1296m},$$ and $a(E,X,B)=a(E,X,tB)-(1-t)\mult_EB\leq 1-\frac{\delta}{1296m}$. We can take $\epsilon=\frac{\delta}{1296m}$ in this case.
\end{proof}

\begin{proof}[Proof of Claim~\ref{claim: bound n for extremaly canonical variety}]
Suppose that $n>3m$. If either $d\geq 4$ or $a\geq 4$, then we can pick positive integers $s_1,s_2$ such that
\begin{itemize}
    \item $s_1+s_2=a'dn$ for some $a'\leq \min\{a,3\}$,
    \item $a'\equiv bs_1\mod n$,
    \item $s_1>n,s_2>n$, and
    \item $\frac{1}{n}(s_1,s_2,a',n)\neq \frac{1}{n}(r_1,r_2,a,n)$.
\end{itemize}
In fact, when $a\geq 4$, we may take $a'=3$. When $d\geq 4$, we may take $a'=1$ and $(s_1,s_2)\neq (r_1,r_2)$. Let $ w':=\frac{1}{n}(s_1,s_2,a',n)$.

Since $a\geq a'$, by Lemma~\ref{lem: irreducibility of exceptional divisors extracted by weighted blow up}, the weighted blow-up with the weight $w'$ at $x\in X$ extracts an analytic prime divisor $E'\neq E$ such that $w'(X\ni x)=\frac{a'}{n}$. By Lemma~\ref{lem: algebraic approximation of weighted blow-up}, we may assume that $E'$ is a prime divisor over $X\ni x$. By our assumption, $a({E'},X,B)=1+{w'}(X\ni x)-{w'}(B)> 1$, thus
$$\frac{1}{m}>\frac{3}{n}\geq \frac{a'}{n}=w'(X\ni x)>  w'(B)=\frac{1}{m}w'(mB),$$
which implies that $w'(h)= w'(mB)<1$.
Since $ w'(x_1)=\frac{s_1}{n}>1$, $w'(x_2)=\frac{s_2}{n}>1$ and $w'(x_4)=1$, up to a scaling of $h$, there exists a positive integer $l$, such that $x_3^{l}\in h$ and ${w'}(mB)=w'(h)= w'(x_3^{l})$. In particular, $w'(h)=l w'(x_3)=\frac{a'l}{n}$ and
$$w'(X\ni x)=\frac{a'}{n}> w'(B)=\frac{a'}{n}\frac{l}{m},$$
this implies that $\frac{l}{m}< 1$. On the other hand,
$$1>a(E,X,B)=1+w(X\ni x)-w(B)\geq 1+\frac{a}{n}- \frac{1}{m}w(x_3^l)=1+\frac{a}{n}(1-\frac{l}{m})>1,$$
a contradiction.
\end{proof}

\begin{proof}[Proof of Theorem~\ref{thm: 1-gap pair}]
Let $(X,B)$ be a threefold pair such that $B\in \Ii$. Possibly replacing $(X,B)$ with a $\Qq$-factorialization and replacing $\Ii$ with $\Ii\cup\{1\}$, we may assume that $(X,B)$ is $\Qq$-factorial dlt. 

If $(X,B)$ is klt, then by Lemma \ref{lem: refduction to extremely canonical}, we may assume that $(X,B)$ is extremely non-canonical. Then either $X$ is not canonical, or $X$ is strictly canonical, or $X$ is terminal. If $X$ is not canonical, then by \cite[Theorem 1.4]{LX21} (see also \cite[Theorem 1.3]{Jia21}), $\mld(X,B)\leq\mld(X)\leq\frac{12}{13}$. If $X$ is strictly canonical, then the theorem follows from Lemma \ref{lem: 1-gam canonical extreme case}. If $X$ is terminal, then the theorem follows from Lemma \ref{lem: 1-gap X terminal extreme case}.

If $(X,B)$ is not klt, then we let $E$ be a prime divisor that is exceptional over $X$ such that $a(E,X,B)=\mld(X,B)$. If $\Center_XE\not\subset\lfloor B\rfloor$, then $$\mld(X,B)=a(E,X,B)=a(E,X,\{B\})\geq\mld(X,\{B\})\geq\mld(X,B)$$
and the theorem follows from the klt case. If $\Center_XE\subset\lfloor B\rfloor$, then there exists a prime divisor $S\subset\lfloor B\rfloor$ such that $\Center_XE\subset S$. We let $K_S+B_S:=(K_X+B)|_S$. By \cite[Corollary 1.4.5]{BCHM10}, $\mld(X,B)=a(E,X,B)$ is equal to the total minimal log discrepancy of $(S,B_S)$. Since $B_S\in D(\Ii)$ which satisfies the DCC, the theorem follows from \cite[Theorem 3.8]{Ale93} (see also \cite{Sho04b} and \cite[Theorem 1.5]{HL20}).
\end{proof}

\begin{proof}[Proof of Theorem \ref{thm: intro global canonical mld acc}]
This follows from Theorems~\ref{thm:  terminal mld acc} and \ref{thm: 1-gap pair}.
\end{proof}

\subsection{Proof of Theorem \ref{thm: mn intro DCC coeff}}

\begin{lem}\label{lem: nak ter 3fold finite coefficient bdd index}
Let $I$ be a positive integer, and $\Ii\subset [0,1]$ a finite set. Then there exists a positive real number $l$ depending only on $I$ and $\Ii$ satisfying the following. Assume that
\begin{enumerate}
    \item $(X\ni x,B)$ is a threefold pair,
    \item $X$ is terminal,
    \item $B\in\Ii$,
    \item $\mld(X\ni x,B)\geq 1$, and
    \item $IK_X$ is Cartier near $x$.
\end{enumerate}
Then there exists a prime divisor $E$ over $X\ni x$, such that $a(E,X,B)=\mld(X\ni x,B)$ and $a(E,X,0)=1+\frac{a}{I}$ for some non-negative integer $a\leq l$.
\end{lem}
\begin{proof}
We will use some ideas of Kawakita \cite[Theorem 4.6]{Kaw21}. Possibly replacing $X$ with a small $\Qq$-factorialization, we may assume that $X$ is $\Qq$-factorial. 

If $\dim x=2$, then the lemma is trivial as we can take $l=0$.

If $\dim x=1$, then $X$ is smooth near $x$. By Lemma~\ref{lem: Terminal blow up}, $\mld(X\ni x, B)=a(E,X,B)$, where $E$ is the exceptional divisor obtained by blowing up $x\in X$. We have $a(E,X,0)=2$, and the theorem holds in this case.

If $\dim x=0$ and suppose that the theorem does not hold, then there exists a sequence of threefold germs $(X_i\ni x_i,B_i)$ satisfying (1-5), and a strictly increasing sequence of positive integers $l_i$, such that for each $i$,
  $$\min\{a(E,X_i,0)\mid \Center_{X_i}E=x_i,a(E,X_i,B_i)=\mld(X_i\ni x_i,B_i)\}=1+\frac{l_i}{I}.$$
By \cite[Lemma 5.1]{Kaw88}, $ID$ is Cartier near $x_i$ for any Weil divisor $D$ on $X_i$. By \cite[Theorem 0.1]{Amb99}, $1\leq\mld(X_i\ni x_i,B_i)\leq 3$ for any $i$. By \cite[Corollary 1.3]{Nak16}, possibly passing to a subsequence, we may assume that there exists a real number $\alpha\geq 1$, such that $\mld(X_i\ni x_i,B_i)=\alpha$ for any $i$. By \cite[Theorem 1.2]{Nak16}, there exists a real number $\alpha'>\alpha$, such that for any $i$ and any prime divisor $F_i$ over $X_i\ni x_i$, if $a(F_i,X_i,B_i)>\alpha$, then $a(F_i,X_i,B_i)>\alpha'$. Therefore,
$$\alpha'\text{-}\lct(X_i\ni x_i,0;B_i)=\frac{1+\frac{l_i}{I}-\alpha'}{1+\frac{l_i}{I}-\alpha}=1-\frac{\alpha'-\alpha}{1+\frac{l_i}{I}-\alpha}<1$$
is strictly increasing, which contradicts Theorem \ref{thm: alct acc terminal threefold}.
\end{proof}

\begin{comment}
\begin{conj}\label{conj: inversion of stalibty for divisors computing MLDs}
	Let $\Ii\subset [0,1]$ be a finite set, $X$ a normal variety, and $x\in X$ a point. Then there exists a positive real number $\tau$ depending only on $\Ii$ and $x\in X$ satisfying the following. 
	
	Assume that $(X\ni x,B)$ and $(X\ni x,B')$ are two lc pairs and $E$ is a prime divisor over $X\ni x$, such that
	\begin{enumerate}
	    \item $B\geq B',||B-B'||<\tau$, and $B\in\Ii$, and
	    \item $a(E,X,B')=\mld(X\ni x,B')$,
	\end{enumerate}
	then $a(E,X,B)=\mld(X\ni x,B)$.
\end{conj}
\end{comment}

The following theorem answers a question of \cite[Conjecture 7.2]{HL20} for terminal threefold pairs. We will use it to prove Theorem \ref{thm: mn intro DCC coeff}.

\begin{thm}\label{thm: divisors compouting mlds}
Let $I$ be a positive integer, and $\Ii\subset [0,1]$ a finite set. Then there exists a positive real number $\tau$ depending only on $I$ and $\Ii$ satisfying the following. Assume that $(X\ni x,B)$ and $(X\ni x,B')$ are two threefold lc pairs and $E'$ is a prime divisor over $X\ni x$, such that
	\begin{enumerate}
	\item $X$ is terminal,
	\item $B\geq B',||B-B'||<\tau$, and $B\in\Ii$, 
	\item $a(E',X,B')=\mld(X\ni x,B')\geq 1$, and
	\item $IK_X$ is Cartier near $x$.
	\end{enumerate}
Then $a(E',X,B)=\mld(X\ni x,B)$.
\end{thm}
\begin{proof}
We may assume that $\{0\}\subsetneq \Ii$, otherwise $B=B'=0$ and the theorem is obvious. 

Since $IK_X$ is Cartier near $x$ and $X$ is terminal, by \cite[Lemma~5.1]{Kaw88}, $ID$ is Cartier for any $\Qq$-Cartier Weil divisor on $X$. By \cite[Theorem 5.6]{HLS19}, there exist positive real numbers $a_1,\dots,a_k\in (0,1]$ depending only on $\Ii$, a positive integer $I'$ depending only on $I$ and $\Ii$, and $\Qq$-divisors $B_1,\dots,B_k\geq 0$ on $X$, such that
\begin{itemize}
\item $\sum_{i=1}^ka_i=1$,
\item $\sum_{i=1}^ka_iB_i=B$, 
\item $(X\ni x,B_i)$ is lc for any $i$. In particular, $K_X+B_i$ is $\Qq$-Cartier for any $i$, and
\item $I'(K_X+B_i)$ is Cartier near $x$ for each $i$.
\end{itemize}
Thus there exists a positive real number $\delta$ depending only on $I$ and $\Ii$, such that for any prime divisor $F$ over $X\ni x$ and $a(F,X,B)>\mld(X\ni x,B)$, we have $a(F,X,B)>\mld(X\ni x,B)+\delta$.

By Lemma \ref{lem: nak ter 3fold finite coefficient bdd index}, there exists a prime divisor $E$ over $X\ni x$, such that $a(E,X,B)=\mld(X\ni x,B)$ and $a(E,X,0)\leq l$ for some positive integer $l$ depending only on $I$ and $\Ii$. In particular, $\mult_E B=a(E,X,0)-a(E,X,B)\le l$.

We show that we may take $\tau:=\frac{\delta}{2l}\cdot\min\{\gamma\in\Ii\mid \gamma>0\}$. In this case, $B'\geq (1-\frac{\delta}{2l})B$. Since
\begin{align*}
    &a(E,X,0)-\mult_EB'=a(E,X,B')\geq a(E',X,B')\geq a(E',X,B)\\
    =&(a(E',X,B)-a(E,X,B))+a(E,X,0)-\mult_E B,
\end{align*}
we have
$$0\leq a(E',X,B)-a(E,X,B)\leq \mult_E(B-B')\leq\frac{\delta}{2l}\mult_E B\le \frac{\delta}{2},$$
which implies that $a(E',X,B)=a(E,X,B)=\mld(X\ni x,B)$.
\end{proof}

%\han{maybe deleted in submitted version.}
\begin{lem}\label{lem: quotient sequences acc}
Let $m_0$ be a positive integer and let $\{a_{i,1}\}_{i=1}^{\infty},\{a_{i,2}\}_{i=1}^{\infty},\dots,\{a_{i,m_0}\}_{i=1}^{\infty}$ be $m_0$ sequences of positive real numbers. Then there exists an integer $1\leq k\leq m_0$, such that possibly passing to a subsequence, $\{\frac{a_{i,j}}{a_{i,k}}\}_{i=1}^{\infty}$ are decreasing (resp. increasing) for all $1\leq j\leq m_0$.
\end{lem}
\begin{proof}
Possibly passing to a subsequence, we may assume that for any $k,j$, $\{\frac{a_{i,j}}{a_{i,k}}\}_{i=1}^{\infty}$ is either decreasing or strictly increasing (resp. either increasing or strictly decreasing). Suppose that the lemma does not hold. Then there exists a function $\pi:\{1,2,\dots,m_0\}\rightarrow\{1,2,\dots,m_0\}$, such that $\{\frac{a_{i,\pi(j)}}{a_{i,j}}\}_{i=1}^{\infty}$ is strictly increasing (resp. strictly decreasing) for any $j$. We may pick $1\leq j_0\leq m$ such that $\pi^{(l)}(j_0)=j_0$ for some positive integer $l$. Then
$$\{1\}_{i=1}^{\infty}=\{\frac{a_{i,\pi(j_0)}}{a_{i,j_0}}\cdot\frac{a_{i,\pi(\pi(j_0))}}{a_{i,\pi(j_0)}}\cdot\dots\cdot \frac{a_{i,\pi^{(l)}(j_0)}}{a_{i,\pi^{(l-1)}(j_0)}}\}_{i=1}^{\infty}$$
is strictly increasing (resp. strictly decreasing), which is absurd.
\end{proof}

\begin{lem}\label{lem: detailed lemma for dcc uniform bdd}
Let $m_0\geq 0, I>0$ be integers, $\Ii_0\subset [0,1]$ a finite set, and $\Ii\subset [0,1]$ a DCC set. Then there exists a positive integer $l$ depending only on $m_0,I,\Ii_0$, and $\Ii$ satisfying the following.

Assume that $\{(X_i\ni x_i,B_i:=\sum_{j=1}^{m_0}b_{i,j}B_{i,j}+B_{i,0})\}_{i=0}^{\infty}$ is a sequence of $\Qq$-factorial threefold germs, such that
\begin{enumerate}
\item $X_i$ is terminal for each $i$,
\item $\{b_{i,j}\}_{i=1}^{\infty}$ is strictly increasing for any fixed $j$, 
\item $b_{i,j}\in\Ii$ and  $B_{i,0}\in\Ii_0$ for each $i$ and $j$,
\item $B_{i,j}\geq 0$ is a Weil divisor on $X_i$ for each $i$ and $j$,
\item $\mld(X_i\ni x_i,B_i)> 1$ for each $i$,
\item $IK_{X_i}$ is Cartier near $x_i$ for each $i$, and
\item $1+\frac{l_i}{I}:=\min\{a(E_i,X_i,0)\mid \Center_{X_i}E_i=x_i,a(E_i,X_i,B_i)=\mld(X_i\ni x_i,B_i)\}.$
\end{enumerate}
Then possibly passing to a subsequence, we have $l_i\leq l$ for each $i$.
\end{lem}
\begin{proof}

\noindent\textbf{Step 1}. 
We prove the lemma by induction on $m_0$. When $m_0=0$, the lemma follows from Lemma \ref{lem: nak ter 3fold finite coefficient bdd index}. Thus we may assume that $m_0\geq 1$. 

Let $\gamma_0:=\min\{1,\gamma\mid \gamma\in\Ii,\gamma>0\}$. Let $b_j:=\lim_{i\rightarrow+\infty}b_{i,j}$, $\bar B_i:=\sum_{j=1}^{m_0}b_jB_{i,j}+B_{i,0}$ for any $i$, and $\Ii_0':=\Ii_0\cup\{b_1,\dots,b_{m_0}\}$. By Theorem \ref{thm: alct acc terminal threefold}, possibly passing to a subsequence, we may assume that $\mld(X_i\ni x_i,\bar B_i)\geq 1$ for each $i$.

By \cite[Lemma~5.1]{Kaw88}, for each $i$, $ID_i$ is Cartier near $x_i$ for any Weil divisor $D_i$ on $X_i$. By \cite[Theorem 0.1]{Amb99}, $1<\mld(X_i\ni x_i,\bar B_i)\leq 3$. By \cite[Theorem 1.2]{Nak16}, possibly passing to a subsequence, we may assume that there exist two real numbers $\alpha\geq 1$ and $\delta>0$, such that for any $i$,  
\begin{itemize}
    \item $\mld(X_i\ni x_i,\bar B_i)=\alpha$, and
    \item for any prime divisor $F_i$ over $X_i\ni x_i$ such that $a(F_i,X_i,\bar B_i)>\mld(X_i\ni x_i,\bar B_i)$, we have $a(F_i,X_i,\bar B_i)>\alpha+\delta$.
\end{itemize}

For each $i$, let $E_i$ be a prime divisor over $X_i\ni x_i$ such that $a(E_i,X_i,B_i)=\mld(X_i\ni x_i,B_i)$ and $a(E_i,X_i,0)=1+\frac{l_i}{I}$. By Theorem \ref{thm: divisors compouting mlds}, possibly passing to a subsequence, we may assume that $a(E_i,X_i,\bar B_i)=\mld(X_i\ni x_i,\bar B_i)=\alpha$.
\medskip

\noindent\textbf{Step 2}. For any $i$ and any $1\leq j\leq m_0$, we define $B'_{i,j}:=\sum_{k\neq j}b_{i,k}B_{i,k}+b_jB_{i,j}+B_{i,0}$.

%for $m-1,I,\Ii_0'$ and $\Ii$,
By the induction for $m_0-1,I,\Ii_0'$ and $\Ii$, possibly passing to a subsequence, we may assume that there exists a positive integer $l'$ depending only on $m_0,I,\Ii_0$ and $\Ii$, such that for any $1\leq j\leq m_0$, there exists a prime divisor $E_{i,j}$ over $X_i\ni x_i$, such that
\begin{itemize}
    \item $a(E_{i,j},X_i,B'_{i,j})=\mld(X_i\ni x_i,B_{i,j}')$, and
    \item $a(E_{i,j},X_i,0)\leq 1+\frac{l'}{I}\le 1+l'$.
\end{itemize}
Since
$$a(E_{i,j},X_i,B_i)\geq\mld(X_i\ni x_i,B_i)=a(E_i,X_i,B_i)$$
and
$$a(E_{i,j},X_i,B'_{i,j})=\mld(X_i\ni x_i,B'_{i,j})\leq a(E_i,X_i,B'_{i,j}),$$
we have $\mult_{E_{i,j}}(B'_{i,j}-B_i)\geq\mult_{E_{i}}(B'_{i,j}-B_i)$. By the construction of $B'_{i,j}$, we have
$$\mult_{E_{i,j}}B_{i,j}\geq\mult_{E_i}B_{i,j}$$
for any $i$ and $1\leq j\leq m_0$.
Since 
\begin{align*}
    1&\leq \mld(X_i\ni x_i,\bar B_i)\leq\mld(X_i\ni x_i,B'_{i,j})\\
    &=a(E_{i,j},X_i,B'_{i,j})\leq a(E_{i,j},X_i,b_{i,j}B_{i,j})\leq a(E_{i,j},X_i,\gamma_0B_{i,j})\\
    &=a(E_{i,j},X_i,0)-\gamma_0\mult_{E_{i,j}}B_{i,j}\leq 1+l'-\gamma_0\mult_{E_{i,j}}B_{i,j},
\end{align*}
for any $i$ and $1\leq j\leq m_0$, we have 
$$\mult_{E_i}B_{i,j}\leq\mult_{E_{i,j}}B_{i,j}\leq\frac{l'}{\gamma_0}.$$

\medskip

\noindent\textbf{Step 3}. Let $a_{i,j}:=b_j-b_{i,j}$ for any $i$ and any $1\leq j\leq m_0$. By Lemma \ref{lem: quotient sequences acc}, possibly re-odering indices and passing to a subsequence, we may assume that $\{\frac{a_{i,j}}{a_{i,1}}\}_{i=1}^{\infty}$ is decreasing for any $1\leq j\leq m_0$. Let $M:=\max\{\frac{a_{1,j}}{a_{1,1}}\mid 1\leq j\leq m_0\}$, $t:=\frac{\delta\gamma_0}{m_0Ml'}$, $\tilde b_{i,j}:=b_j-t\frac{a_{i,j}}{a_{i,1}}$ for any $i,j$, and $\tilde B_{i}:=\sum_{j=1}^{m_0}\tilde b_{i,j}B_{i,j}+B_{i,0}$. Possibly passing to a subsequence, we may assume that $a_{i,1}<t$ for any $i$ as $\lim_{i\to+\infty}a_{i,1}=0$.  

There exist a positive integer $k$ and a finite set $\Lambda\subset\{1,2,\dots,m_0\}$, such that $|\Lambda|=k$, and $\tilde b_{i,j}=\tilde b_{1,j}$ for any $i$ and $j\in\Lambda$ as $\tilde b_{i,1}=b_1-t$. By the induction for $m_0-k,I,\Ii_0\cup\{\tilde b_{1,j}\mid j\in\Lambda\}$, and $\{\tilde b_{i,j}\}_{i\geq 1,1\leq j\leq m_0}$, possibly passing to a subsequence, for any $i$, there exists a prime divisor $\tilde E_i$ over $X_i\ni x_i$, such that
\begin{itemize}
    \item $a(\tilde E_i,X_i,\tilde B_i)=\mld(X_i\ni x_i,\tilde B_i)$, and
    \item $a(\tilde E_i,X_i,0)\leq 1+\frac{l}{I}$ for some positive integer $l$ depending only on $m_0,I,\Ii_0$ and $\Ii$.
\end{itemize}

Since
$$\tilde b_{i,j}\geq b_j-\frac{\delta\gamma_0}{m_0Ml'}\cdot M=b_j-\frac{\delta\gamma_0}{m_0l'},$$
we have
\begin{align*}
&\mld(X_i\ni x_i,\tilde B_i)\leq a(E_i,X_i,\tilde B_i)=a(E_i,X_i,\bar B_i)+\mult_{E_i}(\bar B_i-\tilde B_i)\\
=&\alpha+\sum_{j=1}^{m_0}(b_j-\tilde b_{i,j})\mult_{E_i}B_{i,j}\leq \alpha+\sum_{j=1}^{m_0}\frac{\delta\gamma_0}{m_0l'}\cdot\frac{l'}{\gamma_0}=\alpha+\delta.
\end{align*}

Therefore, $a(\tilde{E_i},X_i,\bar B_i)\leq a(\tilde{E}_i,X_i,\tilde B_i)=\mld(X_i\ni x_i, \tilde{B}_i)\leq \alpha+\delta$, and by our choice $\delta$, we have $a(\tilde{E_i},X_i,\bar B_i)=\mld(X_i\ni x_i,\bar B_i)=\alpha$. By the construction of $\tilde B_i$,
$$B_i=\frac{a_{i,1}}{t}\tilde B_i+(1-\frac{a_{i,1}}{t})\bar B_i.$$
It follows that $a(\tilde E_i,X_i,B_i)=\mld(X_i\ni x_i, B_i)$. Thus $a(E_i,X_i,0)=1+\frac{l_i}{I}\le a(\tilde E_i,X_i,0)\leq 1+\frac{l}{I}$. 
\end{proof}

\begin{thm}\label{thm: bdd mld computing dcc}
Let $\Ii\subset [0,1]$ be a DCC set. Then there exists a positive integer $l$ depending only on $\Ii$ satisfying the following. 

Assume that $(X\ni x,B)$ is a threefold pair such that $X$ is terminal, $B\in\Ii$, and $\mld(X\ni x,B)\ge 1$. Then there exists a prime divisor $E$ over $X\ni x$, such that $a(E,X,B)=\mld(X\ni x,B)$ and $a(E,X,0)\leq 1+\frac{l}{I}$, where $I$ is the index of $X\ni x$. In particular, $a(E,X,0)\leq 1+l$.
\end{thm}
\begin{proof}
Possibly replacing $X$ with a small $\Qq$-factorialization, we may assume that $X$ is $\Qq$-factorial.

Let $\gamma_0:=\min\{\gamma\in\Ii,1\mid \gamma>0\}$. Suppose that the theorem does not hold. Then by Lemmas \ref{lem: Terminal blow up} and \ref{lem: bdd mld computing dcc cA/n type}, Theorems \ref{thm: 3-dim terminal number of coefficients local} and \ref{thm: Nakamura Mustata Conjecture}, there exist a positive integer $I$, an integer $0\leq m\leq\frac{2}{\gamma_0}$, a strictly increasing sequence of positive integers $l_i$, and a sequence of threefold germs $(X_i\ni x_i,B_i=\sum_{j=1}^mb_{i,j}B_{i,j})$, such that
\begin{itemize}
    
    \item $X_i$ is $\Qq$-factorial terminal for each $i$,
    \item $b_{i,j}\in\Ii$, and $B_{i,j}\geq 0$ is a Weil divisor for any $i,j$,
    \item $\mld(X_i\ni x_i,B_i)>1$ for each $i$,
    \item $IK_{X_i}$ is Cartier near $x_i$ for each $i$, and
    \item $1+\frac{l_i}{I}=\min\{a(E_i,X_i,0)\mid \Center_{X_i}E_i=x_i,a(E_i,X_i,B_i)=\mld(X_i\ni x_i,B_i)\}$.
\end{itemize}

Possibly passing to a subsequence, we may assume that $\{b_{i,j}\}_{i=1}^{\infty}$ is increasing for any fixed $j$. We let $b_j:=\lim_{i\rightarrow+\infty}b_{i,j}$ for any $j$, and $\Ii_0:=\{b_1,\dots,b_m\}$. Possibly reordering indices and passing to a subsequence, we may assume that there exists an integer $0\leq m_0\leq m$, such that
\begin{itemize}
    \item $b_{i,j}\neq b_j$ for any $i$ when $j\leq m_0$, and
    \item $b_{i,j}=b_j$ for every $i$ when $j>m_0$.
\end{itemize}
Let $B_{i,0}:=\sum_{j=m_0+1}^{m}b_jB_{i,j}$. Then $B_i=\sum_{j=1}^{m_0}b_{i,j}B_{i,j}+B_{i,0}$. By Lemma \ref{lem: detailed lemma for dcc uniform bdd}, possibly passing to a subsequence, $l_i\le l$ for some positive integer $l$ depending only on $\Ii$, a contradiction.
\end{proof}

\begin{proof}[Proof of Theorem \ref{thm: mn intro DCC coeff}]
This follows from Theorem \ref{thm: bdd mld computing dcc}.
\end{proof}

%\han{maybe deleted in submitted version.}

\begin{thm}\label{thm: uniform bdd mld computing dcc}
Let $\Ii_0=\{b_1,\ldots,b_m\}\subset [0,1]$ be a finite set. Then there exist a positive integer $l$ and a positive real number $\epsilon$ depending only on $\Ii_0$ satisfying the following. 

Assume that $(X\ni x,B'=\sum_{i} b_i'B_i')$ is a threefold pair such that 
\begin{enumerate}
    \item $X$ is terminal,
    \item $b_i-\epsilon<b_i'<b_i$ for each $i$ and $B_i'\ge 0$ are Weil divisors on $X$,
    \item $\mld(X\ni x,B')\ge 1$, and
    \item $IK_X$ is Cartier near $x$ for some positive integer $I$.
\end{enumerate}
Then for any prime divisor $E$ over $X\ni x$ such that $a(E,X,B)=\mld(X\ni x,B)$, we have $a(E,X,0)\leq 1+\frac{l}{I}.$ In particular, $a(E,X,0)\le 1+l$. 
\end{thm}
\begin{proof}
Suppose that the theorem does not hold. Then there exist a strictly increasing sequence of positive integers $l_i$ and a sequence of threefold pairs $(X_i\ni x_i,B_i=\sum_{j=1}^mb_{i,j}B_{i,j})$, such that
\begin{itemize}
    \item $X_i$ is terminal,
    \item $b_{i,j}$ is strictly increasing with $\lim_{i\to +\infty} b_{i,j}=b_j$ for each $j$, and $B_{i,j}\geq 0$ is a Weil divisor for any $i,j$,
    \item $\mld(X_i\ni x_i,B_i)\ge 1$ for any $i$,
    \item $I_iK_{X_i}$ is Cartier near $x_i$ for some positive integer $I_i$, and
    \item there exists a prime divisor $E_i$ over $X_i\ni x_i$ such that $a(E_i,X_i,B_i)=\mld(X_i\ni x_i,B_i)$, and $a(E_i,X_i,0)\geq \frac{l_i}{I_i}$.
\end{itemize}
Possibly replacing each $X_i$ with a small $\Qq$-factorialization, we may assume that $X_i$ is $\Qq$-factorial for each $i$. 

Let $\bar B_i:=\sum_{j=1}^mb_jB_{i,j}$ for any $i$. By Theorem \ref{thm: alct acc terminal threefold}, possibly passing to a subsequence, we may assume that $\mld(X_i\ni x_i,\bar B_i)\ge 1$ for any $i$.  

By Lemma \ref{lem: quotient sequences acc}, possibly reordering the indices and passing to a subsequence, we may assume $\{\frac{b_{i,j}}{b_{i,1}}\}_{i=1}^{\infty}$ is an increasing sequence for each $j$. In particular, $\frac{b_{i,j}}{b_{i,1}}\le \frac{b_j}{b_1}$ as $\lim_{i\to +\infty} \frac{b_{i,j}}{b_{i,1}}=\frac{b_j}{b_1}$. For each $i,j$, let $$b_{i,j}':=b_{i,j}+(b_1-b_{i,1})\frac{b_{i,j}}{b_{i,1}}=\frac{b_{1}b_{i,j}}{b_{i,1}}\le b_{j},$$
and $\bar B_i':=\sum_{j=1}^m b_{i,j}'B_{i,j}\le \bar B_i$. Then $\mld(X_i\ni x_i,\bar B_i')\ge 1$. 
Note that $\Ii':=\{\frac{b_{1}b_{i,j}}{b_{i,1}}\}_{i\in\Zz_{\ge 1}, 1\le j\le m}$ satisfies the DCC. By Theorem \ref{thm: bdd mld computing dcc}, there exists a positive integer $l$ depending only on $\Ii'$ and a prime divisor $E_i'$ over $X_i\ni x_i$ for each $i$, such that $a(E_i',X_i,\bar B_i')=\mld(X_i\ni x_i,\bar B_i')$, and $a(E_i',X_i,0)\leq 1+\frac{l}{I_i}$. Since
\begin{align*}
&a(E_i',X_i,\bar B_i')+\mult _{E_i}(\bar B_i'-B_i)=\mld(X_i\ni x_i,\bar B_i')+\mult _{E_i}(\bar B_i'-B_i)\\
\le &a(E_i,X_i,\bar B_i')+\mult _{E_i}(\bar B_i'-B_i)=a(E_i,X_i,B_i)=\mld(X_i\ni x_i,B_i)\\
\le & a(E_i',X_i,B_i)=a(E_i',X_i,\bar B_i')+\mult _{E_i'}(\bar B_i'-B_i),
\end{align*}
we have $\mult _{E_i}(\bar B_i'-B_i)\le \mult _{E_i'}(\bar B_i'-B_i)$. By the construction of $\bar B_i'$,
$$\bar B_i'-B_i=\sum_{j=1}^m (b_{i,j}'-b_{i,j})B_{i,j}=\sum_{j=1}^m(b_1-b_{i,1})\frac{b_{i,j}}{b_{i,1}}B_{i,j}=\frac{b_1-b_{i,1}}{b_{i,1}}B_i.$$ It follows that $\mult _{E_i} B_i\le \mult _{E_i'} B_i$. Hence
\begin{align*}&a(E_i,X_i,0)=a(E_i,X_i,B_i)+\mult_{E_i}B_i
\le a(E_i',X_i,B_i)+\mult_{E_i'}B_i\\
=&a(E_i',X_i,0)\le 1+\frac{l}{I_i},
\end{align*}
a contradiction.
\end{proof}

%Boundedness of Calabi-Yau type threefolds modulo flops
\subsection{Proof of Theorem \ref{thm: non-canonical klt CY threefold bdd flop}}
\begin{defn}[Log Calabi--Yau pairs]
A log pair $(X, B)$ is called a
\emph{log Calabi--Yau pair} if $K_X+B\sim_\Rr 0$. % by \cite{G}.
\end{defn}

%\han{maybe deleted in submission.}
\begin{defn}[Bounded pairs]\label{sec.bdd}
A collection of varieties $ \mathcal{D}$ is
said to be \emph{bounded} (resp. 
\emph{birationally bounded}, \emph{bounded in codimension one}) if there exists a projective morphism $h\colon \mathcal{Z}\rightarrow S$
of schemes of finite type such that
each $X\in \mathcal{D}$ is isomorphic (resp. birational, isomorphic in codimension one) to $\mathcal{Z}_s$ 
for some closed point $s\in S$.

We say that a collection of log pairs $\mathcal{D}$ is 
{\it log birationally bounded} (resp.  \emph{log bounded},  \emph{log bounded in codimension one})
if there exist a  quasi-projective scheme $\mathcal{Z}$, a 
reduced divisor $\mathcal{E}$ on $\mathcal Z$, and a 
projective morphism $h\colon \mathcal{Z}\to S$, where 
$S$ is of finite type and $\mathcal{E}$ does not contain 
any fiber, such that for every $(X,B)\in \mathcal{D}$, 
there exist a closed point $s \in S$ and a birational
map (resp. isomorphism, isomorphism in codimension one) $f \colon \mathcal{Z}_s \dashrightarrow X$ 
such that $\mathcal{E}_s$ contains the support of $f_*^{-1}B$ 
and any $f$-exceptional divisor (resp. $\mathcal{E}_s$ 
coincides with the support of $f_*^{-1}B$, $\mathcal{E}_s$ 
coincides with the support of $f_*^{-1}B$).

Moreover, if $\mathcal{D}$ is a set of klt Calabi--Yau 
varieties (resp. klt log Calabi--Yau pairs), then it is 
said to be {\it bounded modulo flops} (resp. {\it log 
	bounded modulo flops}) if it is bounded (resp. log bounded) in 
codimension one, each fiber $\mathcal{Z}_{s}$ 
corresponding to a member in $\mathcal{D}$ is normal projective, 
and $K_{\mathcal{Z}_s}$ is $\bQ$-Cartier (resp. $K_{\mathcal{Z}_s}+f_*^{-1}B$ is $\Rr$-Cartier). 
\end{defn}

%Here the name ``modulo flops" comes from the fact that, if we assume that $X$ and $\mathcal{Z}_s$ are both $\bQ$-factorial, then they are connected by flops by running a $(K_X+B+\delta f_*H)$-MMP where $H$ is an ample divisor on $\mathcal{Z}_s$ and $\delta$ is a sufficiently small positive number  (cf. \cite{BCHM10, flops}) .

\begin{proof}[Proof of Theorem \ref{thm: non-canonical klt CY threefold bdd flop}]
We follow the proof of \cite[Theorem 6.1]{Jia21} and \cite[Theorem 5.1]{CDHJS21}. 
By Theorem \ref{thm: intro global canonical mld acc}, there exists a positive real number $\delta<1$ depending only on $\Ii$, such that $\mld(X,B)\le 1-\delta$. By \cite[Corollary 1.4.3]{BCHM10}, there exists a birational morphism $f:Y\to X$ which extracts exactly one exceptional divisor $E$ with $a:=a(E,X,B)\le 1-\delta$. By \cite[Lemma 3.21]{HLS19}, $Y$ is of Fano type over $X$. Possibly replacing $Y$ with the canonical model of $-E$ over $X$, we may assume that $-E$ is ample over $X$, and $\Exc(f)=\Supp E$. We may write
$$K_Y+B_Y+(1-a)E=f^{*}(K_X+B)\equiv 0,$$
where $B_Y$ is the strict transform of $B$ on $Y$. By \cite[Theorem 1.5]{HMX14}, there exists a finite subset $\Ii_0\subset \Ii$ depending only on $\Ii$, such that $B\in\Ii_0$. Possibly replacing $\Ii$ with $\Ii_0$, we may assume that $\Ii$ is finite. By \cite[Theorem 1.5]{HMX14} again (see also the proof of \cite[Lemma 3.12]{CDHJS21}), there exists a positive real number $\epsilon<\frac{1}{2}$ depending only on $\Ii$ such that $(X,B)$ is $(2\epsilon)$-lc. Thus $(Y,B_Y+(1-a)E)$ is a $(2\epsilon)$-lc log Calabi-Yau pair with $1-a\ge \delta>0$. By \cite[Corollary 1.4]{HM07}, each fiber of $f$ is rationally chain connected. Since $\Exc(f)=\Supp E$, $E$ is uniruled. Now by \cite[Proposition 6.4]{Jia21}, the pairs $(Y, B_Y+(1-a)E)$ are log bounded modulo flops. That is, there are finitely many normal varieties $\mathcal{W}_i$, an $\Rr$-divisor $\mathcal{B}_i$ and a reduced divisor $\mathcal{E}_i$ on $\mathcal W_i$, and a projective morphism $\mathcal{W}_i\to S_i$, where $S_i$ is a normal variety of finite type, and $\mathcal{B}_i,\mathcal{E}_i$ do not contain any fiber of $\mathcal{W}_i\to S_i$, such that for every $(Y,B_Y+(1-a)E)$, there is an index $i$, a closed point $s \in S_i$, and a small birational map $g : {\mathcal{W}_{i,s}} \dashrightarrow Y$ such that $\mathcal{B}_{i,s}=g_*^{-1}B_Y$ and $\mathcal{E}_{i,s}=g_*^{-1}E$. We may assume that the set of points $s$ corresponding to such $Y$ is dense in each $S_i$. We may just consider a fixed index $i$ and ignore the index in the following argument.
	
For the point $s$ corresponding to $(Y,B_Y+(1-a)E)$, $$K_{{\mathcal{W}_s}}+g_*^{-1}B_Y+(1-a)g_*^{-1}E\equiv f_*^{-1}(K_Y+B_Y+(1-a)E)\equiv 0$$ and therefore $({\mathcal{W}_s}, g_*^{-1}B_Y+(1-a)g_*^{-1}E)$ is a $(2\epsilon)$-lc log Calabi--Yau pair. 
	
Let $h: \mathcal{W}'\to \mathcal{W}$ be a log resolution of $(\mathcal{W},\mathcal{B}+\mathcal{E})$, $\mathcal{B}'$ the strict transforms of $\mathcal{B}$ on $\mathcal{W}'$, and $\mathcal{E}'$ the sum of all $h$-exceptional reduced divisors and the strict transform of $\mathcal{E}$ on $\mathcal{W}'$. Then there exists an open dense subset $U\subset S$ such that for the point $s\in U$ corresponding to $(Y,B_Y+(1-a)E)$, $h_s: \mathcal{W}'_s\to \mathcal{W}_s$ is a log resolution of $(Y,B_Y+(1-a)E)$. Since $(\mathcal{W}_s,g_*^{-1}B_Y+(1-a)g_*^{-1}E)$ is $(2\epsilon)$-lc, $$K_{\mathcal{W}'_s}+\mathcal{B}_s'+(1-\epsilon)\mathcal{E}'_s-h_s^*(K_{{\mathcal{W}_s}}+g_*^{-1}B_Y+(1-a)g_*^{-1}E)$$
is an $h_s$-exceptional $\Rr$-divisor whose support coincides with $\Supp \mathcal{E}_s'$. Note that $\dim \mathcal{W}_s=3$. By \cite[Lemma 2.10, Theorem 1.1]{HH20}, we may run a $(K_{\mathcal{W}'}+\mathcal{B}'+(1-\epsilon)\mathcal{E}')$-MMP with scaling of an ample divisor over $S$ and reach a relative minimal model $\tilde{\mathcal{W}}$ over $S$. For the point $s\in U$ corresponding to $(Y,B_Y+(1-a)E)$, $\mathcal{E}'_s$ is contracted, and hence $\tilde{\mathcal{W}}_s$ is isomorphic to $X$ in codimension one. This gives a bounded family modulo flops over $U$. Applying Noetherian induction on $S$, the family of all such $X$ is bounded modulo flops.
\end{proof}

\begin{rem}
It is possible to replace Theorem \ref{thm: intro global canonical mld acc} with the uniform lc rational polytopes \cite[Theorem 5.6]{HLS19} and the boundedness of indices of log Calabi-Yau threefolds \cite[Theorem 1.13]{Xu19} to conclude $\mld(X,B)\le 1-\delta$ in the beginning of the proof of Theorem \ref{thm: non-canonical klt CY threefold bdd flop}. We briefly describe the idea here. By \cite[Theorem 1.5]{HMX14}, we may assume that $\Ii$ is a finite set. By \cite[Theorem 5.6]{HLS19}, we may reduce the theorem to the case $\Ii\subset\Qq$. By \cite[Theorem 1.13]{Xu19}, $I(K_X+B)$ is Cartier for some positive integer $I$ which only depends on $\Ii$. In particular, $\mld(X,B)\le 1-\delta$. 
\end{rem}

\section{Questions and conjectures}

It would be interesting to ask if Lemmas \ref{lem: can extract mld place smooth case} \ref{lem: can extract mld place cA/n case} hold for all terminal threefolds.

\begin{conj}\label{conj: compute mld terminal blow-up}
Let $(X\ni x, B)$ be an lc threefold pair, such that $X$ is terminal and $\mld(X\ni x, B)\geq 1$. Then there exists a prime divisor $E$ over $X\ni x$, and a divisorial contraction $f: Y\to X$ of $E$, such that $Y$ is terminal, and $a(E,X,B)=\mld(X\ni x, B)$.
\end{conj}

We remark that the assumption ``$\mld(X\ni x, B)\geq 1$'' is necessary in Conjecture \ref{conj: compute mld terminal blow-up}. Indeed, \cite[Excerise 6.45]{KSC04}, and \cite[Example 5]{Kaw17} show that there exists a $\Qq$-divisor $B$ on $X:=\Cc^3$, such that $\mld(X\ni x,B)=0$, there is exactly one prime divisor $E$ over $X\ni x$ with $a(E,X,B)=\mld(X\ni x,B)$, and $E$ is not obtained by a weighted blow-up. Recall that any divisorial contraction from a terminal threefold to a smooth variety is always a weighted blow-up. 

%\han{delete remark in submitted version}

%We expect the following.
\begin{conj}[{cf. \cite[Introduction]{HL20}}]\label{conj: bdd mld computing divisor}
    Let $d$ be a positive integer and $\Ii\subset[0,1]$ a DCC set. Then there exists a positive real number $l$ depending only on $d$ and $\Ii$ satisfying the following. 
    
    Assume that $(X\ni x,B)$ is an lc pair of dimension $d$ such that $X$ is $\Qq$-Gorenstein and $B\in\Ii$. Then there exists a prime divisor $E$ over $X\ni x$, such that $a(E,X,B)=\mld(X\ni x,B)$ and $a(E,X,0)\leq l$.
\end{conj}
\cite[Conjecture~1.1]{MN18} and \cite[Problem~7.17]{CH21} are exactly Conjecture \ref{conj: bdd mld computing divisor} for the case when $X\ni x$ is a fixed germ and $\Ii$ is a finite set, and when $X\ni x$ is a fixed germ respectively. Conjecture \ref{conj: bdd mld computing divisor} holds when $\dim X=2$ \cite[Theorem~1.2]{HL20}. In this paper, we give a positive answer for terminal threefolds. A much ambitious problem is the following.

%\begin{prob}Let $d$ be a positive integer and $\Ii\in[0,1]$ a DCC set, then there exists a positive real number $l$ depending only on $d$ and $I$ satisfying the following. Let $X$ be a klt variety of dimension $d$, and $(X\ni x, B)$ an lc pair such that $B\in\Ii$, and $\mld(X\ni x,B)\ge 1$. Then there exists a birational morphism $f: Y\to X$ that exactly extracts a prime divisor $E$ over $X\ni x$, such that $-E$ is $f$-ample, $a(E,X,B)=\mld(X\ni x, B)$, and $a(E,X,0)\le l$.\end{prob}

%\han{check here}

\begin{ques}
Let $\Ii\subset [0,1]$ be a DCC set. Assume that $(X\ni x,B)$ is an lc pair such that $X$ is klt near $x$ and $\mld(X\ni x,B)>0$.

\begin{enumerate}
    \item Will there exist a divisorial contraction $f: Y\rightarrow X$ of a prime divisor $E$ over $X\ni x$, such that $a(E,X,B)=\mld(X\ni x,B)$?
    \item Moreover, if $B\in\Ii$, and $\mld(X\ni x,B)\geq 1$, will $a(E,X,B)=\mld(X\ni x,B)$ and $a(E,X,0)\leq l$ for some real number $l$ depending only on $\dim X$ and $\Ii$?
\end{enumerate}
\end{ques}

It was shown in \cite[Theorem~1.1]{HLQ21} that the lc threshold polytopes satisfy the ACC, and a conjecture due to the first author asks whether the volumes of lc threshold polytopes satisfy the ACC. In the same fashion, we ask the following.

\medskip

\begin{ques}[ACC for CT-polytopes]\label{ques: ACC for CT-polytopes}
Let $d,s$ be positive integers, and $\Ii \subseteq \Rr_{\ge0}$ a DCC set. Let $\Ss$ be the set of all $(X, \Delta; D_1,\ldots, D_s)$, where
\begin{enumerate}
\item $\dim X=d$, $(X, \Delta)$ is canonical, and $\Delta\in\Ii$, and
\item $D_1, \ldots, D_s$ are $\Rr$-Cartier divisors, and $D_1,\ldots,D_s\in \Ii$.
\end{enumerate}
Then 
\begin{enumerate}
    \item $\{P(X, \Delta; D_1, \ldots, D_s) \mid (X, \Delta; D_1, \ldots, D_s)\in \Ss\}$ satisfies the ACC (under the inclusion), and
    \item $\{\mathrm{Vol}(P(X, \Delta; D_1, \ldots, D_s)) \mid (X, \Delta; D_1, \ldots, D_s)\in \Ss\}$ satisfies the ACC,
\end{enumerate}
where
$$
P(X, \Delta; D_1,\ldots,D_s)\coloneqq\{(t_1,\ldots,t_s) \in \mathbb{R}_{\geq 0}^s\mid (X,\Delta+t_1 D_1+\ldots+t_s D_s){\rm~is~canonical}\}.
$$
\end{ques}

\begin{comment}
\begin{ques}
Determine accumulation points of $cA/n$. 
\end{ques}
\end{comment}

\appendix

\section{ACC for Newton polytopes}

\begin{defn}\label{defn: newton polytope zd}
Let $n$ be a positive integer. A \emph{Newton polytope} $\mathcal{N}$ is a subset of $\mathbb{Z}_{\geq 0}^n$ satisfying the following: for any point $\bm{x}\in\mathcal{N}$, 
$$\bm{x}+\mathbb{Z}_{\geq 0}^n:=\{\bm{x}+\bm{v}\mid\bm{v}\in\mathbb{Z}_{\geq 0}^n\}\subset\mathcal{N}.$$
\end{defn}

\begin{defn}
Let $n$ be a positive integer, $\bm 0$ the origin of $\Zz^n$ and $\mathcal{N}\subset\mathbb{Z}_{\geq 0}^n$ a Newton polytope. A \emph{vertex} of $\mathcal{N}$ is a point $\bm{u}\in\mathcal{N}$, such that for any $\bm{x}\in\mathcal{N}$ and $\bm{v}\in\mathbb{Z}_{\geq 0}^n$, if $\bm{u}=\bm{x}+\bm{v}$, then $\bm{u}=\bm{x}$ and $\bm{v}=\bm{0}$.
\end{defn}

\begin{lem}\label{lem: weak ACC newton polytopes}
Let $\{\bm{v}_j\}_{j\in\Zz_{\ge1}}$ be a sequence of vectors in $\mathbb{Z}_{\geq 0}^n$. Then the set $\{\mathcal{P}_j=\cup_{i=1}^j (\bm{v}_i+\mathbb{Z}_{\geq 0}^n)\}_{j\in\Zz_{\ge1}}$ satisfies the ACC under the inclusion of polytopes. Furthermore, for any Newton polytope $\mathcal{N}$ in $\Zzzero$, there are only finitely many vertices of  $\mathcal{N}$.
\end{lem}

\begin{proof}
 Suppose that $\{\mathcal{P}_j\}_{j\in\Zz_{\ge1}}$ does not satisfy the ACC, possibly passing to a subsequence, we may assume that $\{\mathcal{P}_j\}_{j\in\Zz_{\ge1}}$ is strictly increasing. As $\mathbb{Z}_{\geq 0}$ satisfies the DCC, we may find a pair $(i, j)$ such that $i<j$ and $\bm{v}_j\in \bm{v}_i+\mathbb{Z}_{\geq 0}^n$. Thus $\mathcal{P}_j= \mathcal{P}_{j-1}$, a contradiction.
 
 Suppose that $\mathcal{N}$ has infinitely many vertices $\bm{v}_i$, $i\in\Zz_{\ge 1}$. Then the set $\{\mathcal{P}_j=\cup_{i=1}^j (\bm{v}_i+\mathbb{Z}_{\geq 0}^n)\}_{j\in\Zz_{\ge1}}$ is strictly increasing, a contradiction. 
\end{proof}

Theorem \ref{thm: acc newton polytope} is proved in \cite{Ste11} based on some results from Russian literature. The following simple proof of Theorem \ref{thm: acc newton polytope} is kindly provided to us by Chen Jiang.

%[Strong ACC for Newton polytopes]
\begin{thm}[ACC for Newton polytopes]\label{thm: acc newton polytope}
Let $n$ be a positive integer, and $\{\mathcal{N}_i\}_{i\in\Zz_{\ge1}}$ a sequence of Newton polytopes in $\Zzzero$. Then there exists a subsequence $\{\mathcal{N}_{i_j}\}_{j\in\Zz_{\ge1}}$ of Newton polytopes, such that $\mathcal{N}_{i_j}\supseteq\mathcal{N}_{i_{j+1}}$ for every positive integer $j$.
\end{thm}
% $1\leq i_1<i_2<\cdots$
%[Proof of Theorem \ref{thm: acc newton polytope}]
%that $\{\mathcal{N}_i\}_{i\geq 1}$ does not satisfy the PROPERTY
\begin{proof}
Suppose that the theorem does not hold.
Then there exists $i_1\in\Zz_{\ge 1}$ such that $\mathcal{N}_{i_1}\not \supseteq \mathcal{N}_k$ for all positive integers $k>i_1$.
%Note that the sequence of Newton polytopes $\{\mathcal{N}_i\}_{i> i_1}$ does not satisfy the PROPERTY.

Inductively, we may construct a sequence of positive integers $i_1<i_2<\cdots$ such that 
$\{\mathcal{M}_j:=\mathcal{N}_{i_j}\}_{j\in\Zz_{\ge1}}$ satisfies that $\mathcal{M}_{l}\not \supseteq \mathcal{M}_m$ for all positive integers $m>l$. Then $\{\mathcal{Q}_j:=\cup_{k\leq j}\mathcal{M}_{k}\}_{j\in \Zz_{\ge1}}$ is a strictly increasing sequence of Newton polytopes. So the statement is equivalent to the ACC for Newton polytopes. 

By Lemma \ref{lem: weak ACC newton polytopes}, any Newton polytope $\mathcal{Q}$ in $\Zzzero$ can be written as a finite union $\mathcal{Q}=\cup_{\bm{v}} (\bm{v}+\mathbb{Z}_{\geq 0}^n)$, where each $\bm{v}$ is a vertex of $\mathcal{Q}$. So we can find a sequence of vectors $\{\bm{v}_i\}_{i\in\Zz_{>0}}$ and a sequence of positive integers $n_1<n_2<\cdots$, such that $\mathcal{Q}_j=\cup_{i=1}^{n_j} (\bm{v}_i+\mathbb{Z}_{\geq 0}^n).$ This contradicts Lemma \ref{lem: weak ACC newton polytopes}.
\end{proof}

\section{Existence of algebraic elephants}

\begin{lem}\label{lem: power series p-adic hensel}
Let $\hat{I}$ be an ideal of $\Cc[[x_1,\ldots,x_n]]$, and $f\in \Cc[[x_1,\ldots,x_n]]$. Suppose that $f\notin\hat{I}$, then there exists a positive integer $c_1$, such that for any $f_{c_1}\in \Cc[[x_1,\ldots,x_n]]$ with $f_{c_1}-f\in \fm^{c_1}$, $f_{c_1}\notin\hat{I}$, where $\fm=(x_1,\ldots,x_n)$.
\end{lem}

\begin{proof}
Let $\bar{\fm}$ be the maximal ideal of $R:=\Cc[[x_1,\ldots,x_n]]/\hat{I}$, and $\bar{f}$ the image of $f$ in $R$. Suppose on the contrary, for any positive integer $i$, there exists $f_{i}\in \hat{I}$, such that $f_{i}-f\in \fm^{i}$. Then $\bar{f}\in \cap_{i=1}^{\infty}\bar{\fm}^i$. Since $R$ is a Noetherian local ring, by \cite[Corollary 10.20]{AM69}, the $\bar\fm$-adic topology of $R$ is Hausdorff. Thus $\bar{f}\in \cap_{i=1}^{\infty}\bar{\fm}^i=0$, or equivalently, $f\in\hat{I}$, a contradiction.
\end{proof}

\begin{rem}
It is also possible to give a constructive proof of Lemma \ref{lem: power series p-adic hensel} by the standard basis introduced by Hironaka \cite{Hir64}, \cite[Theorem (Hironaka Theorem)]{Bec90}.
\end{rem}

\begin{defn}[{\cite[p~322]{HR64}}]\label{defn: equivalence relations}
Let $\tilde{X}$ (resp. $\tilde{X}'$) be a reduced equidimensional analytic space with the analytic structure sheaf $\mathcal{O}_{\tilde{X}}$ (resp. $\mathcal{O}_{\tilde{X}'}$), and $\widetilde{Z}\subset \tilde{X}$ (resp. $\tilde{Z}'\subset \tilde{X}'$) a closed analytic subspace with the ideal sheaf $\mathcal{I}_{\tilde{Z}}$ (resp. $\mathcal{I}_{\tilde{Z}'}$). By $\tilde{Z}_{(v)}$ (resp. $\tilde{Z}_{(v)}'$) we mean the analytic subspace of $\tilde{X}$ (resp. $\tilde{X}'$) with the ideal sheaf $\mathcal{I}_{\tilde{Z}}^v$ (resp. $\mathcal{I}_{\tilde{Z}'}^v$). We say that $\tilde{Z}\subset \tilde{X}$ is \emph{$v$-equivalent} to $\tilde{Z}'\subset \tilde{X}'$ for some $v\in \Zz_{> 0}$ if there exists an analytic isomorphism $\tilde{Z}_{(v)}\cong \tilde{Z}_{(v)}'$. We say that $\tilde{Z}\subset \tilde{X}$ is \emph{formal equivalent} to $\tilde{Z}'\subset \tilde{X}'$ if there exists an isomorphism of locally ringed spaces $(\tilde{Z}, (\hat{\mathcal{O}}_{\tilde{X}})|_{\tilde{Z}})\cong (\tilde{Z}', (\hat{\mathcal{O}}_{\tilde{X}'})|_{\tilde{Z}'})$ which induces an analytic isomorphism $\tilde{Z}\cong {\tilde{Z}}'$, where $\hat{\mathcal{O}}_{\tilde{X}}$ (resp. $\hat{\mathcal{O}}_{\tilde{X}'}$) is the completion of $\mathcal{O}_{\tilde{X}}$ (resp. $\mathcal{O}_{\tilde{X}'}$) along $\mathcal{I}_{\tilde{Z}}$ (resp. $\mathcal{I}_{\tilde{Z}'}$). We say that $\tilde{Z}\subset \tilde{X}$ is \emph{analytic equivalent} to $\tilde{Z}'\subset \tilde{X}'$ if there exists an analytic isomorphism of locally ringed spaces $(\tilde{Z}, (\mathcal{O}_{\tilde{X}})|_{\tilde{Z}})\cong (\tilde{Z}', (\mathcal{O}_{\tilde{X}'})|_{\tilde{Z}'})$ which induces an analytic isomorphism $\tilde{Z}\cong {\tilde{Z}}'$.
\end{defn}

\begin{rem}
Let $i:Z\hookrightarrow X$ be a closed subscheme, and $\mathcal{F}$ a sheaf on $X$. By $\mathcal{F}|_Z$ we mean $i^{-1}\mathcal{F}$, which has nothing to do with the scheme structure of $Z$.
\end{rem}

We collect some results in \cite{HR64} for the reader's convenience.

\begin{lem}\label{lem: basic properties of equivalence relation}
Let $\tilde{X}$ (resp. $\tilde{X}'$) be a reduced equidimensional analytic space with the analytic structure sheaf $\mathcal{O}_{\tilde{X}}$ (resp. $\mathcal{O}_{\tilde{X}'}$), and $\widetilde{Z}\subset \tilde{X}$ (resp. $\tilde{Z}'\subset \tilde{X}'$) a closed analytic subspace with the ideal sheaf $\mathcal{I}_{\tilde{Z}}$ (resp. $\mathcal{I}_{\tilde{Z}'}$). Then the following statements hold.
\begin{enumerate}
    \item If $\tilde{Z}\subset \tilde{X}$ is formal equivalent to $\tilde{Z}'\subset \tilde{X}'$, then $\tilde{Z}\subset \tilde{X}$ is $v$-equivalent to $\tilde{Z}'\subset \tilde{X}'$ for all $v\in \Zz_{>0}$.
    \item If $\tilde{Z}\subset \tilde{X}$ is analytic equivalent to $\tilde{Z}'\subset \tilde{X}'$, then $\tilde{Z}\subset \tilde{X}$ is formal equivalent to $\tilde{Z}'\subset \tilde{X}'$.
    \item If $\tilde{Z}\subset \tilde{X}$ is analytic equivalent to $\tilde{Z}'\subset \tilde{X}'$, then there exists an analytic neighborhood $\tilde{U}$ of $\tilde{Z}$ in $\tilde{X}$ (resp. $\tilde{U}'$ of $\tilde{Z}'$ in $\tilde{X}'$), such that the isomorphism $(\tilde{Z}, (\mathcal{O}_{\tilde{X}})|_{\tilde{Z}})\cong (\tilde{Z}', (\mathcal{O}_{\tilde{X}'})|_{\tilde{Z}'})$ extends to an isomorphism $(\tilde{U}, (\mathcal{O}_{\tilde{X}})|_{\tilde{U}})\cong (\tilde{U}', (\mathcal{O}_{\tilde{X}'})|_{\tilde{U}'})$.
    \item Suppose that $\tilde{Z}$ (resp. $\tilde{Z}'$) is a closed point on $\tilde{X}$ (resp. $\tilde{X}'$) and  $\tilde{X}\setminus \tilde{Z}$ (resp. $\tilde{X}'\setminus \tilde{Z}'$) is smooth. Then there exists a positive integer $c_2$, such that if $\tilde{Z}\subset \tilde{X}$ is $v$-equivalent to $\tilde{Z}'\subset \tilde{X}'$ for some integer $v\geq  c_2$, then $\tilde{Z}\subset \tilde{X}$ is analytic equivalent to $\tilde{Z}'\subset \tilde{X}'$.
\end{enumerate}
\end{lem}

\begin{proof}
Lemma~\ref{lem: basic properties of equivalence relation}(1-3) follows from \cite[p~322~Remark]{HR64}. By \cite[Definition 3]{Gra62}, any closed point is an exceptional subspace. Lemma~\ref{lem: basic properties of equivalence relation}(4) is a special case of \cite[p~314~Theorem]{HR64}.
\end{proof}

\begin{defn}\label{def: truncation}
Let $\bm{\alpha}:=(\alpha_1,\dots, \alpha_n)\in \Zz_{\geq 0}^n$. We define $|\bm{\alpha}|:=\sum_{i=1}^n\alpha_i$ and $\bm{x}^{\bm{\alpha}} :=x_1^{\alpha_1}\cdots x_n^{\alpha_n}$. Let $f:=\sum_{\bm{\alpha}\in \Zz^n_{\geq 0}} a_{\bm{\alpha}} \bm{x}^{\bm{\alpha}}\in \Cc[[x_1,\dots,x_n]]$ be a formal power series. The \emph{$c$-th truncation} of $f$ is defined as the polynomial $t_c(f):=\sum_{|\bm{\alpha}|\leq c} a_{\bm{\alpha}}\bm{x}^{\bm{\alpha}}$.
\end{defn}

\begin{lem}\label{lem: formal approximation divisors}
Let $\tilde{x}\in \tilde{X}$ be a Gorenstein singularity, and $\tilde{y}\in \tilde{Y}$ an analytic singularity, such that there is a $G_{\tilde{Y}}$-action on $\tilde{Y}$ for some finite group $G_{\tilde{Y}}$. Let $B_{\tilde{Y}}$ be an analytic Cartier divisor on $\tilde{Y}$ locally defined by $(f_{\tilde{Y}}=0)$ near $\tilde{y}$, such that $f_{\tilde{Y}}$ is $G_{\tilde{Y}}$-semi-invariant.

Suppose that $\psi:\tilde{x}\in \tilde{X}\to \tilde{y}\in \tilde{Y}$ is an analytic isomorphism. Then for any positive integer $l$, there is a Cartier divisor $C_{\tilde{X}}$ on $\tilde{X}$, such that
\begin{enumerate}
    \item  the image of $C_{\tilde{X}}$ on $\tilde{Y}$, denoted by $C_{\tilde{Y}}$, is locally defined by $(f_{\tilde{Y}}+h_{\tilde{Y}}=0)$ for some $h_{\tilde{Y}}\in (\fm_{\tilde{y}}^{\rm an})^l$, where $\fm_{\tilde{y}}^{\rm an}$ is the maximal ideal of the ring of analytic power series at $\tilde{y}\in \tilde{Y}$,
    \item $f_{\tilde{Y}}+h_{\tilde{Y}}$ is $G_{\tilde{Y}}$-semi-invariant, and
    \item $C_{\tilde{Y}}$ is analytically isomorphic to $B_{\tilde{Y}}$ near $\tilde{y}$.
\end{enumerate}
\end{lem}
\begin{proof}
Locally near $\tilde{x}$, we may write $\tilde{X}=\Spec\Cc[x_1,\ldots,x_n]/I$, where $I$ is an ideal in $\Cc[x_1,\ldots,x_n]$. Let $\hat{I}:=I\otimes \Cc[[x_1,\ldots,x_n]]$, and $\tilde{f}\in \Cc[[x_1,\ldots,x_n]]/\hat{I}$ the image of $f_{\tilde{Y}}$ under the analytic isomorphism $\psi$. Then $\tilde{f}\notin \hat{I}$. Now $\hat{B}_{\tilde{X}}=\Spec\Cc[[x_1,\ldots,x_n]]/(\hat{I},\tilde{f})$ is the completion of the image $B_{\tilde{X}}$ of $B_{\tilde{Y}}$ on $\tilde{X}$ under $\psi$. Let $c_1$ and $c_2$ be the positive integers constructed in Lemmas \ref{lem: power series p-adic hensel} and \ref{lem: basic properties of equivalence relation}(4) depending only on $\Ii,\tilde{f}$ and $\hat{B}_{\tilde{X}}$. Let $c:=\max\{c_1,c_2,l\}$, $t_c(\tilde{f})$ the $c$-th truncation of $\tilde{f}$, and $C_{\tilde{X}}=\Spec\Cc[x_1,\ldots,x_n]/(I,t_c(\tilde{f}))$ the divisor on $\tilde{X}$. By Lemma \ref{lem: power series p-adic hensel}, $t_c(\tilde{f})\notin \hat{I}$, hence $\dim B_{\tilde{X}}=\dim C_{\tilde{X}}$. Since $\tilde{x}\in B_{\tilde{X}}$ and $\tilde{x}\in C_{\tilde{X}}$ are $c$-equivalent, by Lemma~\ref{lem: basic properties of equivalence relation}(4) and our choice of $c$, $\tilde{x}\in B_{\tilde{X}}$ is analytic equivalent to $\tilde{x}\in C_{\tilde{X}}$, by Lemma~\ref{lem: basic properties of equivalence relation}(3), $C_{\tilde{Y}}$ is analytically isomorphic to $B_{\tilde{Y}}$ in a neighborhood of $\tilde{y}$.

Let $G_{\tilde{X}}$ be the action induced by $G_{\tilde{Y}}$ on $\tilde{X}$. Since $\tilde{f}$ is $G_{\tilde{X}}$-semi-invariant, $t_c(\tilde{f})$ is $G_{\tilde{X}}$-semi-invariant. Since $\tilde{f}-t_c(\tilde{f})\in (\fm_{\tilde{y}}^{\rm an})^l$, $C_{\tilde{Y}}$ is locally defined by $(f_{\tilde{Y}}+h_{\tilde{Y}}=0)$ for some $h_{\tilde{Y}}\in (\fm_{\tilde{y}}^{\rm an})^l$, and $f_{\tilde{Y}}+h_{\tilde{Y}}$ is $G_{\tilde{Y}}$-semi-invariant. 
\end{proof}

\begin{thm}
Let $x\in X$ be a threefold terminal singularity which is not smooth.
Then there exists an elephant $C$ of $x\in X$.
\end{thm}
\begin{proof}
By the classification of threefold terminal singularities (cf. \cite[(6.1) Theorem]{Rei87}, \cite[Theorems 12,23,25]{Mor85}), $x\in X$ is analytically isomorphic to a threefold singularity
$y\in Y$ defined by
$$(\phi(x_1,x_2,x_3,x_4)=0)\subset (\mathbb C^4\ni o)/G_{\tilde{Y}}$$ for some cyclic group action $G_{\tilde{Y}}$ on $\Cc^d$ that fixes the origin, and some $G_{\tilde{Y}}$-invariant analytic power series $\phi$.

Let $\tilde{X}\ni \tilde{x}\to X\ni x$ and $\tilde{Y}\ni \tilde{y}\to Y\ni y$ be the index one covers of $x\in X$ and $y\in Y$ respectively. Then we have a natural induced analytic isomorphism $\psi:\tilde{X}\ni \tilde{x}\cong \tilde{Y}\ni \tilde{y}$. Let $B_{\tilde{Y}}$ be an analytic elephant of $\tilde{y}\in \tilde{Y}$ locally defined by $p(x_1,x_2,x_3,x_4)=0$ for some semi-invariant analytic power series $p\in \mathfrak{m}\setminus \mathfrak{m}^2$ (cf. \cite[(6.4)(B)]{Rei87}), where $\mathfrak{m}$ denotes the maximal ideal of $o\in \Cc^4$. Note that $B_Y$ is an analytic elephant of $y\in Y$, where $B_Y$ is the image of $B_{\tilde{Y}}$ on $Y$. By Lemma \ref{lem: formal approximation divisors}, there exists a Cartier divisor $C_{\tilde{X}}$ on $\tilde{X}$ whose image under $\psi$ on $\tilde{Y}$ is an analytic Cartier divisor $C_{\tilde{Y}}$, which is locally defined by $(p+h)$ for some semi-invariant analytic power series $h\in \fm^2$, such that $C_{\tilde{Y}}$ is analytically isomorphic to $B_{\tilde{Y}}$ near $\tilde{y}$. In particular, $C_{\tilde{Y}}$ is canonical near $\tilde{y}$. Hence $C_{\tilde{X}}$ is also canonical, thus normal, near $\tilde{x}$. By \cite[Proposition~5.46]{KM98}, $(\tilde{X},C_{\tilde{X}})$ is canonical near $\tilde{x}$. Let $C$ (resp. $C_Y$) be the image of $C_{\tilde{X}}$ (resp. $C_{\tilde{Y}}$) on $X$ (resp. $Y$).

Since $\frac{p+h}{p}$ is $G_{\tilde{Y}}$-invariant, it defines an analytic Cartier divisor on $Y$. Thus $K_Y+C_Y\sim K_Y+B_Y$ is an analytic Cartier divisor near $y$, by \cite[Lemma~5.1]{Kaw88}, $K_X+C$ is Cartier near $x$. By \cite[Proposition~5.20]{KM98}, $(X,C)$ is plt near $x$ as $(\tilde{X},C_{\tilde{X}})$ is canonical near $\tilde{x}$. We conclude that $(X,C)$ is canonical near $x$ as $K_X+C$ is Cartier near $x$. 
\end{proof}

\section{Admissible weighted blow-ups extracting prime divisors}

Recall that for any vector $\bm{\alpha}=(\alpha_1,\dots, \alpha_d)\in \Zz_{\geq 0}^d$, we define $|\bm{\alpha}|:=\sum_{i=1}^d\alpha_i$ and $\bm{x}^{\bm{\alpha}} :=x_1^{\alpha_1}\dots x_d^{\alpha_d}$.

%\luo{c-truncation and w-weighted leading term, find good notation}

\begin{defn}\label{def: w-weighted leading term} Let $f:=\sum_{\bm{\alpha}\in \Zz^n_{\geq 0}} a_{\bm{\alpha}} \bm{x}^{\bm{\alpha}}\in \Cc[[x_1,\dots,x_d]]$ be a formal power series and $w=(w(x_1),\dots,w(x_d))\in \Qq^d_{> 0}$ a weight. The \emph{$w$-weighted leading term} of $f$ is defined as the polynomial $f_w:=\sum_{w(\bm{x}^{\bm{\alpha}})=w(f)} a_{\bm{\alpha}}\bm{x}^{\bm{\alpha}}$. We call $h\in \Cc[[x_1,\dots,x_d]]$ a \emph{$w$-weighted homogeneous} polynomial if $h=h_w$.
\end{defn}

Let $w:=(w_1,\dots,w_d)\in \Qq_{>0}^d$ be a weight and $\mu\in \Qq_{>0}$. We set $\mu w:=(\mu w_1,\dots,\mu w_d)$. It follows directly from the definition that for any $h\in \Cc[[x_1,\dots,x_d]]$, $h_{w}=h_{\mu w}$.

\begin{lem}\label{lem: product rule for w-weighted leading term}
Let $w=w(x_1,\dots,x_d)\in \Qq_{>0}^n$ be a weight and $h\in \Cc[[x_1,\dots,x_d]]$ a $w$-weighted homogeneous polynomial. If $h=p q$ for some $p,q\in \Cc[[x_1,\dots,x_d]]$, then $h=p_wq_w$.
\end{lem}

\begin{proof}
Possibly replacing $w$ with $\mu w$ for some $\mu\in \Zz_{>0}$, we may assume that $w\in \Zz_{>0}^d$. Let $\lambda\in \Cc$ be a variable. We have $h(\lambda x_1,\dots,\lambda x_d)=\lambda^{w(h)}h(x_1,\dots,x_d)$. Since $h=pq$, $w(h)=w(p)+w(q)$, and
$$    h=\lim_{\lambda\to 0}\frac{h(\lambda x_1,\dots,\lambda x_d)}{\lambda^{w(h)}}=\lim_{\lambda\to 0}\frac{p(\lambda x_1,\dots,\lambda x_d)}{\lambda^{w(p)}}\cdot \frac{q(\lambda x_1,\dots,\lambda x_d)}{\lambda^{w(q)}}=p_wq_w.$$
\end{proof}

Recall that for a weight $w=(w_1,\dots,w_d)\in \Zz^d_{>0}$, there exists an equivalence relation $\sim_w$ on $\Cc^d\setminus \{o\}$: we say that $(x_1,\dots,x_d)\sim_w (y_1,\dots,y_d)$ if $y_i=\lambda^{w_i}x_i$ for all $1\leq i\leq d$ and some complex number $\lambda\neq 0$. We denote by $[x_1,\dots,x_d]_w$ the equivalence class of $(x_1,\dots,x_d)$, and we have a natural quotient morphism $\pi_w: \Cc^d\setminus\{o\}\to \mathbf{P}(w_1,\dots,w_n)$ such that $\pi_w(x_1,\dots,x_d)=[x_1,\dots,x_d]_w$.

\begin{lem}\label{lem: exceptional locus of weighted blow-ups}
Let $X\subset \Cc^d$ be a closed subvariety defined by the ideal $I\subset \Cc[x_1,\dots,x_d]$ near $o\in \Cc^d$, $w=(w_1,\dots,w_d)\in \Zz_{>0}^d$ a weight, $I_w$ the ideal generated by the $w$-weighted leading terms of the elements in $I$, and $X_w\subset \Cc^d$ the subvariety defined by $I_w$. Let $f': W\to \Cc^d$ be the weighted blow-up at $o\in\Cc^d$ with respect to the coordinates $x_1,\dots,x_d$ with the exceptional divisor $E'\cong \mathbf{P}(w_1,\dots,w_d)$, and $f: Y\to X$ the induced weighted blow-up of $X\ni o$ with the exceptional divisor $E:=E'|_Y$, where $Y={f'}_*^{-1} X$. Let $\pi_w: \Cc^d\setminus \{o\}\to \mathbf{P}(w_1,\dots,w_d)$ be the natural quotient morphism. Then $E=\pi_w(X_w)\subset \mathbf{P}(w_1,\dots,w_d)$.
\end{lem}

\begin{proof}
The weighted blow-up at $o\in \Cc^d$ is defined as the closure of the graph $$\Ii_{\pi_w}:=\overline{\{(x_1,\dots,x_d)\times[x_1,\dots,x_d]_w\mid (x_1,\dots,x_d)\in \Cc^d\setminus\{o\}\}}\subset \Cc^d\times \mathbf{P}(w_1,\dots,w_d)\to \Cc^d,$$ where the blow-up morphism $f': \Ii_{\pi_w}=W\to \Cc^d$ is induced by the projection $\Cc^d\times \mathbf{P}(w_1,\dots,w_d)\to \Cc^d$. Under this identification, for any $[y_1,\dots,y_d]_w\in E=E'|_Y$, there exist sequences $\{(y_{1,i},\dots,y_{d,i})\}_{i=1}^{\infty}$ and  $\{\lambda_i\}_{i=1}^\infty$, such that $\lim_{i\to +\infty}y_{k,i}=y_k$ for $1\leq k\leq d$, $\lim_{i\to +\infty} \lambda_i=0$ and $(\lambda_i^{w_1}y_{1,i},\dots,\lambda_i^{w_d}y_{d,i})\in X$ for each $i$. It follows that for any $h\in I$ and any $i$, $h(\lambda_i^{w_1}y_{1,i},\dots,\lambda_i^{w_d}y_{d,i})=0$, hence $$\lim_{i\to \infty} \lambda_i^{-w(h)}h(\lambda_i^{w_1}y_{1,i},\dots,\lambda_i^{w_d}y_{d,i})=h_w(y_1,\dots,y_d)=0.$$ This implies that $E\subset \pi_w(X_w)$.

On the other hand, for any $(y_1,\dots,y_d)\in X_w$, let $h_1,\dots,h_m\in I$ be a set of generators. For each $1\leq k\leq m$, we have $$|\lambda^{-w(h_k)}h_k(\lambda^{w_1}y_1,\dots,\lambda^{w_d}y_d)|=|\lambda^{-w(h_k)}h_k(\lambda^{w_1}y_1,\dots,\lambda^{w_d}y_d)-(h_k)_w(y_1,\dots,y_d)|=O(\lambda^t)$$ for some $t\in \Zz_{>0}$ as $\lambda\to 0$. Hence there exist a sequence of points $\{(y_{1,i},\dots,y_{d,i})\}_{i=1}^{\infty}\subset \Cc^d$ and a sequence of complex numbers $\{\lambda_i\}_{i=1}^\infty$, such that $\lim_{i\to +\infty} \lambda_i= 0$, and for each $i$, $(y_{1,i},\dots,y_{d,i})$ is a common zero of the set of polynomials $$\{\lambda^{-w(h_k)}h_k(\lambda^{w_1}x_1,\dots,\lambda^{w_d}x_d)\mid 1\leq k\leq d\},$$ such that $\lim_{i\to +\infty}y_{k,i}=y_k$ for $1\leq k\leq d$. This implies that $\pi_w(X_w)\subset E$.
\end{proof}

\begin{lem}\label{lem: irr by dominant morphism}
Let $\pi: \tilde{X}\to X$ be a surjective morphism of schemes such that the scheme theoretic image of $\tilde{X}$ under $\pi$ is $X$. If $\tilde{X}$ is integral, then $X$ is integral.
\end{lem}

\begin{proof}
Since $\tilde{X}$ is reduced, $\pi$ factors through $X_{\mathrm{red}}\to X$, hence $X_{\mathrm{red}}=X$ since the scheme theoretic image of $\tilde{X}$ under $\pi$ is $X$, here $X_{\mathrm{red}}$ is the reduced scheme associated to $X$ (cf. \cite[\S2,~Exercise~2.3]{Har77}). Let $\eta_{\tilde{X}}$ be the generic point of $\tilde{X}$. Since $f$ is dominant, $f(\eta_{\tilde{X}})$ is a point whose Zariski closure equals $X$. This implies that $X$ is irreducible.
\end{proof}

\begin{lem}\label{lem: weight leading term irr imply irr exceptional divisor}
Let $X\ni x$ be a germ such that analytically locally,
$$(X\ni x)\cong (\phi_1=\dots =\phi_m=0)\subset (\mathbb C^d\ni o)/\frac{1}{n}(a_1,\dots,a_d)$$ for some positive integers $m,d,n$ and semi-invariant analytic power series $\phi_i\in \fm_o$, where $1\leq i\leq m$, $m<d$ and $\fm_o\subset \Cc\{x_1,\dots,x_d\}$ is the maximal ideal of $o\in \Cc^d$.

Let $w=(w(x_1),\dots,w(x_d))\in \frac{1}{n}\Zz^d_{>0}$ be an admissible weight of $X\ni x$, and $\{\phi_{i,w}\}_{1\leq i\leq m}$ the corresponding $w$-weighted leading terms of $\{\phi_{i}\}_{1\leq i\leq m}$. If the variety defined by $(\phi_{1,w}=\cdots =\phi_{m,w}=0)$ is integral near $o\in \Cc^d$, then the exceptional divisor of the weighted blow-up $f: Y\to X$ with the weight $w$ at $x\in X$ is a prime divisor.
\end{lem}

\begin{proof}
Let $\pi: \tilde{X}\to X$ be the index one cover such that analytically locally, $$(\tilde{X}\ni \tilde{x})\cong (\phi_1=\dots =\phi_m=0)\subset (\mathbb C^d\ni o).$$ Let $\tilde{f}: \tilde{Y}\to \tilde{X}$ be the weighted blow-up with the weight $\tilde{w}:=nw\in \Zz_{>0}$ at $\tilde{x}\in \tilde{X}$. We have the following fiber product commutative diagram:
$$\xymatrix@R=2em{\tilde{E} \ar[d]_{\pi_E} \ar@{^{(}->}[r]
 &\tilde{Y} \ar[d]_{\pi_Y} \ar[r]^{\tilde{f}} &   \tilde{X} \ar[d]_{\pi}\\
 E \ar@{^{(}->}[r]
 & Y \ar[r]_{f} & X
}
$$
where $E$ (resp. $\tilde{E}$) is the exceptional divisor of $f$ (resp. $\tilde{f}$), and $\pi_E,\pi_Y$ are finite surjective morphisms induced by $\pi$. By Lemma~\ref{lem: irr by dominant morphism}, it suffices to prove that $\tilde{E}$ is a prime divisor. Possibly replacing $f: Y\to X,x, E, w$ with $\tilde{f}: \tilde{Y}\to \tilde{X},\tilde{x},\tilde{E}, \tilde{w}$, we may assume that $n=1$.

\medskip

The weighted blow-up at $o\in \Cc^d$ is defined as the closure of the graph $$\Ii_{\pi_w}:=\overline{\{(x_1,\dots,x_d)\times[x_1,\dots,x_d]_w\mid (x_1,\dots,x_d)\in \Cc^d\setminus\{0\}\}}\subset \Cc^d\times \mathbf{P}(w_1,\dots,w_d)\to \Cc^d,$$ and the morphism $f': \Ii_{\pi_w}\to \Cc^d$ is induced by the projection $\Cc^d\times \mathbf{P}(w_1,\dots,w_d)\to \Cc^d$. By Lemma~\ref{lem: exceptional locus of weighted blow-ups}, the exceptional divisor $E$ of $f: Y\to X$ is defined by $(\phi_{1,w}=\dots =\phi_{m,w}=0)\subset \mathbf{P}(w_1,\dots,w_d)$.

Again we consider the quotient morphism $\pi_w: \Cc^d\setminus\{o\}\to \mathbf{P}(w_1,\dots,w_d)$. Let $C_{w}E$ be the variety defined by $(\phi_{1,w}=\dots =\phi_{m,w}=0)\subset \Cc^d$. Note that we have a quotient morphism $C_wE\setminus\{o\}\to E$ induced by $\pi_w$. By our assumption, $C_wE$ is integral near $o$, hence $C_wE\setminus \{o\}$ is integral near $o$, by Lemma~\ref{lem: irr by dominant morphism}, $E$ is a prime divisor. This finishes the proof.
\end{proof}

\begin{lem}\label{lem: algebraic approximation of weighted blow-up}
Let $X\ni x$ be a germ such that analytically locally,
$$\psi:(X\ni x)\cong (\phi_1=\dots =\phi_m=0)\subset (\mathbb C^d\ni o)/\frac{1}{n}(a_1,\dots,a_d)$$ for some positive integers $m,d,n$, some integers $a_1,\dots,a_d$ and some \emph{semi-invariant} analytic power series $\phi_1\dots\phi_m\in\Cc\{x_1,\dots,x_d\}$. Let $B_1,\cdots, B_k$ be $\Qq$-Cartier Weil divisors near $x\in X$.

Let $w=(w(x_1,\dots,w(x_d))\in \frac{1}{n}\Zz^d_{>0}$ be an admissible weight of $X\ni x$. Suppose that the weighted blow-up with the weight $w$ extracts an analytic prime divisor. Then there exists a prime divisor $E$ over $X\ni x$, such that $a(E,X,0)=1+w(X\ni x)$, and $w(B_i)=\mult_E B_i$ for $1\leq i\leq k$. In particular, $w(B)=\mult_E B$ for any $\Rr$-Cartier $\Rr$-divisor $B:=\sum_{i=1}^k b_iB_i$.
\end{lem}

\begin{proof}
Under the analytic isomorphism $\psi$, each $B_i$ is locally defined by $(h_i=0)$ for some $h_i\in \Cc\{x_1,\dots,x_d\}$ that is semi-invariant with respect to the cyclic group action on $\Cc^d\ni o$.

Let $\pi:\tilde{X}\ni \tilde{x}\to X\ni x$ be the index one cover. The analytic isomorphism $\psi$ lifts to an equivariant analytic isomorphism (still we denote it as $\psi$) $(\tilde{X}\ni \tilde{x})\cong (\phi_1=\dots =\phi_m=0)\subset(\Cc^d\ni o)$ with respect to the corresponding cyclic group actions on $\tilde{X}\ni \tilde{x}$ and $\Cc^d\ni o$. By Lemma~\ref{lem: formal approximation divisors}, for each $l\in \Zz_{>0}$, there exist some Cartier divisors $C_{\tilde X,i}$ near $\tilde{x}\in \tilde{X}$ for $1\leq i\leq d$, such that $\psi(C_{\tilde{X},i})$ is locally defined by $(x_i+q_i=0)$ for some $q_i\in \mathfrak{m}_o^l$, and $x_i+q_i$ is semi-invariant with respect to the cyclic group action $\Cc^d\ni o$. Under the analytic coordinate change $$y_1:=x_1+q_1,\dots,y_d:=x_d+q_d,$$ we may write $\phi_i(x_1,\dots,x_d)=\phi_i'(y_1,\dots,y_d)$ and $h_j(x_1,\dots,x_d)=h_j'(y_1,\dots,y_d)$ for $1\leq i\leq m$ and $1\leq j\leq k$. By choosing $l$ big enough, we may assume that $(\phi_i(y_1,\dots,y_d))_{w}=(\phi_i'(y_1,\dots,y_d))_{w}$ and $(h_j(y_1,\dots,y_d))_{w}=(h_j'(y_1,\dots,y_d))_{w}$ for $1\leq i\leq m$ and $1\leq j\leq k$, hence $w(\phi_i')=w(\phi_i)$, $w'(h_j')=w(h_j)$ for each $i,j$. By Lemmas~\ref{lem: exceptional locus of weighted blow-ups} and \ref{lem: weight leading term irr imply irr exceptional divisor}, the exceptional divisor $E'$ of the weighted blow-up with the weight $w$ at $x\in X$ with respect to the local coordinates $y_1,\dots,y_d$ is a prime divisor. Moreover, $a(E',X,0)=1+w(X\ni x)$ and $\mult_{E'}B_i=w(B_i)$ for $1\leq i\leq k$. Let $E:=E'$. We are done.
\end{proof}

\begin{lem}\label{lem: irreducibility of exceptional divisors extracted by weighted blow up}
Let $f: Y\to X$ be a divisorial contraction of a prime divisor over $X\ni x$ and $r_1,r_2,a,d,n$ the corresponding positive integers as in Theorem~\ref{thm: kaw05 1.2 strenghthened}(1). For any weight $w':=\frac{1}{n}(r_1',r_2',a',n)\in\frac{1}{n}\mathbb Z_{>0}^4$ such that $r_1'+r_2'=a'dn$, $a'\leq a$ and $a'\equiv br_1'\mod n$, the weighted blow-up $f':Y'\to X$ with the weight $w'$ extracts an analytic prime divisor, and $w'(X\ni x)=\frac{a'}{n}$. 
\end{lem}

\begin{proof}
By Theorem~\ref{thm: kaw05 1.2 strenghthened}(1), analytically locally, we have $$(X\ni x)\cong (\phi:= x_1x_2+g(x_3^n,x_4)=0)\subset (\mathbb C^4\ni o)/\frac{1}{n}(1,-1,b,0),$$ and $f: Y\to X$ is a weighted blow-up with the weight $w:=\frac{1}{n}(r_1,r_2,a,n)$, such that 
\begin{itemize}
        \item $x_3^{dn}\in g(x_3^n,x_4)$,
        \item $nw(\phi)=r_1+r_2=adn$, and
        \item $a\equiv br_1\mod n$.
\end{itemize}

    Consider the weights $v'=(v'(x_3),v'(x_4))=\frac{1}{n}(a',n)$ and $v=(v(x_3),v(x_4))=\frac{1}{n}(a,n)$. Since $v'\succeq \frac{a'}{a}v$ (see Definition~\ref{defn: compare two weights}), $$a'd=v'(x_3^{nd})\geq v'(g(x_3^n,x_4))\geq \frac{a'}{a}v(g(x_3^n,x_4))=\frac{a'}{a}w(g(x_3^n,x_4))=a'd,$$ and $v'(g(x_3^n,x_4))=a'd$. It follows that $w'(\phi)=\frac{1}{n}(r_1'+r_2')=a'd$ and $$\phi_{w'}=x_1x_2+x_3^{nd}+h(x_3,x_4)$$ for some analytic power series $h(x_3,x_4)=(g(x_3,x_4)-x_3^{nd})_{w'}$. Since $(\phi_{w'}=0)$ defines an integral scheme near $o\in \Cc^n$, by Lemma~\ref{lem: weight leading term irr imply irr exceptional divisor}, the exceptional divisor of $f'$ is a prime divisor. Since $w'(X\ni x)=\frac{1}{n}(r_1'+r_2'+a'+n)-w'(\phi)-1$ and $w'(\phi)=\frac{1}{n}(r_1'+r_2')$, $w'(X\ni x)=\frac{a'}{n}$.
\end{proof}

\begin{lem}\label{lem: cD irreducible case 1}
Let $f: Y\to X$ be a divisorial contraction of a prime divisor over $X\ni x$ and $r,a,d$ the corresponding positive integers as in Theorem~\ref{thm: kaw05 1.2 strenghthened}(2.1).
\begin{enumerate}
    \item For any weight $w':=(r'+1,r',a',1)\in\mathbb Z_{>0}^4$ such that $2r'+1=a'd$ and $2\leq a'<a$, the weighted blow-up $f':Y'\to X$ with the weight $w'$ extracts an analytic prime divisor, and $w'(X\ni x)=a'$.
    \item If $a\geq 4$, then the weighted blow-up $f'':Y''\to X$ with the weight $w'':=(d,d,2,1)$ extracts an analytic prime divisor, and $w''(X\ni x)=2$. 
\end{enumerate}
\end{lem}

\begin{proof}
By Theorem~\ref{thm: kaw05 1.2 strenghthened}(2.1), analytically locally, we have $$(X\ni x)\cong (\phi:= x_1^2+x_1q(x_3,x_4)+x_2^2x_4+\lambda x_2x_3^2+\mu x_3^3+p(x_2,x_3,x_4)=0)\subset (\mathbb C^4\ni o),$$ where $p(x_2,x_3,x_4)\in(x_2,x_3,x_4)^4$, and $f: Y\to X$ is a weighted blow-up with the weight $w:=(r+1,r,a,1)$, such that 
\begin{itemize}
    \item $w(\phi)=w(x_2^2x_4)=w(x_3^d)=2r+1=ad$, where $a,d$ are odd integers,
    \item if $q(x_3,x_4)\not=0$, then $w(x_1q(x_3,x_4))=2r+1$,
    \item $\mu' x_3^d\in \phi$ for some $\mu'\neq 0$ and $d\geq 3$, if $d=3$, then $\mu'=\mu$, and
    \item if $d>3$, then $\mu,\lambda=0$.
\end{itemize}

Consider the weights $v=(v(x_3),v(x_4))=(a,1)$, $v'=(v'(x_3),v'(x_4))=(a',1)$ and $v''=(v''(x_3),v''(x_4))=(2,1)$, we have $v'\succeq \frac{a'}{a}v$ and $v''\succeq \frac{2}{a}v$ (see Definition~\ref{defn: compare two weights}), and for any monomial term $x_3^{k}x_4^l$, $v'(x_3^{k}x_4^l)=\frac{a'}{a}v(x_3^{k}x_4^l)$ (resp. $v''(x_3^{k}x_4^l)=\frac{2}{a}v(x_3^{k}x_4^l)$ when $a\geq 3$) if and only if $l=0$. 

We may write $$p(x_2,x_3,x_4)=g_0(x_3,x_4)+x_2g_1(x_3,x_4)+x_2^2g_2(x_2,x_3,x_4)$$ for some analytic power series $g_0(x_3,x_4),g_1(x_3,x_4)$ and $g_2(x_2,x_3,x_4)$. Since $w(\phi)=2r+1$, $w(g_0(x_3,x_4))=v(g_0(x_3,x_4))\geq 2r+1$ and $w(x_2g_1(x_3,x_4))=r+v(g_1(x_3,x_4))\geq 2r+1$, hence
\begin{align}\label{align: equation for cD case 1}
    v(g_0(x_3,x_4))\geq 2r+1,\ \ \ and \ \ \ v(g_1(x_3,x_4))\geq r+1.
\end{align}
Since $p(x_2,x_3,x_4)\in (x_2,x_3,x_4)^4$, $\mult_o g_1(x_3,x_4)\geq 3$ and $\mult_o g_2(x_2,x_3,x_4)\geq 2$.

\medskip

We now prove part (1) of the Lemma.

We claim that $w'(\mu x_3^3+p(x_2,x_3,x_4))= 2r'+1$ and $(\mu x_3^3+p(x_2,x_3,x_4))_{w'}=\mu' x_3^d$. By (\ref{align: equation for cD case 1}), we have $$w'(g_0(x_3,x_4))=v'(g_0(x_3,x_4))\geq \frac{a'}{a}v(g_0(x_3,x_4))\geq \frac{a'}{a}(2r+1)=a'd=2r'+1,$$ and $x_3^{d}$ is the only possible monomial term in $g_0(x_3,x_4)$ whose weight equals $2r'+1$. Moreover, when $d=3$, $w'(g_0(x_3,x_4))>2r'+1$. We may write $$g_1(x_3,x_4)=h_0(x_3)+x_4h_1(x_3,x_4)$$ for some analytic power series $h_0(x_3)$ and $h_1(x_3,x_4)$. For any monomial term $x_3^{l_1}\in g_1(x_3,x_4)$, by (\ref{align: equation for cD case 1}), we have $v(x_3^{l_1})=al_1\geq r+1=\frac{ad}{2}+\frac{1}{2}$, hence $l_1\geq \lfloor \frac{d}{2}\rfloor+1$, and $$w'(x_3^{l_1})=a'l_1\geq a'(\lfloor \frac{d}{2}\rfloor+1)\geq a'(\frac{d+1}{2})=\frac{a'd-1}{2}+\frac{1+a'}{2}>r'+1,$$ where the last inequality follows from $a'\geq 2>1$. Hence $v'(h_0(x_3))>r'+1$. By (\ref{align: equation for cD case 1}), $v(h_1(x_3,x_4))\geq r$, hence $$v'(h_1(x_3,x_4))\geq \frac{a'}{a}v(h_1(x_3,x_4))\geq \frac{a'}{a}r=r'+\frac{1}{2}(1-\frac{a'}{a})>r'.$$ Now we have
\begin{align*}
    w'(x_2g_1(x_3,x_4))&\geq r'+\min\{v'(h_0(x_3)),1+v'(h_1(x_3,x_4))\}>2r'+1,\ and\\
    w'(x_2^2g_2(x_2,x_3,x_4))&=2r'+v'(g_2(x_2,x_3,x_4))\geq 2r'+2>2r'+1.
\end{align*}
This proves the claim.

If $q(x_3,x_4)\neq 0$, then $w(x_1q(x_3,x_4))=2r+1$, hence $v(q(x_3,x_4))=r$. It follows that
\begin{align*}
    w'(x_1q(x_3,x_4))&=r'+1+v'(q(x_3,x_4))\geq r'+1+\frac{a'}{a}v(q(x_3,x_4))\geq r'+1+\frac{a'}{a}r\\
    &=r'+1+\frac{a'}{a}\frac{ad-1}{2}=2r'+1+\frac{1}{2}(1-\frac{a'}{a})>2r'+1.
\end{align*}

We claim that $w'(\lambda x_2x_3^2)>2r'+1$. If $d>3$, then $\lambda =0$, we are done. We may assume that $d=3$, then $2r'+1=3a'$. We have $$w'(x_2x_3^2)=r'+2a'=r'+\frac{2}{3}(2r'+1)=2r'+\frac{r'+2}{3}>2r'+1,$$ where the last inequality follows from $r'=\frac{3a'-1}{2}>1$.

Since $w'(x_2^2x_4)=2r'+1$ and $w'(x_1^2)=2r'+2>2r'+1$, we have $w'(\phi)=2r'+1$ and $$\phi_{w'}=x_2^2x_4+\mu'x_3^d,\ \ 0\neq\mu'\in \Cc.$$ Since $(\phi_{w'}=0)$ defines an integral scheme near $o\in \Cc^5$, by Lemma~\ref{lem: weight leading term irr imply irr exceptional divisor}, the exceptional divisor of $f'$ is a prime divisor. Since $w'(X\ni x)=r'+1+r'+a'+1-w'(\phi)-1$ and $w'(\phi)=2r'+1$, $w'(X\ni x)=a'$.

\medskip

We now prove part (2) of the Lemma.

We claim that $w''(\mu x_3^3+p(x_2,x_3,x_4))=2d$ and $(\mu x_3^3+p(x_2,x_3,x_4))_{w''}=\mu' x_3^d$. By (\ref{align: equation for cD case 1}), we have $$w''(g_0(x_3,x_4))=v''(g_0(x_3,x_4))\geq \frac{2}{a}v(g_0(x_3,x_4))\geq \frac{2}{a}(2r+1)=2d,$$ and $x_3^{d}$ is the only possible monomial term in $g_0(x_3,x_4)$ whose weight equals $2r'+1$. Moreover, when $d=3$, $w''(g_0(x_3,x_4))>2d$. By (\ref{align: equation for cD case 1}) again, 
\begin{align*}
    w''(x_2g_1(x_3,x_4))&= d+v''(g_1(x_3,x_4))\geq d+\frac{2}{a}(r+1)=d+\frac{ad+1}{a}>2d,\ and\\
    w''(x_2^2g_2(x_2,x_3,x_4))&=2d+v''(g_2(x_2,x_3,x_4))>2d.
\end{align*}
This proves the claim.

If $q(x_3,x_4)\neq 0$, then again we have $v(q(x_3,x_4))=r$. Since $ad=2r+1$ and $\gcd(r,2r+1)=1$, $q(x_3,x_4)=x_4q'(x_3,x_4)$ for some analytic power series $q'(x_3,x_4)$ such that $v(q'(x_3,x_4))=r-1$. It follows that
\begin{align*}
    w''(x_1q(x_3,x_4))&=d+1+v''(q'(x_3,x_4))\geq d+1+\frac{2}{a}v(q'(x_3,x_4))\\
    &=d+1+\frac{2}{a}(\frac{ad-1}{2}-1)=2d+(1-\frac{3}{a})>2d,
\end{align*}
Where the last inequality follows from $a\geq 4$.

We claim that $w''(\lambda x_2x_3^2)>2d$. If $d>3$, then $\lambda =0$, we are done. We may assume that $d=3$. Now $$w''(x_2x_3^2)=7>6=2d,$$ this proves the claim.

Since $w''(x_2^2x_4)=2d+1>2d$ and $w''(x_1^2)=2d$, we have $w''(\phi)=2d$ and $$\phi_{w''}=x_1^2+\mu'x_3^d,\ \ 0\neq \mu'\in \Cc$$ here $d$ is an odd integer. Since $(\phi_{w''}=0)$ defines an integral scheme near $o\in \Cc^5$, by Lemma~\ref{lem: weight leading term irr imply irr exceptional divisor}, the exceptional divisor of $f''$ is a prime divisor. Since $w''(X\ni x)=d+d+2+1-w''(\phi)-1$ and $w''(\phi)=2d$, $w''(X\ni x)=2$.
\end{proof}

\begin{lem}\label{lem: cD irreducible case 2}
Let $f: Y\to X$ be a divisorial contraction of a prime divisor over $X\ni x$ and $r,a,d$ the corresponding positive integers as in Theorem~\ref{thm: kaw05 1.2 strenghthened}(2.2).
\begin{enumerate}
    \item For any weight $w':=(r'+1,r',a',1,r'+2)\in\mathbb Z_{>0}^5$ such that $r'+1=a'd$ and $2\leq a'<a$, the weighted blow-up $f':Y'\to X$ with the weight $w'$ extracts an analytic prime divisor, and $w'(X\ni x)=a'$. 
    \item If $a\geq 2$, then the weighted blow-up $f'':Y''\to X$ with the weight $w'':=(d,d,1,1,d)$ extracts an analytic prime divisor, and $w''(X\ni x)=1$. 
\end{enumerate}
\end{lem}

\begin{proof}
By Theorem~\ref{thm: kaw05 1.2 strenghthened}(2.2), analytically locally, we have  $$(X\ni x)\cong \binom{\phi_1:= x_1^2+x_2x_5+p(x_2,x_3,x_4)=0}{\phi_2:= x_2x_4+x_3^d+q(x_3,x_4)x_4+x_5=0}\subset(\mathbb C^5\ni o),$$where $p(x_2,x_3,x_4)\in(x_2,x_3,x_4)^4$, and $f: Y\to X$ is a weighted blow-up with the weight $w:=(r+1,r,a,1,r+2)$, such that 
\begin{itemize}
    \item $w(\phi_2)=\frac{1}{2}w(\phi_1)=r+1=ad$,
    \item $d\geq 2$, and
    \item if $q(x_3,x_4)\not=0$, then $w(q(x_3,x_4)x_4)=r+1$.
\end{itemize}

Consider the weights $v=(v(x_3),v(x_4))=(a,1)$, $v'=(v'(x_3),v'(x_4))=(a',1)$ and $v''=(v''(x_3),v''(x_4))=(1,1)$, we have $v'\succeq \frac{a'}{a}v$ and $v''\succeq \frac{1}{a}v$ (see Definition~\ref{defn: compare two weights}), and for any monomial term $x_3^{k}x_4^l$, $v'(x_3^{k}x_4^l)=\frac{a'}{a}v(x_3^{k}x_4^l)$ (resp. $v''(x_3^{k}x_4^l)=\frac{1}{a}v(x_3^{k}x_4^l)$) if and only if $l=0$.

We may write $$p(x_2,x_3,x_4)=g_0(x_3,x_4)+x_2g_1(x_3,x_4)+x_2^2g_2(x_3,x_4)+x_2^3g_3(x_2,x_3,x_4)$$ for some analytic power series $g_0(x_3,x_4),g_1(x_3,x_4),g_2(x_3,x_4)$ and $g_3(x_2,x_3,x_4)$. Since $w(p(x_2,x_3,x_4))\geq w(\phi_1)=2(r+1)$, we have
\begin{align}\label{align: equation for cD case 2}
    v(g_0(x_3,x_4))\geq 2(r+1),\ \ v(g_1(x_3,x_4))\geq r+2,\ \ and\ \ v(g_2(x_3,x_4))\geq 2.
\end{align}
Since $p(x_2,x_3,x_4)\in (x_2,x_3,x_4)^4$, we have $\mult_o g_0(x_3,x_4)\geq 4$, $\mult_o g_1(x_3,x_4)\geq 3$, $\mult_o g_2(x_3,x_4)\geq 2$ and $\mult_o g_3(x_2,x_3,x_4)\geq 1$.

\medskip

We now prove part (1) of the Lemma. 

We claim that $(p(x_2,x_3,x_4))_{w'}=\mu x_3^{2d}+\lambda_1 x_2x_3^dx_4+\lambda_2 x_2^2x_4^2$ for some $\mu,\lambda_1,\lambda_2\in \Cc$. First we note that $$w'(x_2^3g(x_2,x_3,x_4))=3r'+v'(g(x_2,x_3,x_4))>2r'+2,$$ where the last inequality follows from $r'=a'd-1\geq 3>2$. By \eqref{align: equation for cD case 2}, we have $$w'(g_0(x_3,x_4))=v'(g_0(x_3,x_4))\geq \frac{a'}{a}v(g_0(x_3,x_4))\geq \frac{a'}{a}(2r+2)=2a'd=2r'+2,$$ and $x_3^{2d}$ is the only possible monomial term in $g_0(x_3,x_4)$ whose weight equals $2r'+2$. We may further write $$g_1(x_3,x_4)=h_0(x_3)+x_4h_1(x_3)+x_4^2h(x_3,x_4)$$ for some analytic power series $h_0(x_3),h_1(x_3)$ and $h(x_3)$. By \eqref{align: equation for cD case 2}, for any monomial term $x_3^{l_1}\in g_1(x_3,x_4)$, $v(x_3^{l_1})=al_1\geq r+2=ad+1$, hence $l_1\geq d+1$, and $w'(x_3^{l_1})=a'l_1\geq a'(d+1)>r'+2$. It follows that $w'(h_0(x_3))>r'+2$. By \eqref{align: equation for cD case 2}, $v(x_4h_1(x_3))=1+v(h(x_3))\geq r+2$, hence $v(h_1(x_3))\geq r+1=ad$. We have $v'(h_1(x_3))\geq \frac{a'}{a}v(h_1(x_3))=a'd=r'+1$, and $x_2x_3^{d}x_4$ is the only possible monomial term in $x_2x_4h_1(x_3)$ whose weight equals $2r'+2$. By \eqref{align: equation for cD case 2}, $v(x_4^2h(x_3,x_4))\geq r+2$, hence $v(h(x_3,x_4))\geq r$, and $$v'(x_4^2h(x_3,x_4))\geq 2+\frac{a'}{a}v(h(x_3,x_4))\geq 2+\frac{a'}{a}r=r'+2+(1-\frac{a'}{a})>r'+2.$$ By \eqref{align: equation for cD case 2} again, $v(g_2(x_3,x_4))\geq 2$, hence $$v'(x_2^2g_2(x_3,x_4))=2r'+v'(g(x_3,x_4))\geq 2r'+v'(x_4^2)=2r'+2,$$ and $x_2^2x_4^2$ is the only possible monomial term in $x_2^2g_1(x_3,x_4)$ whose weight equals $2r'+2$. To sum up, $(g_0(x_3,x_4))_{w'}=\mu x_3^{2d}$, $(x_2g_1(x_3,x_4))_{w'}=\lambda_1 x_2x_3^dx_4$, $(x_2^2g_2(x_3,x_4))_{w'}=\lambda_2x_2^2x_4^2$ and $(x_2^3g(x_2,x_3,x_4))_{w'}=0$ for some $\mu,\lambda_1,\lambda_2\in \Cc$. This proves the claim. Since $w'(x_2x_5)=w'(x_1^2)=2(r'+1)$, we have $w'(\phi_1)=2(r'+1)$ and $$\phi_{1,w''}=x_1^2+x_2x_5+\mu x_3^{2d}+\lambda_1 x_2x_3^dx_4+\lambda_2 x_2^2x_4^2.$$

Since $w'(x_3^d)=w'(x_2x_4)=r'+1$, $w'(x_5)=r'+2>r'+1$ and $$w'(x_4q(x_3,x_4))=v'(x_4q(x_3,x_4))>\frac{a'}{a}w(x_4q(x_3x_4))\geq \frac{a'}{a}(r+1)=a'd=r'+1,$$ we have $w'(\phi_2)=r'+1$ and $$\phi_{2,w'}=x_2x_4+x_3^d.$$ Since $(\phi_{1,w'}=\phi_{2,w'}=0)$ defines an integral scheme near $o\in \Cc^5$, by Lemma~\ref{lem: weight leading term irr imply irr exceptional divisor}, the exceptional divisor of $f'$ is a prime divisor. Since $w'(X\ni x)=r'+1+r'+a'+1+r'+2-w'(\phi_1)-w'(\phi_2)-1$, $w'(\phi_1)=2(r'+1)$ and $w'(\phi_2)=r'+1$, $w'(X\ni x)=a'$.

\medskip

We now prove part (2) of the Lemma. 

By \eqref{align: equation for cD case 2}, we have $$w''(g_0(x_3,x_4))=v''(g_0(x_3,x_4))\geq \frac{1}{a}v(g_0(x_3,x_4))\geq \frac{1}{a}(2r+2)=2d,$$ and $x_3^{2d}$ is the only possible monomial term in $g_0(x_3,x_4)$ whose weight equals $2d$. Since
\begin{align*}
    w''(x_2g_1(x_3,x_4))&\geq d+\frac{1}{a}v(g_1(x_3,x_4))\geq a+\frac{1}{a}(r+2)=d+\frac{1}{a}(ad+1)>2d,\\
    w''(x_2^2g_2(x_3,x_4))&=2d+v''(g_2(x_3,x_4))>2d,\ and\\
    w''(x_2^3g_2(x_2,x_3,x_4))&=3d+v''(g_2(x_2,x_3,x_4))>2d,
\end{align*}
we have $(p(x_2,x_3,x_4))_{w''}=\mu x_3^{2d}$ for some $\mu\in \Cc$. Since $w''(x_2x_5)=w''(x_1^2)=2d$, $w''(\phi_1)=2d$ and $$\phi_{1,w''}=x_1^2+x_2x_5+\mu x_3^{2d}.$$ Since $w''(x_3^d)=w''(x_5)=d$, $w''(x_2x_4)=d+1>d$ and $$w''(x_4q(x_3,x_4))=v''(x_4q(x_3,x_4))>\frac{1}{a}v(x_4q(x_3x_4))\geq \frac{1}{a}(r+1)=d,$$ we have $w''(\phi_2)=d$, and $$\phi_{2,w''}=x_5+x_3^d.$$ Since $(\phi_{1,w''}=\phi_{2,w''}=0)$ defines an integral scheme near $o\in \Cc^5$, by Lemma~\ref{lem: weight leading term irr imply irr exceptional divisor}, the exceptional divisor of $f''$ is a prime divisor. Since $w''(X\ni x)=d+d+1+1+d-w''(\phi_1)-w''(\phi_2)-1$, $w'(\phi_1)=2d$ and $w'(\phi_2)=d$, $w'(X\ni x)=1$.
\end{proof}

\begin{lem}\label{lem: cD/2 irreducible case 1}
Let $f: Y\to X$ be a divisorial contraction of a prime divisor over $X\ni x$ and $r,a,d$ the corresponding positive integers as in Theorem~\ref{thm: kaw05 1.2 strenghthened}(3.1). For any weight $w':=\frac{1}{2}(r'+2,r',a',2)\in\frac{1}{2}\mathbb Z_{>0}^4$ such that $r'+1=a'd$ and $2\leq a'<a$, the weighted blow-up $f':Y'\to X$ with the weight $w'$ extracts an analytic prime divisor, and $w'(X\ni x)=\frac{a'}{2}$. 
\end{lem}

\begin{proof}
By Theorem~\ref{thm: kaw05 1.2 strenghthened}(3.1), analytically locally, we have  $$(X\ni x)\cong (\phi:= x_1^2+x_1x_3q(x_3^2,x_4)+x_2^2x_4+\lambda x_2x_3^{2\alpha-1}+p(x_3^2,x_4)=0)\subset (\mathbb C^4\ni o)/\frac{1}{2}(1,1,1,0),$$ and $f: Y\to X$ is a weighted blow-up with the weight $w:=\frac{1}{2}(r+2,r,a,2)$, such that 
\begin{itemize}
    \item $w(\phi)=w(x_2^2x_4)=r+1=ad$, where $a,r$ is odd,
    \item if $q(x^2_3,x_4)\not=0$, then $w(x_1x_3q(x_3^2,x_4))=r+1$, and
    \item $x_3^{2d}\in p(x_3^2,x_4)$.
\end{itemize}

Consider the weights $v=(v(x_3),v(x_4))=\frac{1}{2}(a,2)$ and $v'=(v'(x_3),v'(x_4))=\frac{1}{2}(a',2)$, we have $v'\succeq \frac{a'}{a}v$ (see Definition~\ref{defn: compare two weights}), and for any monomial term $x_3^{k}x_4^l$, $v'(x_3^{k}x_4^l)=\frac{a'}{a}v(x_3^{k}x_4^l)$ if and only if $l=0$.

If $q(x_3^2,x_4)\neq 0$, then $w(x_1x_3q(x_3^2,x_4))=\frac{1}{2}(r+2+a)+v(q(x_3^2,x_4))= r+1$, and $v(q(x_3^2,x_4))= \frac{r-a}{2}$. It follows that
\begin{align*}
    w'(x_1x_3q(x_3^2x_4))&=\frac{1}{2}(r'+2+a')+v'(q(x_3^2,x_4))\geq \frac{1}{2}(r'+2+a')+\frac{a'}{a}\frac{r-a}{2}\\
    &=\frac{1}{2}(r'+2+a')+\frac{a'}{a}\frac{a(d-1)-1}{2}=r'+1+\frac{1}{2}(1-\frac{a'}{a})> r'+1.
\end{align*}
Since $w(x_2x_3^{2\alpha-1})=\frac{1}{2}(r+(2\alpha-1)a)\geq r+1$, $2\alpha-1\geq \lceil\frac{r+2}{a}\rceil=\lceil\frac{ad+1}{a}\rceil=d+1$, hence
\begin{align*}
    w'(x_2x_3^{2\alpha-1})&=\frac{1}{2}(r'+a'(2\alpha-1))\geq \frac{1}{2}(r'+a'(d+1))=\frac{1}{2}(2r'+1+a')>r'+1,
\end{align*}
where the last inequality follows from $a'\geq 2$. Note that
\begin{align*}
    w'(p(x_3^2,x_4))&=v'(p(x_3^2,x_4))\geq\frac{a'}{a}(r+1)=r'+1=w'(x_3^{2d})\geq w'(p(x_3^2,x_4)),
\end{align*}
and $x_3^{2d}$ is the only monomial term in $p(x_3^2,x_4)$ whose weight equals $r'+1$. Note also that $w'(x_2^2x_4)=r'+1$ and $w'(x_1^2)=r'+2>r'+1$. Hence $w'(\phi)=r'+1$, and the $w'$-weighted leading term $$\phi_{w'}=x_2^2x_4+x_3^{2d}.$$ Since $(\phi_{w'}=0)$ defines an integral scheme near $o\in \Cc^4$, by Lemma~\ref{lem: weight leading term irr imply irr exceptional divisor}, the exceptional divisor of $f'$ is a prime divisor. Since $w'(X\ni x)=\frac{1}{2}(r'+2+r'+a'+2)-w'(\phi)-1$ and $w'(\phi)=r'+1$, $w'(X\ni x)=\frac{a'}{2}$.
\end{proof}

\begin{lem}\label{lem: cD/2 irreducible case 2}
Let $f: Y\to X$ be a divisorial contraction of a prime divisor over $X\ni x$ and $r,a,d$ the corresponding positive integers as in Theorem~\ref{thm: kaw05 1.2 strenghthened}(3.2). For any weight $w':=\frac{1}{2}(r'+2,r',a',2,r'+4)\in\frac{1}{2}\mathbb Z_{>0}^5$ such that $r‘+2=a'd$ and $a'<a$, the weighted blow-up $f':Y'\to X$ with the weight $w'$ extracts an analytic prime divisor, and $w'(X\ni x)=\frac{a'}{2}$. 
\end{lem}

\begin{proof}
By Theorem~\ref{thm: kaw05 1.2 strenghthened}(3.2), analytically locally, we have $$(X\ni x)\cong \binom{\phi_1:= x_1^2+x_2x_5+p(x_3^2,x_4)=0}{\phi_2:= x_2x_4+x_3^d+q(x_3^2,x_4)x_3x_4+x_5=0}\subset(\mathbb C^5\ni o)/\frac{1}{2}(1,1,1,0,1),$$ and $f: Y\to X$ is a weighted blow-up with the weight  $w:=\frac{1}{2}(r+2,r,a,2,r+4)$, such that 
\begin{itemize}
    \item $2w(\phi_2)=w(\phi_1)=r+2=ad$,
    \item $d$ is odd, and
    \item if $q(x_3^2,x_4)\not=0$, then $w(q(x_3^2,x_4)x_3x_4)=\frac{r+2}{2}$.
\end{itemize}

Consider the weights $v=(v(x_3),v(x_4))=\frac{1}{2}(a,2)$ and $v'=(v'(x_3),v'(x_4))=\frac{1}{2}(a',2)$, we have $v'\succeq \frac{a'}{a}v$ (see Definition~\ref{defn: compare two weights}), and for any monomial terms $x_3^kx_4^l$, $v'(x_3^{k}x_4^l)=\frac{a'}{a}v(x_3^{k}x_4^l)$ if and only if $l=0$. 

Since $w(p(x_3^2,x_4))=v(p(x_3^2,x_4))\geq r+2=ad$, we have $$w'(p(x_3^2,x_4))=v'(p(x_3^2,x_4))\geq \frac{a'}{a}v(p(x_3^2,x_4))\geq a'd=r'+2,$$ and $x_3^{2d}$ is the only possible monomial term in $p(x_3^2,x_4)$ whose weight equals $r'+2$. Note that $w'(x_2)=w'(x_2x_5)=r'+2$, hence $w'(\phi_1)=r'+2$ and $$\phi_{1,w'}=x_1^2+x_2x_5+\mu x_3^{2d}$$ for some $\mu\in \Cc$. Since $w'(x_2x_4)=w'(x_3^d)=\frac{1}{2}(r'+2)$, $w'(x_5)=\frac{1}{2}(r'+4)>\frac{1}{2}(r'+2)$ and $$w'(q(x_3^2,x_4)x_3x_4)=v'(q(x_3^2,x_4)x_3x_4)> \frac{a'}{a}v(q(x_3^2,x_4)x_3x_4)=\frac{1}{2}a'd=\frac{1}{2}(r'+2),$$ we have $w'(\phi_2)=\frac{1}{2}(r'+2)$ and $$\phi_{2,w'}=x_2x_4+x_3^d.$$ Since $(\phi_{1,w'}=\phi_{2,w'}=0)$ defines an integral scheme near $o\in \Cc^5$, by Lemma~\ref{lem: weight leading term irr imply irr exceptional divisor}, the exceptional divisor of $f'$ is a prime divisor. Since $w'(X\ni x)=\frac{1}{2}(r'+2+r'+a'+2+r'+4)-w'(\phi_1)-w'(\phi_2)-1$, $w'(\phi_1)=r'+2$ and $w'(\phi_2)=\frac{1}{2}(r'+2)$, $w'(X\ni x)=\frac{a'}{2}$.
\end{proof}

\section{Ideal-adic versions}\label{sec: Ideal-adic versions}

In some scenarios, especially in the study of fixed germs and smooth varieties, people may care more about the structure of singularities of ``ideal-adic pairs" (or ``ideal-adic triples"). That is, pairs of the form $(X\ni x,\aaa)$ (or triples of the form $(X,B,\aaa)$), where $\aaa$ is an $\Rr$-ideal instead of a boundary divisor $B$ on $X$. This kind of pairs also appear in some literature in the study of mlds (cf. \cite{Kaw14,Kaw15b,Nak16,MN18,NS20,Kaw21,NS21}). For future researchers' convenience, we prove the ideal-adic version of some of our main theorems, and apply them to some conjectures on singularities of fixed germs.

We will adopt the definitions and notation as in \cite[Section 2]{Kaw21}, except that we will use $\mld(X\ni x,B,\aaa)$ instead of $\mld_x(X,B,\aaa)$. First we prove a useful lemma.

\begin{lem}\label{lem: transform pair to ideal adic}
Let $\Ii\subset [0,+\infty)$ be a DCC (resp. finite) set. Then there exists a positive real number $\delta$ and a DCC (resp. finite) set $\Ii'\subset [0,+\infty)$ depending only on $\Ii$ satisfying the following. Assume that $(X\ni x,B,\aaa)$ is a $\Qq$-factorial threefold triple such that $X$ is terminal, $B,\aaa\in\Ii$ and $\mld(X\ni x,B,\aaa)\geq 1-\delta$. Suppose that $\aaa\not=\mathfrak{0}$. Then
\begin{enumerate}
    \item there exists an $\Rr$-ideal $\aaa'\in\Ii'$, such that $a(E,X,B,\aaa)=a(E,X,\aaa')$ for any prime divisor $E$ over $X\ni x$, and
    \item there exists an $\Rr$-divisor $B'\in\Ii'$, such that $\mld(X\ni x,B')=\mld(X\ni x,B,\aaa)$.
\end{enumerate}
\end{lem}
\begin{proof}
Let $\gamma_0:=\min\{1,\gamma\in\Ii\mid \gamma>0\}$ and $I:=\lceil\frac{2}{\gamma_0}\rceil!$. We show that we can take $\delta=\frac{\gamma_0}{2}$ and $\Ii':=\{\sum n_i\gamma_i\mid n_i\in\Zz_{\geq 0},\gamma_i\in\frac{1}{I}\Ii'\}$. Since $\aaa\not=\mathfrak{0}$,
$$\mld(X\ni x)\geq\mld(X\ni x,B,\aaa)+\gamma_0\geq 1+\frac{\gamma_0}{2}.$$
By the classification of terminal threefold singularities and \cite[Lemma~5.1]{Kaw88}, $ID$ is Cartier near $x$ for any $\Qq$-Cartier Weil divisor on $X$.  Assume that $B=\sum b_iB_i$ where $B_i$ are the distinct prime divisors of $B$, then we may take $\aaa':=\aaa\cdot\prod_i\bbb_i^{\frac{b_i}{I}}$, where $\bbb_i:=\mathcal{O}_X(-IB_i)$. This proves (1).

To prove (2), possibly replacing $(X\ni x,B,\aaa)$ with $(X,\aaa')$ we may assume that $B=0$. Assume that $\aaa=\prod_i\aaa_i^{\gamma_i}$ where each $\gamma_i\in\Ii$ and each $\aaa_i$ is an ideal. By \cite[Proposition 2.2 iv)]{MN18}, possibly replacing each $\aaa_i$ with $\aaa_i+\mmm_x^l$ for some $l\gg 0$ where $\mmm_x$ is the maximal ideal of $\mathcal{O}_{X,x}$, we may assume that each $\aaa_i$ is $\mmm_x$-primary. We take a log resolution $f: Y\rightarrow X$ of $(X\ni x,\aaa)$. For each $i$, we have $\aaa_i\mathcal{O}_Y=\mathcal{O}_X(-A_i)$ for some $A_i\geq 0$ such that $-A_i$ is free over $X$. Pick general Weil divisors $H_i\in |-A_i|/X$ for each $i$. Then $f$ is also a log resolution of $(X\ni x,\sum_i\gamma_iH_i)$, and we may take $B'=\sum_i\gamma_iH_i$.
\end{proof}

Now we prove the ideal-adic ACC conjecture for mlds for canonical threefold triples. Note that in the following theorem, we do not require $X\ni x$ to be a fixed germ.

\begin{thm}[Ideal-adic version of Theorem {\ref{thm: intro global canonical mld acc}}]\label{thm: ideal-adic acc mld [1,3]}
Let $\Ii\subset [0,+\infty)$ be a DCC set. Then there exists a positive real number $\delta$ depending only on $\Ii$, such that
$$\{\mld(X,B,\aaa)\mid \dim X=3, B,\aaa\in\Ii\}\cap [1-\delta,+\infty)$$
satisfies the ACC.
\end{thm}
\begin{proof}
Possibly replacing $X$ with a small $\Qq$-factorialization, we may assume that $X$ is $\Qq$-factorial.  If $X$ is not terminal and $\aaa\not=\mathfrak{0}$, then the theorem is trivial. If $X$ is not terminal and $\aaa=\mathfrak{0}$, then the theorem follows from Theorem \ref{thm: intro global canonical mld acc}. If $X$ is terminal, then the theorem follows from Theorems \ref{thm:  terminal mld acc} and \ref{thm: 1-gap pair} and Lemma \ref{lem: transform pair to ideal adic}.
\end{proof}

We prove the ideal-adic ACC conjecture for threefold $a$-lc thresholds when $a\geq 1$. Notice that when $a=1$, we get the ideal-adic ACC for canonical thresholds for threefolds. Moreover, we do not require $X\ni x$ to be a fixed germ.

\begin{thm}[Ideal-adic version of Theorem \ref{thm: alct acc terminal threefold}]\label{thm: ideal-adic alct acc terminal threefold}
Let $a\geq 1$ be a real number and $\Ii,\Ii'\subset [0,+\infty)$ two DCC sets. Then
$$\{a\text{-}\lct(X\ni x,B,\aaa;D,\bbb)\mid \dim X=3, X\text{ is terminal, }B,\aaa\in\Ii,D,\bbb\in\Ii'\}$$
satisfies the ACC.
\end{thm}
\begin{proof}
Suppose that the theorem does not hold, then there exist a strictly increasing sequence of positive real number $t_i$, a sequence of lc triples $(X_i\ni x_i,B_i,\aaa_i)$ of dimension $3$ such that $X_i$ is terminal and $B_i,\aaa_i\in\Ii$, a sequence of $\Rr$-Cartier $\Rr$-divisors $D_i\in\Ii'$ and a sequence of $\Rr$-ideals $\bbb_i\in\Ii'$, such that $a\text{-}\lct(X_i\ni x_i,B_i,\aaa_i;D_i,\bbb_i)=t_i>0$. 
Possibly replacing each $X_i$ with a small $\Qq$-factorialization, we may assume that each $X_i$ is $\Qq$-factorial. 

Suppose that $\mld(X_i\ni x_i,B_i+t_iD_i,\aaa_i\bbb^{t_i}_i)>a$ for infinitely many $i$. Possibly passing to a subsequence, we have that $(X_i,B_i+t_iD_i)$ is lc but not klt near $x_i$, and $t_i=\lct(X_i\ni x_i,B_i;D_i)$, which contradicts \cite[Theorem 1.1]{HMX14}. Thus we may assume that $\mld(X_i\ni x_i,B_i+t_iD_i,\aaa_i\bbb^{t_i}_i)=a$ for each $i$.

By Theorem \ref{thm: alct acc terminal threefold}, we may assume that either $\aaa_i\not=\mathfrak{0}$ or $\bbb_i\not=\mathfrak{0}$ for each $i$. Let $\gamma_0:=\min\{\gamma\mid\gamma>0, \gamma\in\Ii\text{ or }\gamma\in t_1\Ii'\}$. Then 
$$\mld(X_i\ni x_i)\geq \mld(X_i\ni x_i,B_i+t_iD_i,\aaa_i\bbb^{t_i}_i)+\gamma_0=a+\gamma_0\geq 1+\gamma_0.$$
By the classification of terminal threefold singularities and \cite[Lemma~5.1]{Kaw88}, there exists a positive integer $I$ depending only on $\Ii$ such that $ID$ is Cartier near $x$ for any $\Qq$-Cartier Weil divisor $D$ on $X$. Assume that $B_i=\sum b_{i,j}B_{i,j}$ and $D_i=\sum d_{i,j}D_{i,j}$ where $B_{i,j}$ and $D_{i,j}$ are the irreducible components of $B_i$ and $D_i$ respectively. Possibly replacing $\Ii$ with $\Ii\cup\frac{1}{I}\Ii$, $\Ii'$ with $\Ii'\cup\frac{1}{I}\Ii'$, $\aaa_i$ with $\aaa_i\prod_j\rrr_{i,j}^{\frac{b_{i,j}}{I}}$, and $\bbb_i$ with $\bbb_i\prod_j\ddd_{i,j}^{\frac{d_{i,j}}{I}}$, where $\rrr_i:=\mathcal{O}_X(-IB_i)$ and $\ddd_i=\mathcal{O}_X(-ID_i)$, we may assume that $B=D=0$.

Assume that $\aaa_i=\prod_j\aaa_{i,j}^{\gamma_{i,j}}$ and $\bbb_i=\prod_j\bbb_{i,j}^{\gamma'_{i,j}}$ where each $\gamma_i\in\Ii$, $\gamma_i'\in\Ii'$ and each $\aaa_{i,j},\bbb_{i,j}$ is an ideal. By \cite[Proposition 2.2 iv)]{MN18}, possibly replacing each $\aaa_{i,j}$ with $\aaa_{i,j}+\mmm_{x_i}^{l_j}$ and $\bbb_{i,j}$ with $\bbb_{i,j}+\mmm_{x_i}^{l_j}$ for some $l_j\gg 0$, where $\mmm_{x_i}$ is the maximal ideal of $\mathcal{O}_{X_i,x_i}$, we may assume that each $\aaa_{i,j}$ and $\bbb_{i,j}$ are $\mmm_{x_i}$-primary. We take a log resolution $f_i: Y_i\rightarrow X_i$ of $(X_i\ni x_i,\aaa_i\bbb^{t_i}_i)$. For each $i,j$, we have $\aaa_{i,j}\mathcal{O}_{Y_i}=\mathcal{O}_{X_i}(-A_{i,j})$ and $\bbb_{i,j}\mathcal{O}_{Y_i}=\mathcal{O}_{X_i}(-R_{i,j})$ for some $A_{i,j},R_{i,j}\geq 0$ such that $-A_{i,j},-R_{i,j}$ are free over $X_i$. Thus we may pick general Weil divisors $H_{i,j}\in |-A_{i,j}|/X_i$ and $S_{i,j}\in |-R_{i,j}|/X_i$. Then $f_i$ is also a log resolution of $(X_i\ni x_i,\sum_j\gamma_{i,j}H_{i,j}+t_i\sum_j\gamma_{i,j}'S_{i,j})$, and by our construction, $\mld(X_i\ni x_i,\sum_j\gamma_{i,j}H_{i,j}+t_i\sum_j\gamma_{i,j}'S_{i,j})=\mld(X_i\ni x_i,\aaa_i\bbb_i^{t_i})=a$. But this contradicts Theorem \ref{thm: index conjecture DCC threefold} as $t_i$ is strictly increasing.
\end{proof}

Now we prove the ideal-adic version of Theorem \ref{thm: mn intro DCC coeff}. Note that we do not require $X\ni x$ to be a fixed germ.

\begin{thm}[Ideal-adic version of Theorem \ref{thm: mn intro DCC coeff}]\label{thm: ideal-adic mn DCC coeff}
Let $\Ii\subset [0,+\infty)$ be a DCC set. Then there exists a positive integer $l$ depending only on $\Ii$ satisfying the following. Assume that $(X\ni x,B,\aaa)$ is a threefold lc triple such that $X$ is terminal, $B,\aaa\in\Ii$, and $\mld(X\ni x,B,\aaa)\geq 1$. Then there exists a prime divisor $E$ over $X\ni x$, such that $a(E,X,B,\aaa)=\mld(X\ni x,B,\aaa)$ and $a(E,X,0)\leq l$.
\end{thm}
\begin{proof}
Possibly replacing $X$ with a small $\Qq$-factorialization, we may assume that $X$ is $\Qq$-factorial. If $\aaa=\mathfrak{0}$, then the theorem follows Theorem \ref{thm: mn intro DCC coeff}. Thus we may assume that $\aaa\not=\mathfrak{0}$. By Lemma \ref{lem: transform pair to ideal adic}(1), we may assume that $B=0$. 

Assume that $\aaa=\prod_i\aaa_i^{\gamma_i}$ where each $\gamma_i\in\Ii$ and each $\aaa_i$ is an ideal. Since $\mld(X\ni,\aaa)\geq 1$, there are only finitely many prime divisors $E_1,\dots,E_m$ over $X\ni x$ which compute $\mld(X\ni x,\aaa)$. By \cite[Proposition 2.2 iv]{MN18}, possibly replacing each $\aaa_i$ with $\aaa_i+\mmm_x^l$ for some $l\gg 0$, where $\mmm_x$ is the maximal ideal of $\mathcal{O}_{X,x}$, we may assume that each $\aaa_i$ is $\mmm_x$-primary. We take a log resolution $f: Y\rightarrow X$ of $(X\ni x,\aaa)$ such that $E_1,\dots,E_m$ are on $Y$. For each $i$, we have $\aaa_i\mathcal{O}_Y=\mathcal{O}_X(-A_i)$ for some $A_i\geq 0$ such that $-A_i$ is free over $X$. Pick general divisors $H_i\in|-A_i|/X$, then $f$ is also a log resolution of $(X\ni x,\sum_i\gamma_iH_i)$. By our construction, 
$$a(E_k,X,\sum_i\gamma_iH_i)=\mld(X\ni x,\sum_i\gamma_iH_i)=\mld(X\ni x,\aaa)=a(E_k,X,\aaa)$$
for any $1\leq k\leq m$. Theorem \ref{thm: ideal-adic acc mld [1,3]} follows from Theorem \ref{thm: mn intro DCC coeff}.
\end{proof}

There are two conjectures that are known to be equivalent to the ACC conjecture for mld for fixed germs: the \emph{uniform $\mmm$-adic semicontinuity conjecture} and the \emph{generic limit conjecture}. We can also investigate these two conjectures for terminal threefold pairs.

First we have the uniform $\mmm$-adic semi-continuity conjecture (cf. \cite[Conjecture 7.1]{MN18}, \cite[Conjecture 4.2(iii)]{Kaw21})  for canonical threefold triples. We remark that in the theorem below, $l$ is independent of the choice of $X$, which is not expected in the original statement of the uniform ideal-semicontinuity conjecture. Moreover, we can choose $\Ii$ to be any DCC set, which is also not expected in the original statement of the conjecture where only the finite coefficient case is expected.

\begin{thm}[Uniform $\mmm$-adic semicontinuity for canonical threefold triples]\label{thm: ideal adic semicontinuity termial threefolds}
Let $\Ii$ be a DCC set of real numbers. Then there exists a positive real integer $l$ depending only on $\Ii$ satisfying the following. Assume that $(X\ni x,B,\aaa)$ and $(X\ni x,B,\bbb)$ are two threefold triples, such that
\begin{itemize}
    \item $X$ is terminal,
    \item $\mld(X\ni x,B,\aaa)\geq 1$ and $\mld(X\ni x,B,\bbb)\geq 1$,
\item $\aaa=\prod^m_i\aaa_i^{\gamma_i}$ and $\bbb=\prod^m_i\bbb_i^{\gamma_i}$, where $m>0$, $\aaa_i$ and $\bbb_i$ are ideals, $\gamma_i\in\Ii$, and
    \item $\aaa_i+\mmm^l=\bbb_i+\mmm^l$ for each $i$, where $\mmm$ is the maximal ideal of $\mathcal{O}_{X,x}$.
\end{itemize}
Then $\mld(X\ni x,B,\aaa)=\mld(X\ni x,B,\bbb)$.
\end{thm}
\begin{proof}
Possibly replacing $X$ with a small $\Qq$-factorialization, we may assume that $X$ is $\Qq$-factorial. If either $\aaa$ or $\bbb$ is $\mathfrak{0}$, then since $\aaa_i+\mmm^l=\bbb_i+\mmm^l$ for each $i$, we have $\aaa=\bbb=\{0\}$, and the Theorem is obvious. Thus we may assume that $\aaa,\bbb\not=\mathfrak{0}$. By Lemma \ref{lem: transform pair to ideal adic}, we may assume that $B=0$.

We only need to show that $\mld(X\ni x,\aaa)\geq\mld(X\ni x,\bbb)$. Let $\gamma_0:=\min\{\gamma\in\Ii, 1\mid \gamma>0\}$. By Theorem \ref{thm: ideal-adic mn DCC coeff}, there exists a positive real number $M$ depending only on $\Ii$ and a prime divisor $E$ over $X\ni x$, such that $a(E,X,\aaa)=\mld(X\ni x,\aaa)$ and $a(E,X,0)\leq M$. We show that $l:=\frac{M}{\gamma_0}+1$ satisfies our requirements. We have
$$1\leq\mld(X\ni x,\aaa)=a(E,X,\aaa)=a(E,X,0)-\mult_E\aaa\leq M-\sum_{i=1}^m\gamma_i\mult_E\aaa_i.$$
Thus we have $\sum_{i=1}^m\gamma_i\mult_E\aaa_i\leq M-1$, which implies that $\mult_{E}\aaa_i\leq\frac{M}{\gamma_i}<l\leq l\mult_E\mmm_x$ for each $i$, where $\mmm_x$ is the maximal ideal of $\mathcal{O}_{X,x}$. By \cite[Proposition 2.2 iv)]{MN18}, we have
\begin{align*}
    \mld(X\ni x,\aaa)&=\mld(X\ni x,\prod_i\aaa_i^{\gamma_i})=\mld(X\ni x,\prod_{i=1}^m(\aaa_i+\mmm^l)^{\gamma_i})\\
    &=\mld(X\ni x,\prod_{i=1}^m(\bbb_i+\mmm^l)^{\gamma_i})\geq\mld(X\ni x,\prod_{i=1}^m\bbb_i^{\gamma_i})=\mld(X\ni x,\bbb),
\end{align*}
and the theorem follows.
\end{proof}

In the end, we prove the generic limit conjecture (cf. \cite[Conjecture 5.7]{Kaw14}) for canonical threefold pairs.

\begin{thm}[Generic limit conjecture for canonical threefold pairs]\label{thm: gl conjecture terminal threefolds}
Let $\gamma_1,\dots,\gamma_m\in [0,1]$ be real numbers. Let $x\in X$ be a terminal threefold singularity and $\{\aaa_i=\prod_{i=1}^m\aaa_{i,j}^{\gamma_j}\}_{i=1}^{\infty}$ a sequence of $\Rr$-ideals on $X$, such that $\mld(X\ni x,\aaa_i)\geq 1$ for each $i$ and each $\aaa_{i,j}$ is an ideal. Let $(\widehat X\ni \widehat x,\widehat\aaa)$ be the generic limit of $\{(X\ni x,\aaa_i)\}_{i=1}^{\infty}$ and let $\widehat\aaa_{j}$ be the generic limit of $\{\widehat\aaa_{i,j}\}_{i=1}^{\infty}$ for each $j$. Then there exists an infinite subset $\Lambda\subset\Zz_{>0}$, such that
\begin{enumerate}
    \item $\widehat\aaa_j$ is the generic limit of $\{\aaa_{i,j}\}_{i\in\Lambda}$ for each $j$, and
    \item $\mld(\widehat X\ni\widehat x,\widehat\aaa)=\mld(X\ni x,\aaa_i)$ for every $i\in\Lambda$.
\end{enumerate}
\end{thm}
\begin{proof}
By Theorem \ref{thm: ideal adic semicontinuity termial threefolds} and \cite[Proposition 2.2 iv)]{MN18}, there exists a positive integer $l$, such that for any integer $s\geq l$ and any $i$, we have
\begin{equation}\label{equ: first equation gl conj}
\mld(X\ni x,\prod_{i=1}^m(\aaa_{i,j}+\mmm_x^s)^{\gamma_j})=\mld(X\ni x,\prod_{i=1}^m\aaa_{i,j}^{\gamma_j})=\mld(X\ni x,\aaa_i),  
\end{equation}
where $\mmm_x$ is the maximal ideal of $\mathcal{O}_{X,x}$. 

Let $E$ be a prime divisor over $\widehat X\ni\widehat x$ which computes $\mld(\widehat X\ni\widehat x,\widehat\aaa)$. Let $s_0\geq l$ be an integer, such that $s_0\mult_E\mmm_{\widehat x}\geq\mult_E\widehat a_j$ for each $j$, where $\mmm_{\widehat x}$ is the maximal ideal of $\mathcal{O}_{\widehat X,\widehat x}$. By \cite[Proposition 2.2 iv)]{MN18}, we have 
\begin{equation}\label{equ: second equation gl conj}
\mld(\widehat X\ni\widehat x,\widehat\aaa)=\mld(\widehat X\ni\widehat x,\prod_{i=1}^m\widehat\aaa_j^{\gamma_j})=\mld(\widehat X\ni\widehat x,\prod_{i=1}^m(\widehat\aaa_j+\mmm_{\widehat x}^{s_0})^{\gamma_j}).
\end{equation}
By \cite[Remark 3.2]{MN18}, there exists an infinite subset $\Lambda\subset\Zz_{>0}$, such that
\begin{equation}\label{equ: third equation gl conjecture}
\mld(\widehat X\ni\widehat x,\prod_{i=1}^m(\widehat\aaa_j+\mmm_{\widehat x}^{s_0})^{\gamma_j})=\mld(X\ni x,\prod_{i=1}^m(\aaa_{i,j}+\mmm_x^{s_0})^{\gamma_j})
\end{equation}
for any $i\in\Lambda$. We get the desired result by combining \eqref{equ: first equation gl conj}, \eqref{equ: second equation gl conj}, and \eqref{equ: third equation gl conjecture}.
\end{proof}

\section{Uniform boundedness of divisors computing mlds for toric varieties}
\begin{center}
	Joint work with\\
	Lingyao Xie and Qingyuan Xue\\
Department of Mathematics, The University of Utah, Salt Lake City, UT 84112, USA\\  \email{lingyao@math.utah.edu}, \email{xue@math.utah.edu}\\
	\footnote[1]{We thank Guodu Chen, Junpeng Jiao, and Lu Qi for useful discussions. LX and QX are partially supported by NSF research grants no: DMS-1801851, DMS-1952522 and by a grant from the Simons Foundation; Award Number: 256202. }
\end{center}

In this section, we adopt the standard notation and definitions in \cite{Oda88} and \cite{Amb06} and will freely use them. 

Let $\Delta$ be a fan in a lattice $N$. We denote by $T_N \rm{emb}(\Delta)$ the toric variety with the normal fan $\Delta$. Each lattice point of $N$ determines a vector on $N_{\Rr}$ starting from the origin. By saying the lattice points $e_1,\dots,e_k\subset N$ are linearly independent, we mean the corresponding vectors in $N_{\Rr}$ are linearly independent. We denote by $\Delta(1)$ the set of primitive lattice points that generate the one-dimensional cones of $\Delta$. Let $\sigma$ be a strongly convex cone in a lattice $N$ and $o\in T_N\rm{emb}(\sigma)$ the corresponding torus invariant closed point. Then each torus invariant prime divisor over $T_N\rm{emb}(\sigma)\ni o$ corresponds to a lattice point in $N\cap \rm{relin}(\sigma)$.

\begin{thm}\label{thm: toric mld bdd}
Let $(X\ni x, B)$ be a $\Qq$-Gorenstein toric lc pair such that $\dim X=d$. Then there exists a torus invariant prime divisor $E$ over $X\ni x$, such that $a(E,X,B)=\mld(X\ni x, B)$, and $a(E,X,0)\leq d$. Moreover, if $(X\ni x, B)$ is klt, then for any torus invariant prime divisor $E$ over $X\ni x$ such that $a(E,X,B)=\mld(X\ni x, B)$, we have $a(E,X,0)\leq d$.
\end{thm}

As a consequence of Theorem~\ref{thm: toric mld bdd}, Conjecture~\ref{conj: bdd mld computing divisor} holds for all toric lc pairs.

\begin{proof}[Proof of Theorem~\ref{thm: toric mld bdd}]
The proof follows directly from the definition of toric mlds.

Let $\Delta$ be a fan in a lattice $N$ such that $X=T_N\rm{emb}(\Delta)$, and $B=\sum_{i} b_i V(e_i)$ for some $e_i\in\Delta(1)$, where $V(e_i)\subset X$ is the torus invariant prime divisor corresponding to $e_i$. Recall that $K_X=-\sum_{e_i\in \Delta(1)} V(e_i)$, and the
$\Rr$-Cartier property of $K_X+B$ (resp. $K_X$) means that there exists a function $\psi: |\Delta|\to \Rr$ (resp. $\psi_0: |\Delta|\to \Rr$) such that
$\psi(e_i)=1-b_i$ (resp. $\psi_0(e_i)=1$) for every $i$, and $\psi$ (resp. $\psi_0$) is a linear function restricted on each cone of $\Delta$.

There exists a unique convex cone $\sigma\subset |\Delta|$ such that $x\in \rm{orb}(\sigma)$. Let $c$ be the
codimension of $\rm{orb}(\sigma)$ in $X$. Consider the toric pair $$(X'\ni x', B')\times \mathbb{A}_{\Cc}^{d-c}:=(T_{N}\cap(\sigma-\sigma)\rm{emb}(\sigma)\ni o,\sum_{e\in \sigma(1)}\mult_{V(e)} V(e))\times \mathbb{A}_{\Cc}^{d-c},$$ where $o$ is a torus invariant closed point, and $\sigma-\sigma:=\{x-y\mid x,y\in \sigma\}$. For any torus invariant prime divisor $E$ over $X\ni x$, there exists a torus invariant prime divisor $E'$ over $X'\ni x'$, such that $a(E,X,B)=a(E',X'B')$, and $a(E,X,0)=a(E',X',0)$. Then possibly replacing $(X\ni x, B)$ with $(X'\ni x',B')$ and $d$ with $d-c$, we may assume that $\dim \sigma=d$, $\sigma$ is strongly convex, $X=T_N\rm{emb}(\sigma)$, and $x$ is a torus invariant closed point.

Recall that $$\mld(X\ni x, B)=\min(\psi|_{N\cap\rm{relint}(\sigma)}),$$ there exists a lattice point $e\in\rm{relint}(\sigma)$ such that $\psi(e)=\mld(X\ni x, B)$. Since $\psi_0(e)=a(E,X,0)$, it suffices to prove that $\psi_0(e)\leq d$. By Carath\'a{e}odory’s Theorem
(see \cite[Theorem~A.15]{Oda88}), there exist a positive integer $s\leq d$ and a subset $\{e_1,\dots, e_s\} \subset \sigma(1)$, such
that $e_1,\dots,e_s$ are linearly independent and $e$ belongs to the relative interior of the cone
spanned by $e_1,\dots,e_s$. We may write $e=\sum_{i=1}^s \lambda_ie_i$. Since $(X\ni x, B)$ is lc, $\psi(e_i)\geq 0$ for $1\leq i\leq s$, hence possibly replacing $e$ with $\sum_{i=1}^s (\lambda_i+1-\lceil \lambda_i\rceil)e_i$, we may assume that $\lambda_i\in (0,1]$, it follows that $$a(E,X,0)=\psi_0(e)=\sum_{i=1}^s \lambda_i\psi_0(e_i)=\sum_{i=1}^s \lambda_i\leq s\leq d.$$ When $(X\ni x, B)$ is klt, $\psi(e_i)>0$, hence for any $e=\sum_i \lambda_ie_i$ such that $\psi(e)=\min(\psi|_{N\cap\rm{relint}(\sigma)})$, $0<\lambda_i\leq 1$, and $$a(E,X,0)=\psi_0(e)=\sum_{i=1}^s \lambda_i\psi_0(e_i)=\sum_{i=1}^s \lambda_i\leq s\leq d.$$ This finishes the proof.
\end{proof}

\begin{rem}
It is shown in \cite[Proposition~4.2]{Amb99} that if $(X\ni x,B)$ is a toric germ of dimension $d$ such that $\mld(X\ni x, B)>d-1$, then $X$ is smooth near $x$.
\end{rem}

\end{document}